\documentclass{amsart}
\usepackage{amsmath}
\usepackage{amssymb}
\usepackage{amsfonts}
\usepackage{amsthm}
\usepackage{graphicx}
\DeclareMathOperator{\curl}{curl}

\DeclareMathOperator{\diverg}{div}

\DeclareMathOperator{\im}{Im}

\DeclareMathOperator{\re}{Re}

\newtheorem{lemma}{Lemma}[section]

\theoremstyle{definition}

\theoremstyle{definition}
\newtheorem{remark}{Remark}[section]
\theoremstyle{definition}
\newtheorem{example}{Example}[section]
{\catcode `\@=11 \global\let\AddToReset=\@addtoreset}
\AddToReset{equation}{section}

\newcommand{\N}{{\mathbb N }}
\newcommand{\R}{{\mathbb R}}
\newcommand{\C}{{\mathbb C}}
\newcommand{\e}{{\varepsilon }}
\newcommand{\de}{{\delta}}
\newcommand{\ie}{{\sl i.e.\/ }}
\newcommand{\cf}{{\sl cf.\/ }}
\newcommand{\eg}{{\sl e.g.\/}}

\newcommand{\Eb}{\boldsymbol{\rm{E}}}
\newcommand{\Ab}{\boldsymbol{\rm{A}}}
\newcommand{\Bb}{\boldsymbol{\rm{B}}}
\newcommand{\Jb}{\boldsymbol{\rm{J}}}
\newcommand{\alb}{\boldsymbol{\alpha}}
\newcommand{\xb}{\boldsymbol{\rm{x}}}

\newcommand{\p}{\partial}

\newcommand{\Og}{\Omega}
\newcommand{\f}[2]{\frac{#1}{#2}}
\newcommand{\dt}{\delta}

\newcommand{\nn}{\nonumber}
\newcommand{ \ap}{\alpha}

\newcommand{\Ld}{\Lambda}
\newcommand{\ld}{\lambda}

\newcommand{\vp}{\varphi}

\newcommand{\ift}{ \infty}
\newcommand{\vep}{\varepsilon}

\newcommand{\Dt}{\Delta}

\newcommand{\btd}{\nabla}
\newcommand{\btu}{\Delta}
\newcommand{\tg}{\triangle}

\newcommand{\be}{\begin{equation}}
\newcommand{\ee}{\end{equation}}
\newcommand{\ba}{\begin{array}}
\newcommand{\ea}{\end{array}}
\newcommand{\bea}{\begin{eqnarray}}
\newcommand{\eea}{\end{eqnarray}}
\newcommand{\beas}{\begin{eqnarray*}}
\newcommand{\eeas}{\end{eqnarray*}}
\newcommand{\dpm}{\displaystyle}

\def\({\left(}
\def\){\right)}
\def\<{\left\langle}
\def\>{\right\rangle}
\def\O{\mathcal O}
\newcommand{\newpar}{\par}\parindent =0pt\parskip=3pt\textheight = 615pt

\newcommand{\Id}[1]{{\rm I\kern-2pt I_{#1}}}
\renewcommand{\hbar}{{\displaystyle\bar{\phantom{x}}\kern-6pt h}}

\numberwithin{equation}{section}

\begin{document}
\title[A time-splitting spectral method for Maxwell-Dirac system]
{A time-splitting spectral scheme for the Maxwell-Dirac system}
\author[Z. Huang]{Zhongyi Huang}
\author[S. Jin]{Shi Jin}
\author[P. A. Markowich]{Peter A. Markowich}
\author[C. Sparber]{Christof Sparber}
\author[C. Zheng]{Chunxiong Zheng}
\address[Z. Huang and C. Zheng]{Dept. of Mathematical Sciences, Tsinghua University,
Beijing 100084, China}
\email{zhuang@math.tsinghua.edu.cn}
\email{czheng@math.tsinghua.edu.cn}
\address[S. Jin]{Dept. of Mathematics,
University of Wisconsin, Madison, WI 53706, USA and
Dept. of Mathematical Sciences, Tsinghua University,
Beijing 100084, China}
\email{jin@math.wisc.edu}
\address[P. A. Markowich and C. Sparber]{Fakult\"at f\"ur Mathematik der
Universit\"at Wien\\ Nordbergstra\ss e 15\\ A-1090 Vienna\\ Austria}
\email{peter.markowich@univie.ac.at}
\email{christof.sparber@univie.ac.at}
\begin{abstract}
We present a time-splitting spectral scheme for the Maxwell-Dirac system and
 similar time-splitting methods for the corresponding asymptotic problems
in the semi-classical and the non-relativistic regimes.
The scheme for the Maxwell-Dirac system conserves
the Lorentz gauge condition, is unconditionally stable
and highly efficient as our numerical examples show. In particular we focus in
our examples on the creation of positronic modes in the semi-classical regime and on the
electron-positron interaction in the non-relativistic regime. Furthermore, in the non-relativistic regime,
our numerical method exhibits uniform convergence in the small parameter $\dt$, which is
the ratio of the characteristic speed and the speed of light.
\end{abstract}
\thanks{This work was partially supported by the Wittgenstein Award 2000 of
P.~A.~M., NSF grant No. DMS-0305080, the NSFC Projects no. 10301017 and
10228101,  Basic Research Projects of Tsinghua
University number JC 2002010, SRF for ROCS, SEM
and the Austrian-Chinese Technical-Scientific
Cooperation Agreement.
C.S. has been supported by the APART grant of the Austrian Academy of
Science.}
\maketitle
\begin{center}
version: \today
\end{center}
\section{Introduction and asymptotic scaling}
The \emph{Maxwell-Dirac system} (MD) describes
the time-evolution of fast, \ie \emph{relativistic spin-$1/2$ particles}, say electrons and positrons,
within external and \emph{self-consistent} electromagnetic fields. In \emph{Lorentz gauge}
it is given by the following set of equations:
\begin{equation}\label{dm0}
\left \{
\begin{aligned}
& i \hbar \partial_t  \psi   = \sum_{k=1}^3 \alpha^k
\left(\frac{\hbar c}{i}\partial_k -q (A_k+A_k^{ex})\right) \psi +
q(V+V^{ex})\psi +mc^2 \beta \psi ,\\
& \left(\frac{1}{c^{2}} \partial_{tt} - \Delta \right) V = \frac{1}{4\pi\epsilon_0} \rho,
\quad \left(\frac{1}{c^{2}} \partial_{tt} - \Delta\right) \Ab = \frac{1}{4\pi \epsilon_0c}
\,\Jb,\qquad \xb \in \R^3, \ t\in \R,
\end{aligned}
\right.
\end{equation}
subject to \emph{Cauchy initial data}:
\begin{equation}
\label{dm1}
\left \{
\begin{aligned}
& V \big| _{t=0}=V^{(0)}(\xb), \quad \, \partial _{t}V \big| _{t=0}= V^{(1)} (\xb) ,  \\
&\Ab \big| _{t=0}=\Ab^{(0)}(\xb)  ,
\quad  \partial _{t}\Ab \big| _{t=0}= \Ab^{(1)} (\xb), \quad \psi  \big| _{t=0}= \psi^{(0)}(\xb).
\end{aligned}
\right.
\end{equation}
The \emph{particle-} and \emph{current-densities} $\rho$ and $\Jb=(j_1,j_2,j_3)$ are defined by:
\begin{equation}\label{dm2}
\rho:= q|\psi|^2, \quad j_k:= qc \langle \psi , \alpha^k \psi
\rangle_{\C^4}\equiv qc\,\bar\psi \cdot \ap^k\,\psi, \quad \ k=1,2,3,
\end{equation}
where the spinor field $\psi = \psi (t,\xb)=(\psi_1,\psi_2,\psi_3,\psi_4)^T\in \mathbb C^4$ is normalized s.t.
\begin{equation}\label{nor}
\int_{\mathbb R^3} |\psi(t,\xb) |^2 d\xb =1,
\end{equation}
with $t$, $\xb\equiv (x_1,x_2,x_3)$, denoting the time - resp. spatial coordinates.
Further, $V(t,\xb)$ and $V^{ex}(\xb)\in \mathbb R$ are the \emph{self-consistent} resp.
\emph{external electric potential}
and $A_k(t,\xb)\in \mathbb R$, resp. $A^{ex}_k(\xb)\in \mathbb R$, represents the $k$th-components
of the self-consistent, resp. external, \emph{magnetic potential}, \ie $\rm \Ab=(A_1,A_2,A_3)$.
Here and in the following we shall only consider
\emph{static} external fields.
The complex-valued, Hermitian \emph{Dirac matrices}, \ie $\beta, \alpha^k$, are explicitly given by:
\begin{equation}
\beta:=
\begin{pmatrix}
\Id{2} & 0\\
0 & -\Id{2}
\end{pmatrix}
,\quad
\alpha^k:=
\begin{pmatrix}
0 & \sigma ^k\\
\sigma ^k & 0
\end{pmatrix},
\end{equation}
with $\Id{2}$, the $2\times 2$ identity matrix and $\sigma ^k$ the $2\times 2$ \emph{Pauli matrices}, \ie
\begin{equation}
\sigma  ^1:=
\begin{pmatrix}
0 & 1\\
1 & 0
\end{pmatrix}
,\quad
\sigma  ^2:=
\begin{pmatrix}
0 & -i\\
i & 0
\end{pmatrix}
,\quad
\sigma  ^3:=
\begin{pmatrix}
1 & 0\\
0 & -1
\end{pmatrix}.
\end{equation}
Finally, the physical constants, appearing in \eqref{dm0}-\eqref{dm2}, are the normalized
Planck's constant $\hbar$, the speed of light $c$, the permittivity
of the vacuum $\epsilon_0$, the particle mass $m$ and its charge $q$.
\newpar
Additionally to (\ref{dm0}), we impose the \emph{Lorentz gauge condition}
\begin{equation}
\label{lor}
\partial _t V(t,\xb) + c\diverg \Ab(t,\xb)=0,
\end{equation}
for the initial potentials $V^{(0)}(\xb), V^{(1)}(\xb)$, and $\Ab^{(0)}(\xb), \Ab^{(1)}(\xb)$. That means
\[\f{1}{c}V^{(1)}(\xb)+\btd\cdot\Ab^{(0)}(\xb)=0,\quad \btu V^{(0)}(\xb)+\f{q}{4\pi\epsilon_0}|\psi^{(0)}|^2
+\f{1}{c}\btd\cdot\Ab^{(1)}(\xb)=0.
\]
Then the gauge is henceforth conserved during the time-evolution.
This ensures that the corresponding \emph{electromagnetic fields} $\Eb$, $\Bb$ are uniquely determined by
\begin{equation}
\Eb(t,\xb):=- \frac{1}{c} \partial _t \Ab(t,\xb)-\nabla V(t,\xb), \quad \Bb(t,\xb):=\curl \Ab(t,\xb).
\end{equation}
Also it is easily seen that multiplying the Dirac equation with $\overline{\psi}$ implies
the following conservation law
\begin{equation}
\partial_t \rho + \diverg \Jb =0.
\end{equation}
The MD equations are the underlying field equations of relativistic \emph{quantum electro-dynamics}, \cf \cite{Sc},
where one considers the system within the formalism of \emph{second quantization}. Nevertheless,
in order to obtain a deeper understanding for the interaction of matter and radiation, there is a growing
interest in the MD system also  for classical fields, since one can expect at least qualitative results, \cf \cite{EsSe}.
Analytical results concerning local and global well-posedness of \eqref{dm0}-\eqref{dm2}, have been obtained in
\cite{BoRa, Ch, FST, Gro}. Also the rigorous study of asymptotic descriptions for the MD system has been
a field of recent research. In particular the \emph{non-relativistic limit} and the \emph{semi-classical} asymptotic
behavior (in the weakly coupled regime) have been discussed in \cite{BMS2, SpMa1}. For the former case
a numerical study can be found in \cite{BaLi}. Since our numerical simulations shall deal with both asymptotic
regimes, let us discuss now more precisely the corresponding scaling for these physical situations.
\subsection{The MD system in the (weakly coupled) semi-classical regime}
First, we consider the semi-classical or \emph{high-frequency regime} of fast (relativistic) particles,
\ie particles which have a reference speed $v\approx c$. (Of course for particles with mass $m>0$ we always have $v<c$.)
To do so we rewrite the MD system in dimensionless form, such that there remains only one positive real parameter
\begin{equation}
\kappa_0= \frac{4\pi \hbar c \epsilon_0 } {q^{2}}.
\end{equation}
As described in \cite{SpMa1}, we obtain the following rescaled MD system:
\begin{equation}\label{dmd}
\left \{
\begin{aligned}
& i \kappa_0\, \partial_ t  \psi   = -i \kappa_0 \, \alb \cdot \nabla \psi - \alb \cdot (\Ab+\Ab^{ex})\psi
+ (V+V^{ex})\psi + \beta \psi ,\\
& (\partial_{tt}-\Delta) V =  \rho,\\
& (\partial_{tt}-\Delta) \Ab =  \Jb,
\end{aligned}
\right.
\end{equation}
where from now on we shall also use the shorthand notation $\alb \cdot \nabla:=\sum \alpha^k\partial_k$.
Notice that if $|q|=e$, \ie in the case of electrons or positrons where $q$ equals the elementary charge $\pm e$, the
parameter $\kappa_0\approx 137$ is nothing but the reciprocal of the famous \emph{fine structure constant}.
Thus for fast (relativistic) particles which are not too heavily charged, $\kappa_0$ in general is not small and
therefore asymptotic expansions as $\kappa_0 \rightarrow 0$ do not make sense.
In order to describe the semi-classical regime we therefore suppose that the given external electromagnetic potentials
are \emph{slowly varying} w.r.t. the microscopic scales,
\ie $V^{ext}=V^{ext} (\xb \e/\kappa_0)$ and likewise $A^{ext}=A^{ext}(\xb \e/\kappa_0)$,
where from now on $0<\varepsilon\ll 1$ denotes the small \emph{semi-classical parameter}.
Here we fix $\kappa_0$ and include it in the scaling which conveniently eliminates this factor
from the resulting equations. Finally, observing the time-evolution on macroscopic scales we are led to
\begin{equation}
\tilde \xb= \frac{\varepsilon }{\kappa_0} \xb , \quad \tilde t = \frac{\varepsilon}{\kappa_0} t.
\end{equation}
and we set
\begin{equation}
\psi^\e (\tilde t,\tilde \xb)=\left(\frac{\e}{\kappa_0}\right)^{-3/2} \psi\left(\tilde t \frac{\kappa_0}{\e},
\tilde \xb \frac{\kappa_0}{\e}\right)
\equiv \left(\frac{\e}{\kappa_0}\right)^{-3/2} \psi(t,\xb),
\end{equation}
in order to satisfy the normalization condition (\ref{nor}). Plugging this into \eqref{dmd} and
omitting all ``$ \ \tilde { } \ $'' we obtain the following \emph{semi-classically scaled MD system}:
\begin{equation}\label{dmsc}
\left \{
\begin{aligned}
& i \e \, \partial_ t \psi^\e   = -i  \e \, \alb \cdot \nabla \psi^\e - \alb \cdot (\Ab^\e+\Ab^{ex})\psi^\e
+ (V^\e+V^{ex})\psi^\e
+ \beta \psi^\e ,\\
& (\partial_{tt}-\Delta) V^\e =  \e|\psi^\e|^2,\\
& (\partial _{tt}-\Delta) A_k^\e = \e \langle \psi^\e , \alpha^k \psi^\e \rangle_{\C^4}, \quad \ k=1,2,3,
\end{aligned}
\right.
\end{equation}
with $0<\e\ll1$. Note the additional factor $\e$ in the source terms appearing
on the right hand side of the wave equations governing $V^\e$ and $\Ab^\e$,
which implies that we are dealing with a \emph{weak nonlinearity} in the sense of \cite{DoRa, Ra}.
The scaled particle-density in this case is $\rho^\e:=|\psi^\e|^2$ and we also
have $\Jb^\e:=(\langle \psi^\e , \alpha^k \psi^\e \rangle_{\C^4})_{k=1,2,3}$.
\begin{remark} Note that, \emph{equivalently}, we could consider small
asymptotic solutions $\psi^\e\sim\O(\sqrt{\e})$ which again satisfy the semi-classical scaled MD system
\eqref{dmsc} but with source terms of order $\O(1)$ in the wave equations. This point of view is adopted in \cite{SpMa1}.
\end{remark}
\subsection{The MD system in the non-relativistic regime}
We shall also deal with the non-relativistic regime for the MD
system, \ie we consider particles which have a reference speed
$v\ll c$. Introducing a reference length $L$,
time $T$ and writing $v = L/T $, we rescale the time and the spatial
coordinates in \eqref{dm0} by
\begin{equation}\label{sc1}
\tilde \xb = \frac{\xb}{L}, \quad \tilde t = \frac{t}{T}.
\end{equation}
Moreover we set $\widetilde \psi (\tilde t,\tilde \xb)= L^{3/2} \psi(t,\xb)$, such that (\ref{nor}) is satisfied,
and we also rescale the electromagnetic potentials by
\begin{equation}\label{sc2}
\widetilde \Ab^{(ex)}(\tilde t,\widetilde \xb)= \lambda \Ab^{(ex)}(t,\xb), \quad
\widetilde V^{(ex)}(\tilde t,\widetilde \xb)= \lambda V^{(ex)}(t,\xb),
\end{equation}
where $\lambda = q/(4\pi L\e_0)$, \cf \cite{BaLi, BMP}.
In this case we have again two important dimensionless parameters, namely
\begin{equation}
\de = \frac{v}{c}\ll1 , \quad \kappa=\frac{4\pi \hbar v \varepsilon_0 } {q^{2}}.
\end{equation}
Note that for $v\approx c$ we get $\kappa \approx \kappa_0$. Choosing for convenience $v=
q^2/(4\pi\hbar \e_0)$ and $L=q/4\pi \e_0$, we shall from now on denote by $\psi^\delta(\tilde t,\tilde \xb)$
the rescaled wave function $\widetilde \psi (\tilde t,\tilde \xb)$, which is obtained for this particular choice of
$v=L/T$. Then, similarly as before, $\psi^\de$ satisfies
a dimensionless one-parameter model (again omitting all ``$ \ \tilde { } \ $''), given by
\begin{equation}\label{dmnr}
\left \{
\begin{aligned}
& i \partial_ t  \psi^\de   = - \frac{i}{\de} \, \alb \cdot \nabla \psi^\de -
 \alb \cdot (\Ab^\de+\Ab^{ex})\psi^\de +
(V^\de+V^{ex})\psi^\de + \frac{1}{\delta^2}\beta \psi^\de ,\\
& (\delta^2 \partial_{tt}-\Delta) V^\de =  |\psi^\de|^2 ,\\
& (\de^2 \partial_{tt}-\Delta) A_k^\de = \langle \psi^\de , \alpha^k \psi^\de \rangle_{\C^4}, \quad \ k=1,2,3.
\end{aligned}
\right.
\end{equation}
In analogy to the semi-classical case, this system will be called the \emph{non-relativistically scaled MD system}.
In this case the scaled particle density is $\rho^\de := |\psi^\de|^2$, whereas $\Jb^\de:=
\de^{-1} (\langle \psi^\de , \alpha^k \psi^\de \rangle_{\C^4})_{k=1,2,3}$.
Note that in this scaling $\Jb^\de\sim \O(1)$, $\Ab^\de\sim \O(\de)$, (due to a rather complex cancellation
mechanism already known in
the linear case \cf \cite{BMP}) such that
the magnetic field is a relativistic effect which does not appear in the zeroth order
approximation of the MD system, \cf \cite{BMP, BMS2, MaMa} (see also \cite{BMS1} for a similar study).
\newpar
As in the corresponding numerical simulations for semi-classical nonlinear Schr\"odinger equations,
\cf \cite{BJM1}, the main difficulty is  to find an efficient and convenient numerical scheme
with best possible properties in the limiting regimes $\e \to 0$ and $\delta \to 0$, \ie
in particular with uniform convergence properties in $\de$.
\newpar
In the following we present a time-splitting spectral method for the MD system,
and its semi-classical and non-relativistic limiting systems.
The time-splitting spectral methods have been proved to be the best numerical approach to solve linear and
nonlinear Schr\"odinger type systems in the semi-classical regime, \cf \cite{BJM1, BJM2}.
Besides the usual properties of the time-splitting spectral method, such as
the conservation of the Lorentz gauge condition and the unconditionally stability property,
here we shall pay special attention to its performance in both the semi-classical and non-relativistic
regimes. Note that in particular the semi-classical asymptotics has not been studied in \cite{BaLi}.
The method proposed here is similar to the one used for the Zakharov
system in \cite{JMZ}. A distinguished feature of the scheme developed in \cite{JMZ} is that it can be used,
in the sub-sonic regime, with mesh size and time step independent of the
subsonic parameter, a possibility not shared by works before \cite{BS, BSW}.
For the MD system, our time splitting spectral method allows the use
of mesh size and time steps {\it independent} of the relativistic parameter
$\delta$, allowing {\it coarse} grid computations in this asymptotic
regime. This is achieved by the Crank-Nicolson time discretization for
the Maxwell equations, a scheme shown to perform better for wave equations
in the subsonic regime than the exact time integration, as studied in \cite{JMZ}.
For the same reason,  the previously proposed time-splitting spectral
method for the MD system in \cite{BaLi}
does \emph{not} possess this property since it uses the exact time integration
for the Maxwell equations.
\newpar
The paper is now organized as follows: In section \ref{ssc}, we give the
time-splitting spectral method for the MD system and one simple example to show the
reliability, efficiency  and the convergent rate of our method. Our method has spectral convergence
for space discretization and second order convergence for time discretization.
In section \ref{sec:sc} and \ref{sec:nr}, we discuss the time-splitting methods for the
asymptotic systems (the semi-classical regime and non-relativistic regime) and give some examples
for them respectively. We conclude the paper in section 5.

\section{A Time-splitting spectral method for the Maxwell-Dirac system}\label{ssc}
\subsection{A time-splitting method}
Before we describe our time-splitting spectral method, we combine
the rescaled MD system \eqref{dmsc} and \eqref{dmnr}, using two parameters:
\bea\label{dmscnr}
\left \{
\begin{aligned}
& i \e \, \partial_ t \psi   = - \f{i\e}{\dt} \, \alb \cdot \nabla \psi - \alb \cdot (\Ab^\e+\Ab^{ex})\psi
+ (V^\e+V^{ex})\psi
+ \f{1}{\dt^2} \beta\psi ,\\
& (\dt^2\partial_{tt}-\Delta) V =  \e|\psi|^2,\\
& (\dt^2\partial_{tt}-\Delta) A_k = \e \langle \psi , \alpha^k \psi \rangle_{\C^4}, \quad \ k=1,2,3.
\end{aligned}
\right.
\eea
In the following we shall denote by
\begin{align}
\label{do}
\mathcal D_{\Ab}(D)\psi := \alb \cdot(-i\nabla-\Ab^{ex}(\xb))\psi + \beta\psi + V^{ex}(\xb)\psi,
\end{align}
the standard Dirac operator with (external) electromagnetic fields, $D:=-i\nabla$. The corresponding
$4\times4$ matrix-valued symbol is given by
\begin{align}
\mathcal D_{\Ab}(\xi ) = \alb\cdot (\xi - \Ab^{ex}(\xb))+\beta \psi +
V^{ex}(\xb),
\quad
\end{align}
where $\xb,\xi \in \mathbb R^{3}$. Likewise the \emph{free Dirac operator} will be written as
\begin{align}
\label{fdo}
\mathcal D_0(D_x)\psi := -i \alb \cdot \nabla\psi +\beta \psi
\end{align}
Its symbol admits a simple orthogonal decomposition given by
\begin{align}
\label{sym}
\mathcal D_0(\xi ) \equiv \alpha \cdot  \xi  +\beta =  \lambda_0(\xi) \Pi^+_{0} (\xi )-\lambda_0(\xi)\Pi_{0}^-(\xi),
\end{align}
where
\begin{equation}\label{lam}
\lambda_0(\xi):=\sqrt{|\xi|^2+1},
\end{equation}
and
\begin{align}
\label{frpr}
\Pi_{0}^\pm(\xi):=\frac{1}{2}\left(\Id{4}\pm \frac{1}{\lambda_{0}(\xi )}\mathcal D_{0} (\xi )\right).
\end{align}
The time-splitting scheme we propose is then as follows:
\newpar
\textbf{Step 1.} Solve the system
\begin{equation}
\label{st1}
\left \{
\begin{aligned}
&\, i\e \partial _t\psi - \frac{1}{\de^2}\mathcal D_0(\de \e D_x)\psi=0, \\
&\,(\de^2\partial_{tt}-\Delta) V =  \e |\psi|^2,\\
& \,(\de^2\partial_{tt}-\Delta) A_k =   \e \left<\psi, \alpha^k \psi\right>,\quad k=1,2,3,
\end{aligned}
\right.
\end{equation}
on a fixed time-interval $\tg t$, using the spectral decomposition \eqref{sym}.
\newpar
\textbf{Step 2.} Then, in a second step we solve
\begin{equation}
\label{st2}
\begin{aligned}
& \,i\e \partial _t\psi  + \alb\cdot(\Ab+\Ab^{ex})\psi-(V+V^{ex})\psi=0, \\
\end{aligned}
\end{equation}
on the same time-interval, where the solution obtained in step 1 serves as
initial condition for step 2. Also the fields $\Ab$, $V$ are taken from step 1.
It is then easy to see that this
scheme conserves the particle density and the Lorentz gauge.
\subsection{The numerical algorithm}\label{sec:num}
In the following, for the convenience
of computation, we shall deal with the system
(\ref{dmscnr}) on a bounded domain, for example, on
the cubic domain
\bea
\Og=\{\xb=(x_1,x_2,x_3)\ |\ a_j\le x_j \le b_j, j=1,2,3\},
\eea
imposing \emph{periodic boundary conditions}. We choose the time step $\tg t=T/M$ and spatial mesh
size $\tg x_j=(b_j-a_j)/N_j$, $j=1,2,3$, in $x_j$-direction, with given $M, N_j\in \N$ and
$[0,T]$ denoting the computational time interval. Further we denote the time grid points by
\bea
t_n=n\triangle t,\quad t_{n+1/2}=\left(n+\f{1}{2}\right)\triangle t, \quad t=0,1,\dots,M
\eea
and the spatial grid points by
\bea
\mathbf{x}_\mathbf{m}=(x_{1,m_1},\ x_{2,m_2},\ x_{3,m_3}), \quad
\mbox{where $x_{j,m_j}:=  \, a_j+m_j \tg x_j$}, \quad j=1,2,3,
\eea
and $\mathbf{m}=(m_1,m_2,m_3)\in \mathcal{M}$, with
\begin{align}
\mathcal{M}=  \, \left\{(m_1,m_2,m_3)\,\Big|\ 0\le m_j\le N_j, \ j=1,2,3\right\}.
\end{align}
In the following let $\Psi^n_{\mathbf{m}}$, $V^n_{\mathbf{m}}$,
and $\mathbf{A}^n_{\mathbf{m}}$ be the numerical approximations of
$\psi(t_n,\mathbf{x}_{\mathbf{m}})$,
$V(t_n,\mathbf{x}_{\mathbf{m}})$, and
$\mathbf{A}(t_n,\mathbf{x}_{\mathbf{m}})$, respectively. Suppose
that we are given $\Psi^n$, $V^n$, and $\mathbf{A}^n$, then we obtain
$\Psi^{n+1}$, $V^{n+1}$ and $\mathbf{A}^{n+1}$ as follows:
\newpar

\textbf{Step 1.} For the first step we denote the value of $\Psi$ at time $t$ by
$\Phi(t)$. Then we approximate the spatial derivative in (\ref{st1}) by the spectral differential
operator. More precisely we first take a \emph{discrete Fourier transform}
(DFT) of (\ref{st1}):
\bea\label{eq:numst11}
\left \{
\begin{aligned}
& \, \p_t\hat{\Phi}=  -\f{i}{\e \dt^2}
\left(\e \dt{\alb}\cdot\xi+\beta\right)\hat\Phi\equiv \mathbb M_1\hat\Phi,\\
& \,(\de^2\p_{tt}+|\xi|^2)\hat{V} = \e \widehat{\ |\Phi|^2},\\
& \,(\de^2 \p_{tt}+|\xi|^2)\hat A_k= \e
\widehat{\left<\Phi,\ap^k\Phi\right>},\quad \mbox{for }k=1,2,3,
\end{aligned}
\right.
\eea
where $\hat f$ is the DFT of function $f$.
As the matrix $\mathbb M_1\in \C^{4\times4}$ is diagonalizable, \ie
there exists a Hermitian matrix $D_1$ such that
\begin{equation}
\bar D_1^{T}\mathbb M_1 D_1 = \mbox{diag}\,[\ld,\ld,-\ld,-\ld] \equiv \Lambda
\end{equation}
is a purely imaginary diagonal matrix with entries
\begin{equation}
\lambda =\frac{i}{\e \dt^2} \sqrt{1+\e^2\dt^2|\xi|^2}.
\end{equation}
Then the value of $\hat\Phi$ at
time $t_{n+1}$ is given by
\bea \hat\Phi^{n+1}&=& D_1\,
\exp\left(\Ld \tg t\right)\, \bar  D_1^{T}\hat\Psi^n \nn\\
&=&\begin{pmatrix}
c_\lambda-i s_\lambda & 0 & -i\e\dt s_\lambda \xi_3  & -\e\dt s_\lambda(\xi_2+i\xi_1)\\
0 & c_\lambda-is_\lambda & \e\dt s_\lambda(\xi_2-i\xi_1) &i\e\dt s_\lambda\xi_3\\
-i\e\dt s_\lambda\xi_3  & -\e\dt s(\xi_2+i\xi_1)& c_\lambda+is_\lambda & 0  \\
\e\dt s_\lambda(\xi_2-i\xi_1) &i\e\dt s_\lambda\xi_3 &0 & c_\lambda+is_\lambda\\
\end{pmatrix}\hat\Psi^n, \label{eq:numst14}\eea
where
\begin{equation}\label{cos}
c_\lambda :=\cos (-i \lambda \tg t), \quad
s_\lambda:=\sin (-i \lambda \tg t)(1+|\e\dt\xi|^2)^{-1/2}.
\end{equation}
Then we obtain the value of $\Phi^{n+1}$ by an \emph{inverse discrete Fourier
transform} (IDFT). Hence from (\ref{eq:numst11}),
we can find the values of $\hat V$ and $\hat{\mathbf A}$ by the Crank-Nicolson scheme, \ie
\begin{align}\label{eq:numst15}
\left(1+\f{\tg t^2|\xi|^2}{4\dt^2}\right)\left(\ba{c}\hat V ^{n+1}\\
\partial_t \hat V^{n+1}\ea\right) =\,
\begin{pmatrix}
1-\f{\tg t^2|\xi|^2}{4\dt^2} & \tg t\\
-\f{\tg t|\xi|^2}{\dt^2} &1-\f{\tg t^2|\xi|^2}{4\dt^2}
\end{pmatrix}\left(\ba{c}\hat V ^{n}\\
\partial_t\hat V^{n}\ea \right)+\e\left(\ba{c}\f{\tg t^2}{4\dt^2} \\ \f{\tg t}{2\dt^2}\ea \right)
\left(\hat \rho^{n} +\hat \rho^{n+1}\right)
\end{align}
and
\begin{align}
\left(1+\f{\tg t^2|\xi|^2}{4\dt^2}\right)\left(\ba{c}\hat {\mathbf A}^{n+1} \\
\partial_t \hat \Ab^{n+1}\ea\!\!\!\right) =\,
\begin{pmatrix}
1-\f{\tg t^2|\xi|^2}{4\dt^2} & \tg t\\
-\f{\tg t|\xi|^2}{\dt^2} &1-\f{\tg t^2|\xi|^2}{4\dt^2}
\end{pmatrix}
\left(\!\!\ba{c}\hat \Ab ^{n}\\
\partial_t \hat \Ab^{n}\ea\right)+\e\de\left(\!\!\ba{c}\f{\tg t^2}{4\dt^2} \\ \f{\tg t}{2\dt^2}\ea \right)
\left( \hat \Jb^n+ \hat \Jb^{n+1}\right), \label{eq:numst16}
\end{align}
where for $k=1,2,3$, we denote
\begin{align}
\rho^n=|\Psi^n|^2,\ \rho^{n+1}=|\Phi^{n+1}|^2, \
\Jb_k^n=\de^{-1} \left<\Psi^n,\ap^k\Psi^n\right>,\ \Jb_k^{n+1}=\de^{-1}\left<\Phi^{n+1},
\ap^k\Phi^{n+1}\right>.
\end{align}
Performing an IDFT of $\hat V^{n+1}$ and $\hat{\mathbf A}^{n+1}$,
we finally obtain $V^{n+1}$ and $\mathbf{A}^{n+1}$.
\newpar
\textbf{Step 2.} Since $V$ and $A_k$ do not change in Step 2, we only have to update
$\Psi$. First we shall rewrite the equation (\ref{st2}) in the following form:
\be\label{eq:numst21}
\p_t\Psi = \f{i}{\e} \, \alb\cdot(\Ab+\Ab^{ex})\Psi-(V+V^{ex})
\Psi \equiv\mathbb M_2\Psi.
\ee
Then there exists again a Hermitian matrix $ D_2$ such that
\begin{equation}
\bar D_2^{T}\mathbb M_2  D_2= \mbox{diag} \, [\mu_1,\mu_1,\mu_2,\mu_2]\equiv
\Theta,
\end{equation}
where $\Theta$ is a purely imaginary diagonal matrix with
\begin{equation}
\mu_1=-\frac{i}{\e}\left((V+V^{ex})-|\mathbf{A}+\mathbf{A}^{ex}|\right),\quad
\mu_2=-\frac{i}{\e}\left((V+V^{ex})+|\mathbf{A}+\mathbf{A}^{ex}|\right).
\end{equation}
Hence, the value of $\Psi$ at time $t_{n+1}$ is given by
\bea\Psi^{n+1}&=&D_2\, \exp\left(\Theta \tg t\right)\, \bar
{D}_2^{T}\Phi^{n+1}\nn\\
&=&\begin{pmatrix}
\f{c_1+c_2-i(s_1+s_2)}{2} & 0 & \f{A_3}{|\Ab|}(c_0-is_0) &\f{A_1-iA_2}{|\Ab|}(c_0-is_0)\\
0 & \f{c_1+c_2-i(s_1+s_2)}{2} & \f{A_1+iA_2}{|\Ab|}(c_0-is_0) &-\f{A_3}{|\Ab|}(c_0-is_0)\\
 \f{A_3}{|\Ab|}(c_0-is_0) &\f{A_1-iA_2}{|\Ab|}(c_0-is_0)& \f{c_1+c_2-i(s_1+s_2)}{2} & 0  \\
\f{A_1+iA_2}{|\Ab|}(c_0-is_0) &-\f{A_3}{|\Ab|}(c_0-is_0) &0 & \f{c_1+c_2-i(s_1+s_2)}{2}\\
\end{pmatrix}\Phi^{n+1}, \label{eq:numst22}\eea
where we use a notation analogous to \eqref{cos} and write
\begin{equation}
\exp(\mu_1\tg t)\equiv c_1-is_1, \ \exp(\mu_2\tg t)\equiv c_2-is_2, \ c_0:=c_1-c_2, \ s_0:=s_1-s_2.
\end{equation}
Clearly, the algorithm given above is first order in time.
We can get a second order scheme by the Strang splitting method, which means that
we use Step 1 with time-step $\tg t/2$, then Step 2 with time-step
$\tg t$, and finally integrate Step 1 again with $\tg t/2$.
Our algorithm given above is an `explicit' and unconditional stable scheme.
The main costs are DFT and IDFT.

\begin{figure} 
\begin{center}
\resizebox{2.4in}{!}{\includegraphics{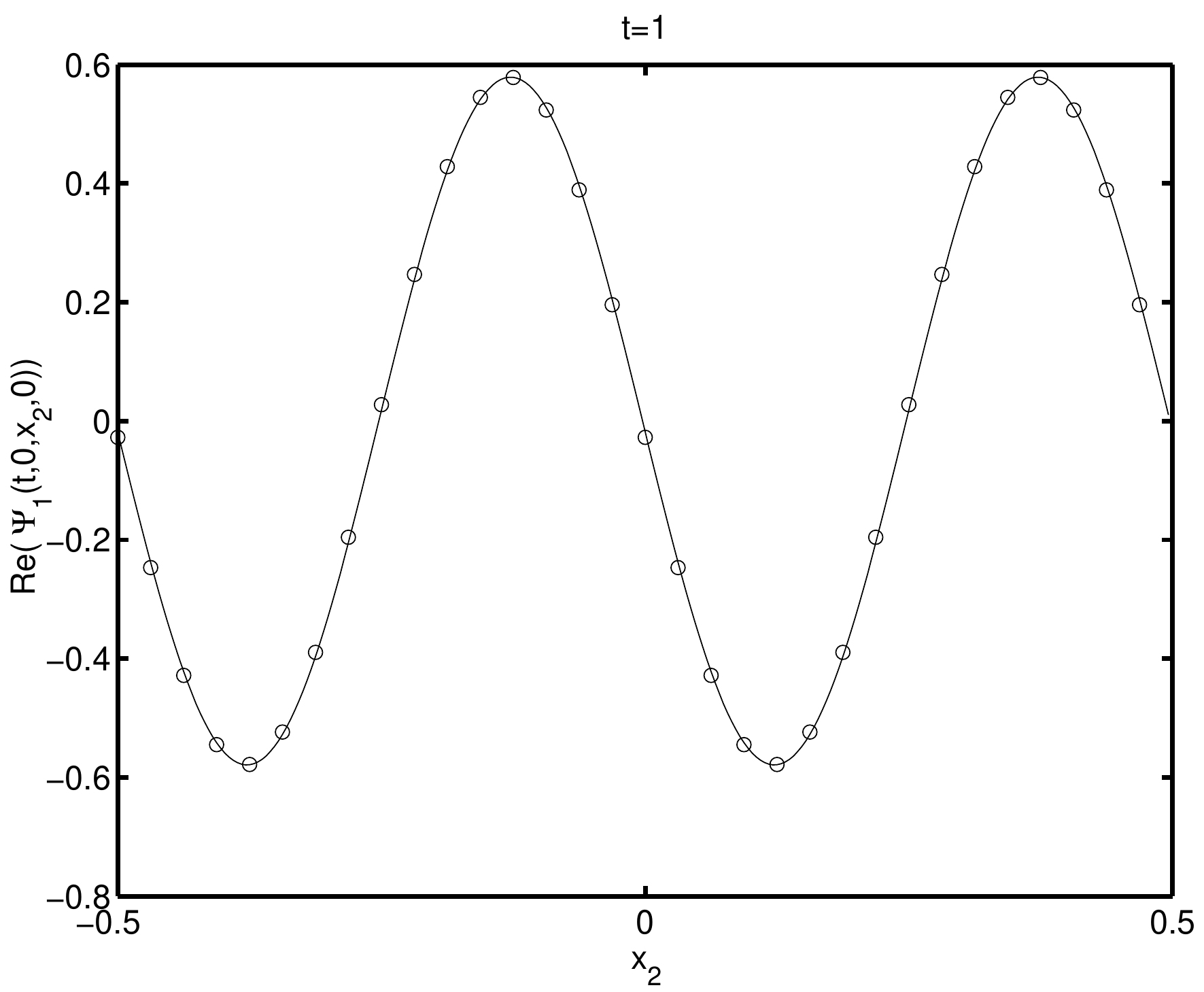}}
\resizebox{2.4in}{!}{\includegraphics{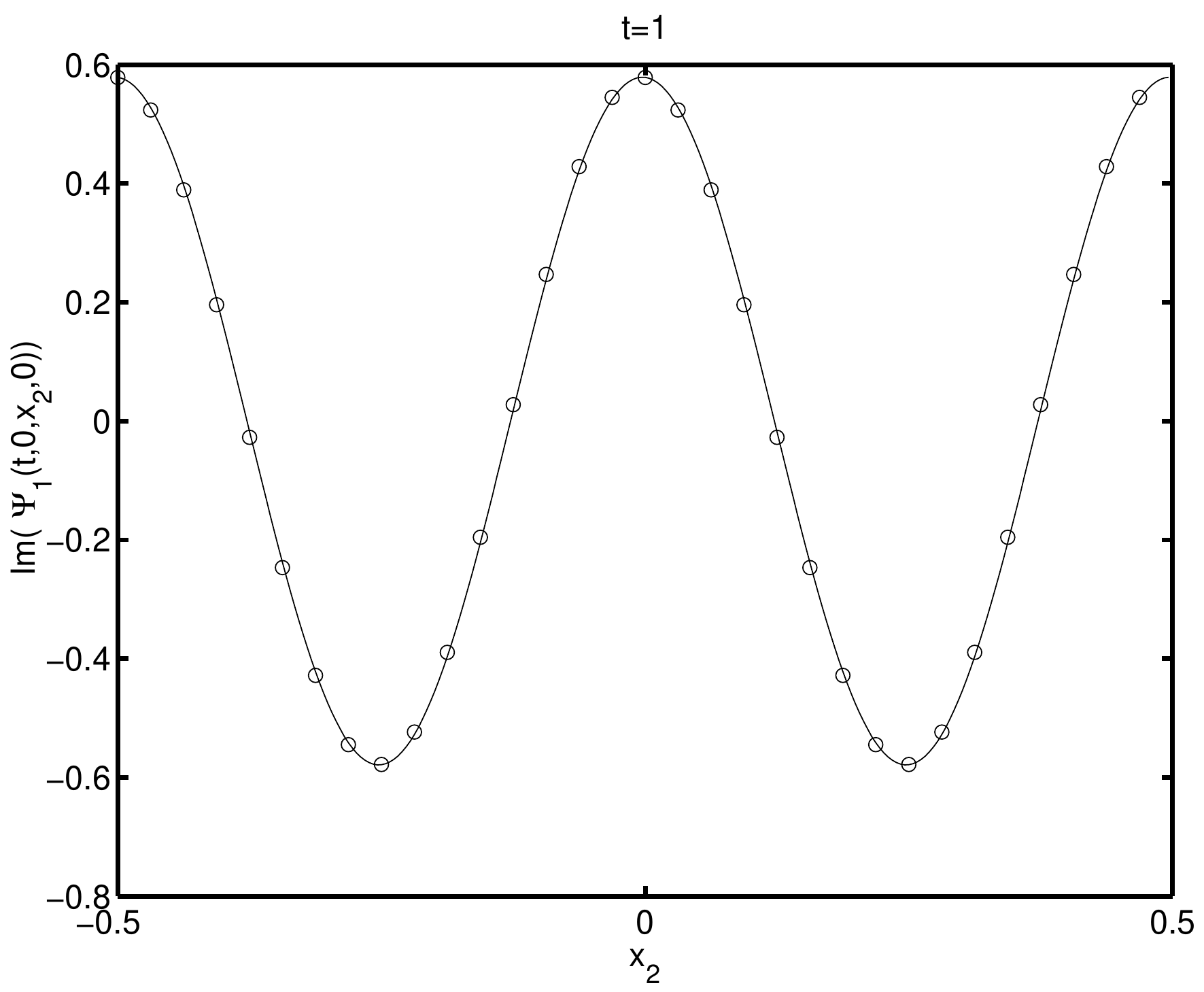}}
\resizebox{2.4in}{!}{\includegraphics{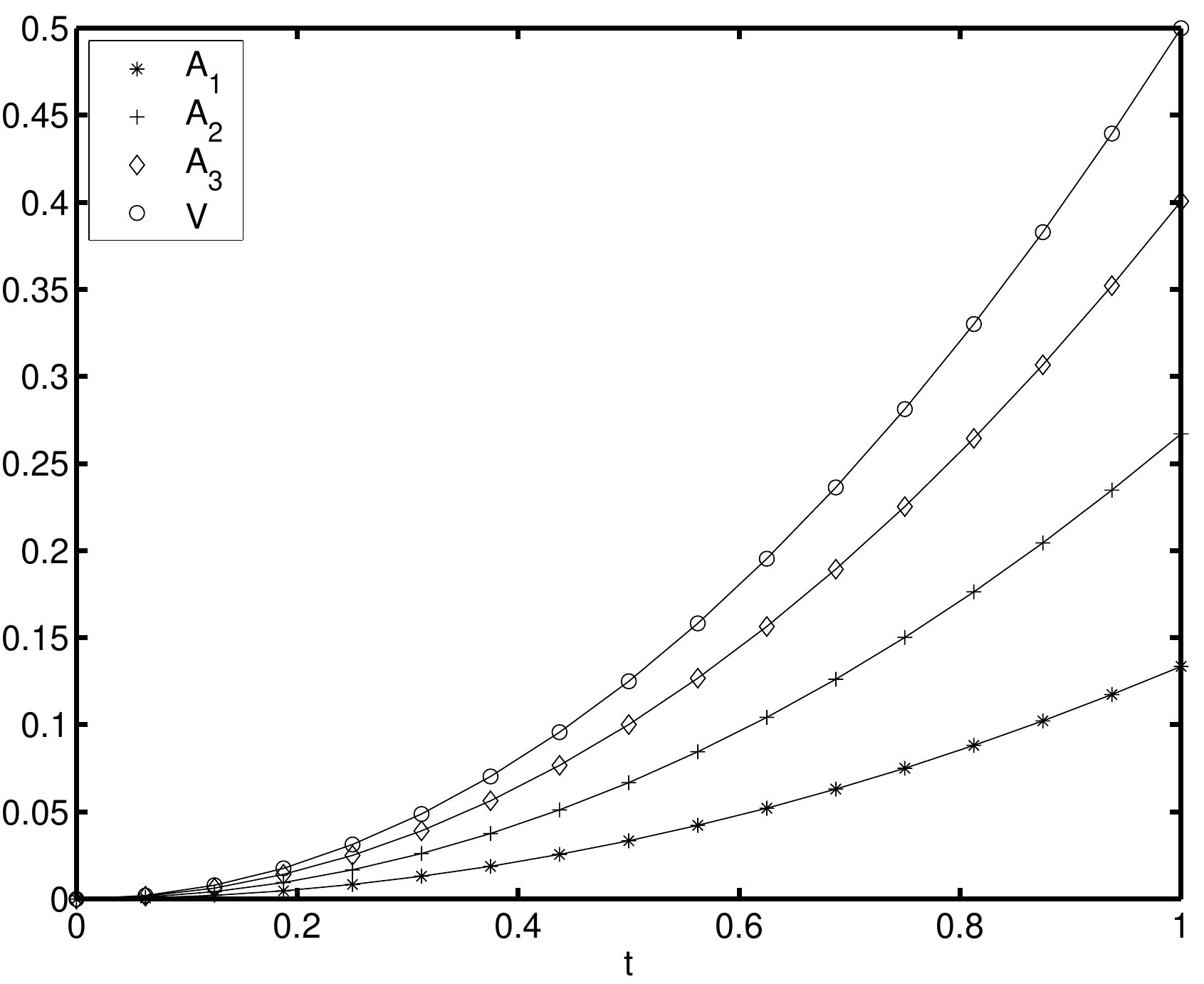}}
\resizebox{2.4in}{!}{\includegraphics{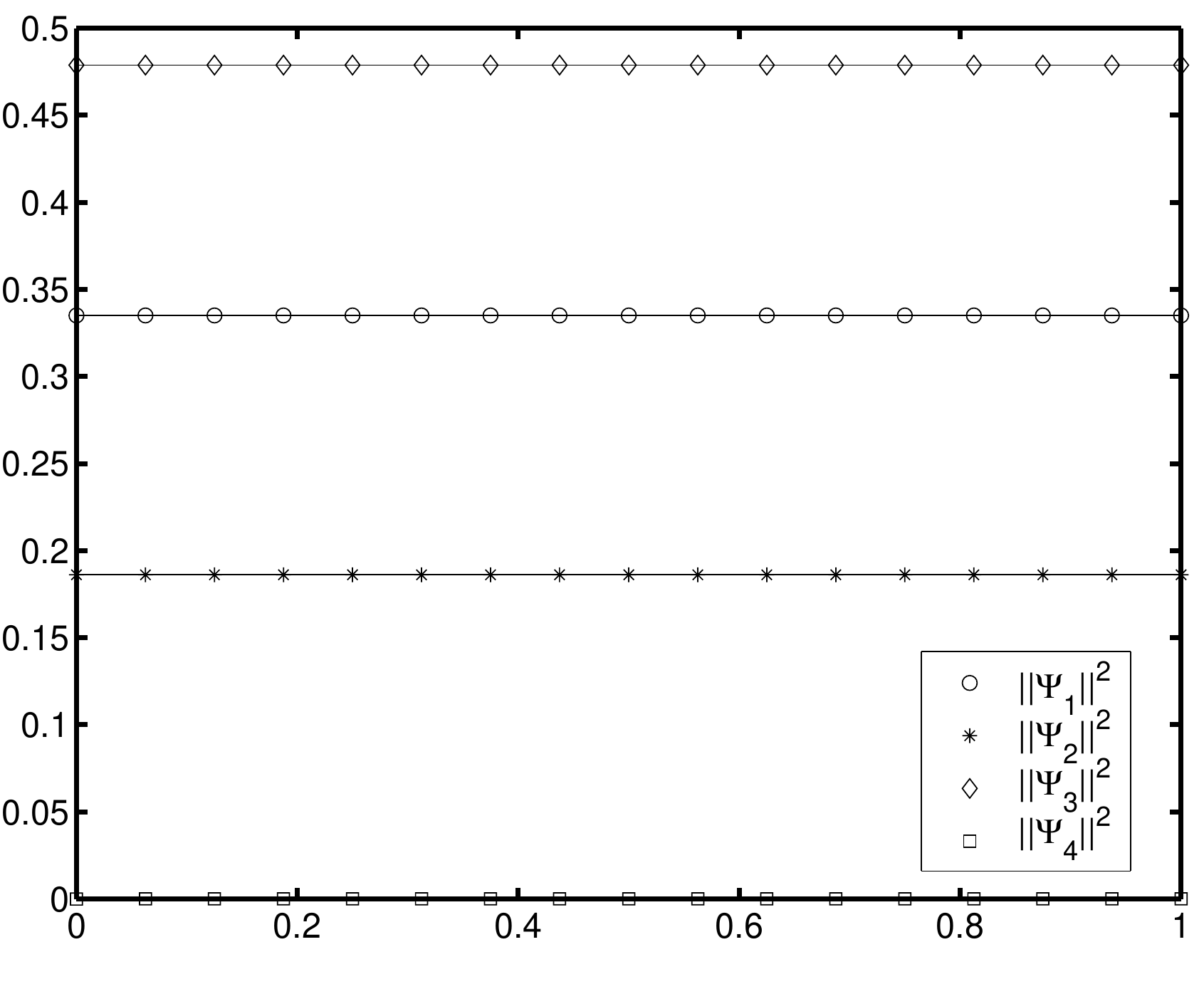}}
\end{center}
\caption{The numerical solutions of example \ref{examFDM}.
Here $\tg t=\f{1}{128}$, $\tg x=\f{1}{32}$.
The top two graphs are the real and imaginary parts of $\psi^\e_1(t,0,x_2,0)|_{t=1}$.
The bottom two graphs are electromagnetic potentials and the norms of $\psi^\e_k$, $k=1,2,3,4$.
The solid lines are the exact solutions, `ooo', `***', `$\diamond\diamond\diamond$'
and `$\square\square\square$' are numerical solutions.} \label{fig11}
\end{figure}

\begin{lemma}\label{lm2}
{Our numerical scheme  conserves the particle density in the discrete
$l^2$ norm (discrete total charge) and the Lorentz gauge.}
\end{lemma}
{\bf Proof: }
From \eqref{eq:numst14} and \eqref{eq:numst22},
it is easy to check that the discrete total charge is conserved. From the initial conditions and \eqref{eq:numst15},
we have
\[\dt\partial_t\hat V^{0}+i\xi\cdot\hat \Ab ^{0}=0,\quad
\vep\hat\rho^0=|\xi|^2\hat V^0-i\dt\xi\cdot\p_t\hat\Ab^0,\quad
\hat \rho^{1} =\hat \rho^{0}-\f{i\tg t}{2\dt}\xi\cdot\left( \hat \Jb^0+ \hat \Jb^{1}\right).\]
From \eqref{eq:numst15} and \eqref{eq:numst16}, we obtain
\beas\left(1+\f{\tg t^2|\xi|^2}{4\dt^2}\right)\left(
\dt\partial_t \hat V^{n+1}+i\xi\cdot\hat \Ab^{n+1}\right) &=&
\left(1-\f{\tg t^2|\xi|^2}{4\dt^2}\right)\left(\dt\partial_t\hat V^{n}+i\xi\cdot\hat \Ab ^{n}\right)\\
&+&\f{\tg t}{\dt}\left(-|\xi|^2\hat V ^{n}+i\dt\xi\cdot\p_t\hat \Ab^{n}\right)\\
&+&\f{\vep\tg t}{2\dt}
\left(\hat \rho^{n} +\hat \rho^{n+1}+\f{i\tg t}{2\dt}\xi\cdot\left(\hat \Jb^n+ \hat \Jb^{n+1}\right)\right).
\eeas
Then it is clear that for all $n$, we have
\[\dt\partial_t\hat V^{n}+i\xi\cdot\hat \Ab ^{n}=0,\quad
\vep\hat\rho^n=|\xi|^2\hat V^n-i\dt\xi\cdot\p_t\hat\Ab^n,\quad
\hat \rho^{n+1} =\hat \rho^{n}-\f{i\tg t}{2\dt}\xi\cdot\left( \hat \Jb^n+ \hat \Jb^{n+1}\right).\quad \Box\]
In order to test the numerical scheme we consider the example of an exact solution for the
full MD system, \cf\cite{DaKa}. In all of the following examples, we take the computational domain
$\Og$ to be the unit cubic $[-0.5,0.5]^3$.
\setcounter{lemma}{0}
\begin{example}[Exact solution for the MD system]\label{examFDM}
Let us consider the MD system for $\e=\dt=1$ with initial data
\begin{equation}\label{eq:e05}
\left \{
\begin{aligned}
& \psi^{(0)}(\xb)=
\dpm\f{\exp(i\xi\cdot \xb)}{\sqrt{2(1+|\xi|^2-\sqrt{1+|\xi|^2}) }}\ \chi,
\quad \chi= (\xi_3,\xi_1+i\xi_2,\sqrt{1+|\xi|^2}-1,0),
\\
& V^{(0)}(\xb)= V^{(1)}(\xb) =0, \quad \Ab^{(0)}(\xb)=\Ab^{(1)}
(\xb)=0,
\end{aligned}
\right.
\end{equation}
and external fields given by
\begin{equation}
V^{ex}=\dpm-\f{ t^2}{2},\qquad \Ab^{ex}=-\f{ t^2\xi}{2\sqrt{1+|\xi|^2}},\quad \xi=(2\pi,4\pi,6\pi)\in \R^3.
\end{equation}
In this case, there is an exact plane wave solutions for the MD system in the following form, \cf \cite{DaKa}:
\begin{equation}\label{eq:exact}
\left \{
\begin{aligned}
\psi(t,\xb)&=\dpm\f{\exp\big(i(\xi\cdot \xb-t\sqrt{1+|\xi|^2})\big)}{\sqrt{2(1+|\xi|^2-\sqrt{1+|\xi|^2}) }}\ \chi,\\
V(t,\xb)&=\dpm\f{t^2}{2},\quad \Ab(t,\xb)=\f{t^2\xi}{2\sqrt{1+|\xi|^2}}.
\end{aligned}
\right.
\end{equation}
In Figure \ref{fig11}, we see that our method gives a very good agreement with the exact result.
\end{example}
To test the accuracy of our time-splitting method for the MD system,
we did the spatial and temporal discretization error tests (see Table \ref{tb1} and \ref{tb2}).
Table \ref{tb1} shows  the {\it spectral convergence} for spatial discretization.
Table \ref{tb2} shows the convergence rate for temporal discretization is about 2.0.
Here $\psi^{\tg x,\tg t}(t,\cdot)$ is the numerical solution for mesh size $\tg x$ and time step $\tg t$,
and $\psi(t,\cdot)$ is the exact solution given by \eqref{eq:exact}.
In the following also show the charge conservation test (see Table \ref{tb3}):
\begin{table}[ht]
\begin{center}
\caption{Spatial discretization error test: at time t=0.25 under $\tg t=1/1024$ ($\e=\dt=1$).}\label{tb1}
\begin{tabular}{c|cccc}\hline
mesh size & $\tg x=1/4$ &   $\tg x=1/8$ & $\tg x=1/16$ & $\tg x=1/32$\\ \hline 
$\ba{c}\vspace*{-4mm} \\
\dpm\f{{\left\|\psi^{\tg x,\tg t}(t,\cdot)-\psi(t,\cdot)\right\|}_{l^2}}{{\|\psi(t,\cdot)\|}_{l^2}}\\ \vspace{-4mm} \ea$
 & 8.40E-2 &  2.68E-3 & 6.95E-5 &  5.00E-8\\ \hline
convergence order & & 4.9& 5.3 & 10.4 \\ \hline 
\end{tabular}
\end{center}
\end{table}
\begin{table}[ht]
\begin{center}
\caption{Temporal discretization error test: at time t=0.25 under $\tg x=1/32$ ($\e=\dt=1$).}\label{tb2}
\begin{tabular}{c|cccc}\hline
time step &          $\tg t=\f{1}{16}$ &   $\tg t=\f{1}{32}$ & $\tg t=\f{1}{64}$ & $\tg t=\f{1}{128}$ \\ \hline 
$\ba{c}\vspace*{-4mm} \\
\dpm\f{{\left\|\psi^{\tg x,\tg t}(t,\cdot)-\psi(t,\cdot)\right\|}_{l^2}}{{\|\psi(t,\cdot)\|}_{l^2}}\\ \vspace{-4mm} \ea$
 &    2.59E-4 &  5.14E-5  & 1.29E-5 &  3.21E-6
\\ \hline
convergence order & & 2.3& 2.0 & 2.0 \\ \hline
\end{tabular}
\end{center}
\end{table}
\begin{table}[ht]
\begin{center}
\caption{Charge conservation test: under $\tg x=1/32$, $\tg t=1/1024$ ($\e=\dt=1$).}\label{tb3}
\begin{tabular}{cccccc}\hline
\hspace*{5mm} &time &          t=0 &   t=0.5 & t=1.0 &\hspace*{5mm}\\ \hline 
&${\|\psi^{\tg x,\tg t}(t,\cdot)\|_{l^2}}$
&    1.00000000 &  0.99999998  & 0.99999997 &
\\ \hline 
\end{tabular}
\end{center}
\end{table}
\section{The semi-classical regime}\label{sec:sc}
We shall consider in the following the semi-classically scaled MD system \eqref{dmsc}. First we shall discuss the
(formal) asymptotic description as $\e\rightarrow 0$ and then consider some particular numerical test cases.
\subsection{Formal asymptotic description}
To describe the limiting behavior of $\psi^\e$ as $\e\rightarrow 0$ we introduce the following notations:
\newpar
Analogously to the free Dirac operator, the matrix-valued symbol $\mathcal D_{\Ab}(\xi)$ can be
(orthogonally) decomposed into
\begin{align}
\label{dec}
\mathcal D_{\Ab}(\xi ) = h_{\Ab}^+ (\xi) \Pi^+_{\Ab} (\xi ) +h_{\Ab}^-(\xi)\Pi_{\Ab}^-(\xi ),\quad \xi \in\mathbb R^3,
\end{align}
where
\begin{align}
\label{ham}
h^\pm_{\Ab}(\xi ):=  \pm \lambda_{\Ab}(\xi)+V(\xb),
\end{align}
with
\begin{align}
\lambda_{\Ab}(\xi):= \sqrt{{1+|\xi-\Ab^{ex}(\xb)|}^2}+V^{ex}(\xb).
\end{align}
The corresponding (orthogonal) projectors $\Pi_{\Ab}^\pm(\xi)$ are then given by
\begin{align}
\label{pro}
\Pi_{\Ab}^\pm(\xi):=\frac{1}{2}\left(\Id{4}\pm \frac{1}{\lambda_{\Ab}(\xi )}\
\left(\mathcal D_{\Ab} (\xi )-V^{ex}(\xb)\Id{4}\right)\right).
\end{align}
Clearly, we obtain the corresponding decomposition of the free Dirac operator \eqref{sym}, \eqref{frpr},
by setting $\Ab^{ex}(\xb)=0$ and $V^{ex}(\xb)=0$ in the above formulas.
Note that $h^{\pm}_{\Ab}(\xi)$ is nothing but the \emph{classical relativistic Hamiltonian}
(corresponding to positive resp. negative energies) for a particle with momentum $\xi$. These
particles can be interpreted as \emph{positrons and electrons}, resp., at least in the limit $\e\rightarrow 0$,
as we shall see below. Finally, we also define the relativistic \emph{group-velocity} by
\begin{equation}
\label{vel}
\omega_{\Ab} ^\pm(\xi):=\nabla _\xi  h^\pm_{\Ab} (\xi ).
\end{equation}
The group velocity for \emph{free} relativistic particles is then $\omega_{0} ^\pm(\xi)=\xi/\lambda_0(\xi)$.
\newpar
The semi-classical limit for solution of the weakly nonlinear MD system \eqref{dmsc}
can now be described by means of \emph{WKB-techniques} as given in \cite{SpMa1} (see also \cite{SpMa2}) .
To do so we assume (well prepared) \emph{highly oscillatory initial data} for $\psi^\e$, \ie
\begin{equation}
\psi^{(0)}(\xb) \sim u^+_I(\xb)e^{i \phi^+_I(\xb)/\e}+
u^-_I(\xb)e^{i \phi^-_I(\xb)/\e}+\O(\e).
\end{equation}
We then expect that $\psi^\e(t,\xb)$ can be described
in leading order (as $\e\rightarrow 0$) by a WKB-approximation of the following form
\begin{align}
\label{u0}
\psi^\e(t,\xb)\sim  u^+(t,\xb)e^{i \phi^+ (t,\xb)/\varepsilon }+u^-(t,\xb)e^{i \phi^- (t,\xb)/\varepsilon }+\O(\e).
\end{align}
Here, the \emph{phase functions} $\phi^\pm (t,\xb)\in \R$, resp. satisfy the
\emph{electronic} or \emph{positronic eiconal equation}
\begin{equation}
\label{eic}
\partial_t \phi^\pm (t,\xb)+ h^\pm_{\Ab}(\nabla \phi^\pm(t,\xb))=0, \quad \phi^\pm (0,\xb)=\phi^\pm_I (\xb).
\end{equation}
As usual in WKB-analysis we can expect an approximation of the form \eqref{u0} to be valid only locally in time,
\ie for $|t|<t_c$, where $t_c$ denotes the time at which the first \emph{caustic} appears in the solution of \eqref{eic}.
\begin{remark}
We want to stress that the self-consistent fields $\Ab^\e$, $V^\e$ do \emph{not} enter in \eqref{eic}, \ie
the eiconal equation is found to be the same as in the linear case.
This is due to the weakly nonlinear scaling described in the introduction. In particular,
\ie for the Dirac equation without Maxwell coupling,
this setting allows us to compute the rays of geometrical optics, \ie the characteristics for \eqref{eic},
\emph{independently} of $\Ab^\e$, $V^\e$.
\end{remark}
It is shown in \cite{SpMa1}, for the simplified case where $\Ab^{ex}(\xb)=V^{ex}(\xb)=0$, that the
\emph{principal-amplitudes} $u^\pm(t,\xb)\in \C^4$ solve a \emph{nonlinear} first order system, given by
\begin{equation}\label{trans}
\left \{
\begin{aligned}
& \, \left( \partial _t  + (\omega_0^+ (\nabla \phi^+ ) \cdot \nabla) \right) u^+
+ \frac{1}{2}\diverg (\omega_{0}^+(\nabla \phi^+))u^+
= i\mathcal N^+[u]\,u^+,\\
& \, \left( \partial _t  + (\omega_0^- (\nabla \phi^- ) \cdot \nabla) \right) u^-
+ \frac{1}{2}\diverg (\omega_{0}^-(\nabla \phi^- ))u^-
= i\mathcal N^-[u]\, u^-,\\
\end{aligned}
\right.
\end{equation}
with initial condition
\begin{equation}
u^\pm(0,\xb):=  \Pi_0^\pm(\nabla \phi_I^\pm) u_I(\xb).
\end{equation}
The nonlinearity on the r.h.s. of \eqref{trans} is given by
\begin{equation}\label{nonl}
\mathcal N^\pm[u] := \mathcal A \cdot  \omega_0^{\pm} (\nabla \phi^\pm) - \mathcal V,
\end{equation}
where the fields $\mathcal V$, $\mathcal A$ are computed self-consistently through
\begin{align}\label{zi}
(\partial_{tt}-\Delta) \mathcal V = \rho^0  , \quad  (\partial_{tt}-\Delta) \mathcal A  = \Jb^0.
\end{align}
with source terms
\begin{align}
\rho^0:=|u^+|^2+|u^-|^2, \quad \Jb^0:=\omega_0^+(\nabla \phi^+)|u^+|^2  \ + \omega_0^- (\nabla \phi^-) {|u^-|}^2.
\end{align}
The polarization of $u^\pm$ is henceforth preserved, \ie
\begin{align}
u^\pm(t,\xb) =\Pi_0^\pm(\nabla \phi^\pm)u^\pm(t,\xb), \quad \mbox{for all $|t|<t_c$,}
\end{align}
and we  call $u^+$ the \emph{(semi-classical) electronic amplitude} and
$u^-$ the \emph{(semi-classical) positronic amplitude}.
Note that in this case, \ie without external fields, we have the simplified
relation
\begin{align}
\phi^+(t,\xb) = - \phi^-(t,\xb),
\end{align}
if this holds initially, which we will henceforth assume.
The fact that \eqref{trans} conserves the polarization of $u^\pm$, is crucial. It allows us to justify the
interpretation in terms of electrons and positrons. In other words, the WKB-analysis given above shows that
the energy-subspaces, defined via \eqref{pro}, remain \emph{almost invariant} in time,
\ie up to error terms of order $\mathcal O(\e)$. This, so called, \emph{adiabatic decoupling phenomena}
is already known from the linear semi-classical scaled Dirac equation \cite{BoKe, Te, Sp}.
However we want to stress the fact that in our non-linear setting rigorous proofs so far are
only valid locally in time \cite{SpMa1}. More precisely, it holds
\begin{align}\label{eq:scerr}
\sup_{0\leq |t|< t_c-\tau}\left\|\, \psi^\e(t) -
\sum_{\pm}u^\pm(t)e^{i \phi^\pm(t)/\varepsilon } \,
\right\|_{L^2(\R^3)\otimes \C^4}&=\O(\e),\quad \mbox{for every $0<\tau<t_c$.}
\end{align}
On the other hand we want to remark that in the case of the linear Dirac equation,
global-in-time results are available which also confirm the
adiabatic decoupling for all $t\in \R$, \cf \cite{Te, Sp}.
\newpar
Note that the nonlinearity in \eqref{trans} is purely imaginary. Hence for the densities
$\rho^\pm:=|u^\pm|^2$ we find
\begin{equation}
\label{conl}
\partial_t \rho^\pm + \diverg\left(\omega_0^\pm (\nabla \phi^\pm )\rho^\pm\right)=0,
\end{equation}
which clearly implies the important property of charge-conservation:
\begin{equation}\label{cons}
\int_{\mathbb R^3} \left(\rho^+ (t,\xb) +\rho^-(t,\xb)\right) d\xb =  \mbox {const.}
\end{equation}
In the case of non-vanishing external fields, \ie $\Ab^{ex}(\xb)\not=0$, $V^{ex}(\xb)\not=0$ the system
\eqref{trans} becomes much more complicated. First $\omega_0^\pm$ has to be replaced by $\omega_{\Ab}^\pm$
in the above given formulas and second, an additional matrix-valued potential
has to be added, the, so called, \emph{spin-transport term}, \cf \cite{BoKe, Te, Sp},
which mixes the components of each $4$-vector $u^\pm$ (\cf \cite{Te} for a broad discussion on this).
We shall not go into further details here since in
our (semi-classical) numerical examples below we shall always assume $\Ab^{ex}(\xb)=0$
and $V^{ex}(\xb)=0$, since we are mainly interested in studying the influence of the self-consistent fields.
The only exception is Example \ref{exha} below, where we treat the harmonic oscillator case with
$V^{ex}(\xb)=|\xb|^2$.
\begin{remark}
Strictly speaking, the results obtained in \cite{SpMa1} do not include the
most general case of non-vanishing external fields \emph{and} mixed initial data, \ie $u^\pm(0,\xb)\not=0$.
Rather, the given results only hold in one of the following two (simplified) cases: Either
$\Ab^{ex}(\xb)=V^{ex}(\xb)=0$ and $u^\pm(0,\xb)\not=0$, or: $\Ab^{ex}(\xb)\not=0$, $V^{ex}(\xb)\not=0$,
but then one needs to assume $u^+(0,\xb)=0$, or $u^-(0,\xb)=0$, respectively.
The reason for this is that the analysis given in \cite{SpMa1}
heavily relies on a \emph{one-phase} WKB-ansatz, which is needed (already on a formal level) to
control the additional oscillations induced for example through the, so called, \emph{Zitterbewegung} \cite{Sc} of
$\Jb^\e$, \cf \cite{SpMa1}, \cite{SpMa2}, for more details.
\end{remark}
\subsection{Numerical methods for the WKB-system}
In order to solve the Hamilton-Jacobi equation (\ref{eic}) numerically we shall rely on a relaxation method
as presented in \cite{JiXi}.
Then we can solve the system of transport equations \eqref{trans} by a time-splitting spectral scheme,
similar to the one proposed for the full MD system (\cf Section \ref{sec:num}). Using similar
notations, suppose that we know the values $u^{\pm,n}$, $V^{n}$ and $\Ab^{n}$.
\newpar
\textbf{Step 1.} First, we solve the following problem:
\begin{equation}\label{sst1}
\left \{
\begin{aligned}
\partial _t  u^\pm + \nabla\cdot \mathbf{v} (u^\pm) = &\, \eta(u^\pm),\\
(\p_{tt}-\btu) \mathcal V = & \,  \rho^0,\\
(\p_{tt}-\btu)\mathcal A = & \,\Jb^0,
\end{aligned}
\right.
\end{equation}
by a pseudo-spectral method, where we use the shorthanded notations
\begin{equation}
\mathbf{v}(u^\pm):=\omega_0^\pm (\nabla \phi^\pm )\otimes  u^\pm, \quad
\eta(u^\pm):= \frac{1}{2}\diverg (\omega_{0}^\pm(\nabla \phi^\pm))u^\pm.
\end{equation}
First, we take a DFT of \eqref{sst1}, \ie
\bea\label{eq:sc11}
\left \{
\begin{aligned}
\p_t\hat{u}^\pm+i\xi\cdot\hat {\mathbf{v}} (u^\pm)= & \, \hat \eta(u^\pm),\\
(\p_{tt}+|\xi|^2)\hat{\mathcal V}= & \,\hat{\rho}^0, \\
(\p_{tt}+|\xi|^2)\hat {\mathcal A}= & \,\hat{\Jb}^0.
\end{aligned}
\right.
\eea
Let us denote by $u^{\pm,n}$, the value of $u^\pm$ at time $t_{n}$ in Step 1.
Then we can find the values of $\hat u^{\pm,n+1}$, $\hat {\mathcal V}^{n+1}$, and $\hat{\mathcal A}^{n+1}$ by
the Crank-Nicolson scheme. After an IDFT, we obtain the values of
$u^{\pm,n+1}$, $\mathcal V^{n+1}$, and $\mathcal A^{n+1}$.
\newpar
\textbf{Step 2.} It remains to solve the ordinary differential equation
\begin{equation}\label{sst2}
\begin{aligned}
& \p_t u^{\pm} = i \mathcal N^\pm[u] \, u^\pm,
\end{aligned}
\end{equation}
with $\mathcal N$ given by \eqref{nonl}.
Because $\mathcal N^\pm[u]$ does not change in step 2, we have
\[u^{\pm,n+1}=\exp\left(i\mathcal N^\pm[u] \tg t\right) u^{\pm,n}.\]
\begin{remark}
We can also use the \emph{Strang-splitting method} to obtain a second order scheme in time.
Again, it is easy to see that this algorithm conserves \eqref{cons}.
\end{remark}
The solution of the Hamilton-Jacobi equation (\ref{eic}) may develop singularities at caustic manifolds, also
the group velocities $\omega_0^\pm (\nabla \phi^\pm )$ and the principal amplitudes
become singular. This makes the numerical approximation
of the transport equations (\ref{trans}) a difficult task. Actually, we are
not aware of a previous numerical study on such transport equations with caustic type singularities. Our
computational experience indicates that it is important to
conserve the density in the
transport problem \eqref{trans}, which relies on an accurate (high-order)
numerical approximation of
the terms $\omega_0^\pm (\nabla \phi^\pm )$
and $\diverg (\omega_{0}^\pm(\nabla \phi^\pm))$. However, the
Hamilton-Jacobi equation is typically solved by a shock capturing type
method, which reduces to first order at singularities.
In order to get a better numerical approximation, we still use a shock
capturing method, namely \emph{the relaxation scheme} developed in \cite{JiXi},
spatially for the Hamilton-Jacobi equation (\ref{eic}), but use
the \emph{fourth order Runge-Kutta method} temporally.
For the transport problem \eqref{trans} we found that the pseudo-spectral method behaves better
than finite difference schemes.

\subsection{Numerical examples in the semi-classical regime}
In all of the following examples we shall assume for simplicity
\begin{equation}
V^{(0)}(\xb)= V^{(1)}(\xb) =0, \quad \Ab^{(0)}(\xb)=\Ab^{(1)} (\xb)=0,
\end{equation}
since different, \ie non-zero, initial conditions would only add to the homogeneous solution of
the corresponding wave equation.
\begin{remark} Remark that in the following numerical examples $\phi_I$ has to be chosen such that it
satisfies the periodic boundary conditions.
\end{remark}
\begin{example}[\textbf{Self-consistent steady state}]\label{exsc1}
Consider the system (\ref{dmsc}) with initial
condition
\be
\label{eq:e09} \psi^{\e}\big|_{t=0}= \chi \,
\exp{\left(-\frac{|\xb|^2}{4d^2}\right)},
\quad \chi=(1,0,0,0),\ d=1/16,
\ee
and zero external potentials, \ie
$A^{ex}_k(\xb)=V^{ex}(\xb)=0$. This example models a wave packet
with initial width $d$ and zero initial speed, propagating
only under its self-interaction. Note that in this case
$\phi^\pm_I(\xb)\equiv0$ and $u^+(0,\xb)$ is simply given by
\eqref{eq:e09}, whereas $u^-(0,\xb)\equiv 0$, hence $u^-(t,\xb) =
0$, for $t>0$. First, we choose $\e=10^{-2}$ and compare the solution
of the full MD system with the numerical solution obtained by
solving the asymptotic WKB-system \eqref{eic}, \eqref{trans}. From
Figure \ref{fig21} we see that the two numerical solutions agree
very well for such a small $\e$. In particular the creation of
positrons in the full MD system is small, \ie $O(\e)$ as one
expects from the semi-classical analysis. This is clearly visible
in \cf Figure \ref{fig22}, which shows that the projectors
$\Pi^\pm_0(\btd\phi)$ are indeed good approximations of
$\Pi^\pm_0(-i\e\btd)$ for $\e$ is small. However for $\e=1$ this
is no longer true. Furthermore, because in this case the WKB-phase is found to be simply given by
$\phi^+(t,x)=-t$, we thus have $\btd\phi^+\equiv 0$ and $\btd\cdot\omega^+_0=0$, and
hence the transport equation \eqref{trans} simplifies to
$$\partial_t u^+ + i \mathcal V u^+=0,$$
which implies $|u^+(t,\xb)|^2$ to be constant. In this particular case,
we can use a very coarse mesh to get satisfactory results (\cf Table \ref{tb101}).
Remark that the results in Table \ref{tb101} also illustrate the validity of \eqref{eq:scerr}.
\begin{table}[ht]
\begin{center}
\caption{Difference between the asymptotic solution and the full MD system
for example \ref{exsc1} ($\tg t=1/128$, $\tg x=1/32$):}\label{tb101}
\begin{tabular}{cccc}\hline
$\e$ &         0.0001  &   0.001 & 0.01 \\ \hline 
$\ba{c}\vspace{-4mm}\\ \dpm\sup_{0\le t\le 0.25}\left\|\, \psi^\e -
\sum_{\pm}u^\pm e^{i \phi^\pm/\varepsilon } \,
\right\|_{L^2(\Og)\otimes \C^4}\ea$
 & 3.20E-3 &  3.34E-2 & 2.98E-1 \vspace*{0.5mm}\\ \hline
$\ba{c}\vspace{-4mm}\\ \dpm\sup_{0\le t\le 0.25}\left\|\, \psi^\e -
\sum_{\pm}u^\pm e^{i \phi^\pm/\varepsilon } \,
\right\|_{L^\ift(\Og)\otimes \C^4}\ea$
 & 4.90E-3 &  5.01E-2 & 4.40E-1 \vspace*{0.5mm}\\ \hline 
\end{tabular}
\end{center}
\end{table}
\begin{figure} 
\begin{center}
\resizebox{2in}{!} {\includegraphics{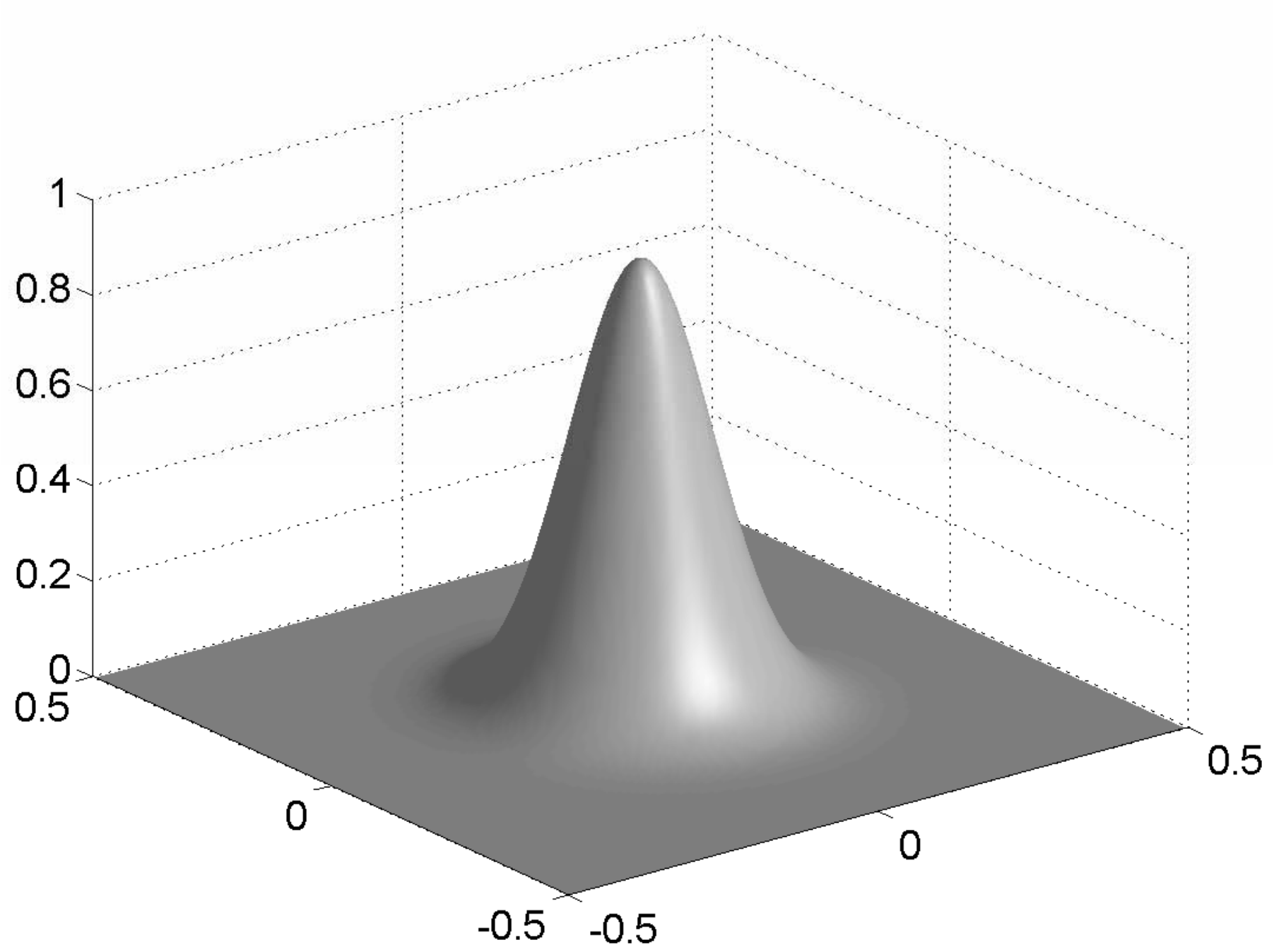}}
\resizebox{2in}{!} {\includegraphics{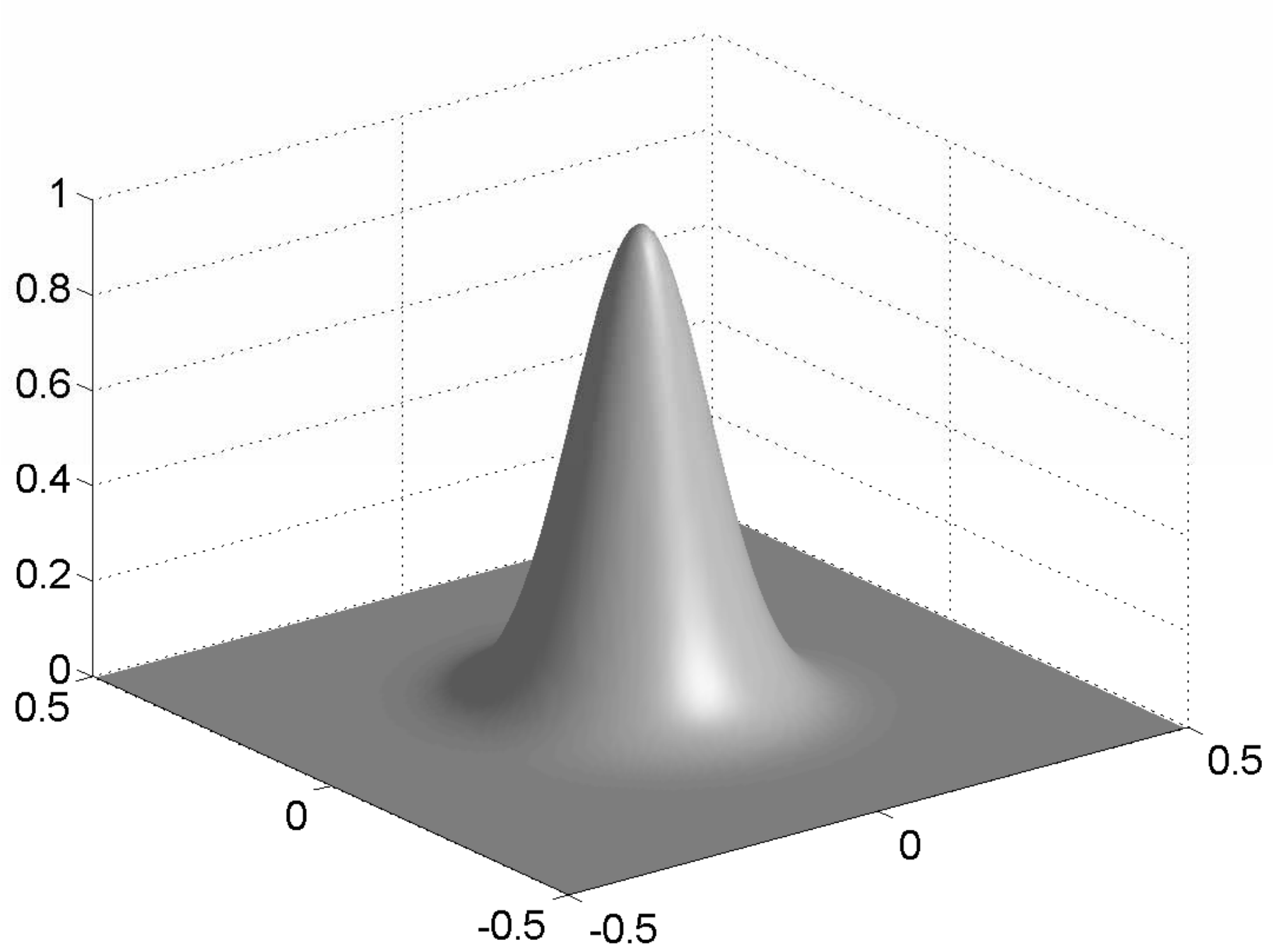}}

{$\left|\psi^\e(t,\xb)|_{x_3=0}\right|^2$ and $\left|\sum_\pm
u^\pm(t,\xb)e^{i\phi^\pm/\e}|_{x_3=0}\right|^2$ at
$t=0.25$}\vspace*{2mm}

\resizebox{2in}{!}{\includegraphics{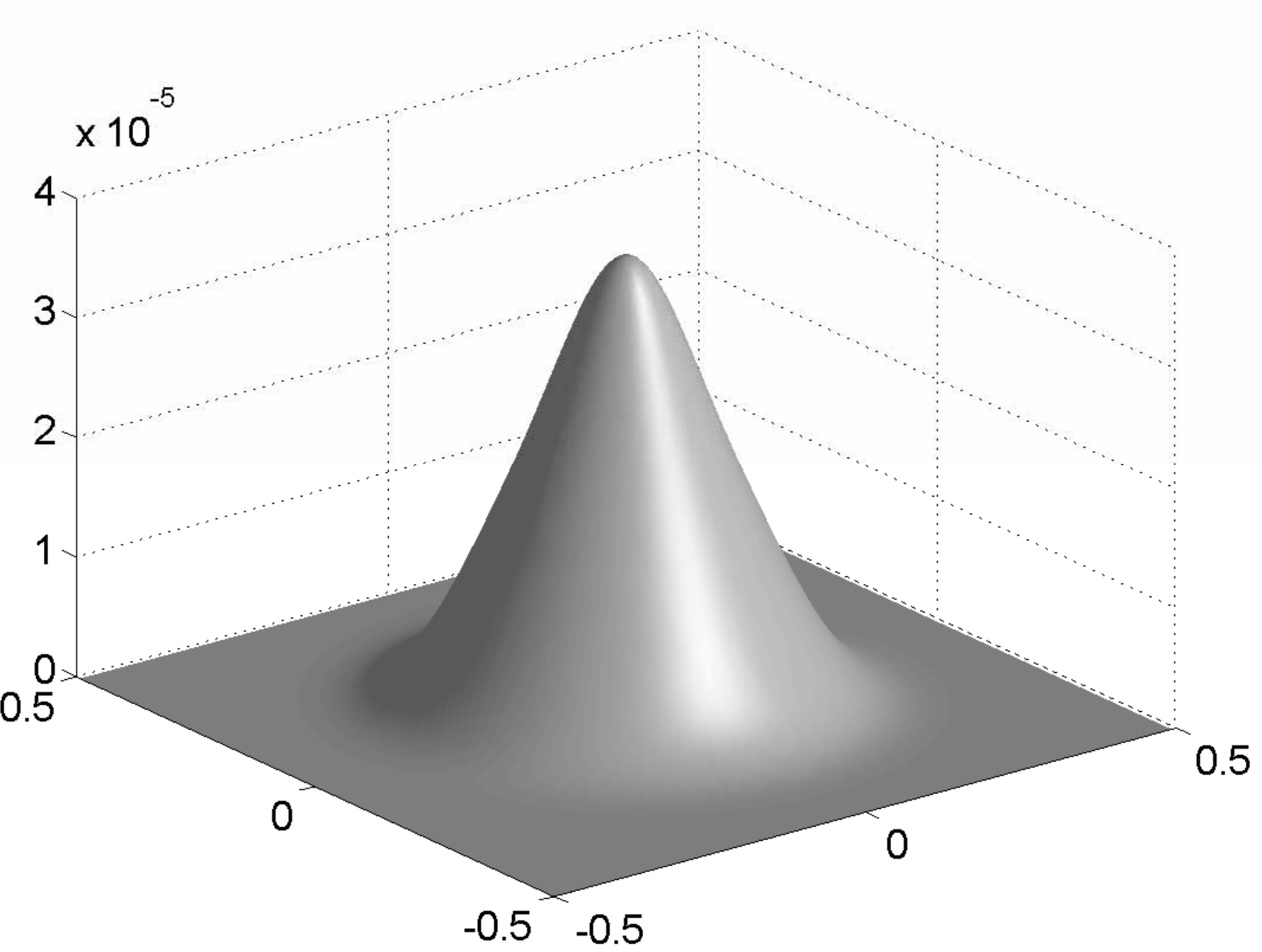}}
\resizebox{2in}{!}{\includegraphics{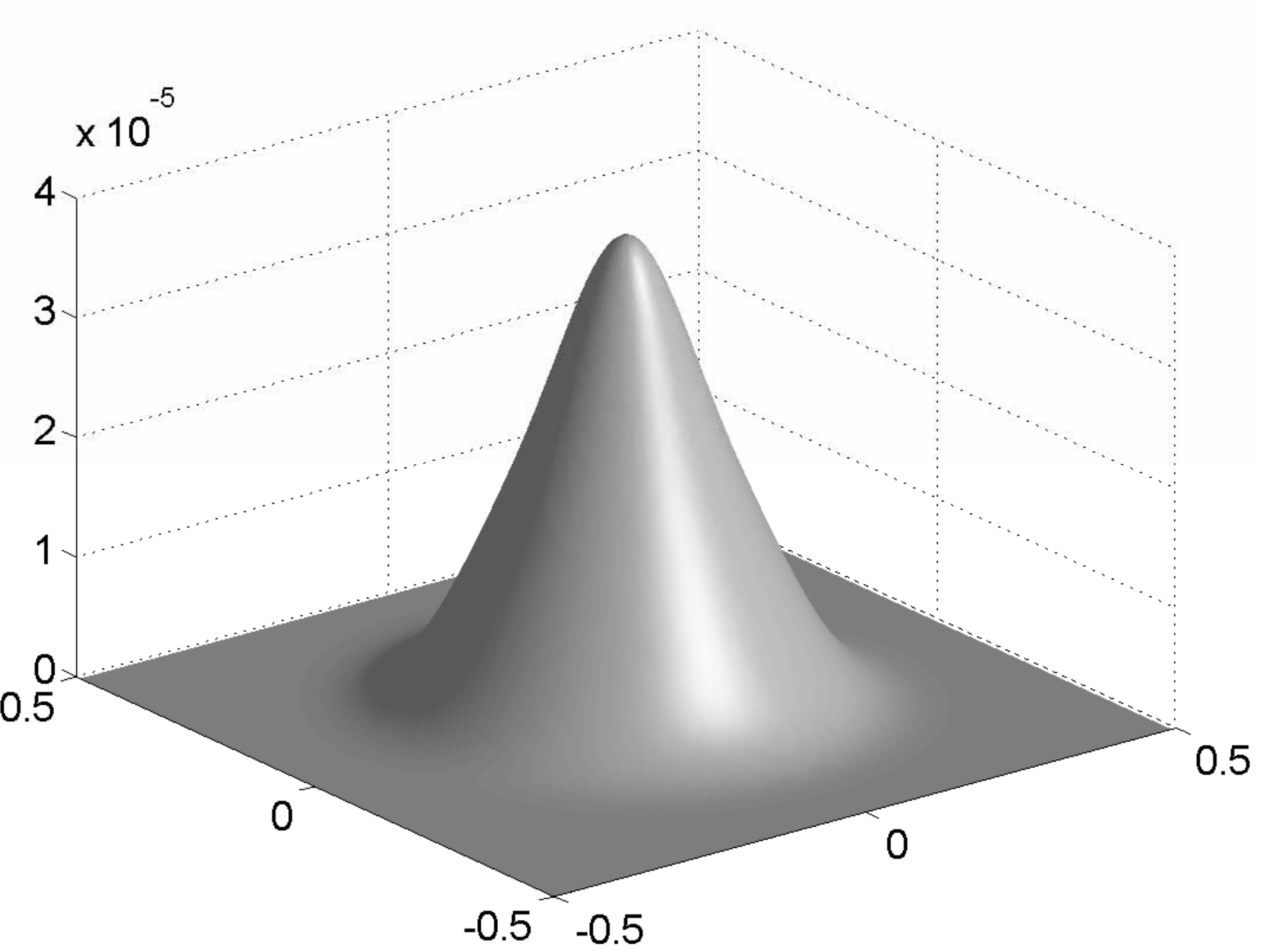}}

{$V^\e(t,\xb)|_{x_3=0}$ and $V(t,\xb)|_{x_3=0}$ at
$t=0.25$}\vspace*{5mm}

\resizebox{2in}{!}{\includegraphics{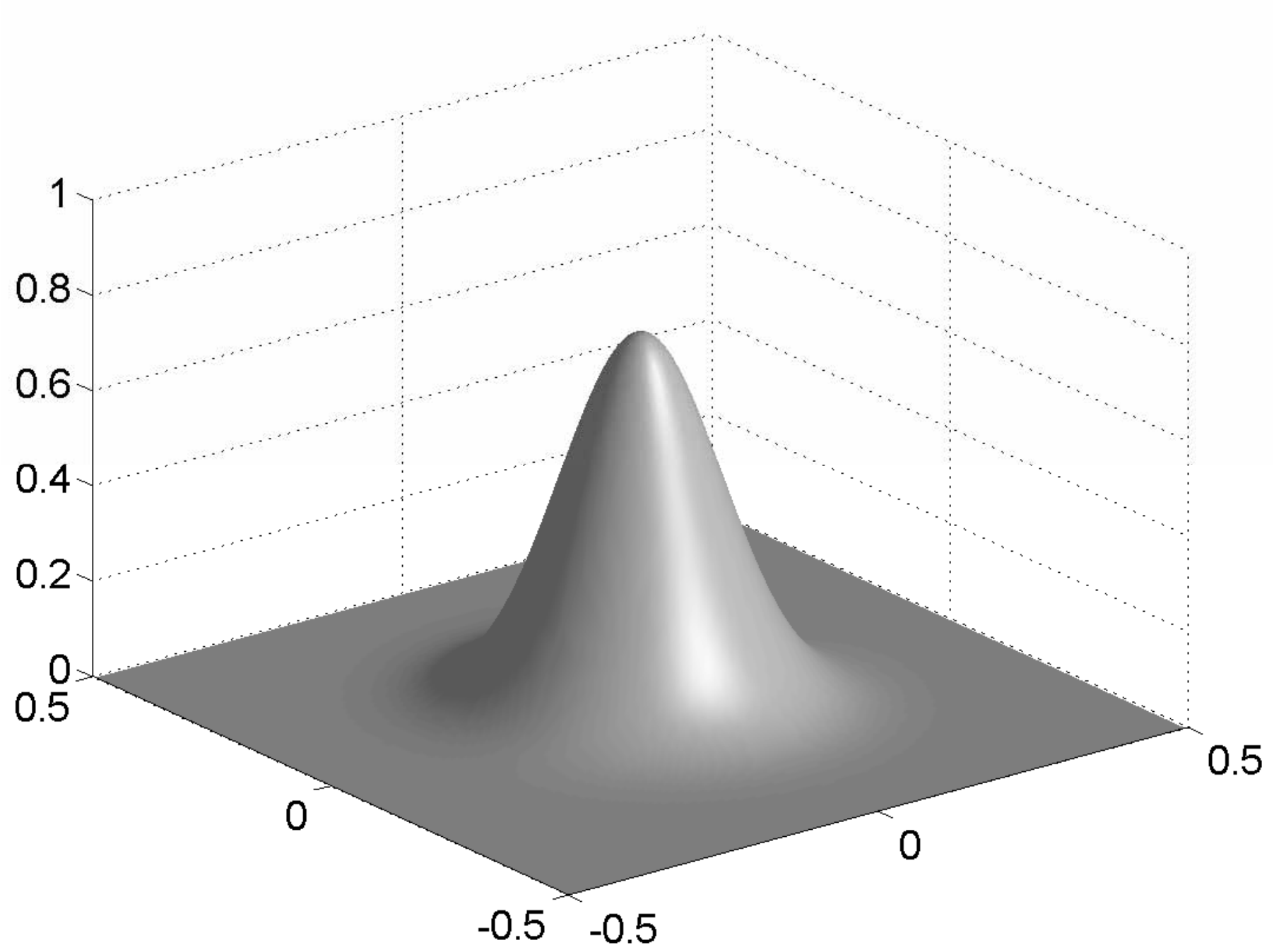}}
\resizebox{2in}{!}{\includegraphics{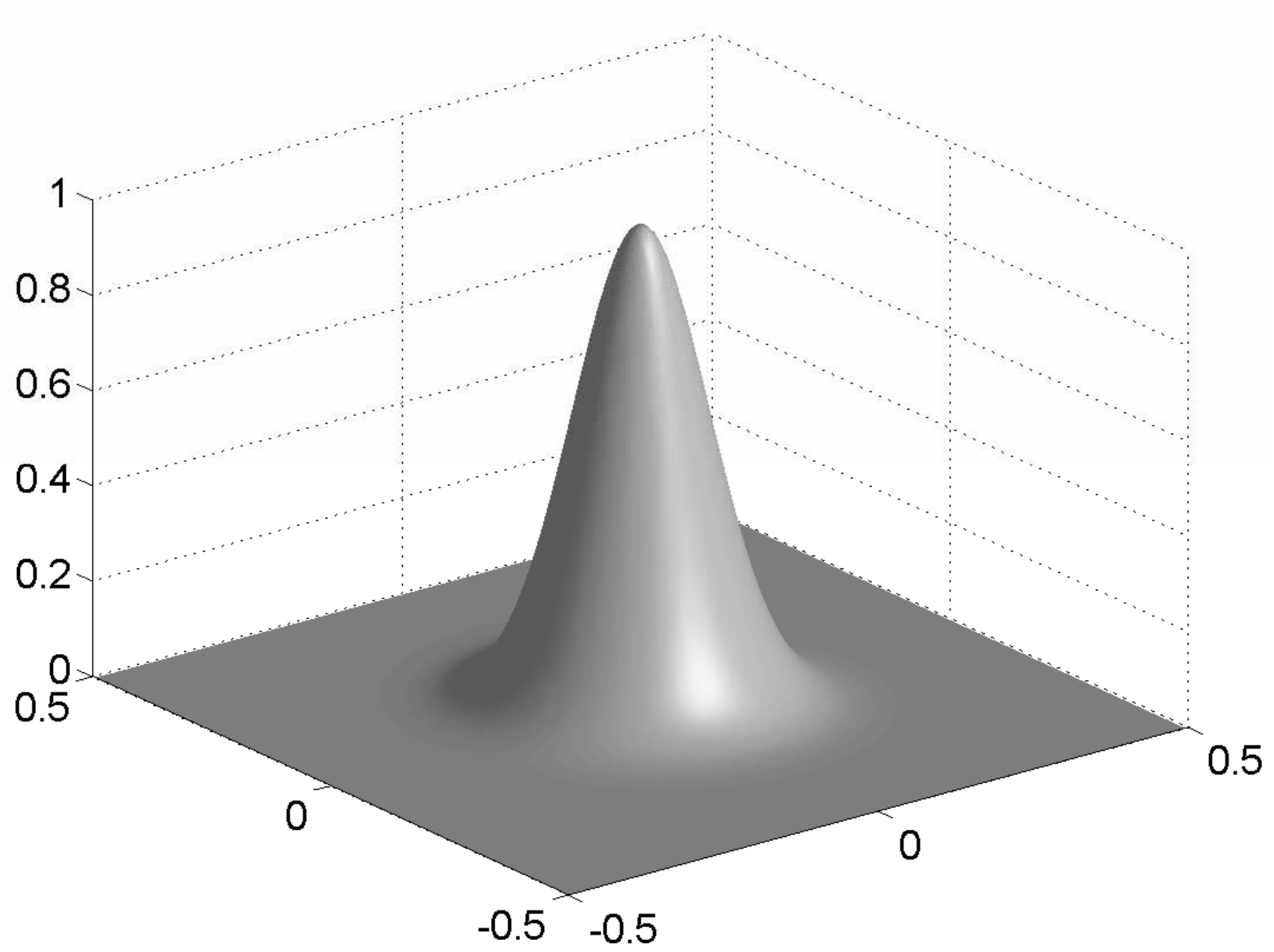}}

{$\big|\psi^\e(t,\xb)|_{x_3=0}\big|^2$ and $\left|\sum_\pm
u^\pm(t,\xb)e^{i\phi^\pm/\e}|_{x_3=0}\right|^2$
 at $t=0.5$}\vspace*{2mm}

\resizebox{2in}{!}{\includegraphics{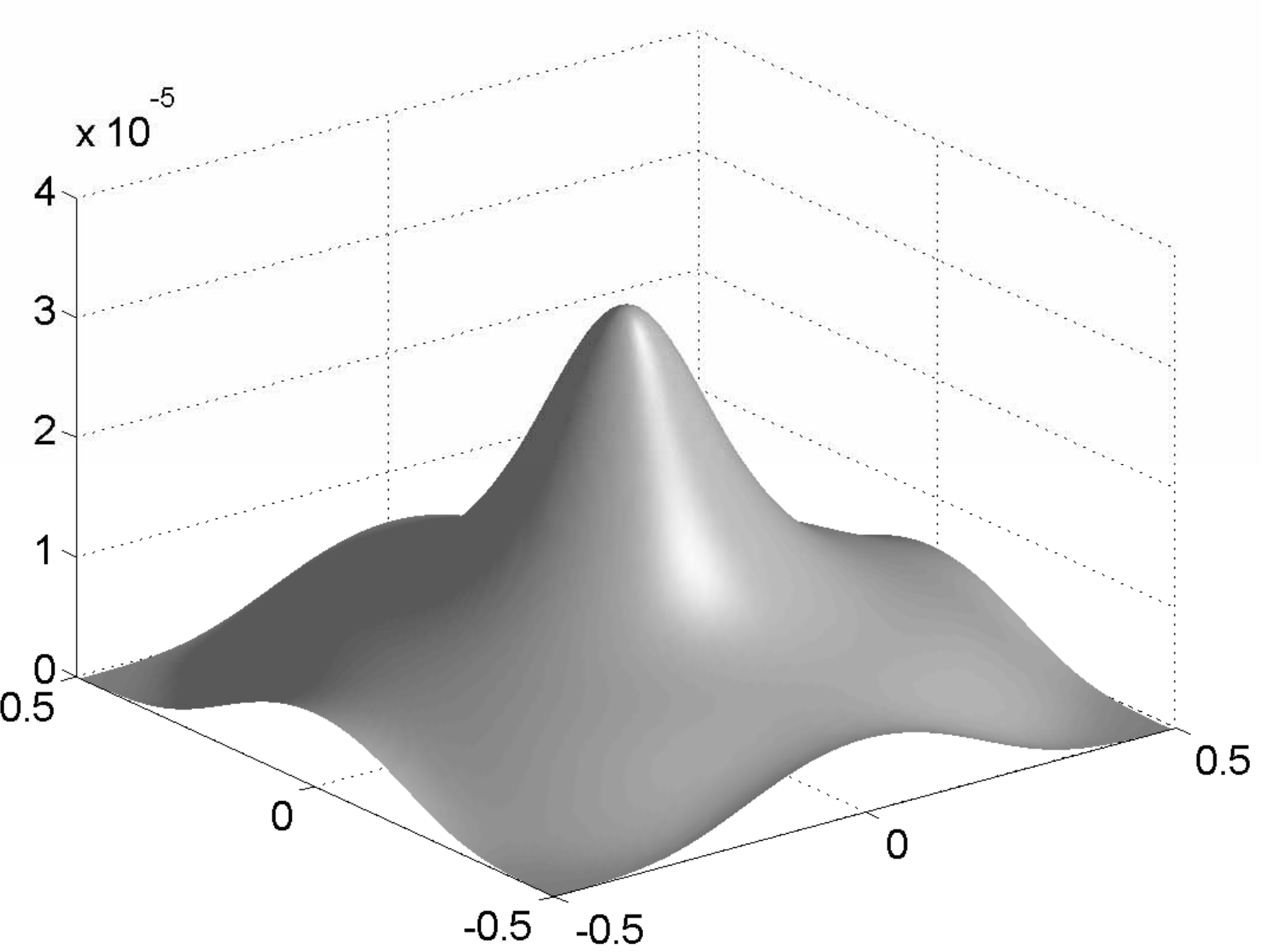}}
\resizebox{2in}{!}{\includegraphics{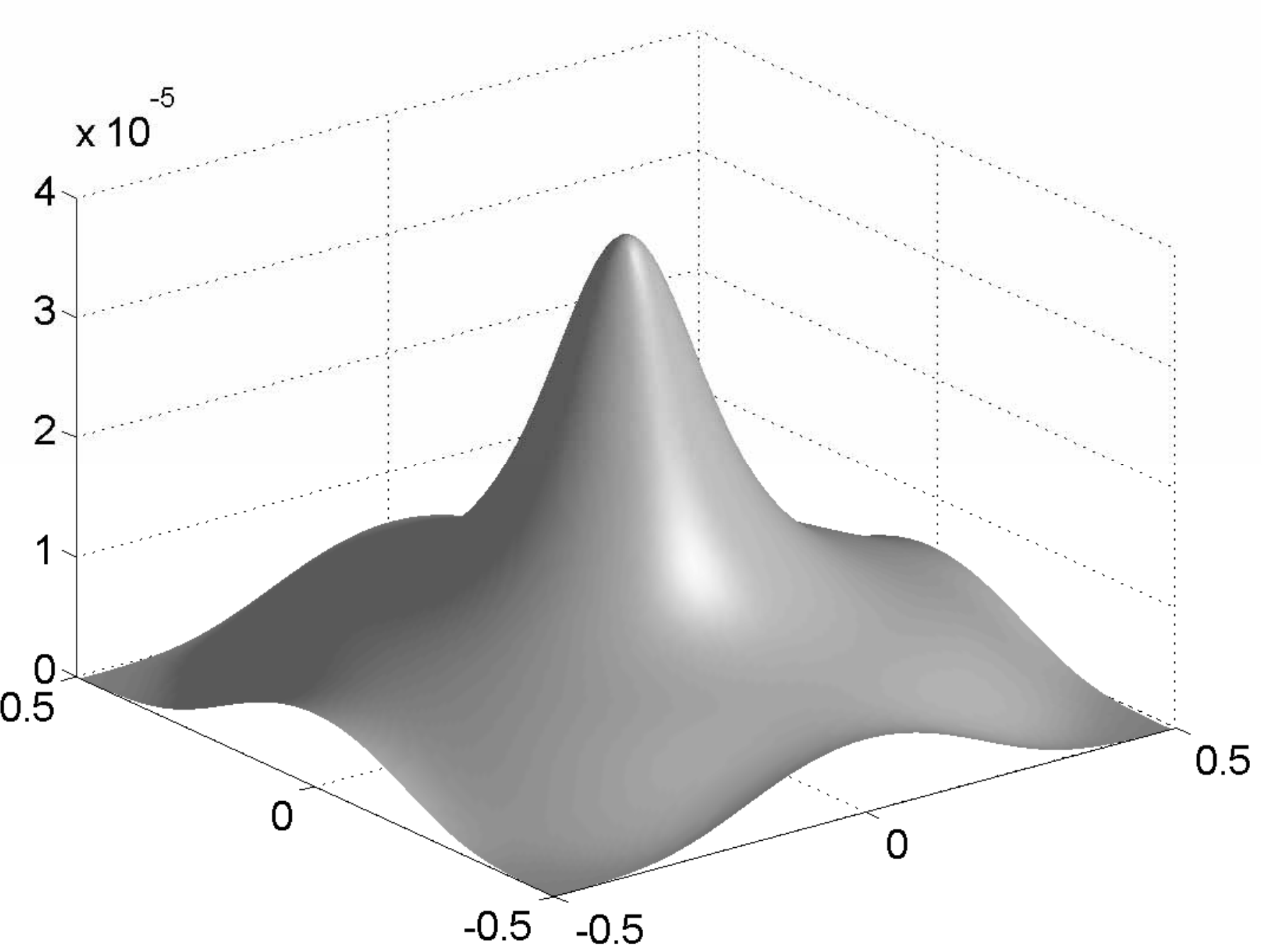}}

{$V^\e(t,\xb)|_{x_3=0}$ and $V(t,\xb)|_{x_3=0}$ at $t=0.5$}
\end{center}
\caption{Numerical results for example \ref{exsc1}.
The left column shows the graphs of the solution of the MD system,
the right column shows the graphs of the solution of the asymptotic problem.
Here  $\e=0.01$, $\tg t=\frac{1}{128}$, $\tg x=\frac{1}{32}$.}\label{fig21}
\end{figure}
\begin{figure} 
\begin{center}
\resizebox{2in}{!} {\includegraphics{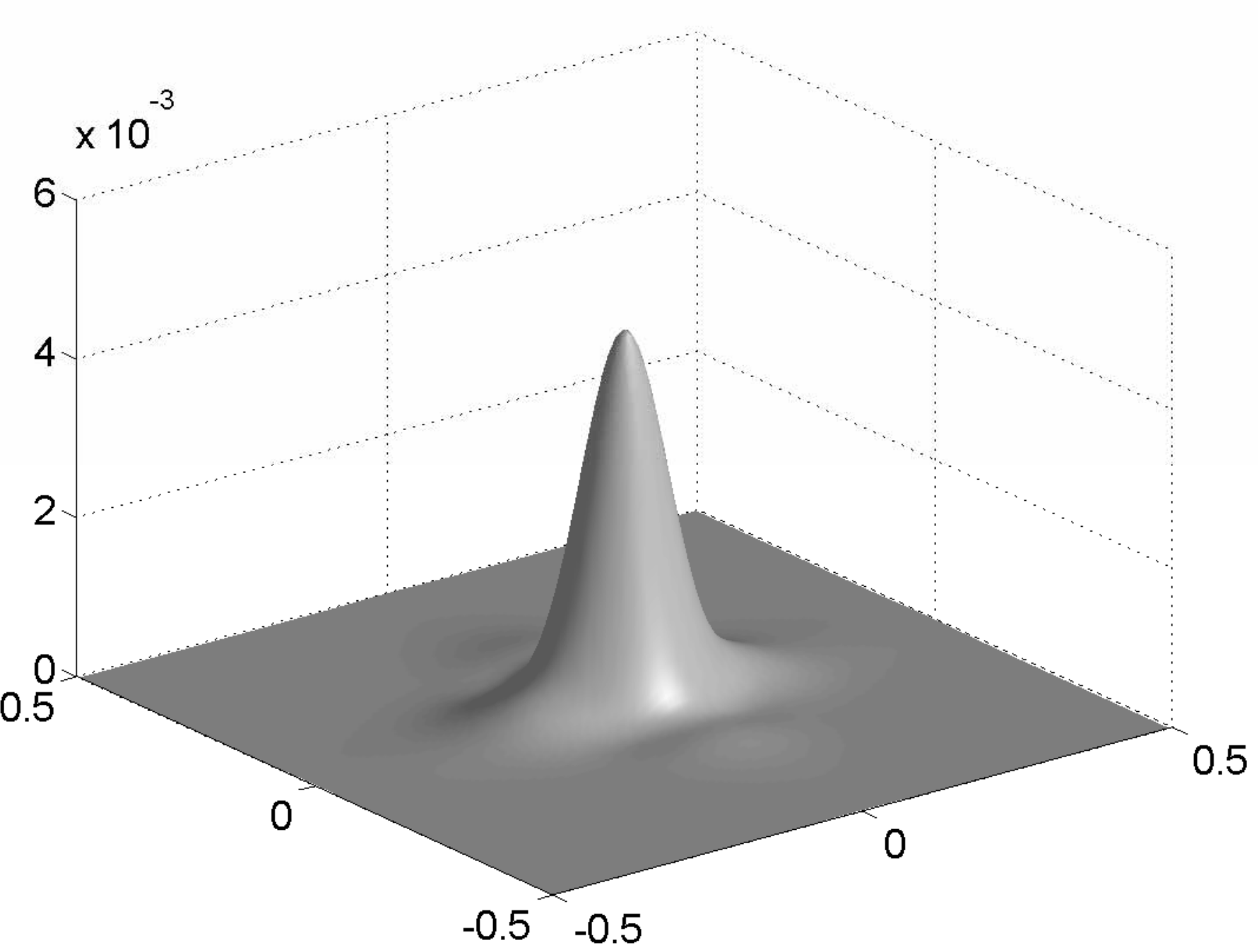}}
\resizebox{2in}{!} {\includegraphics{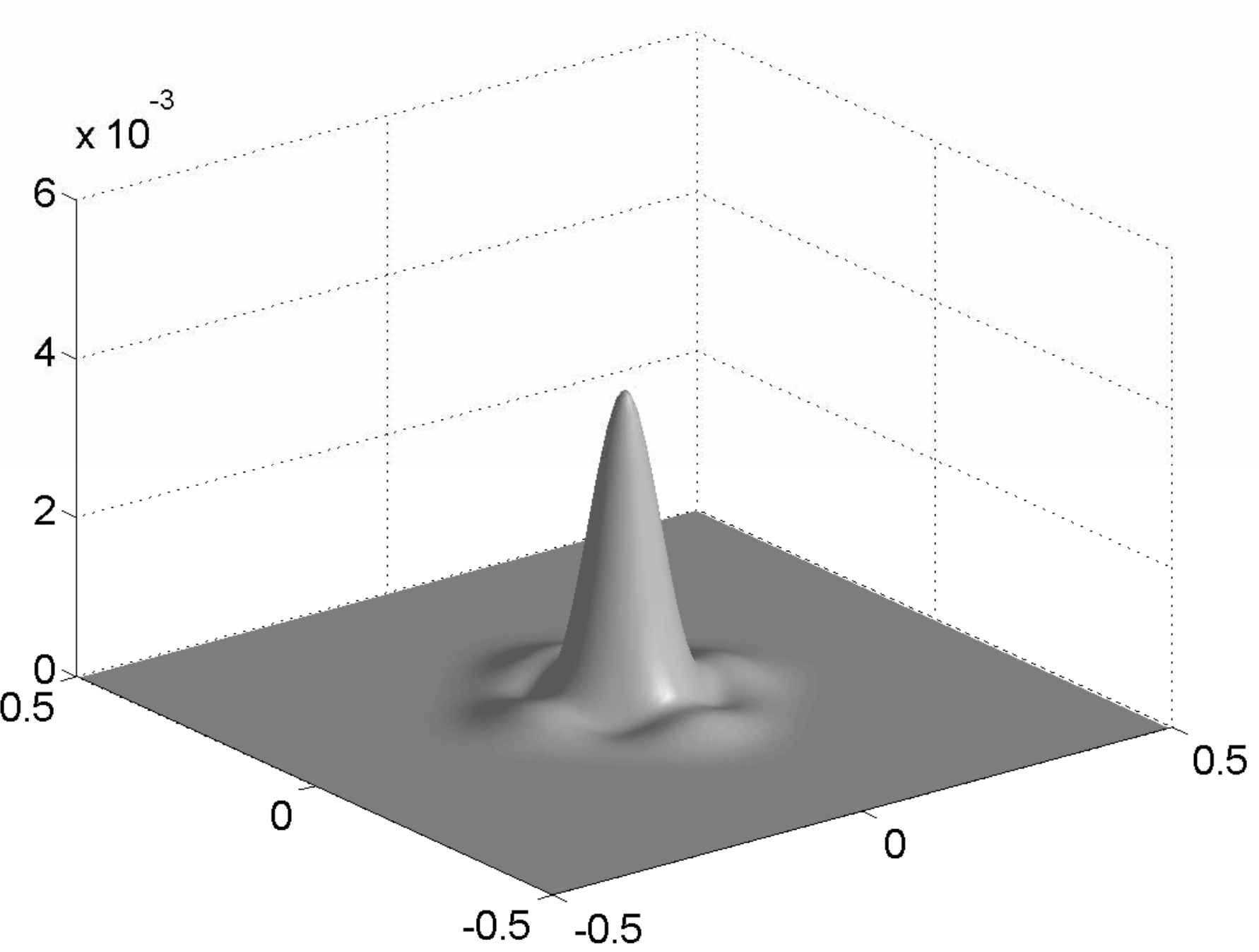}}

{$\e=0.01$, $\big|\Pi^-_0(-i\e\btd)\psi^\e(t,\xb)|_{x_3=0}\big|^2$ and
$\big|\Pi^-_0(\btd\phi)\psi^\e(t,\xb)|_{x_3=0}\big|^2$}\vspace*{2mm}

\resizebox{2in}{!} {\includegraphics{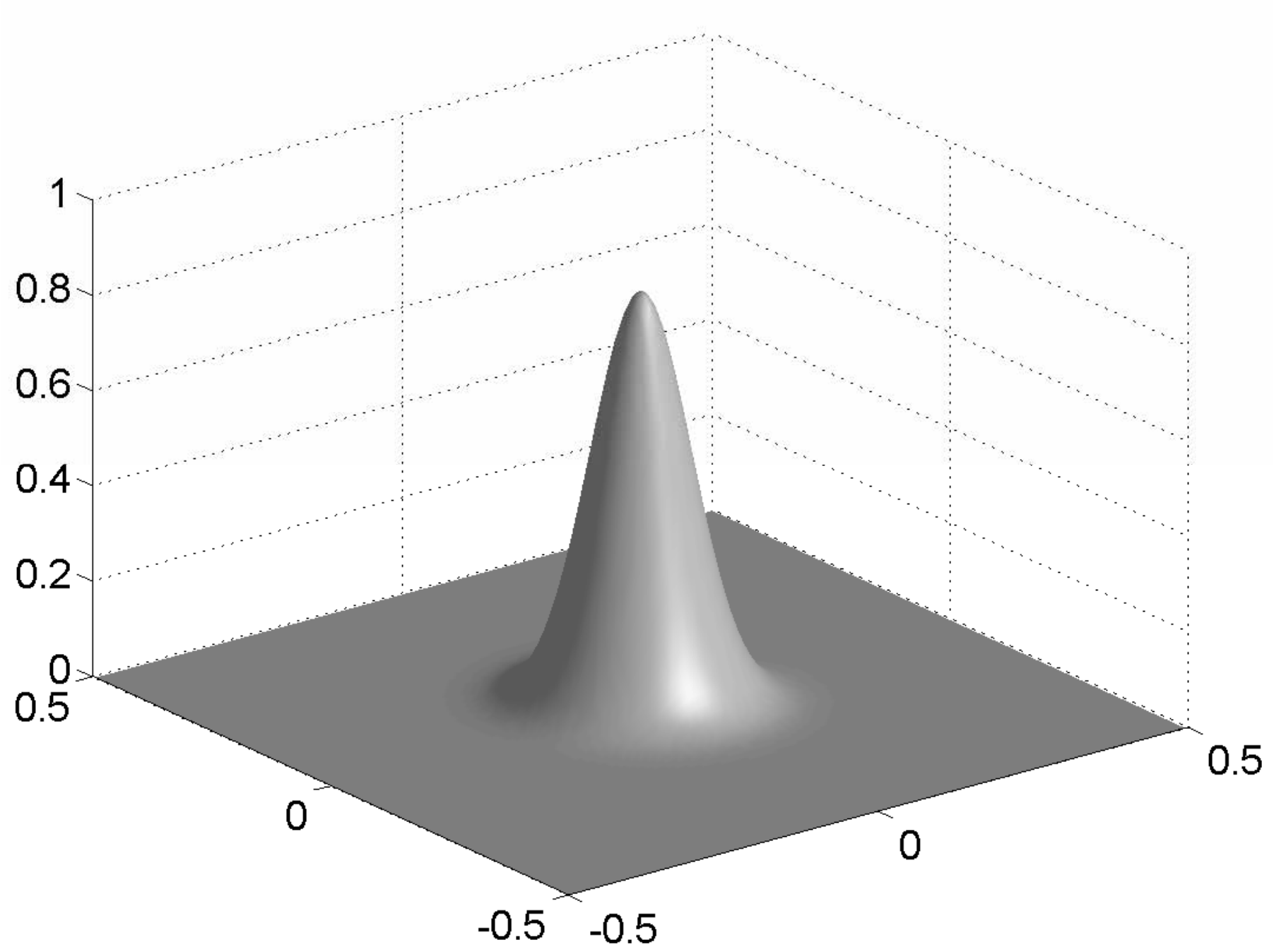}}
\resizebox{2in}{!} {\includegraphics{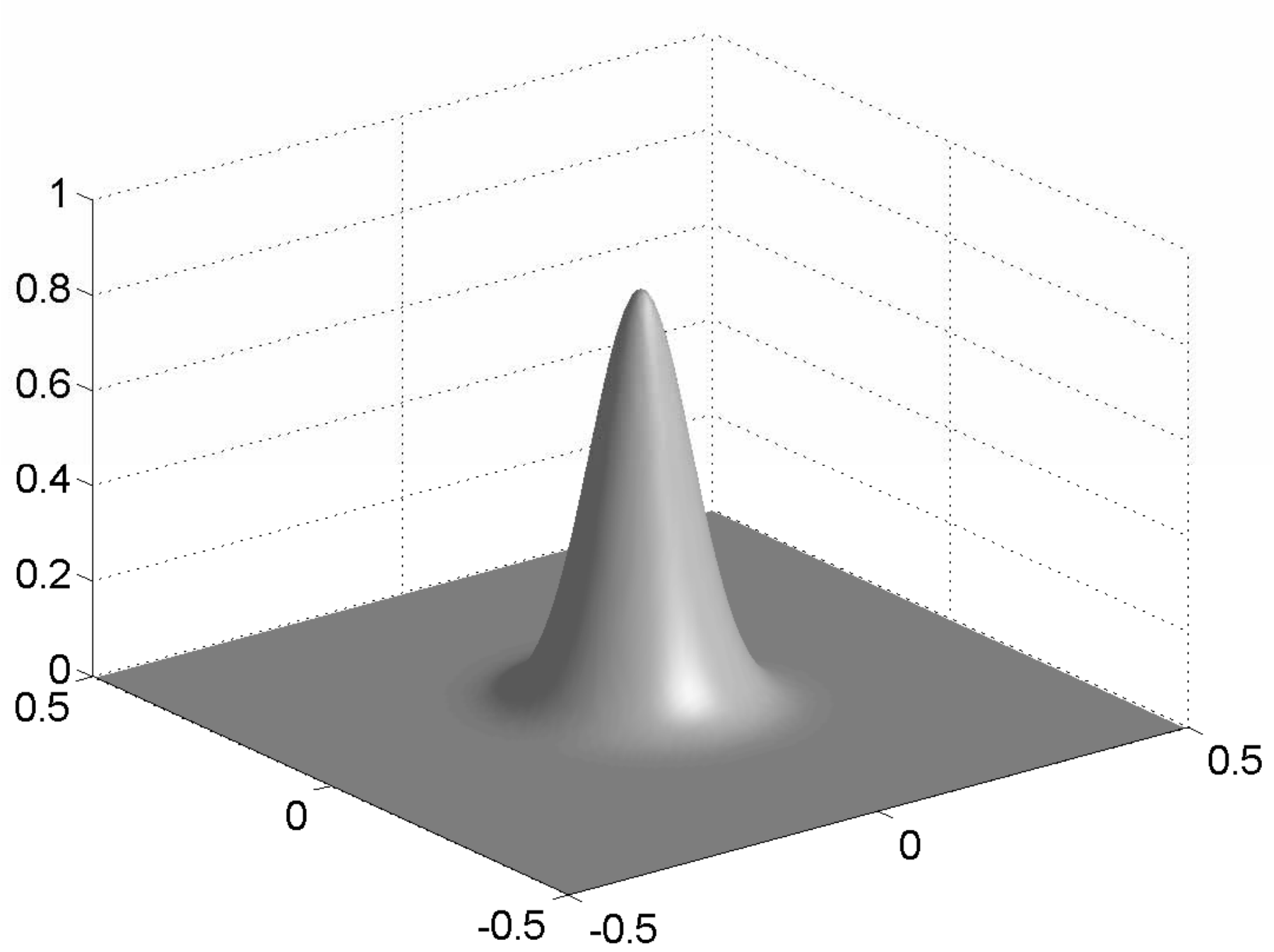}}

{$\e=0.01$, $\big|\Pi^+_0(-i\e\btd)\psi^\e(t,\xb)|_{x_3=0}\big|^2$ and
$\big|\Pi^+_0(\btd\phi)\psi^\e(t,\xb)|_{x_3=0}\big|^2$}\vspace*{2mm}

\resizebox{2in}{!} {\includegraphics{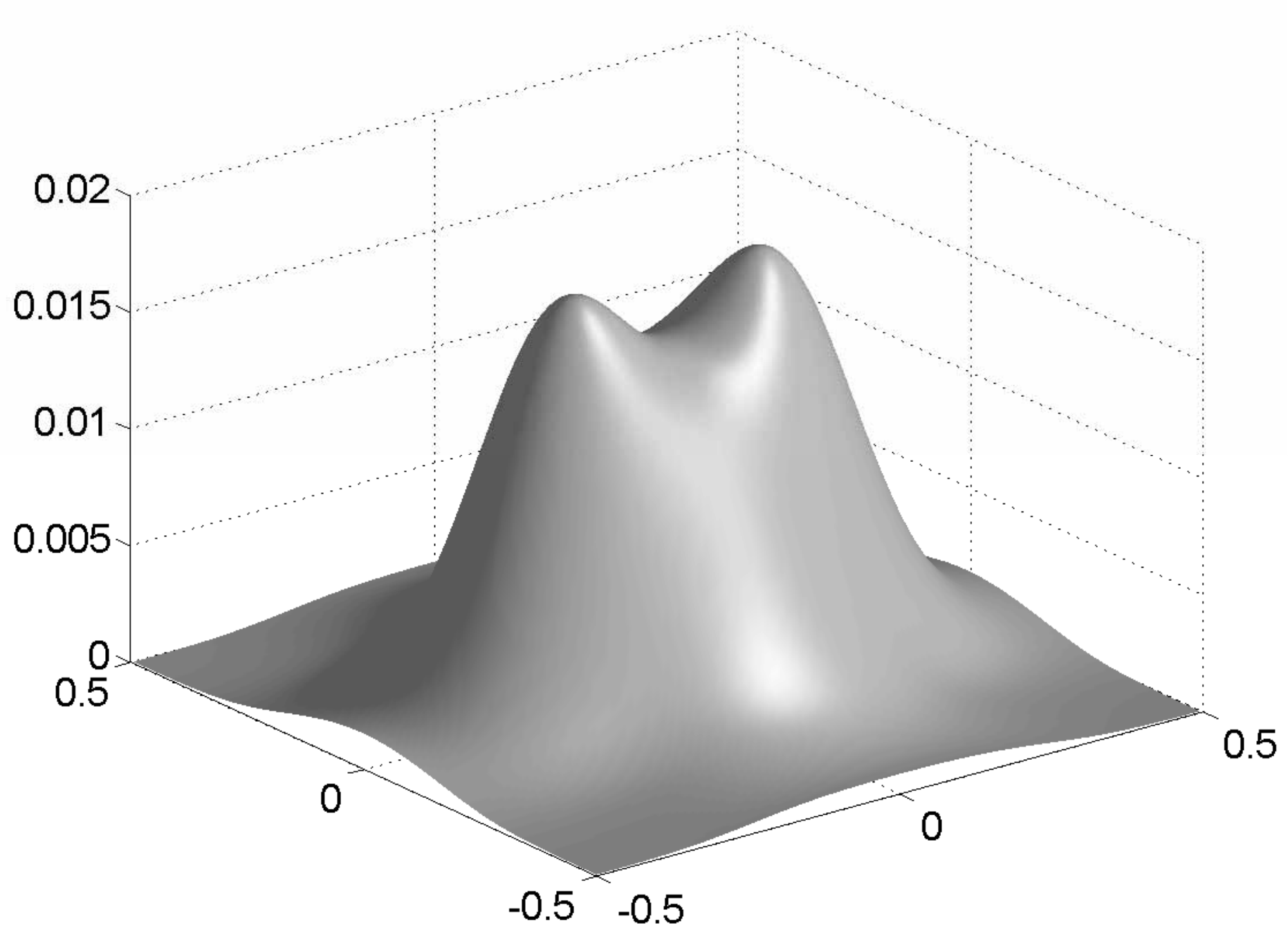}}
\resizebox{2in}{!} {\includegraphics{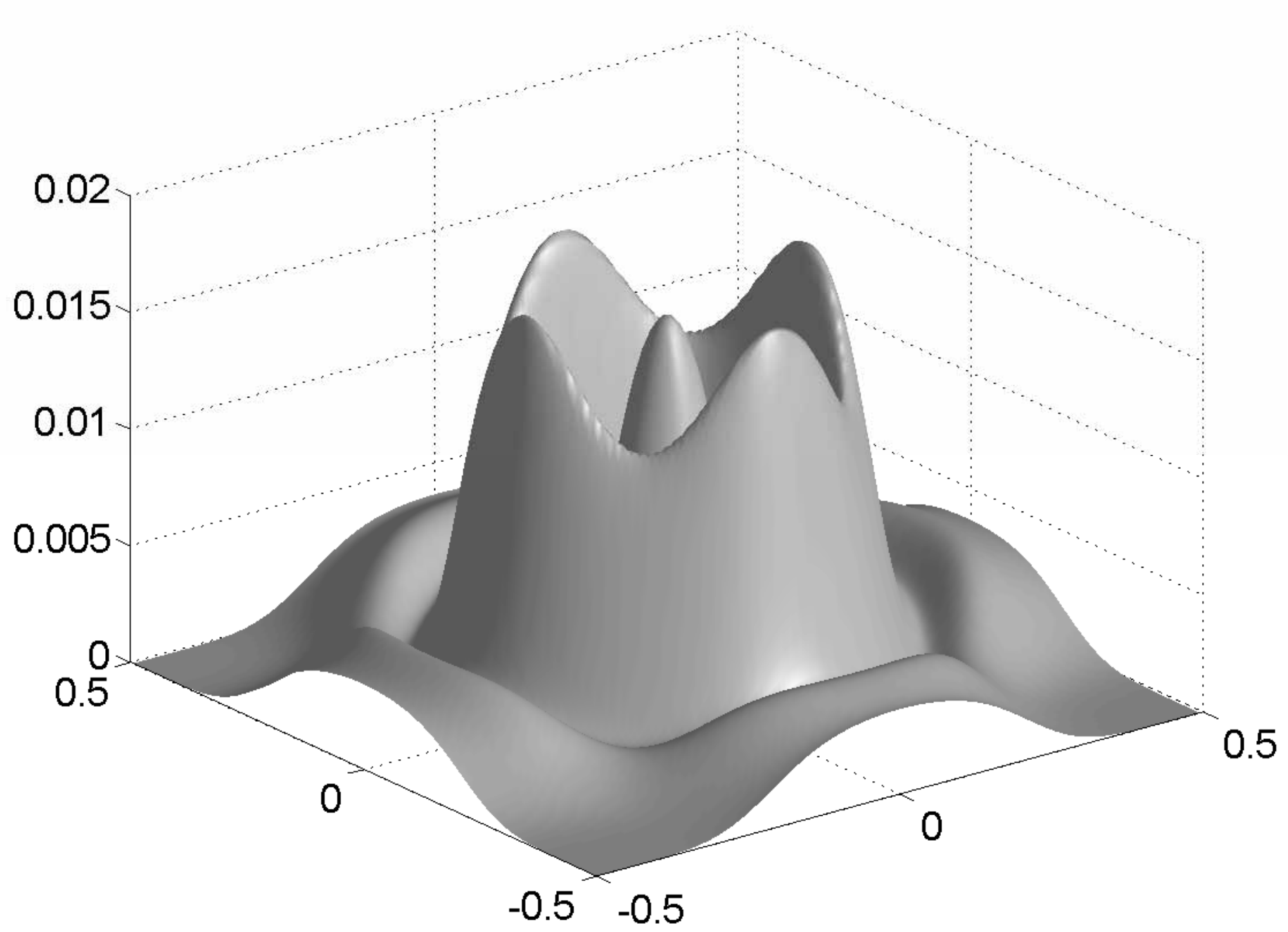}}

{$\e=1.0$, $\big|\Pi^-_0(-i\e\btd)\psi^\e(t,\xb)|_{x_3=0}\big|^2$ and
$\big|\Pi^-_0(\btd\phi)\psi^\e(t,\xb)|_{x_3=0}\big|^2$}\vspace*{2mm}

\resizebox{2in}{!} {\includegraphics{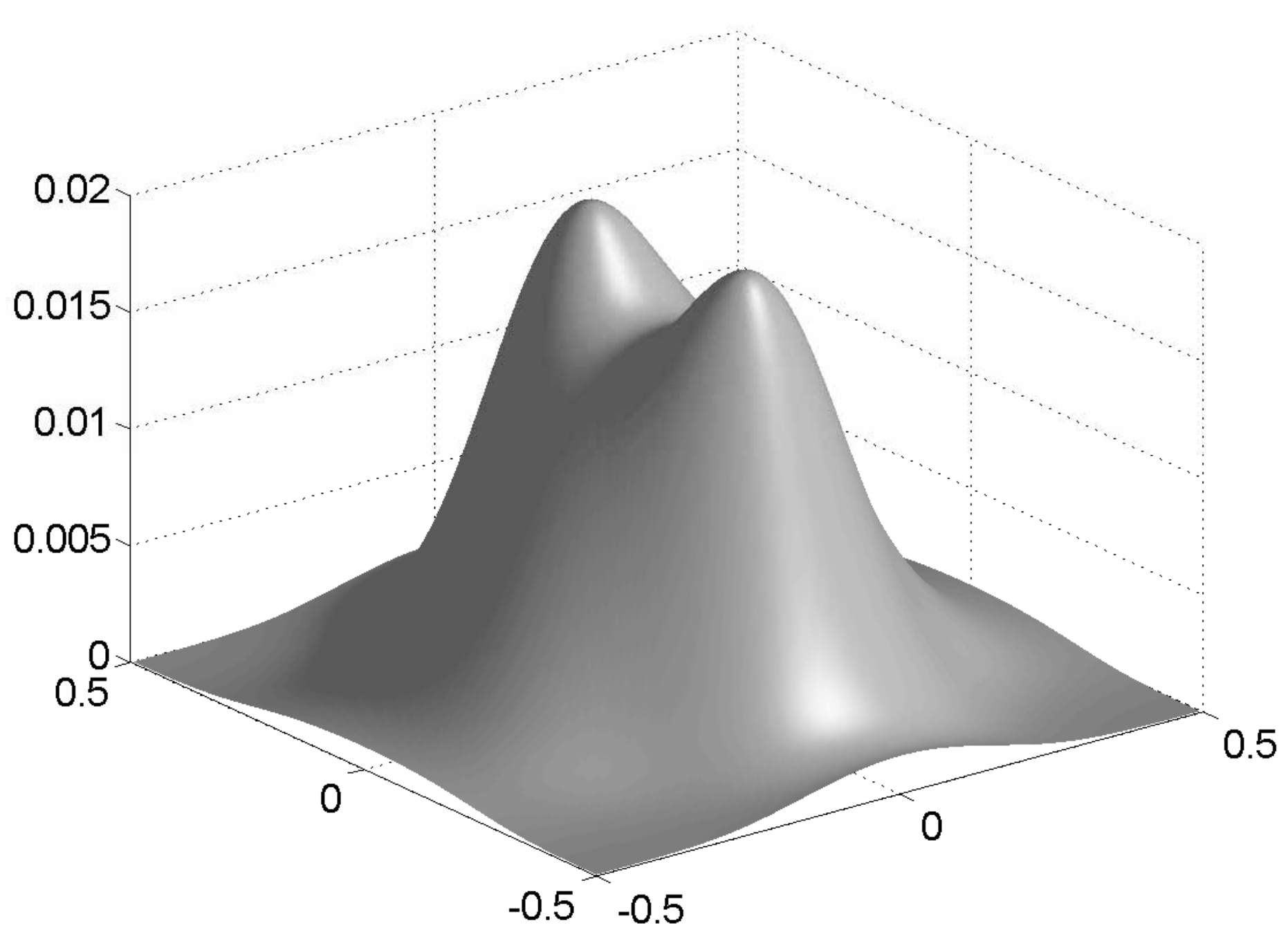}}
\resizebox{2in}{!} {\includegraphics{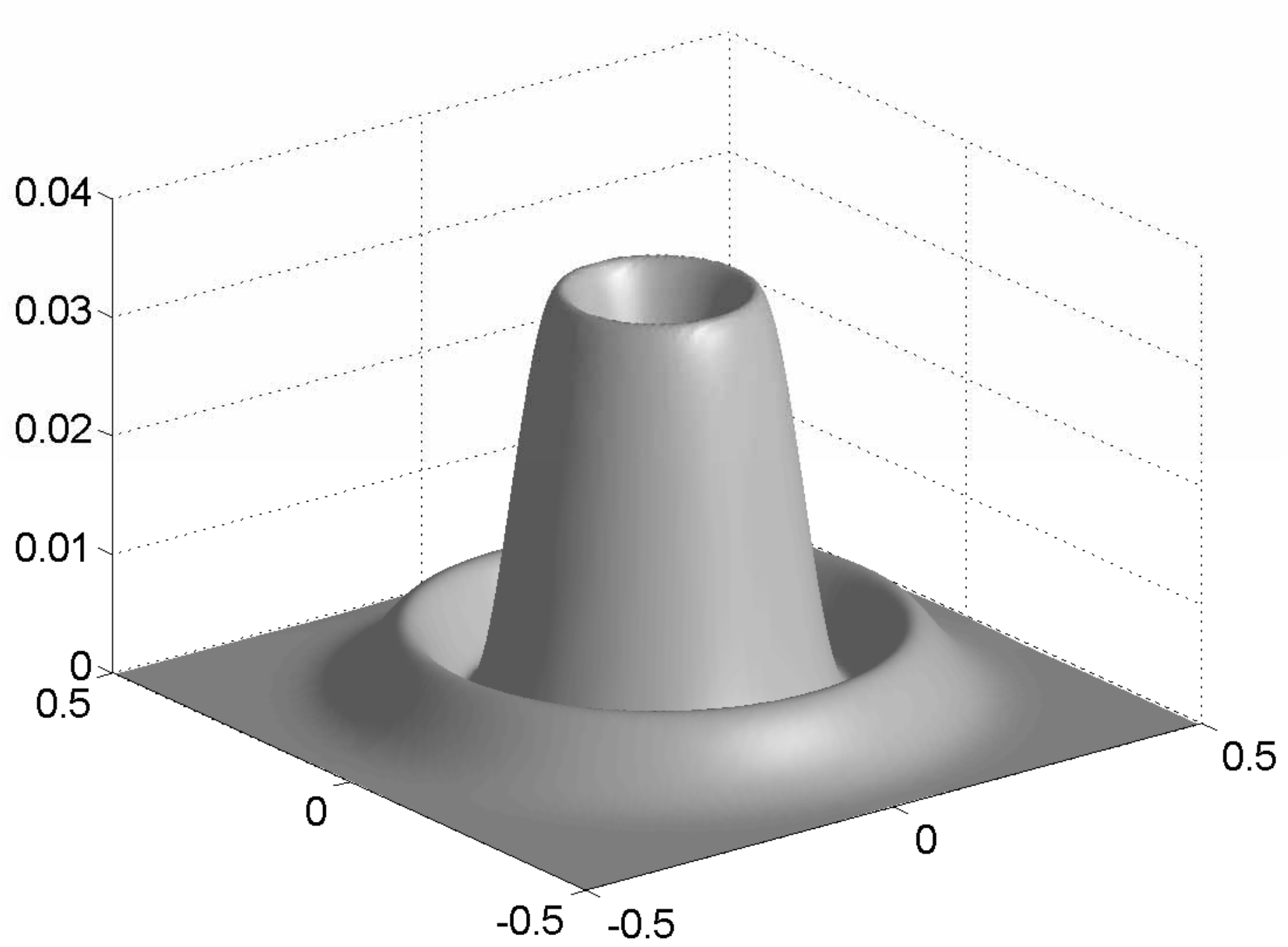}}

{$\e=1.0$, $\big|\Pi^+_0(-i\e\btd)\psi^\e(t,\xb)|_{x_3=0}\big|^2$ and
$\big|\Pi^+_0(\btd\phi)\psi^\e(t,\xb)|_{x_3=0}\big|^2$}
\end{center}
\caption{Numerical results of the densities of electron/positron projectors for example \ref{exsc1}.
The left column is
$\big|\Pi^\pm_0(-i\e\btd)\psi^\e(t,\xb)|_{x_3=0}\big|^2$,
the right column is $\big|\Pi^\pm_0(\btd\phi)\psi^\e(t,\xb)|_{x_3=0}\big|^2$.
Here $t=0.25$, $\tg x=1/32$, $\tg t=1/128$.}
\label{fig22}
\end{figure}
\end{example}
\newpage
\begin{example}[\textbf{Purely self-consistent motion}]\label{exsc2}
In this example, again zero external fields are assumed, but we modify the initial condition for
$\psi^\e$ as follows:
\bea\label{eq:e10}
\psi^\e\big|_{t=0}= \chi(\xb)\,\exp{\left(-\frac{|\xb|^2}{4d^2}+ i\frac{\phi_I(\xb)}{\e}\right)},\quad d=1/16,
\eea
where the phase function describing the $\e$-oscillations is given by
\bea
\phi_I(\xb)  =   \frac{1}{40}\, (1+\cos 2\pi x_1)(1 + \cos2\pi x_2)
\eea
and we choose the initial amplitude such that $\Pi_0^+(\nabla \phi_I(\xb)) \chi(\xb) = \chi(\xb)$, \ie
\begin{equation}
\label{chi}
\chi(\xb) = \left(\dpm \f{\xi^2_{1}(\xb)+\xi^2_{2}(\xb)}{2(\sqrt{1+|\xi|^2}-1)},
-\f{\xi_{3}(\xb)(\xi_1(\xb)+i\xi_2(\xb))}{2(\sqrt{1+|\xi|^2}-1)},0,
\dpm\f{\xi_{1}(\xb)+i\xi_{2}(\xb)}{2}\right),\quad \xi = \nabla \phi_I(\xb).
\end{equation}
As in the above example we thus have $u^-(t,\xb)\equiv 0$.
Note that for $\phi_I =0$, \eqref{eq:e10} reduces to \eqref{eq:e09}.
The numerical solution of the eiconal equation \eqref{eic} \cite{JiXi}
indicates
a \emph{kink-type singularity} in the phase of our asymptotic description
at about $t\simeq 0.56$, \cf Figure \ref{fig32}.
Hence the asymptotic WKB-type approximation for the spinor field
is no longer correct for $t>0.56$,
\newpar
The numerical results for both the MD system and the semi-classical
limit for $\e=0.01$ are given in Figure \ref{fig311}.
Table \ref{tb21} attempts to show the validity of \eqref{eq:scerr}.
Compared to Table \ref{tb101}, the difference between two systems is somewhat larger than $O(\e)$. Our
experience indicates that this has to do with the numerical difficulties
mentioned before and with the fact that discretization errors
``pollute'' the solution of the semi-classical system as time evolves,
preventing a more accurate comparison at later time.
Due to our computing capacity, we are unable to conduct more refined calculation,
which would have provided a better justification of the
ansatz \eqref{eq:scerr} for this problem.
For the same problem, we also present the numerical solutions of
the Maxwell-Dirac system at later time in Figures \ref{fig312}.
We also present a numerical simulation of the case $\e=1.0$, \ie
away from the semi-classical regime, see
Figure \ref{fig33}. From the plots it becomes clear that the ``exact'' spinor field and the solution of the asymptotic
WKB-problem are qualitatively ``close'' for small values of $\e$ and before caustics, while they are even
qualitatively different away from the semi-classical regime.
\begin{figure} 
\begin{center}
\resizebox{2.5in}{!} {\includegraphics{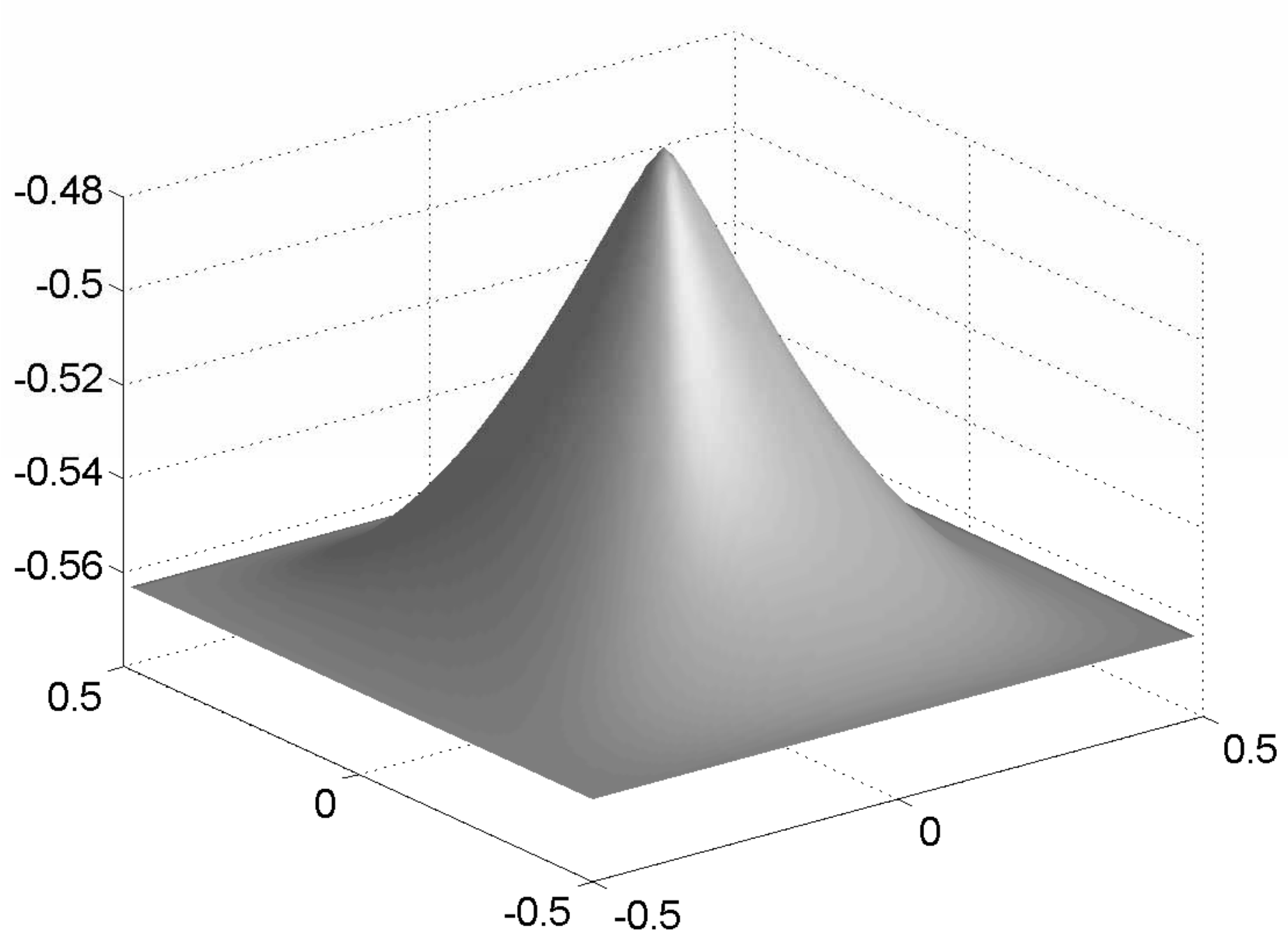}}
\end{center}
\caption{The graph of the phase $\phi(t,\xb)|_{x_3=0}$ at
$t=0.5625$ for example \ref{exsc2}. It shows the phase becomes singular
at the tip.} \label{fig32}
\end{figure}
\begin{table}[ht]
\begin{center}
\caption{Difference between the asymptotic solution and the full
MD system for example \ref{exsc2} ($\tg t=1/128$,
$\tg x=1/64$):}\label{tb21}
\begin{tabular}{cccc}\hline
$\e$ &         0.01  &   0.1 \\ \hline
$\ba{c}\vspace{-4mm}\\ \dpm\sup_{0\le t\le 0.125}\left\|\, \psi^\e -
\sum_{\pm}u^\pm e^{i \phi^\pm/\varepsilon } \,
\right\|_{L^2(\Og)\otimes \C^4}\ea$
 & 0.196 &  0.926 \vspace*{0.5mm}\\ \hline
$\ba{c}\vspace{-4mm}\\ \dpm\sup_{0\le t\le 0.125}\left\|\, \psi^\e -
\sum_{\pm}u^\pm e^{i \phi^\pm/\varepsilon } \,
\right\|_{L^\ift(\Og)\otimes \C^4}\ea$
 & 0.115 &  0.646 \vspace*{0.5mm}\\ \hline
\end{tabular}
\end{center}
\end{table}
\begin{figure} 
\begin{center}
\resizebox{2in}{!}{\includegraphics{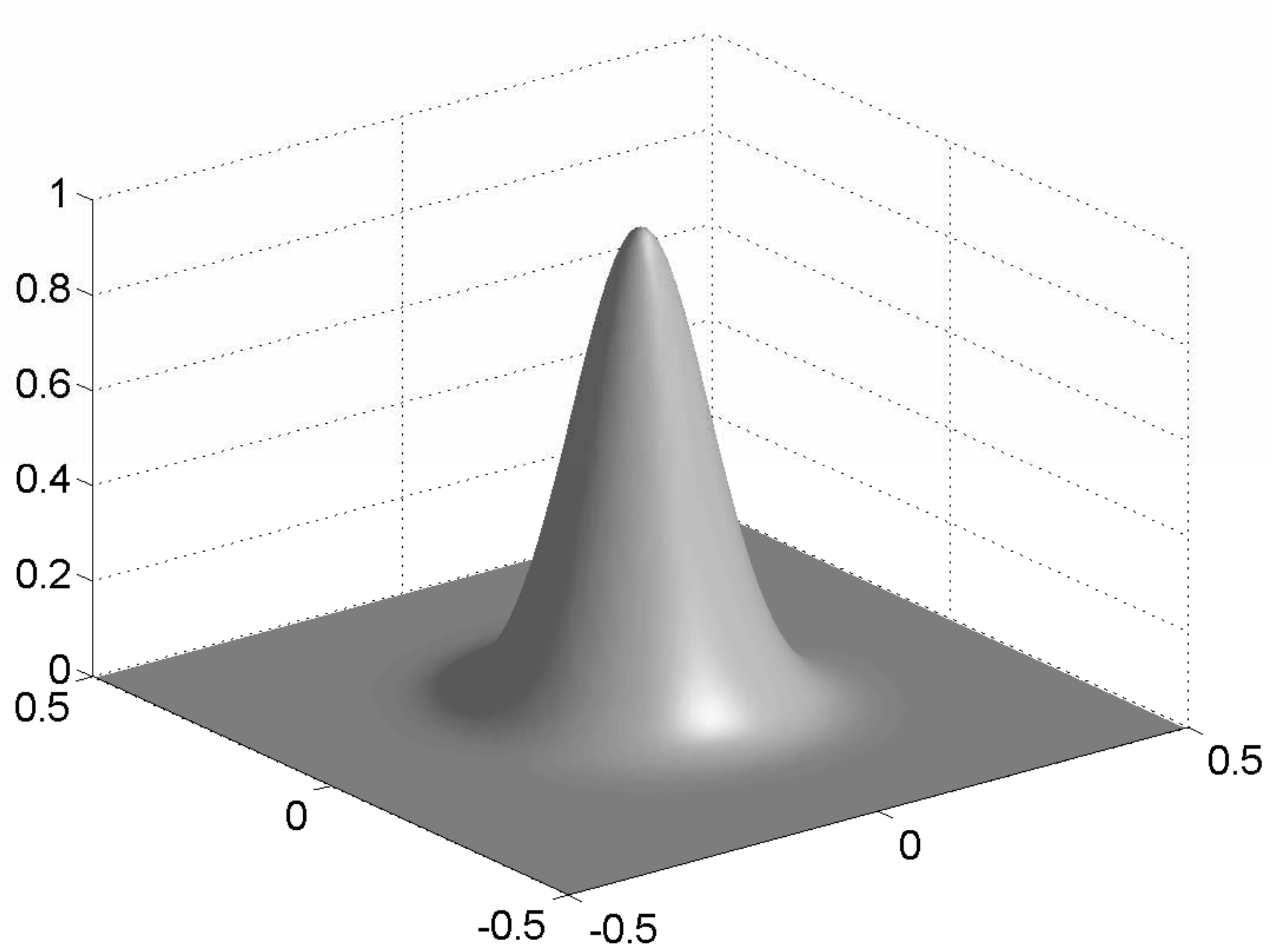}}
\resizebox{2in}{!}{\includegraphics{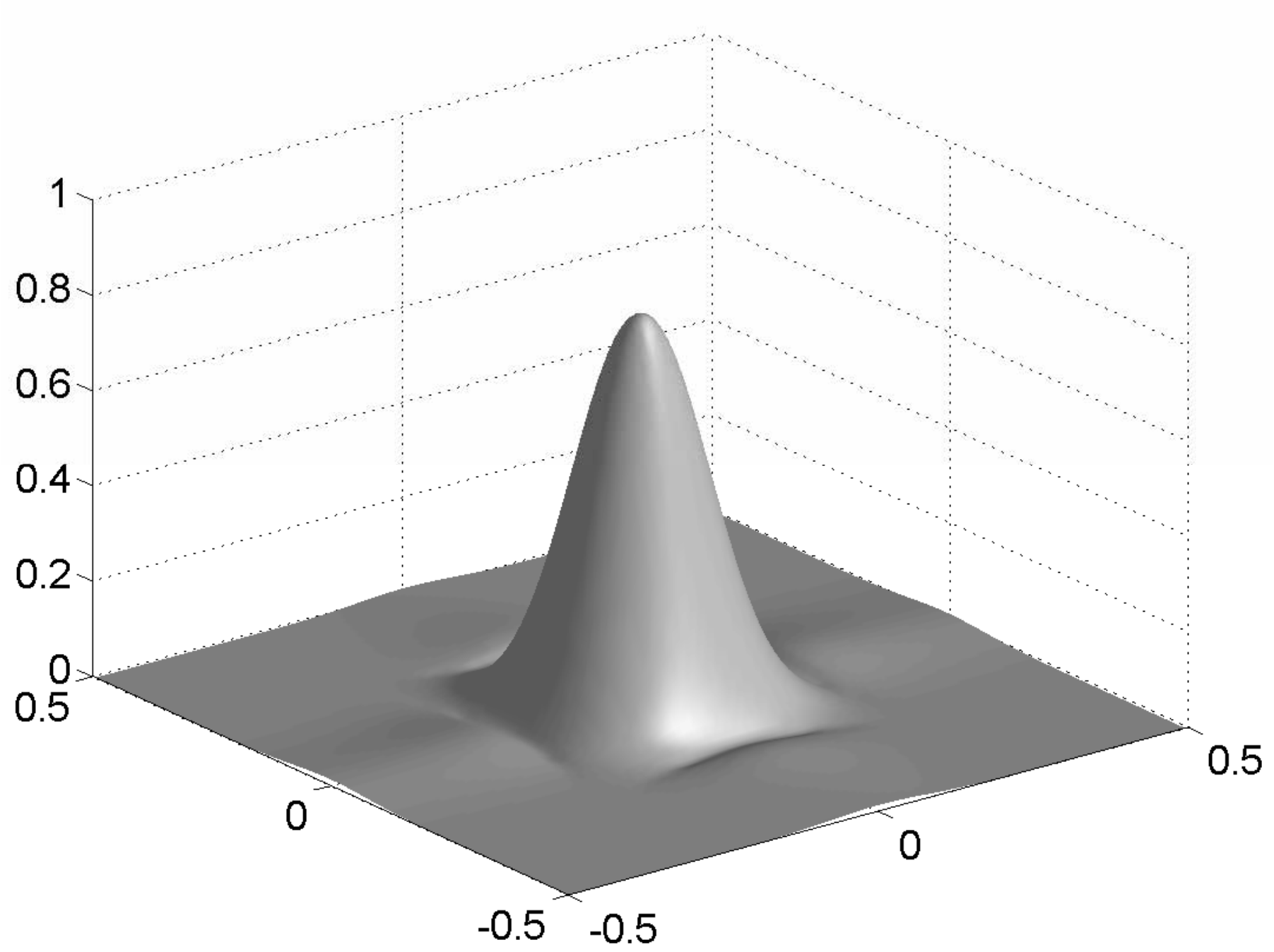}}

{$\big|\psi^\e(t,\xb)|_{x_3=0}\big|^2$ and $\left|\sum_\pm u^\pm(t,\xb)e^{i\phi^\pm/\e}|_{x_3=0}\right|^2$
at $t=0.25$}\vspace*{1mm}

\resizebox{2in}{!}{\includegraphics{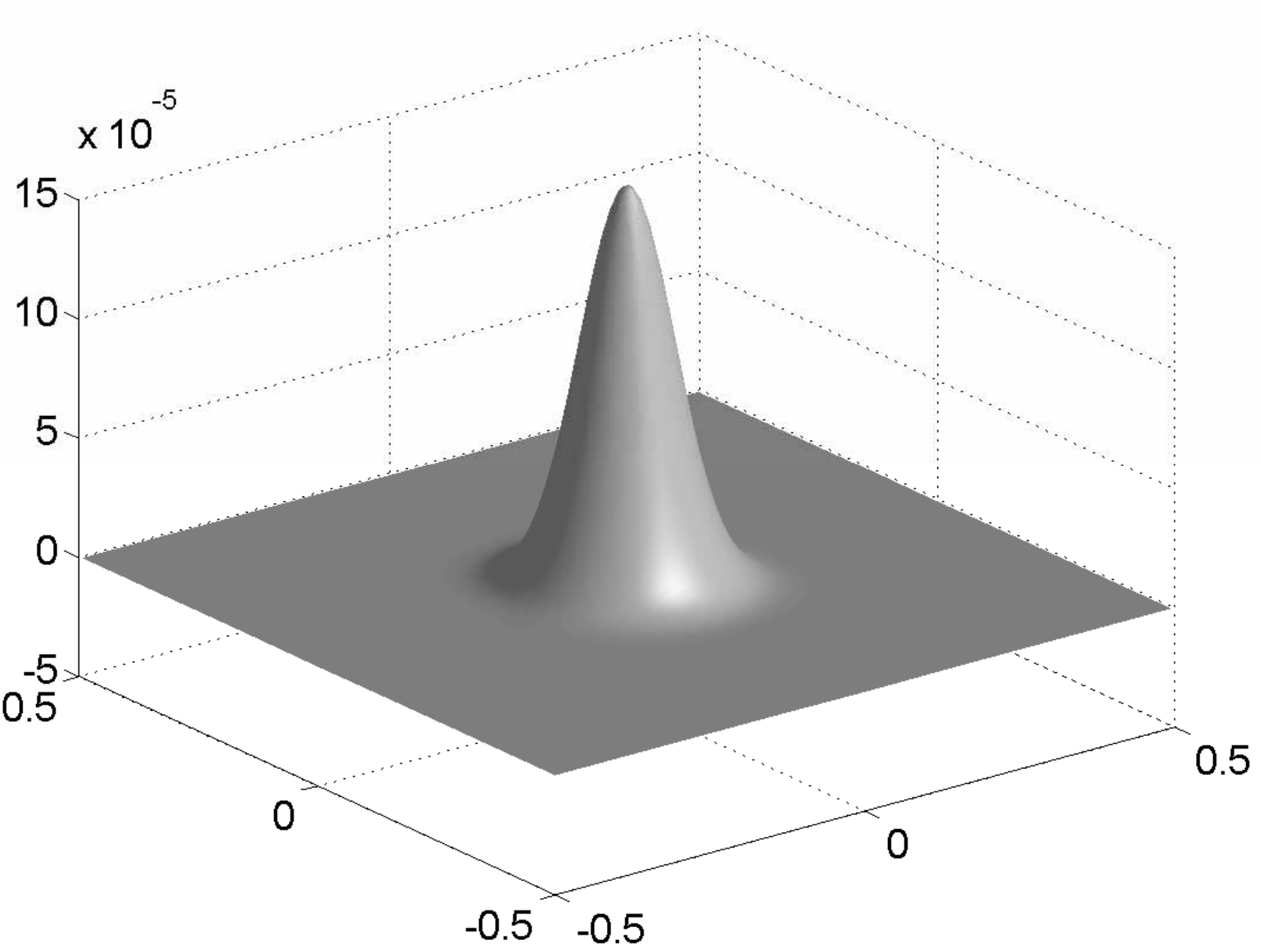}}
\resizebox{2in}{!}{\includegraphics{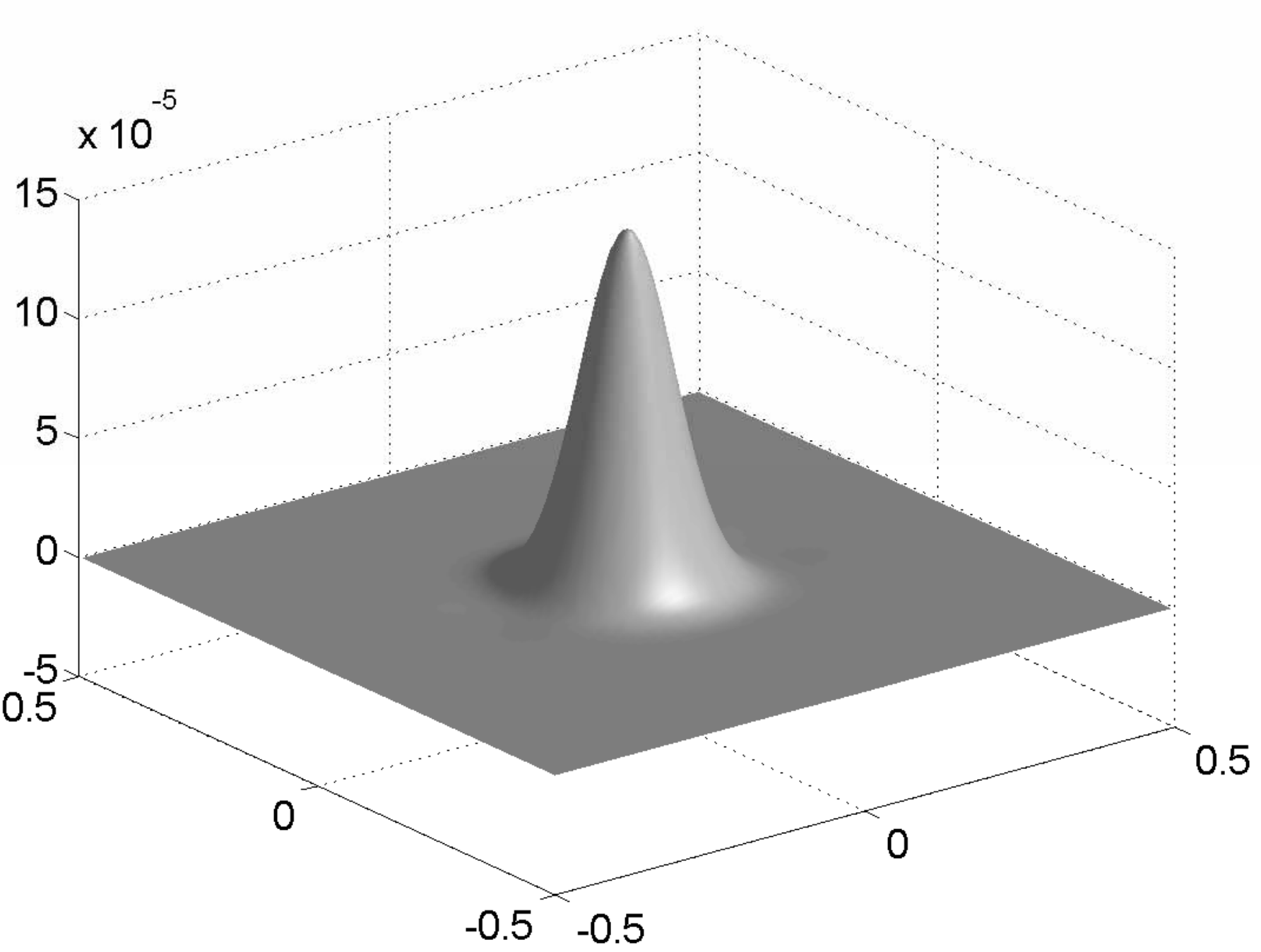}}

{$V^\e(t,\xb)|_{x_3=0}$ and $V(t,\xb)|_{x_3=0}$  at $t=0.25$}\vspace*{1mm}

\resizebox{2in}{!}{\includegraphics{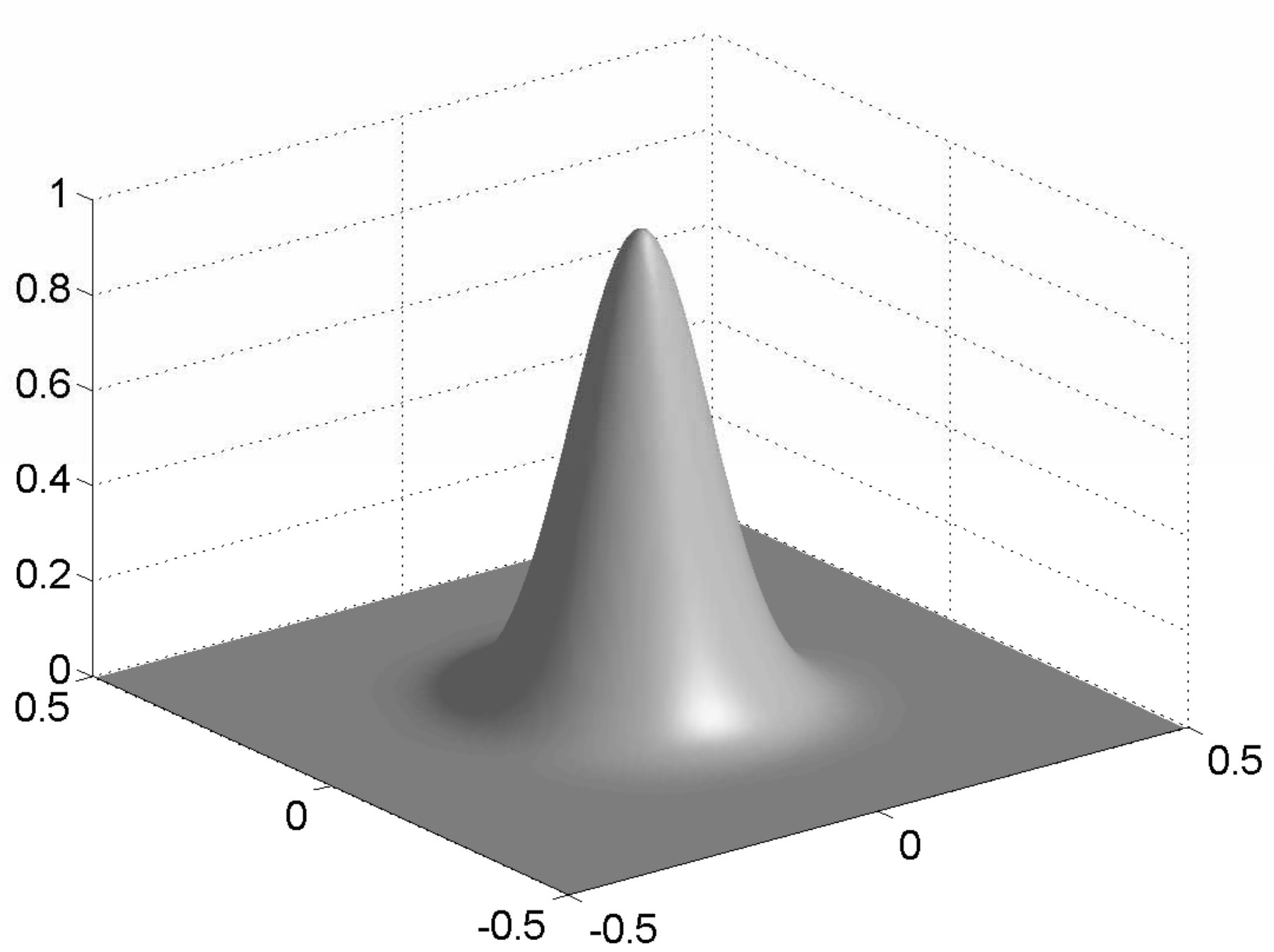}}
\resizebox{2in}{!}{\includegraphics{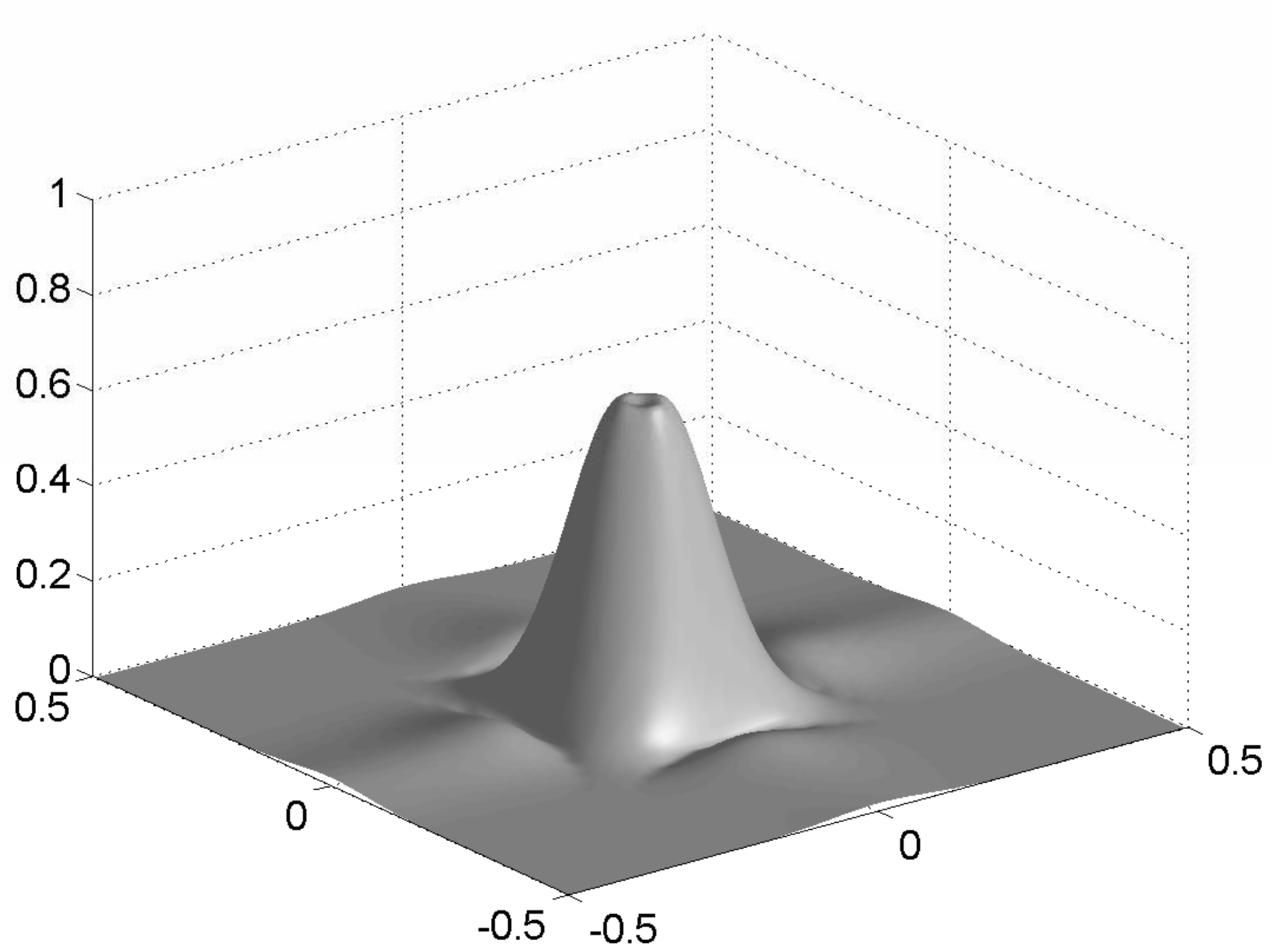}}

{$\big|\psi^\e(t,\xb)|_{x_3=0}\big|^2$ and $\left|\sum_\pm u^\pm(t,\xb)e^{i\phi^\pm/\e}|_{x_3=0}\right|^2$
at $t=0.375$}\vspace*{1mm}

\resizebox{2in}{!}{\includegraphics{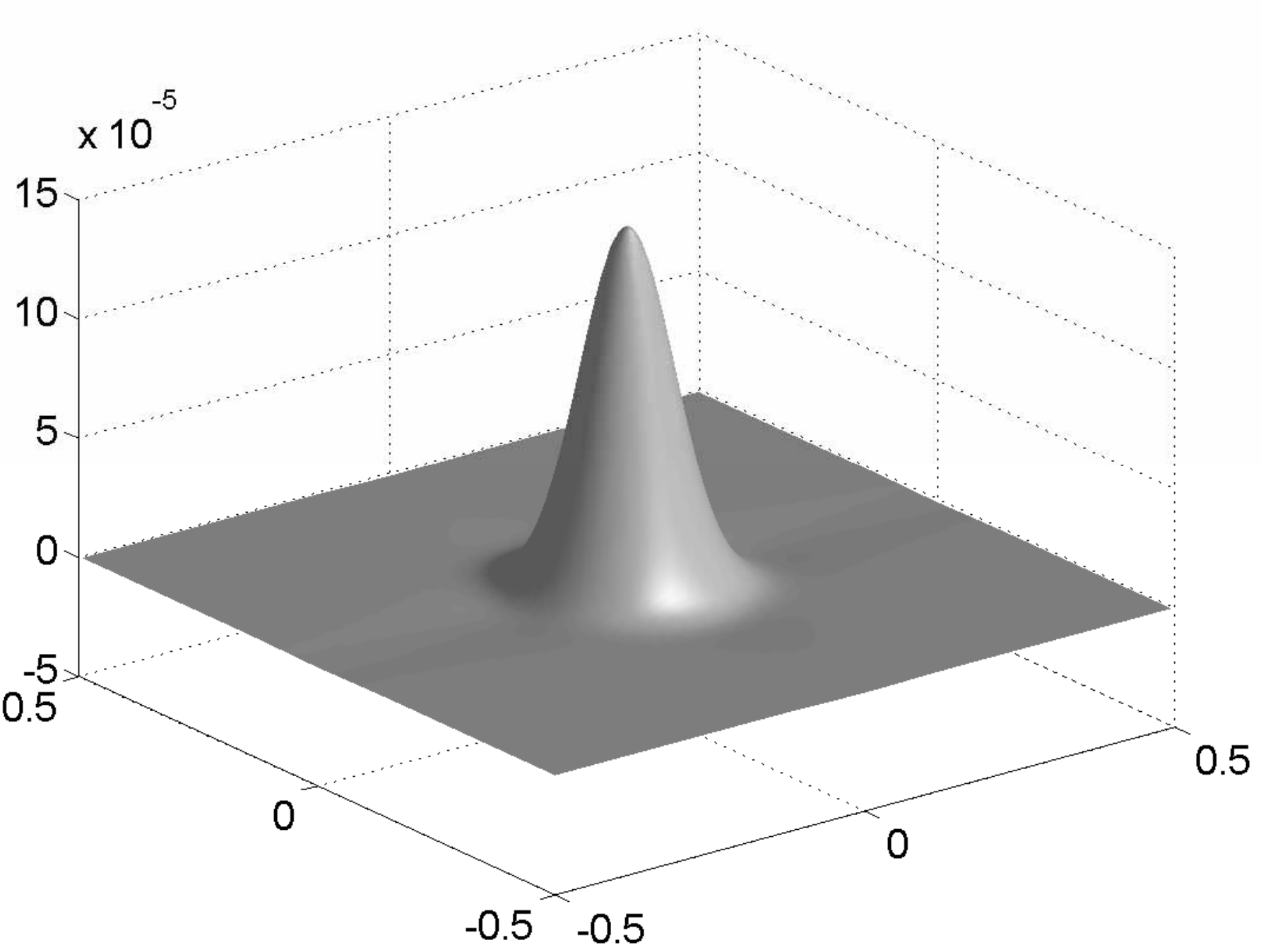}}
\resizebox{2in}{!}{\includegraphics{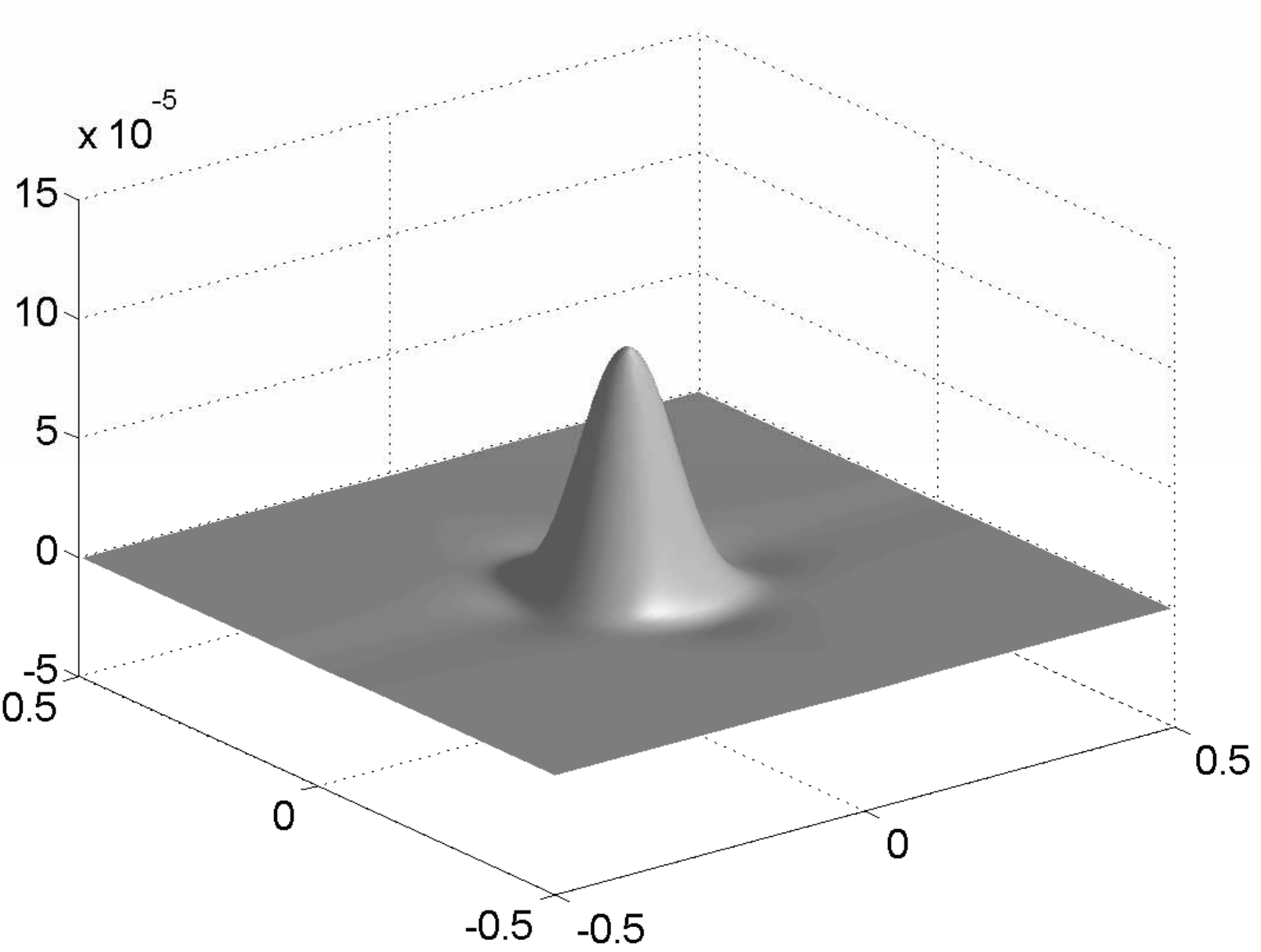}}

{$V^\e(t,\xb)|_{x_3=0}$ and $V(t,\xb)|_{x_3=0}$  at $t=0.375$}
\end{center}
\caption{Numerical results for example \ref{exsc2}.
The left column shows the graphs of the solution of the MD system,
the right column shows the graphs of the solution of the asymptotic problem.
Here  $\e=0.01$, $\tg t=\frac{1}{128}$, $\tg x=\frac{1}{32}$.}\label{fig311}
\end{figure}

\begin{figure} 
\begin{center}
\resizebox{2in}{!}{\includegraphics{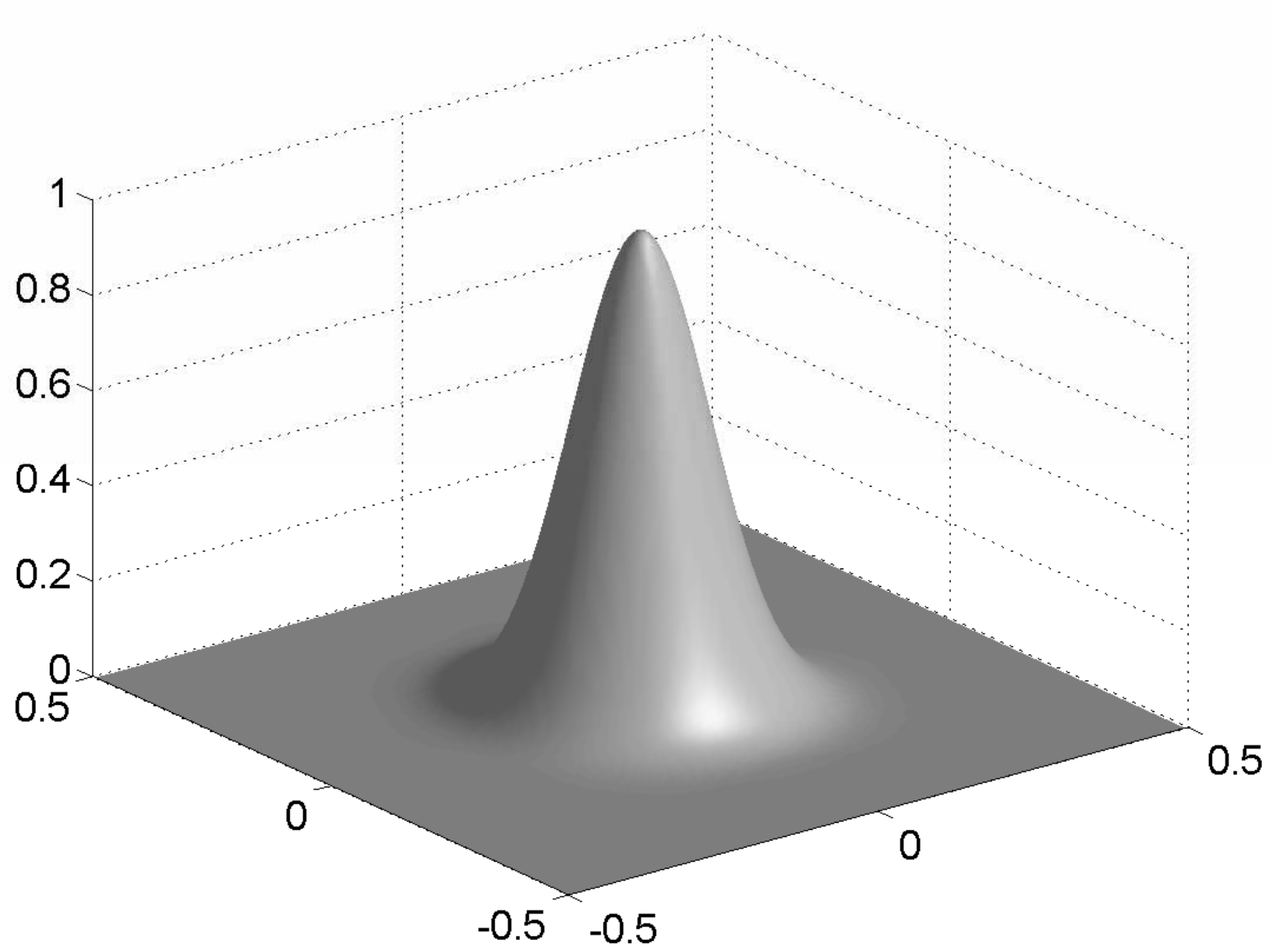}}
\resizebox{2in}{!}{\includegraphics{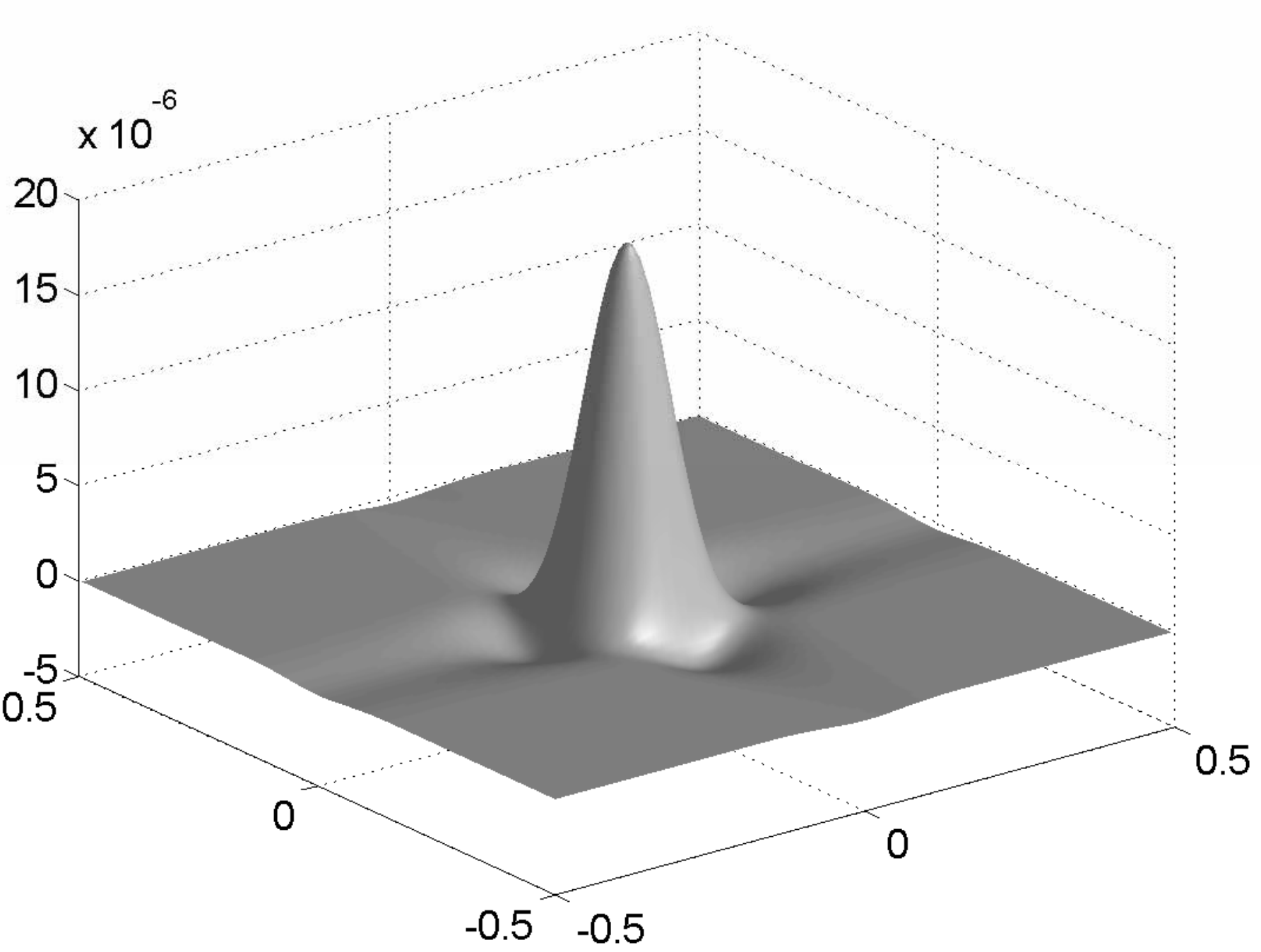}}

{$\big|\psi^\e(t,\xb)|_{x_3=0}\big|^2$ and $V^\e(t,\xb)|_{x_3=0}$
at $t=0.53$}\vspace*{2mm}

\resizebox{2in}{!}{\includegraphics{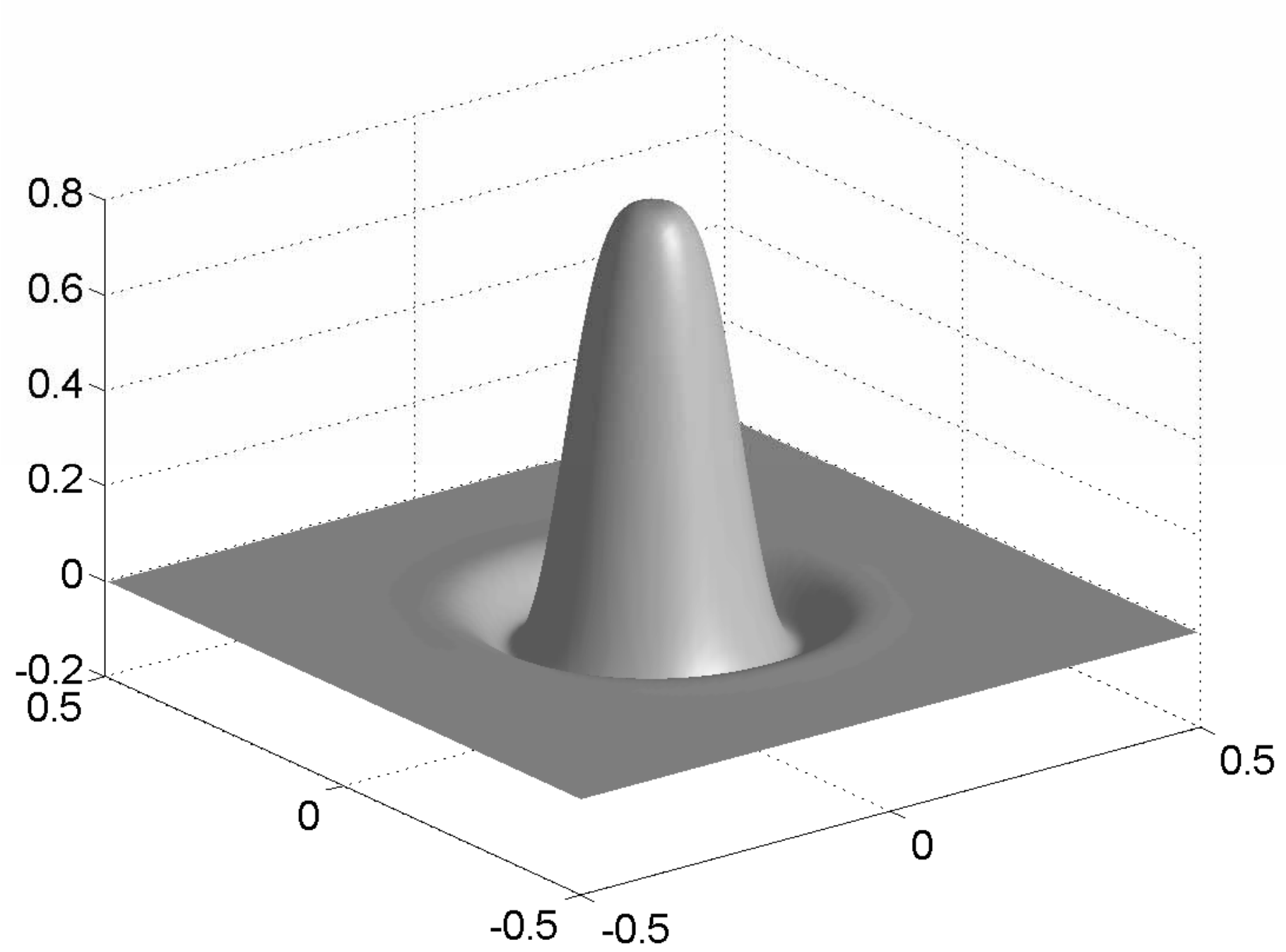}}
\resizebox{2in}{!}{\includegraphics{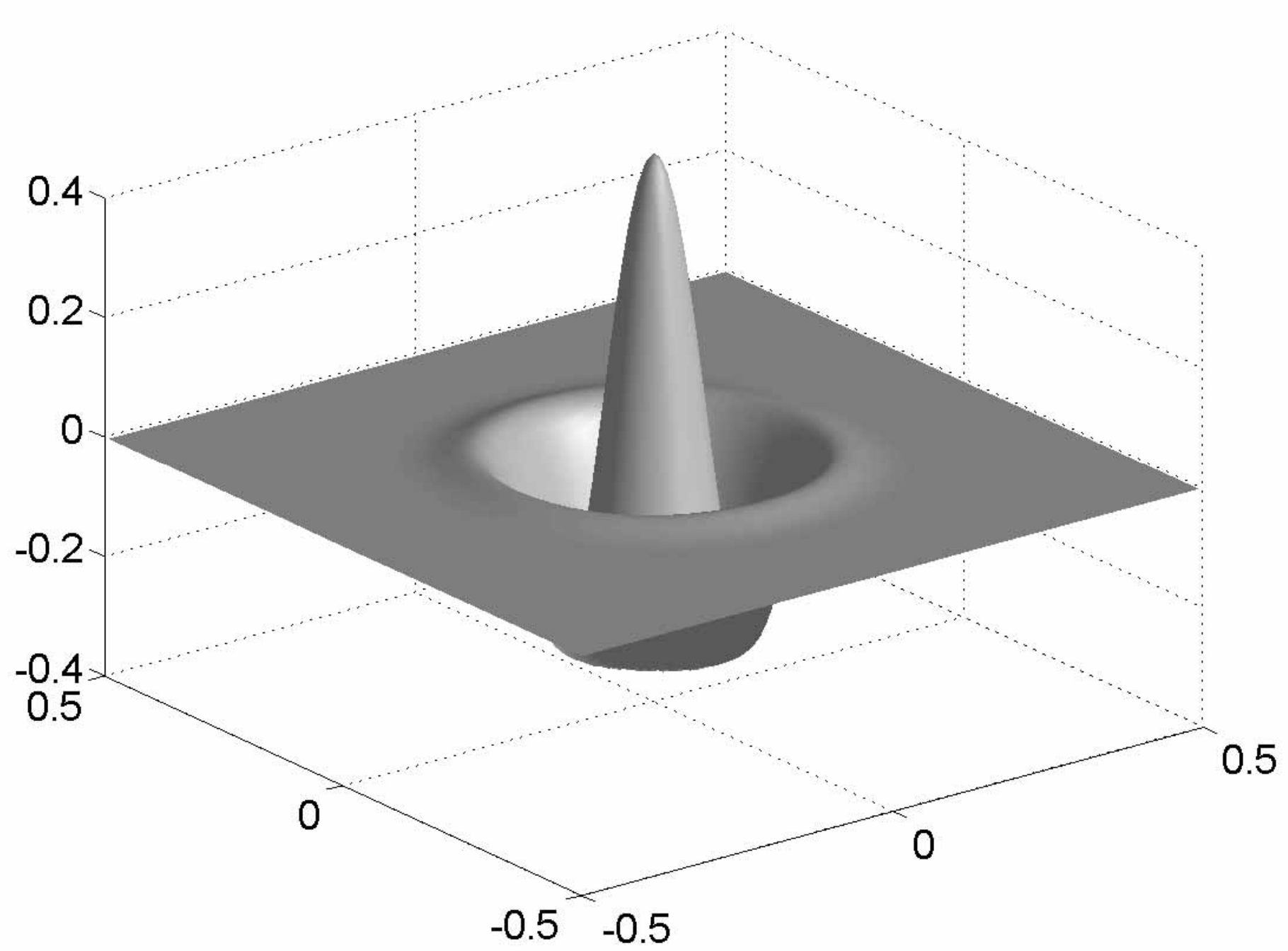}}

{Re$\big(\psi^\e_1(t,\xb)\big)|_{x_3=0}$
and Im$\big(\psi^\e_1(t,\xb)\big)|_{x_3=0}$  at $t=0.53$.}\vspace*{2mm}

\resizebox{2in}{!}{\includegraphics{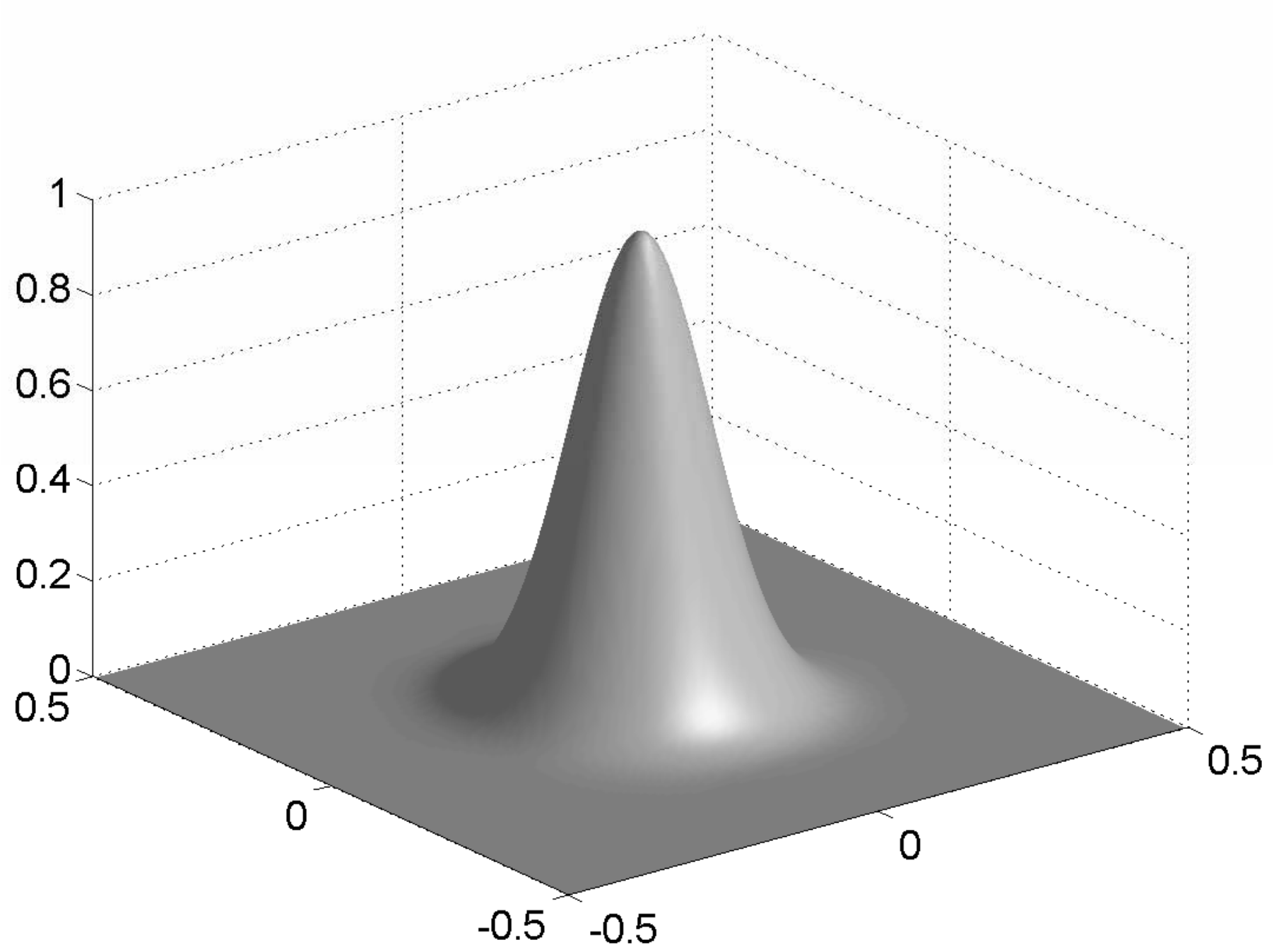}}
\resizebox{2in}{!}{\includegraphics{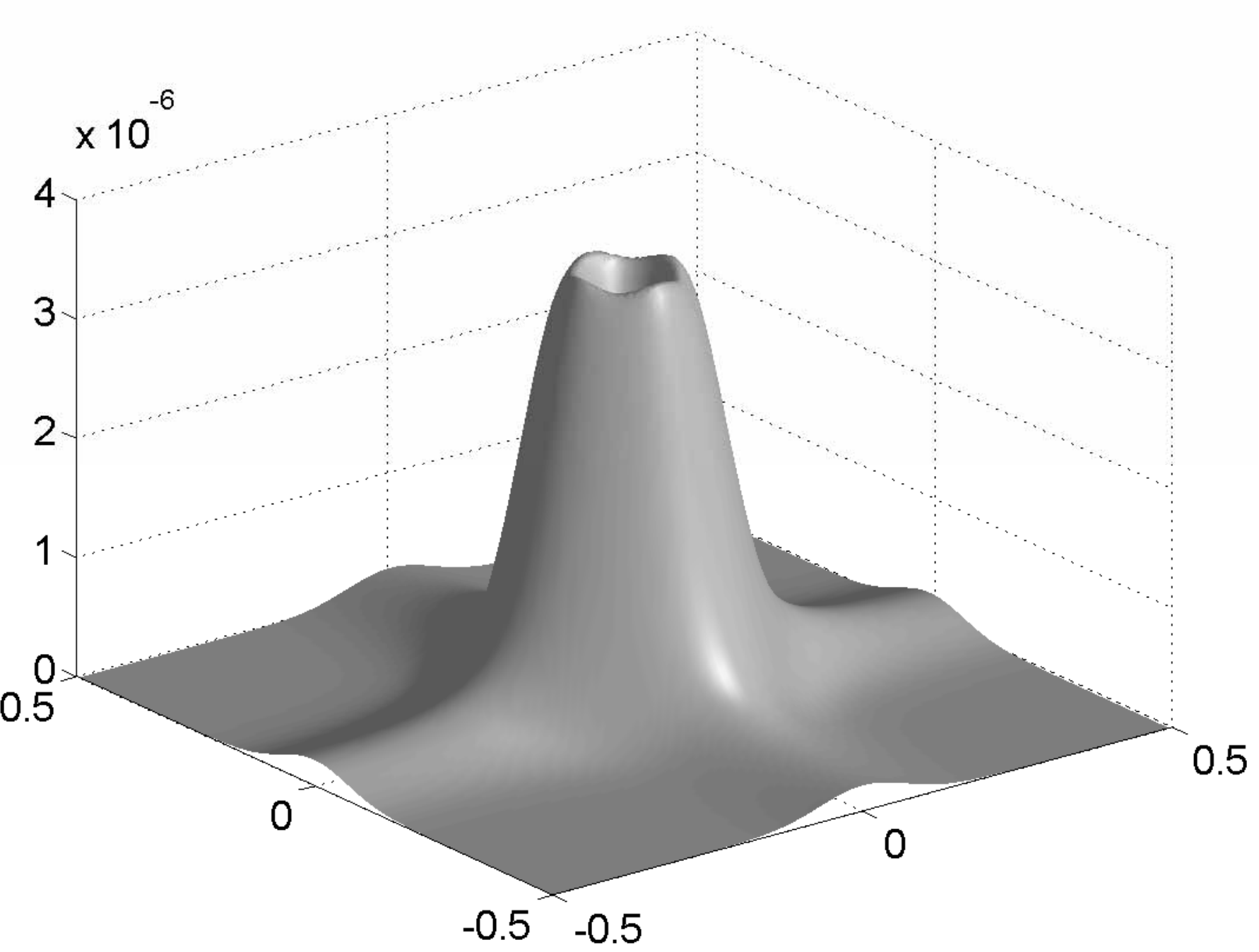}}

{$\big|\psi^\e(t,\xb)|_{x_3=0}\big|^2$ and $V^\e(t,\xb)|_{x_3=0}$ at $t=0.625$}\vspace*{2mm}

\resizebox{2in}{!} {\includegraphics{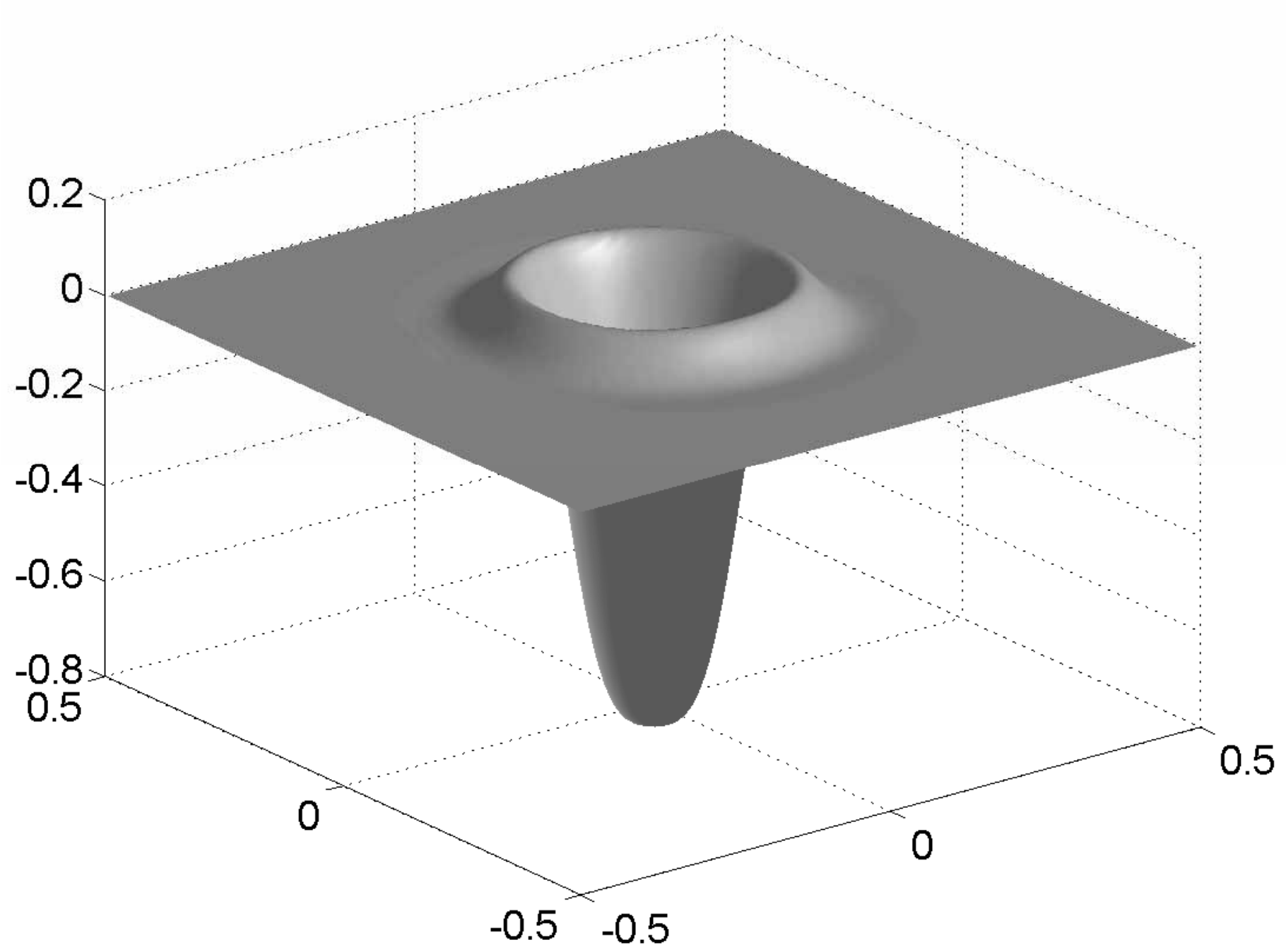}}
\resizebox{2in}{!} {\includegraphics{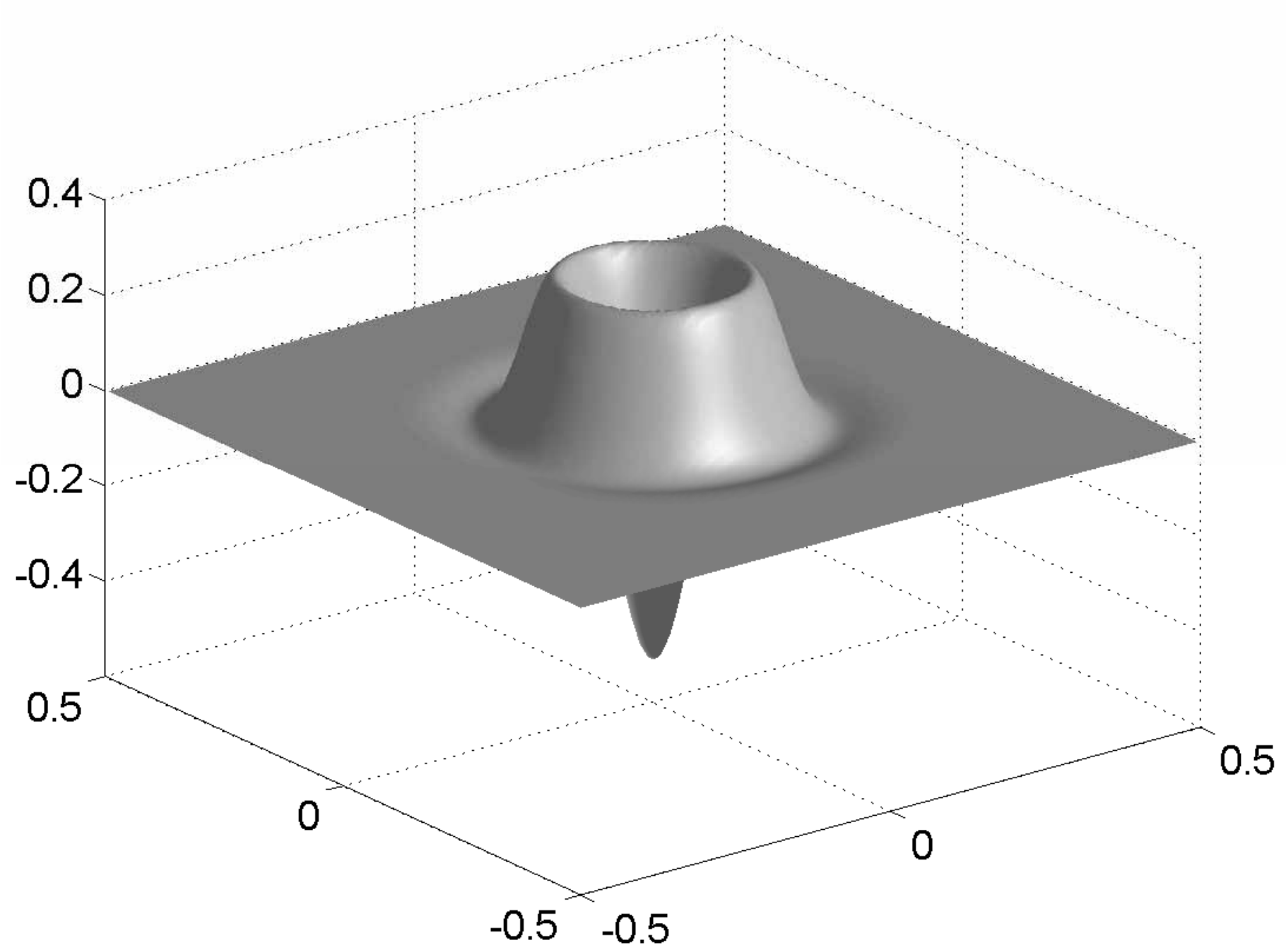}}

{Re$\big(\psi^\e_1(t,\xb)\big)|_{x_3=0}$ and Im$\big(\psi^\e_1(t,\xb)\big)|_{x_3=0}$  at $t=0.625$.}
\end{center}
\caption{Numerical results of the MD system for example \ref{exsc2}.
Here  $\e=0.01$, $\tg t=\frac{1}{128}$, $\tg x=\frac{1}{32}$.}\label{fig312}
\end{figure}

\begin{figure} 
\begin{center}
\resizebox{2in}{!} {\includegraphics{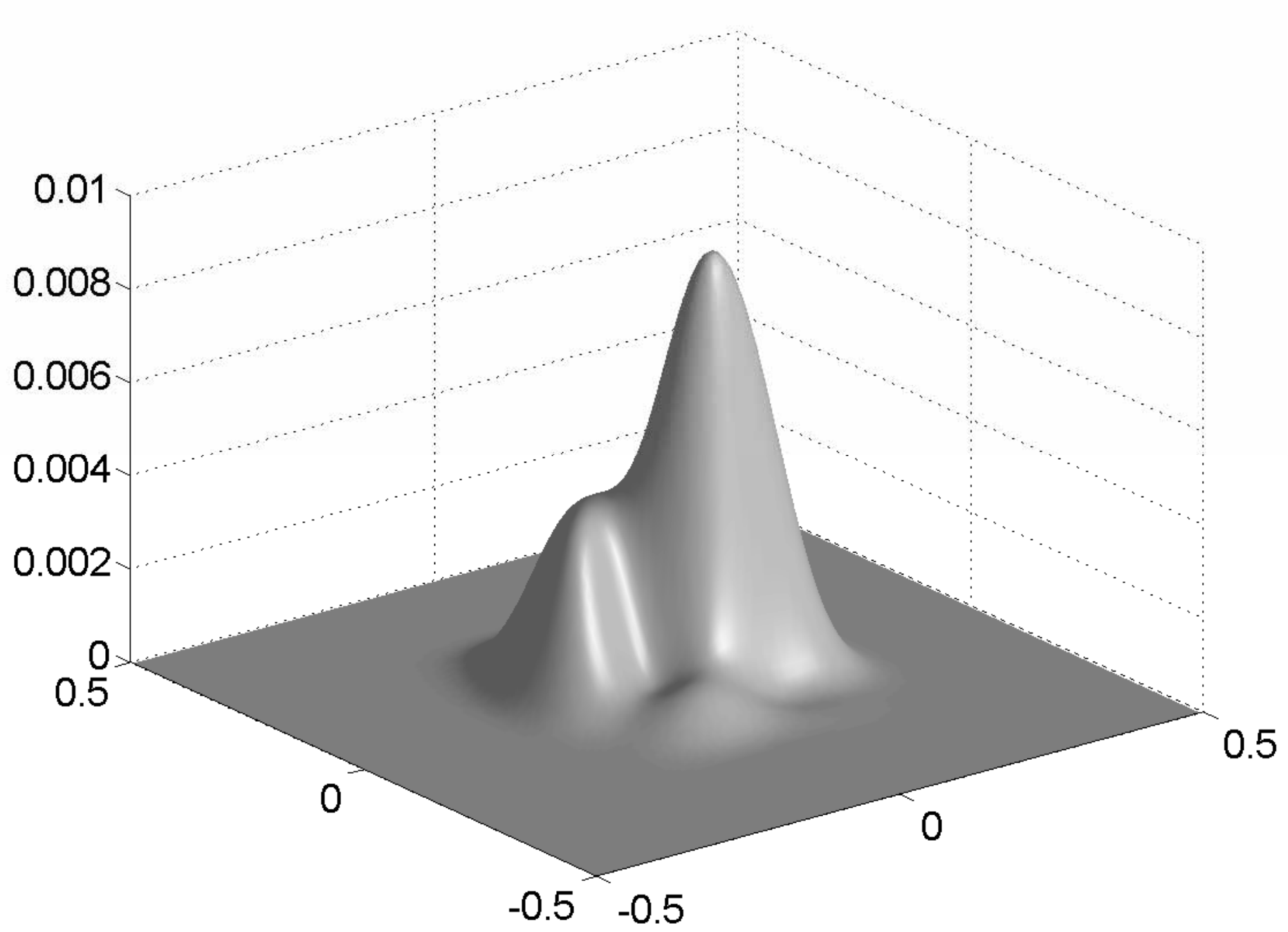}}
\resizebox{2in}{!} {\includegraphics{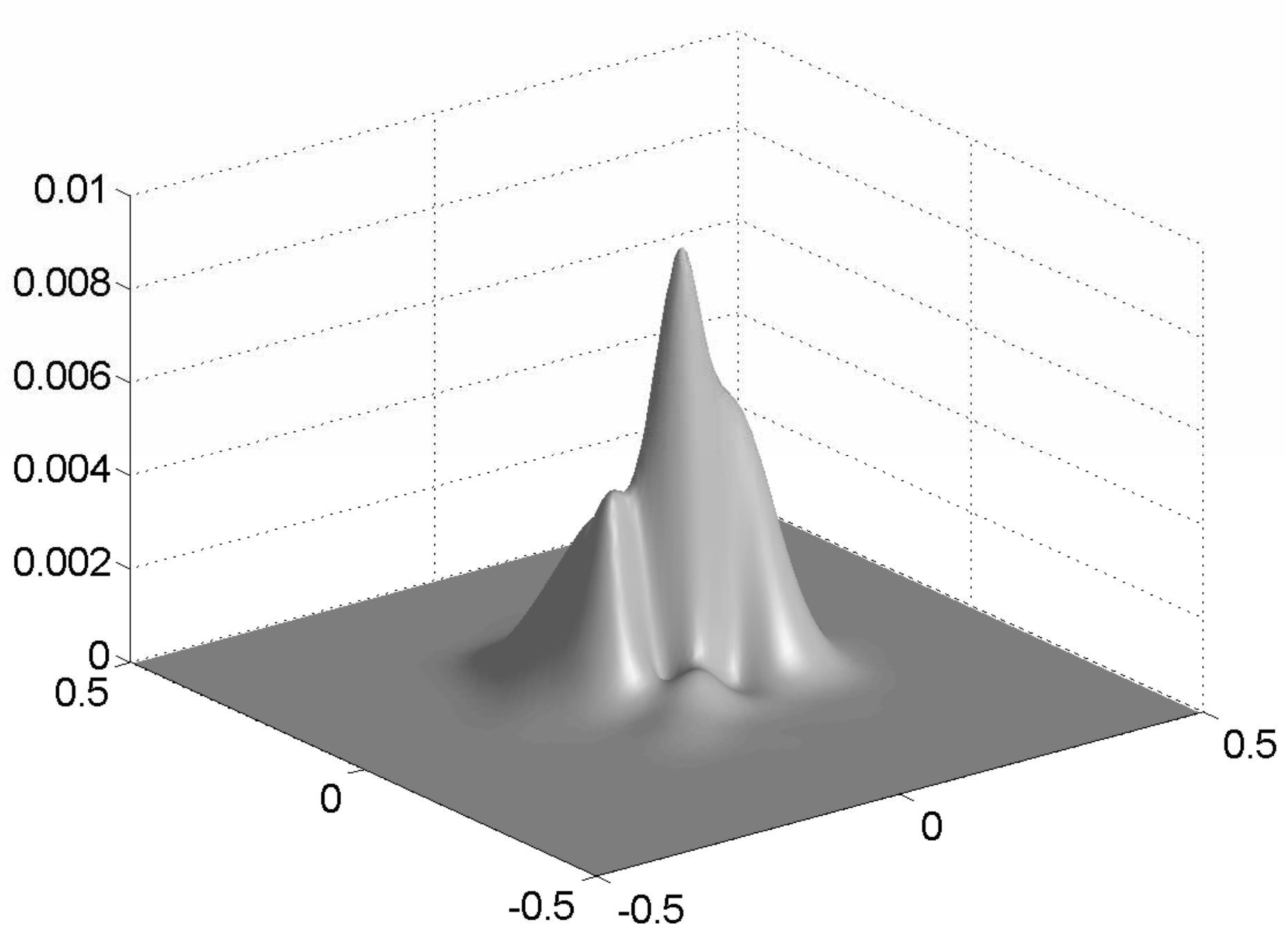}}

{$t=0.25$, $\big|\Pi^-_0(-i\e\btd)\psi^\e(t,\xb)|_{x_3=0}\big|^2$ and
$\big|\Pi^-_0(\btd\phi)\psi^\e(t,\xb)|_{x_3=0}\big|^2$}\vspace*{2mm}

\resizebox{2in}{!} {\includegraphics{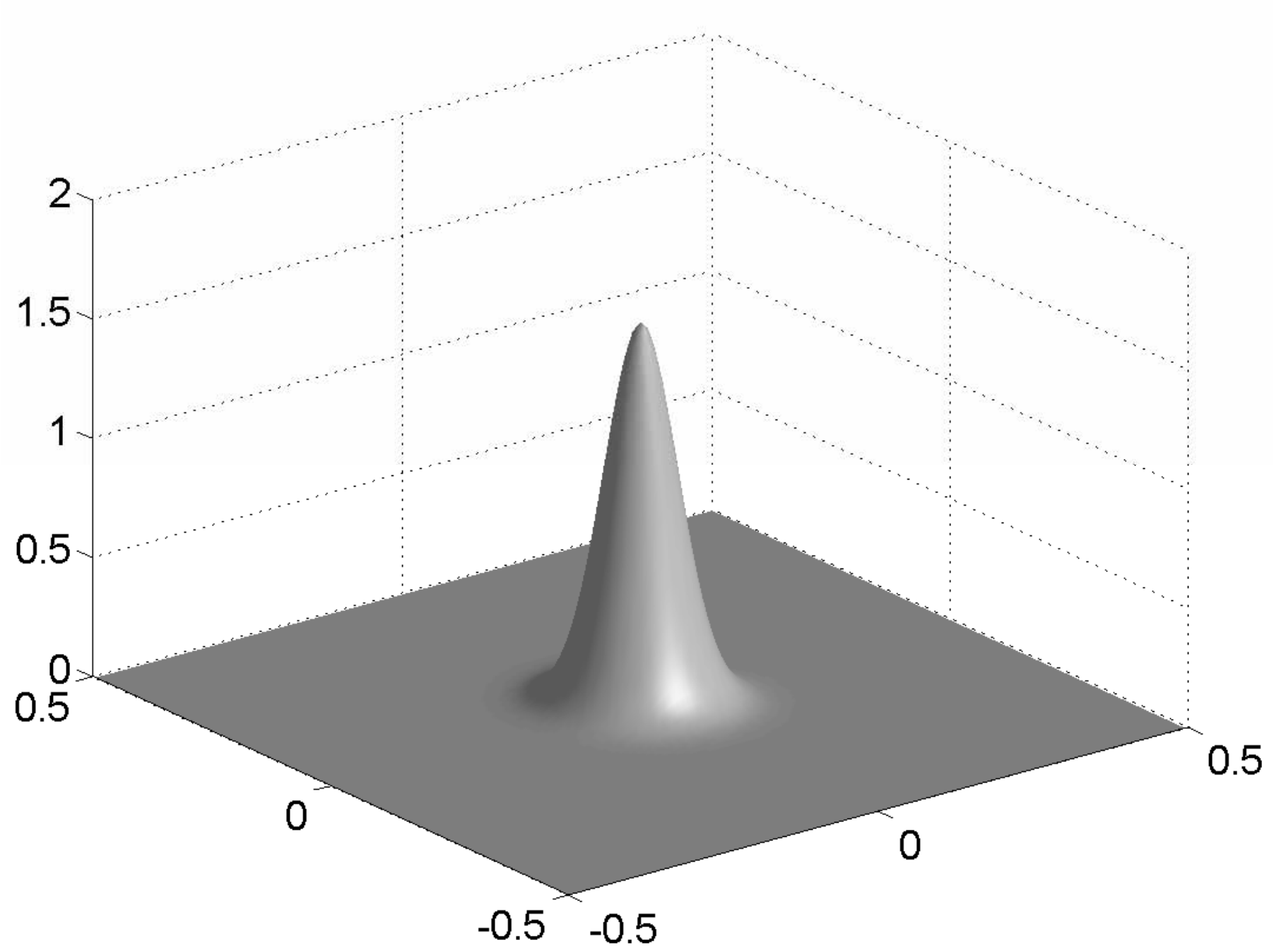}}
\resizebox{2in}{!} {\includegraphics{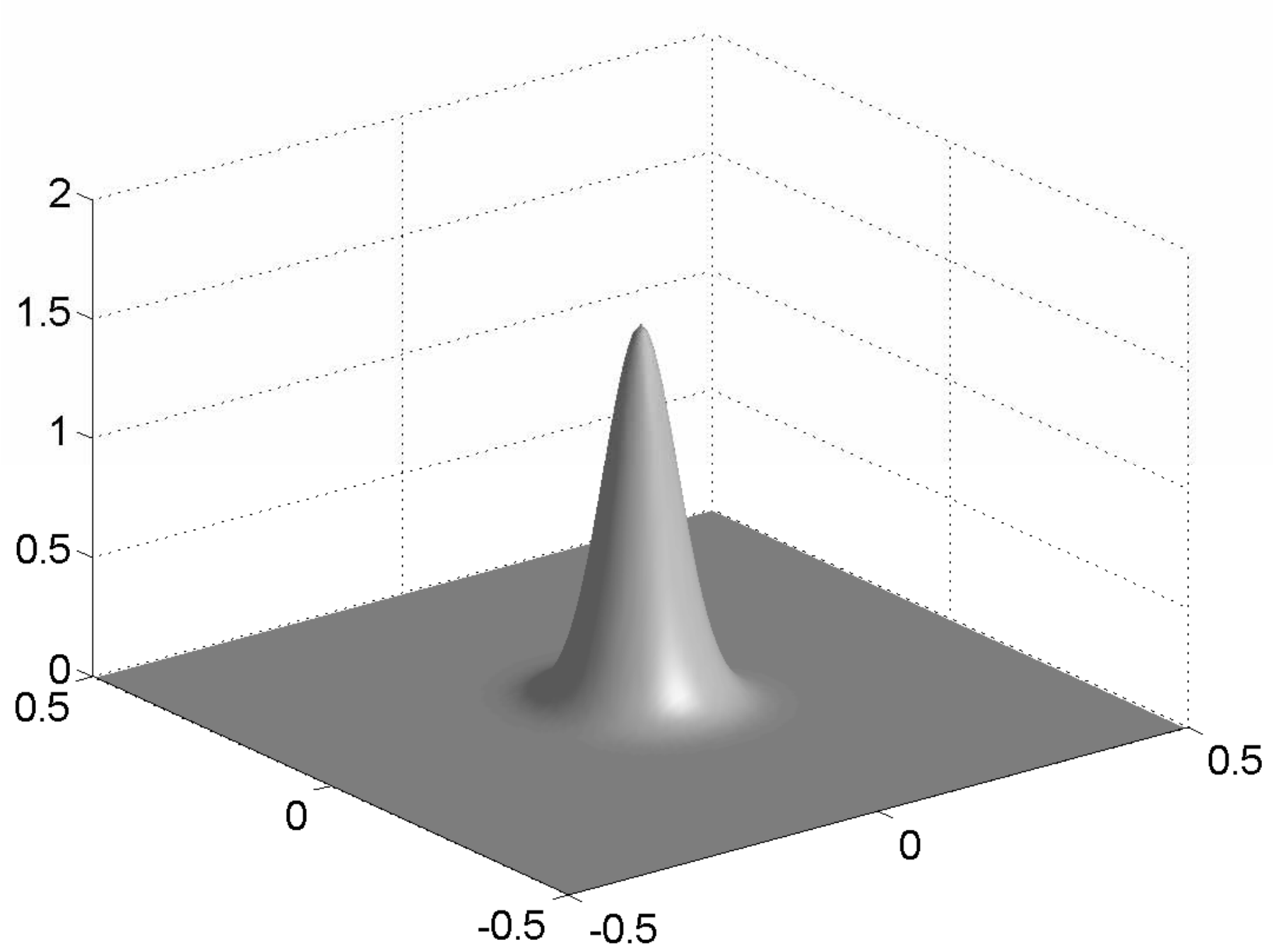}}

{$t=0.25$, $\big|\Pi^+_0(-i\e\btd)\psi^\e(t,\xb)|_{x_3=0}\big|^2$ and
$\big|\Pi^+_0(\btd\phi)\psi^\e(t,\xb)|_{x_3=0}\big|^2$}\vspace*{2mm}

\resizebox{2in}{!}{\includegraphics{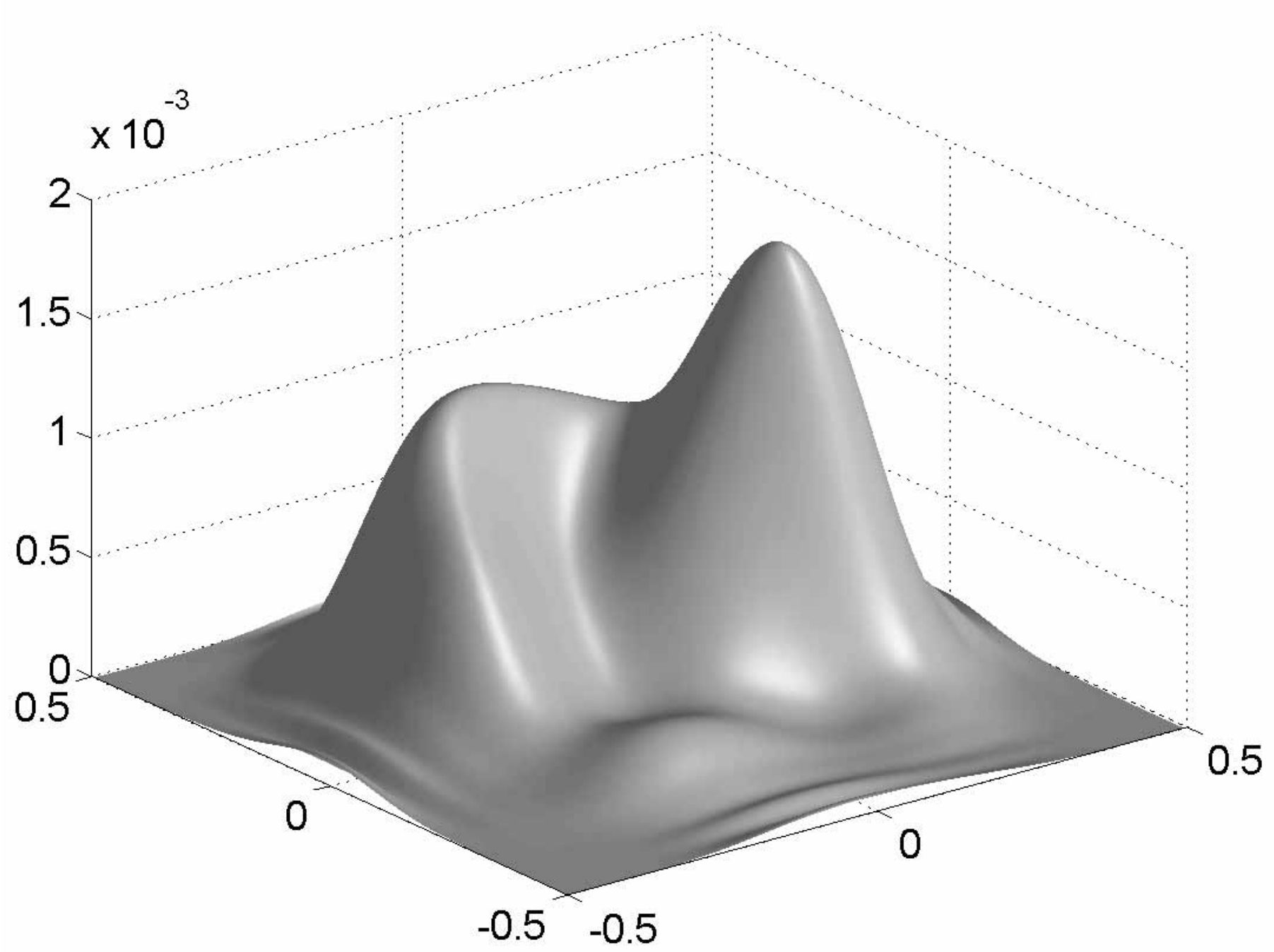}}
\resizebox{2in}{!}{\includegraphics{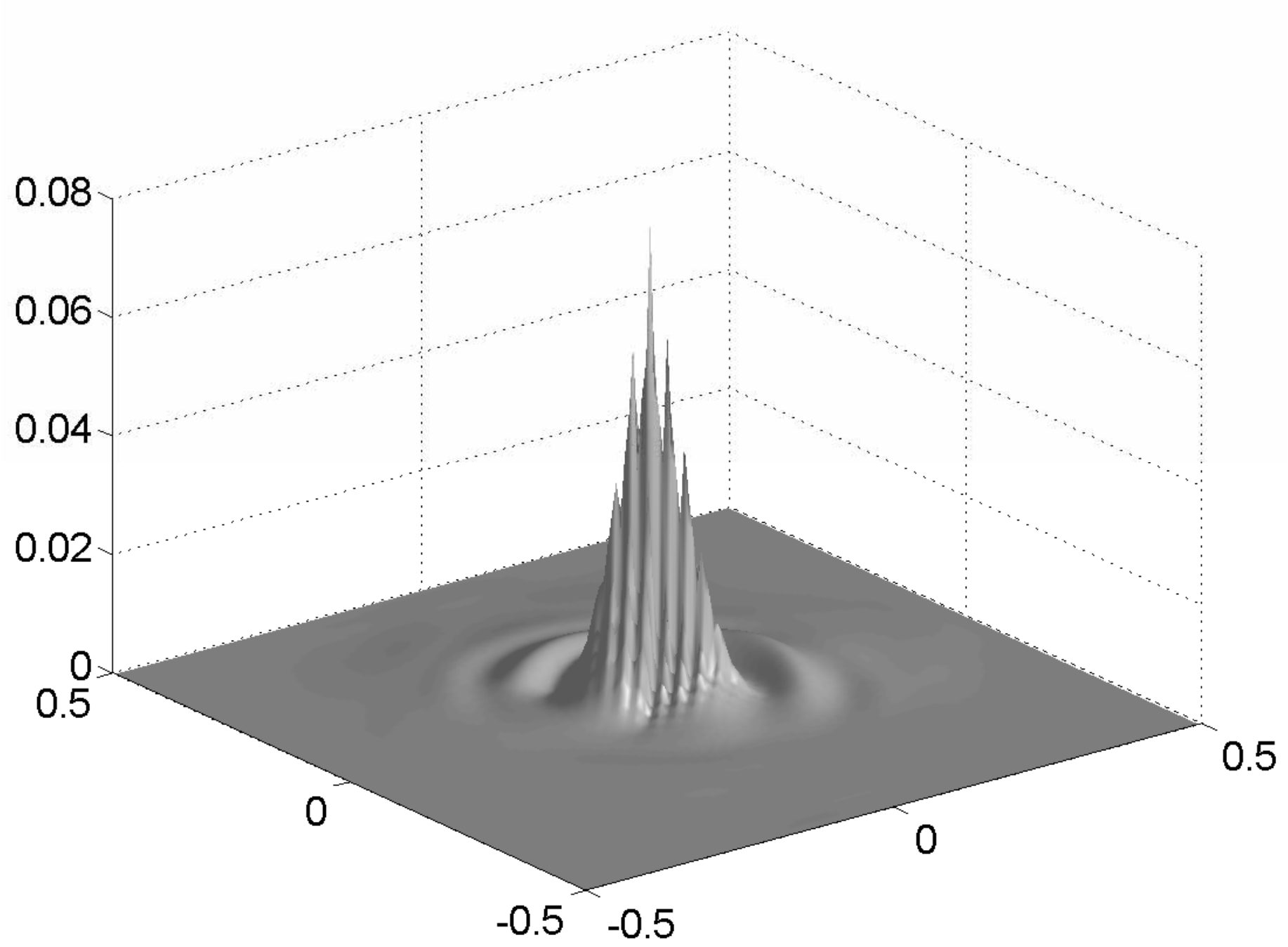}}

{$t=0.75$, $\big|\Pi^-_0(-i\e\btd)\psi^\e(t,\xb)|_{x_3=0}\big|^2$ and
$\big|\Pi^-_0(\btd\phi)\psi^\e(t,\xb)|_{x_3=0}\big|^2$}\vspace*{2mm}

\resizebox{2in}{!} {\includegraphics{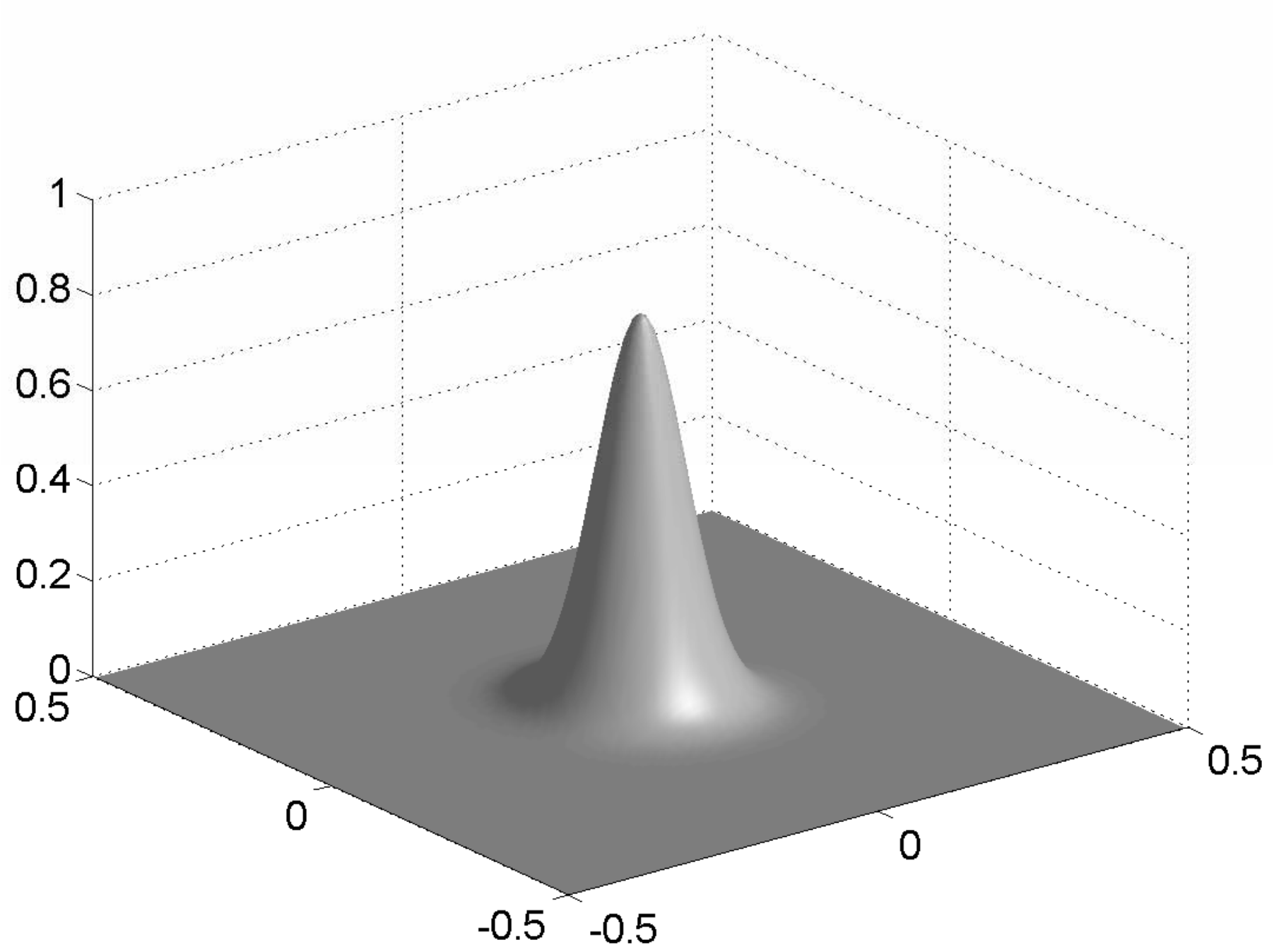}}
\resizebox{2in}{!} {\includegraphics{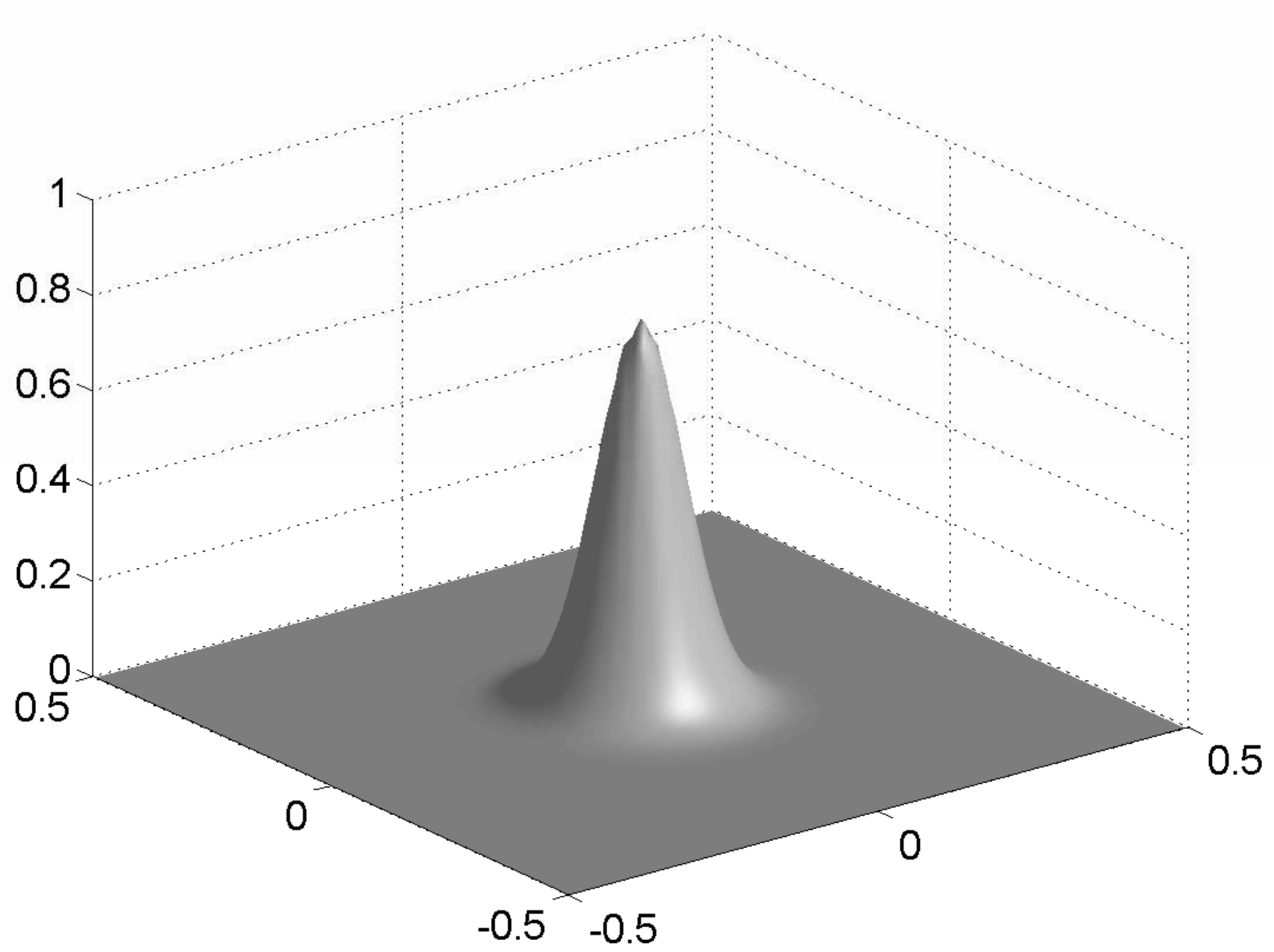}}

{$t=0.75$, $\big|\Pi^+_0(-i\e\btd)\psi^\e(t,\xb)|_{x_3=0}\big|^2$ and
$\big|\Pi^+_0(\btd\phi)\psi^\e(t,\xb)|_{x_3=0}\big|^2$}
\end{center}
\caption{Numerical results  of the densities of electron/positron projectors for example \ref{exsc2}.
The left column is
$\big|\Pi^\pm_0(-i\e\btd)\psi^\e(t,\xb)|_{x_3=0}\big|^2$,
the right column is $\big|\Pi^\pm_0(\btd\phi)\psi^\e(t)|_{x_3=0}\big|^2$.
The graphs show that the
matrices $\Pi^-_0(\btd\phi)$ do not mimic $\Pi^-_0(-i\e\btd)$ after the caustic point.
Here $\e=0.01$, $\tg x=1/32$, $\tg t=1/128$.}
\label{fig34}
\end{figure}
\begin{figure} 
\begin{center}
\resizebox{1.6in}{!} {\includegraphics{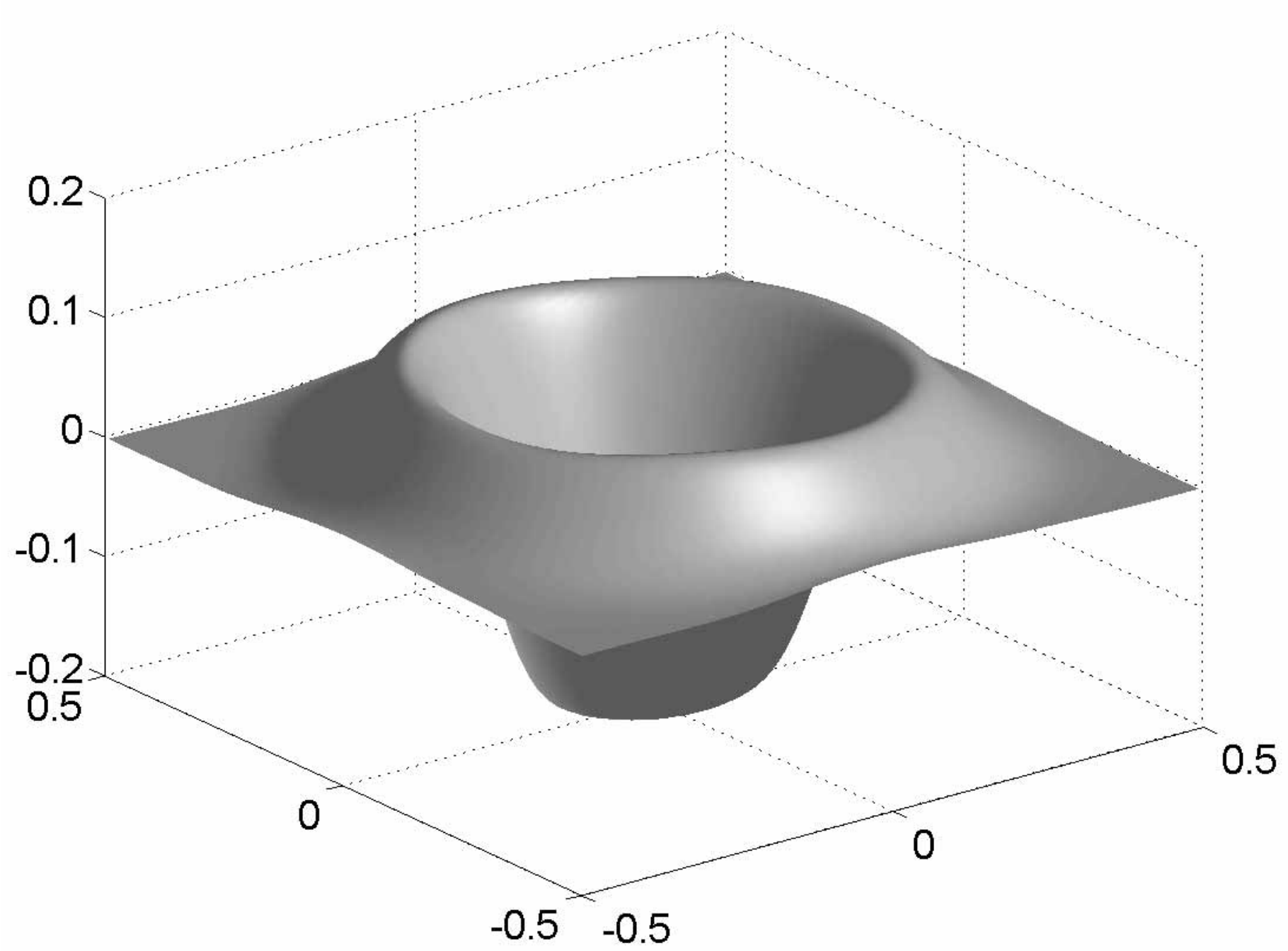}}
\resizebox{1.6in}{!} {\includegraphics{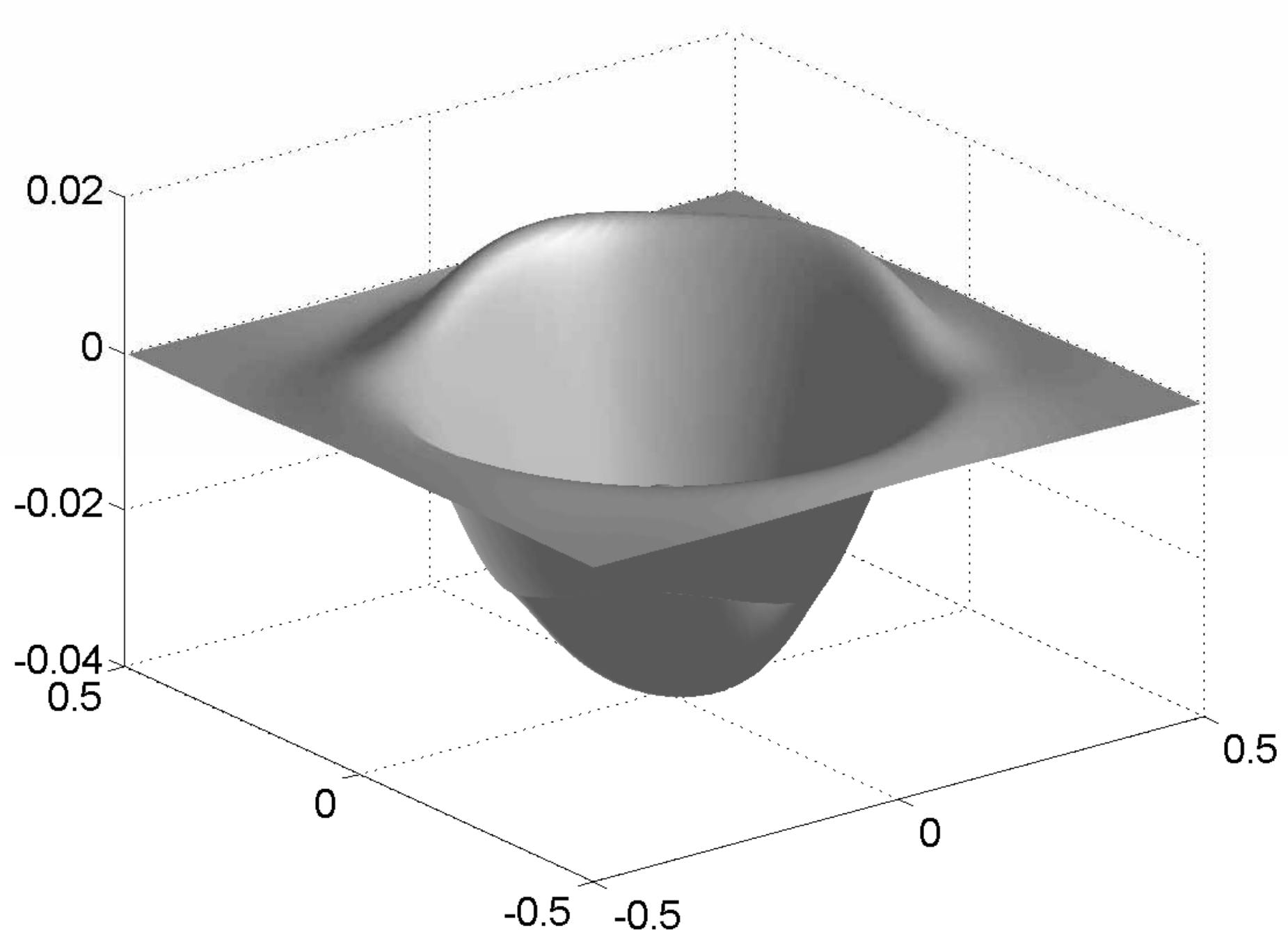}}
\resizebox{1.6in}{!}{\includegraphics{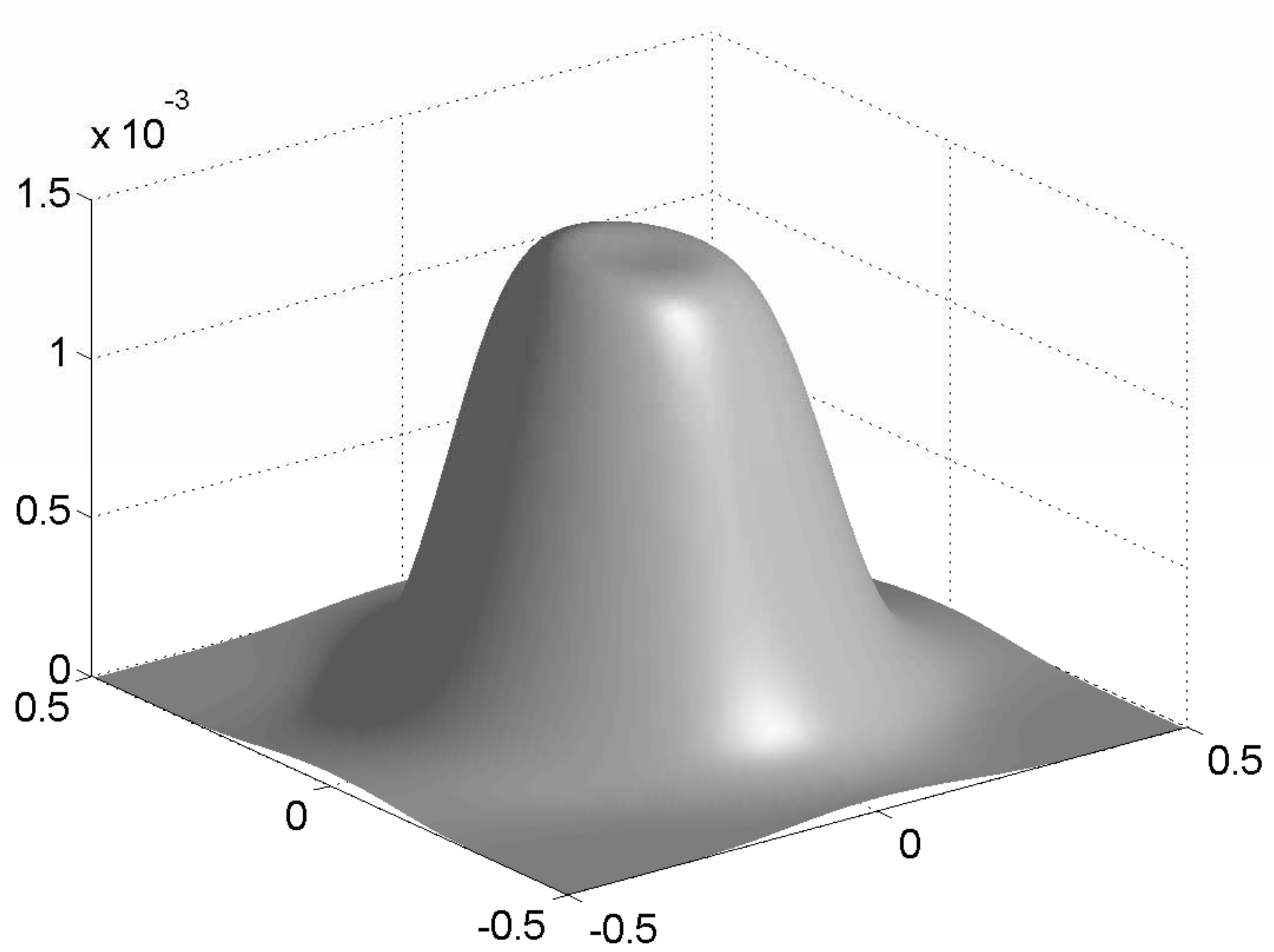}}

{Re$\big(\psi^\e_1(t,\xb)\big)|_{x_3=0}$,
Im$\big(\psi^\e_1(t,\xb)\big)|_{x_3=0}$ and $V(t,\xb)|_{x_3=0}$
at $t=0.25$.\vspace*{3mm}}

\resizebox{1.6in}{!} {\includegraphics{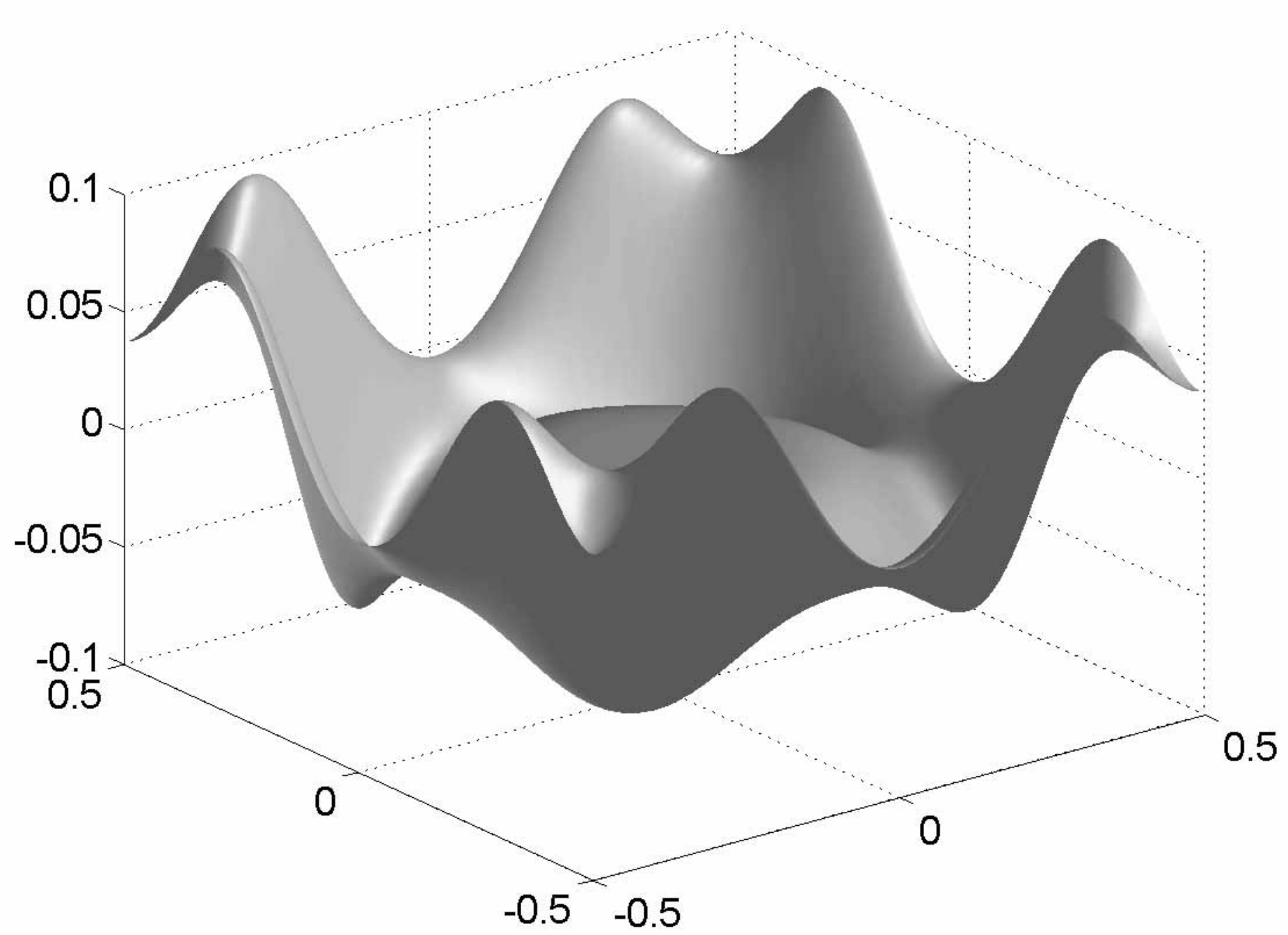}}
\resizebox{1.6in}{!} {\includegraphics{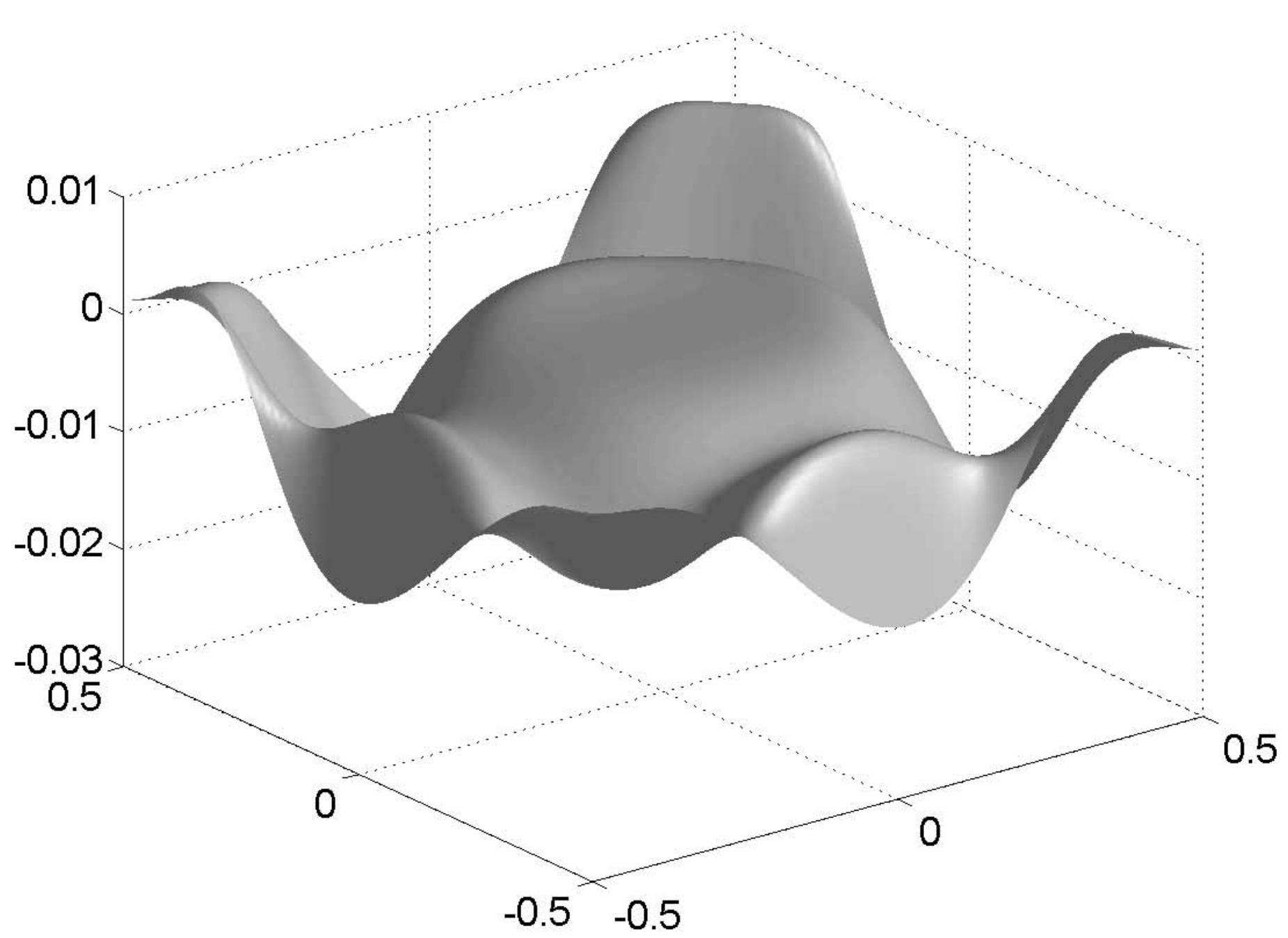}}
\resizebox{1.6in}{!}{\includegraphics{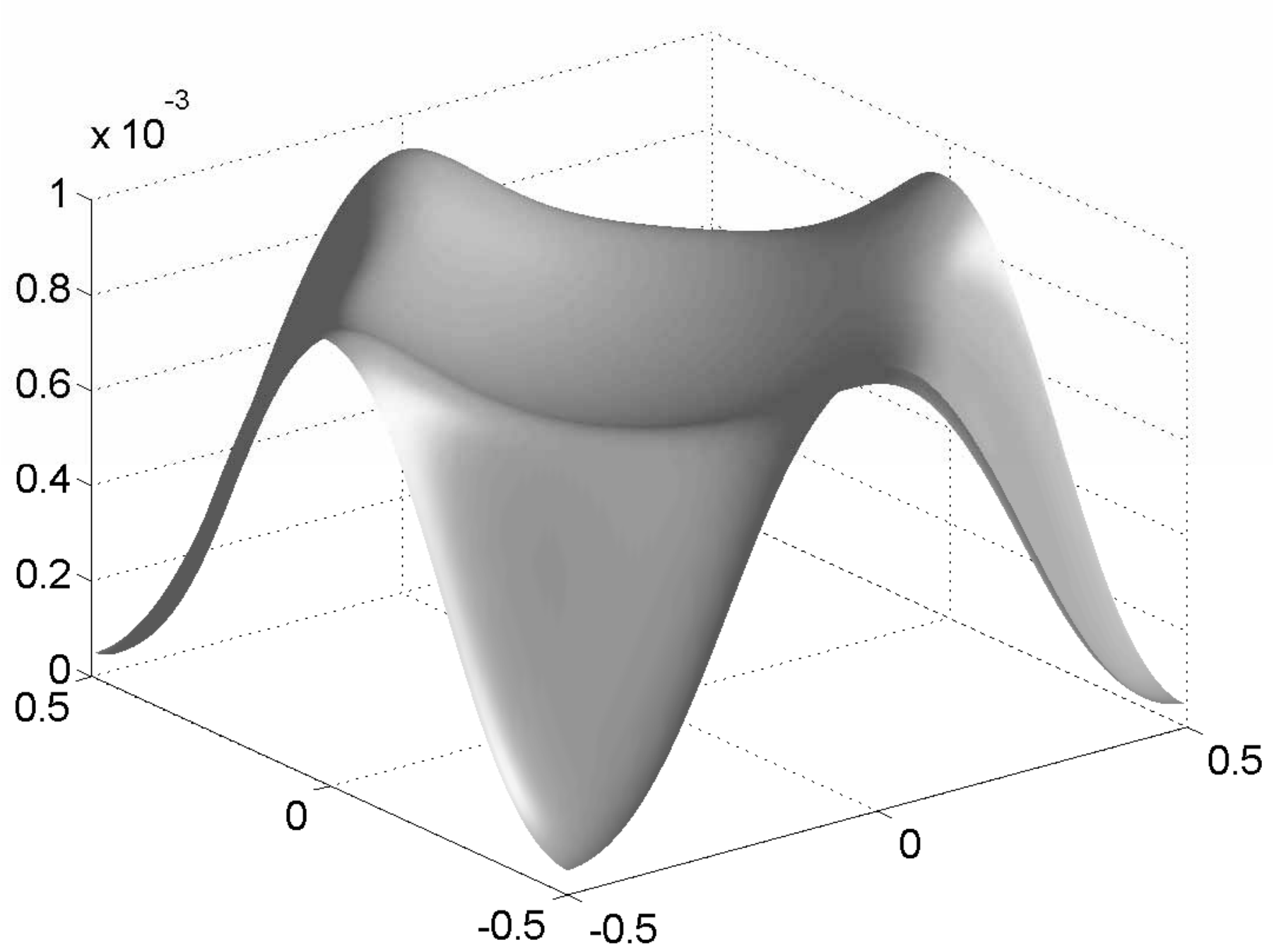}}

{Re$\big(\psi^\e_1(t,\xb)\big)|_{x_3=0}$,
Im$\big(\psi^\e_1(t,\xb)\big)|_{x_3=0}$ and $V(t,\xb)|_{x_3=0}$
at $t=0.5$.\vspace*{3mm}}

\resizebox{1.6in}{!} {\includegraphics{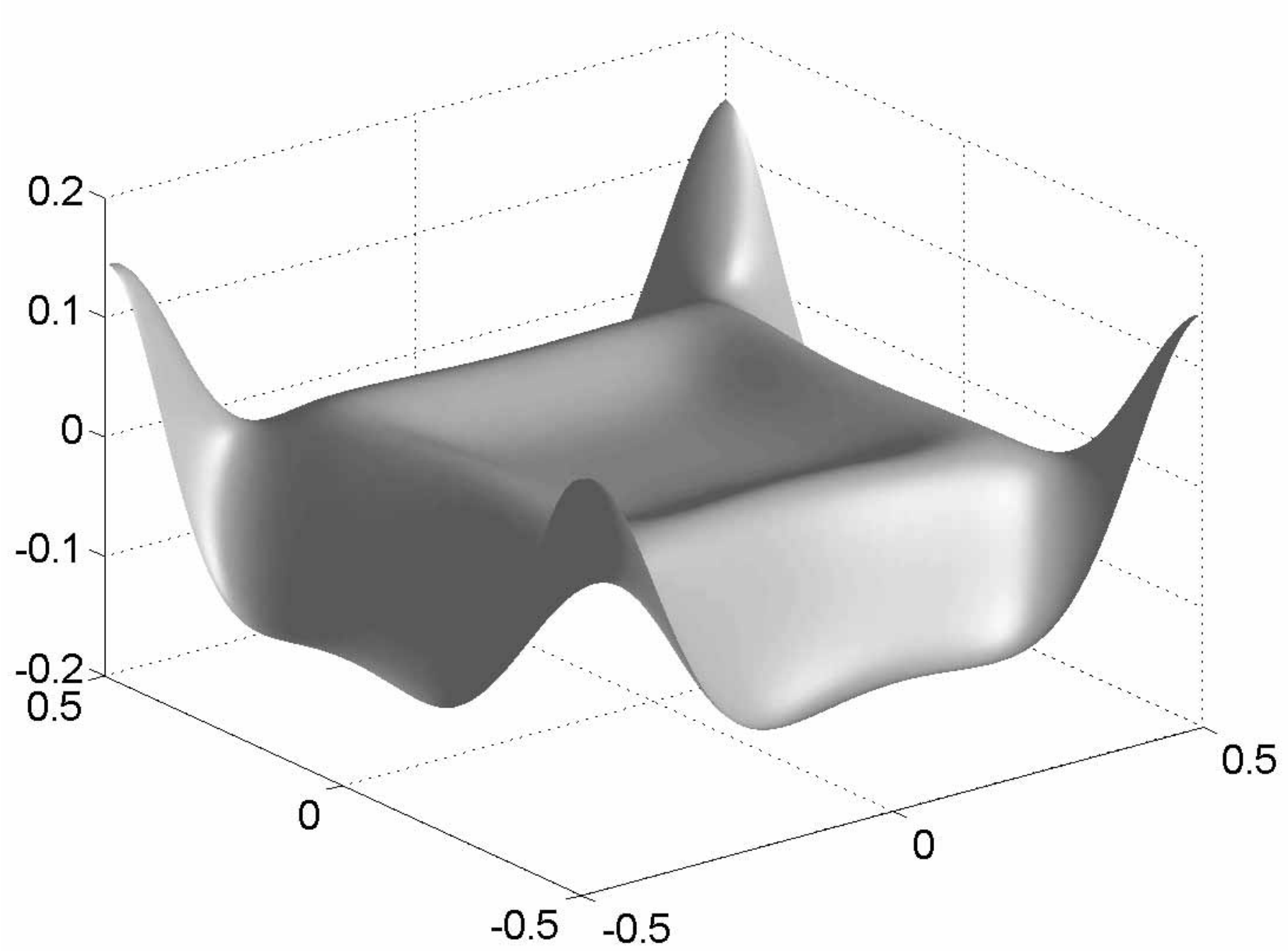}}
\resizebox{1.6in}{!} {\includegraphics{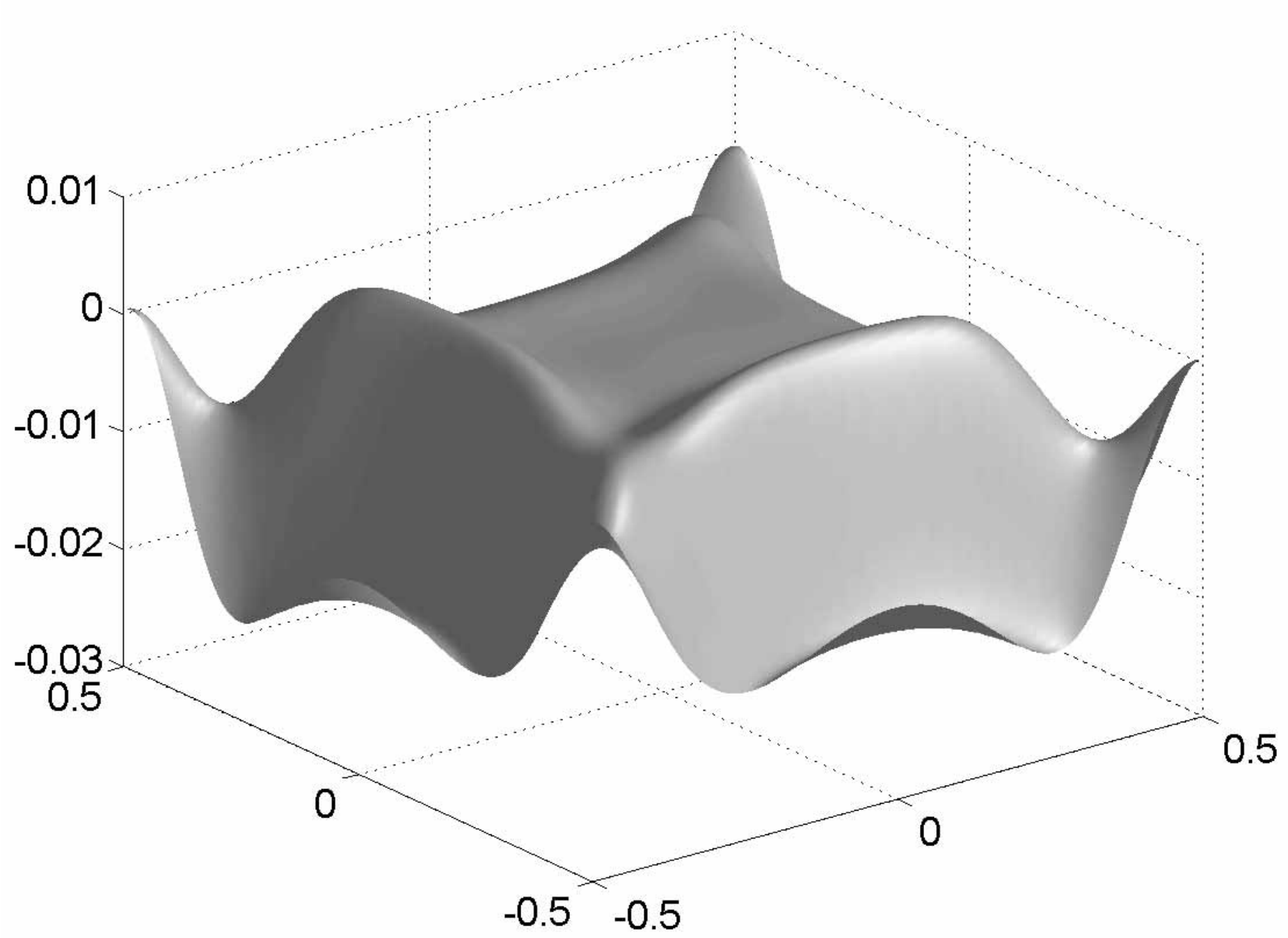}}
\resizebox{1.6in}{!}{\includegraphics{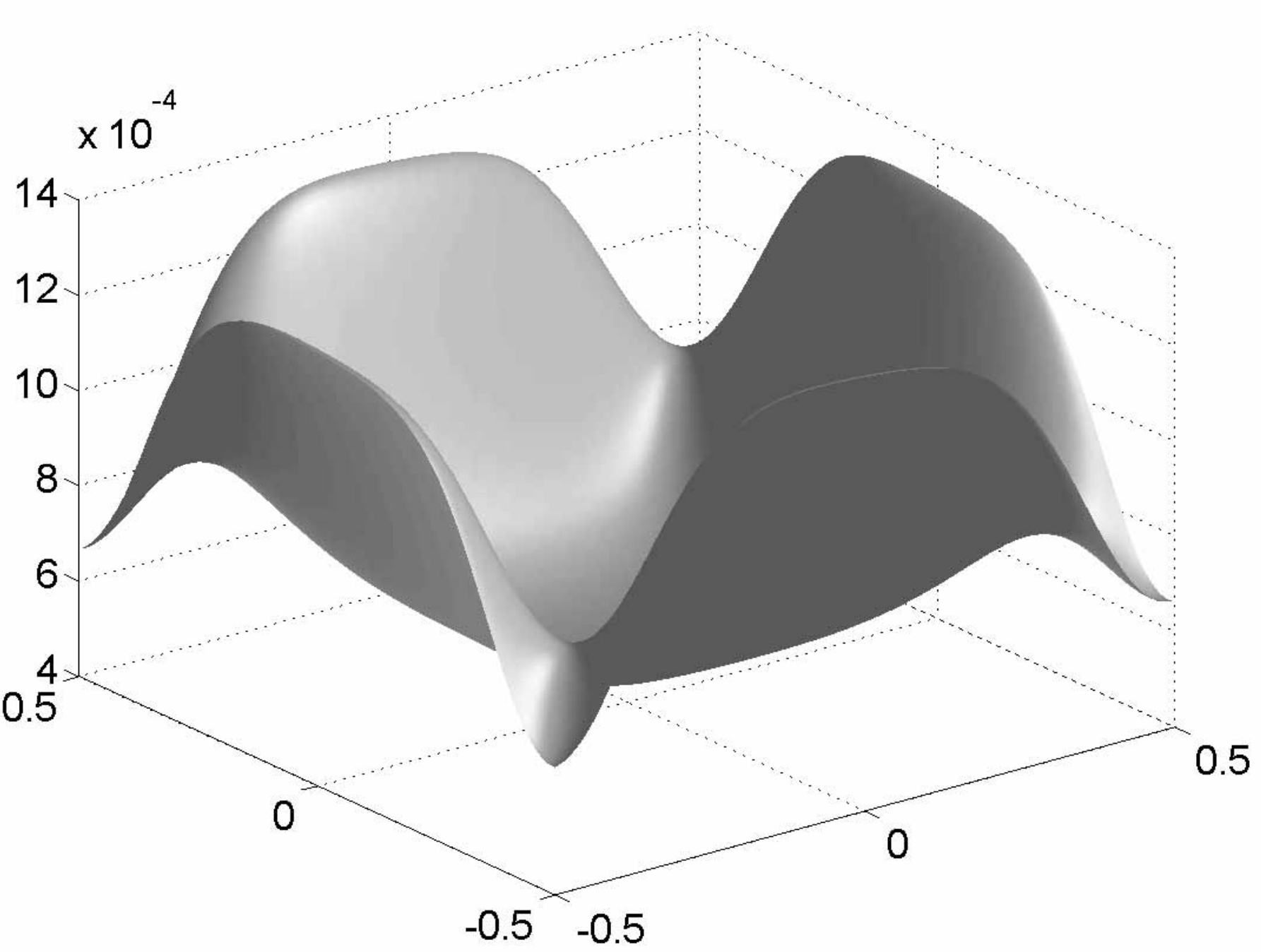}}

{Re$\big(\psi^\e_1(t,\xb)\big)|_{x_3=0}$,
Im$\big(\psi^\e_1(t,\xb)\big)|_{x_3=0}$ and $V(t,\xb)|_{x_3=0}$
at $t=0.625$.}
\end{center}
\caption{Numerical results of the MD system for example \ref{exsc2}.
Here $\e=1.0$, $\tg t=\f{1}{128}$, $\tg x=\f{1}{32}$. } \label{fig33}
\end{figure}
\end{example}
\newpage
\begin{example}[\textbf{Harmonic oscillator}]\label{exha}
Finally, we take $\Ab^{ex}(\xb)=0$ and include a confining electric
potential of harmonic oscillator type, \ie $V^{ex}(\xb)=|\xb|^2$. Hence $\phi^\pm$ satisfies
\begin{equation}
\label{eicnr}
\partial_t \phi^\pm (t,\xb)+ \sqrt{|\nabla \phi^\pm|^2+1}+|\xb|^2=0, \quad \phi^\pm (0,\xb)=\phi_I (\xb),
\end{equation}
which implies $\omega_{\Ab}^\pm(\xi)=\omega^\pm_0(\xi)$ in this case.
Due to the presence of the external potential, the semi-classical transport equations
\eqref{trans} have to be generalized by including a spin-transport term, \cf \cite{Sp},
which however only enters in the phase of $u^\pm$. Thus the conservation law for the densities $\rho^\pm$
is the same as in \eqref{conl}.
\newpar
Let us consider the system (\ref{dmsc}) with initial condition
\begin{equation}
\psi^{\e}\big|_{t=0}= \chi \, \exp{\left(-\frac{(x_1-0.1)^2+(x_2+0.1)^2+x_3^2}{4d^2}\right)}, \quad
\chi=(1,0,0,0),\ d=1/16,
\end{equation}
In this case we choose $\e=10^{-2}$, $\tg t=1/32$, $\tg x=1/32$. The numerical results are give in Figure \ref{fig41}.
We see that the wave packet moves in circles due to its interaction with the harmonic potential.
\begin{figure} 
\begin{center}
\resizebox{1.8in}{!} {\includegraphics{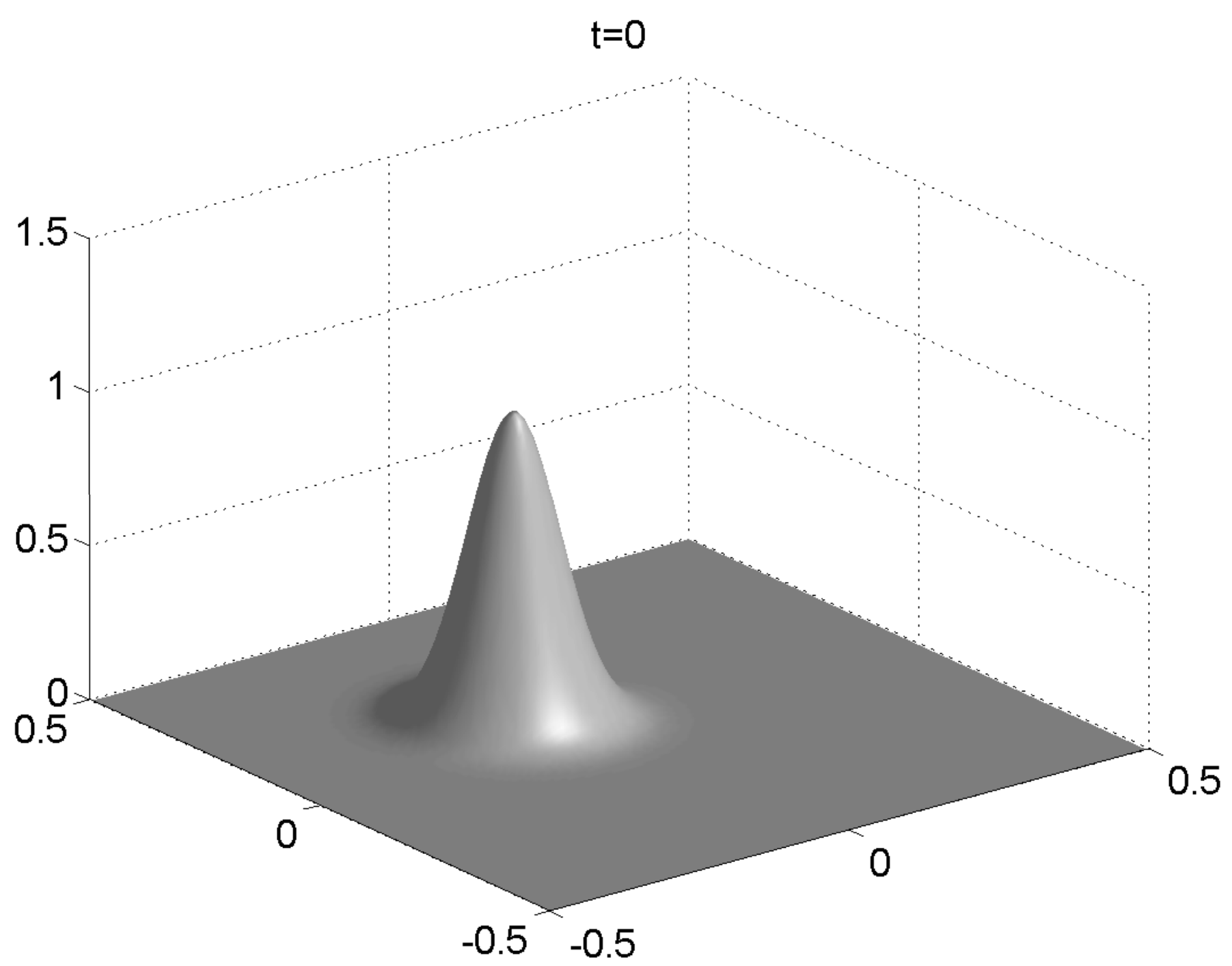}}
\resizebox{1.8in}{!} {\includegraphics{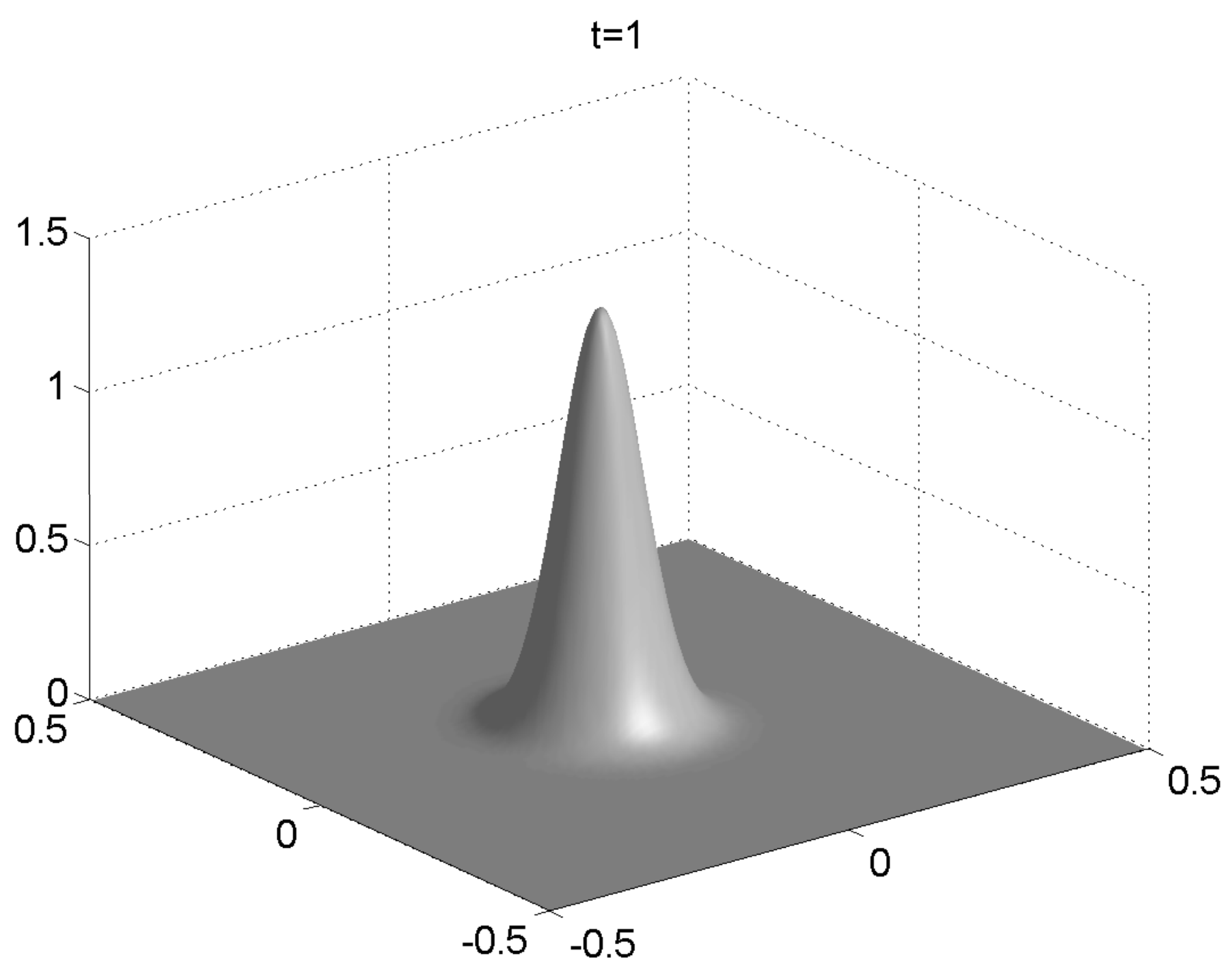}}
\resizebox{1.8in}{!} {\includegraphics{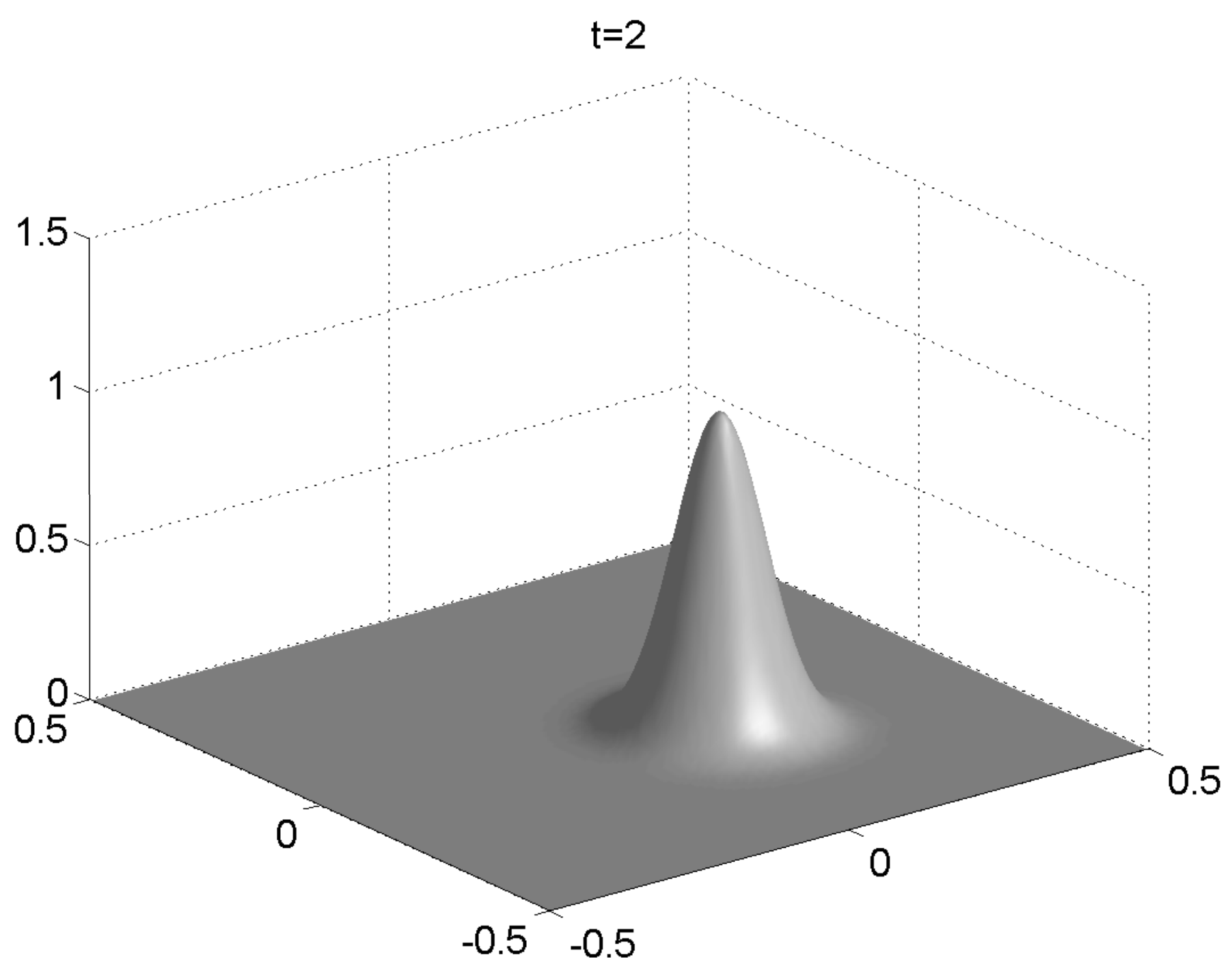}}
\vspace{4mm}

\resizebox{1.8in}{!} {\includegraphics{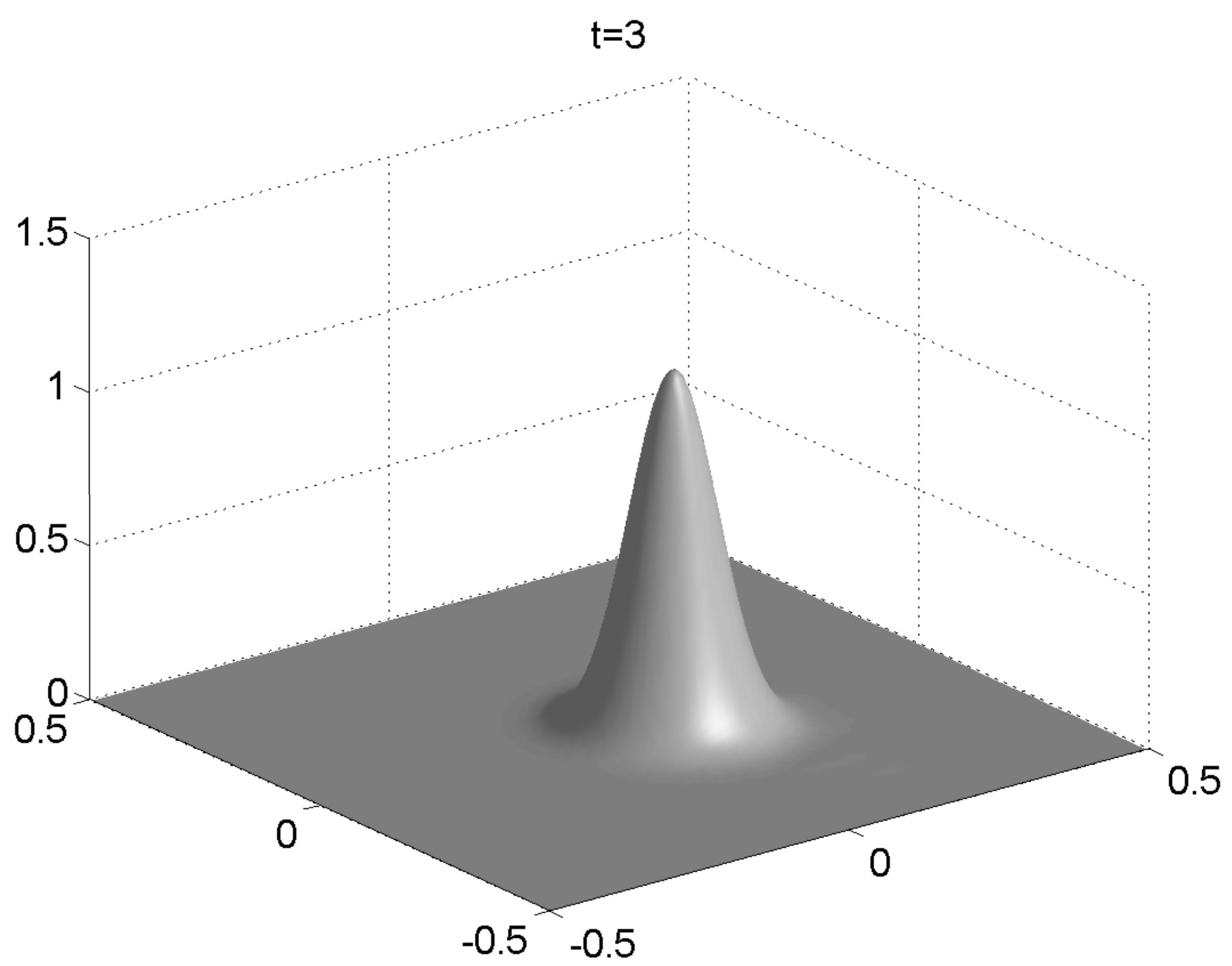}}
\resizebox{1.8in}{!} {\includegraphics{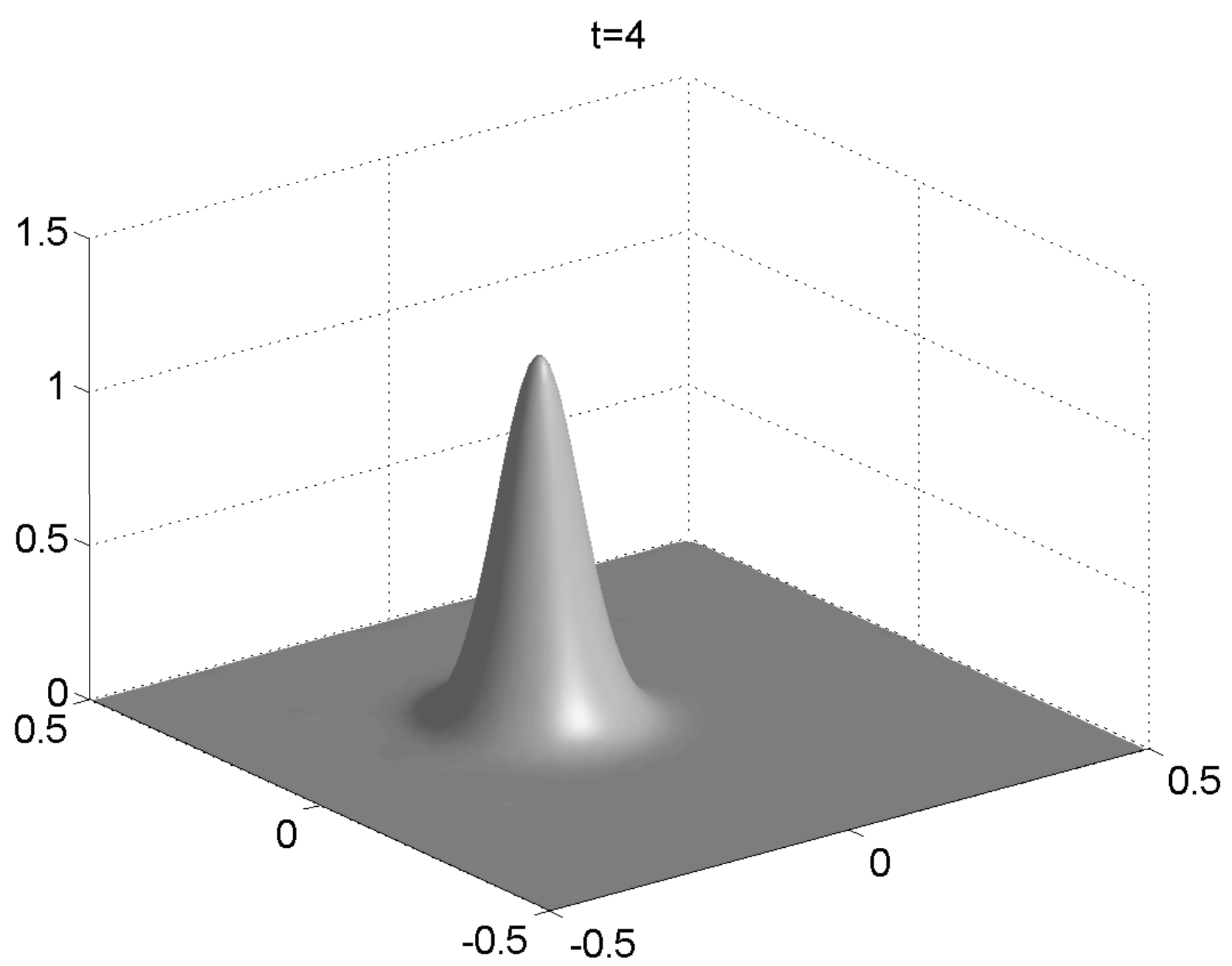}}
\resizebox{1.8in}{!} {\includegraphics{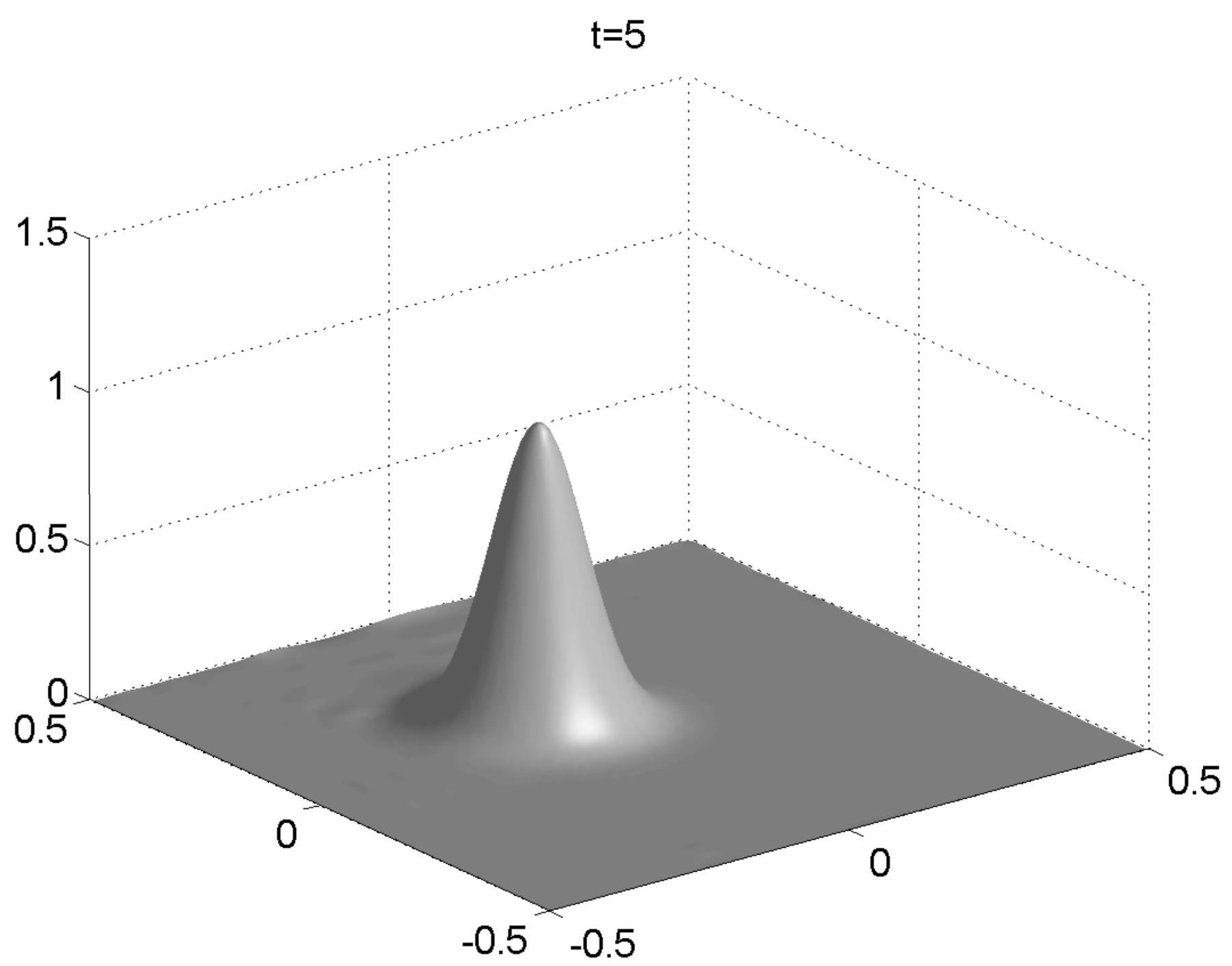}}
\vspace{4mm}

\resizebox{1.8in}{!} {\includegraphics{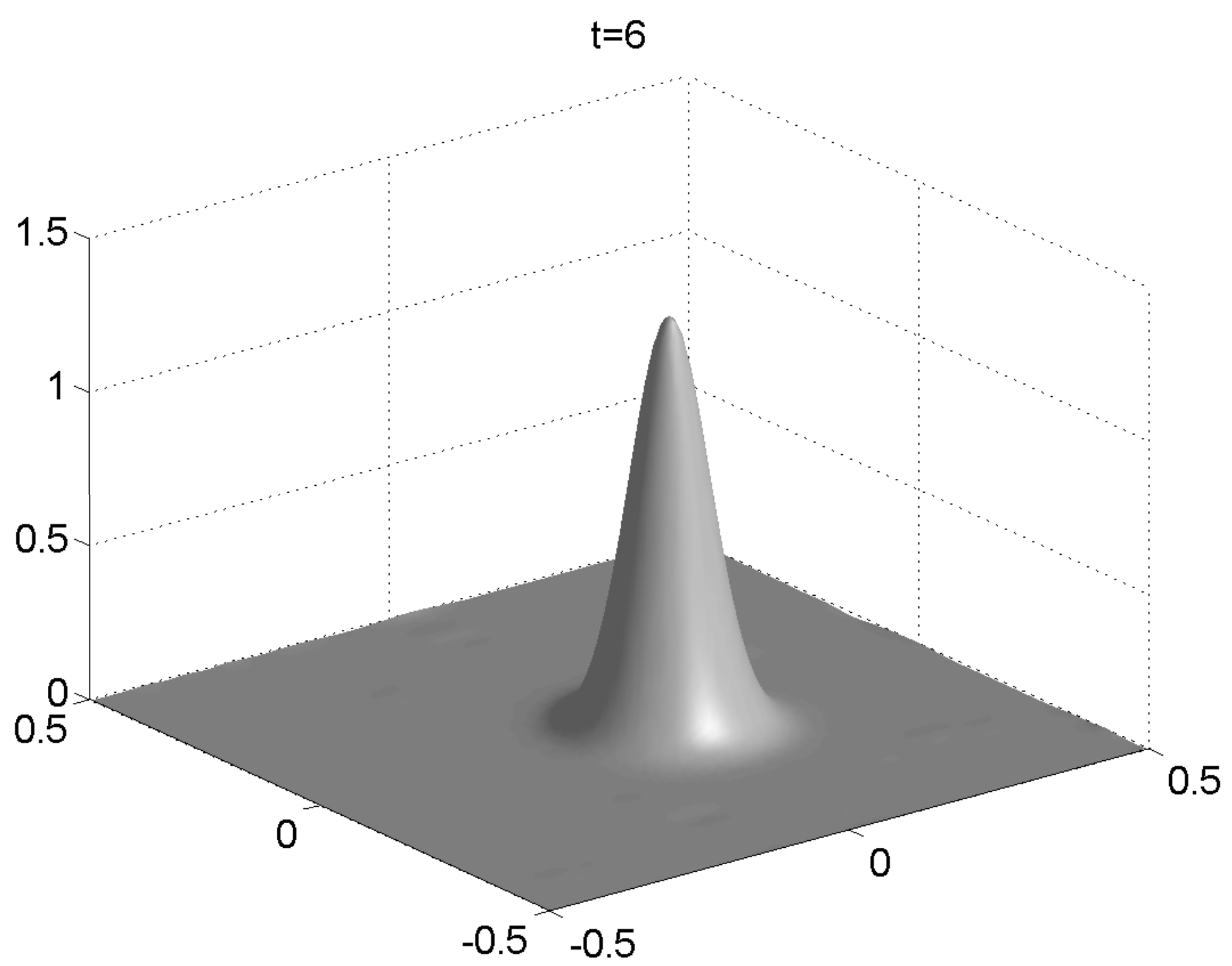}}
\resizebox{1.8in}{!} {\includegraphics{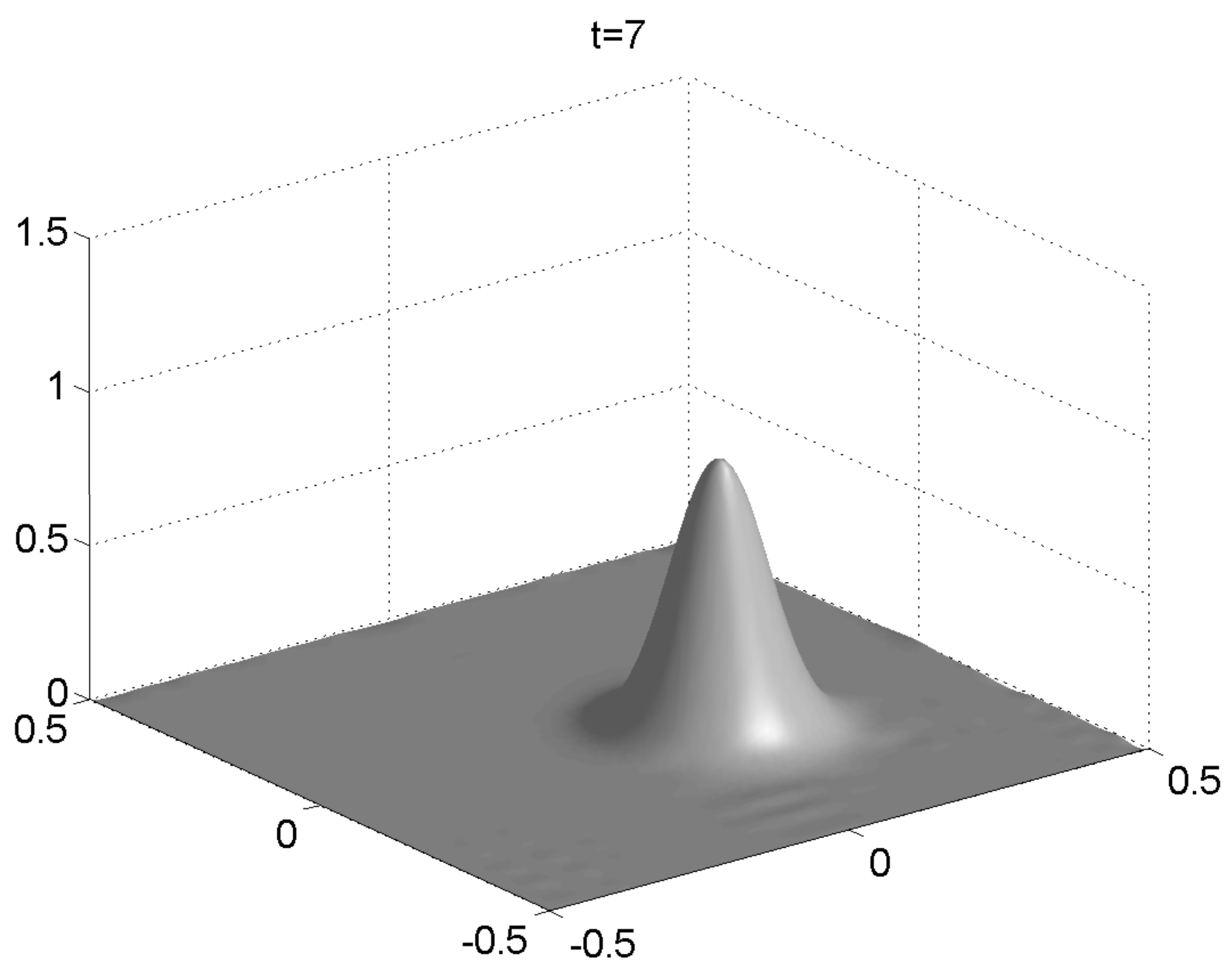}}
\resizebox{1.8in}{!} {\includegraphics{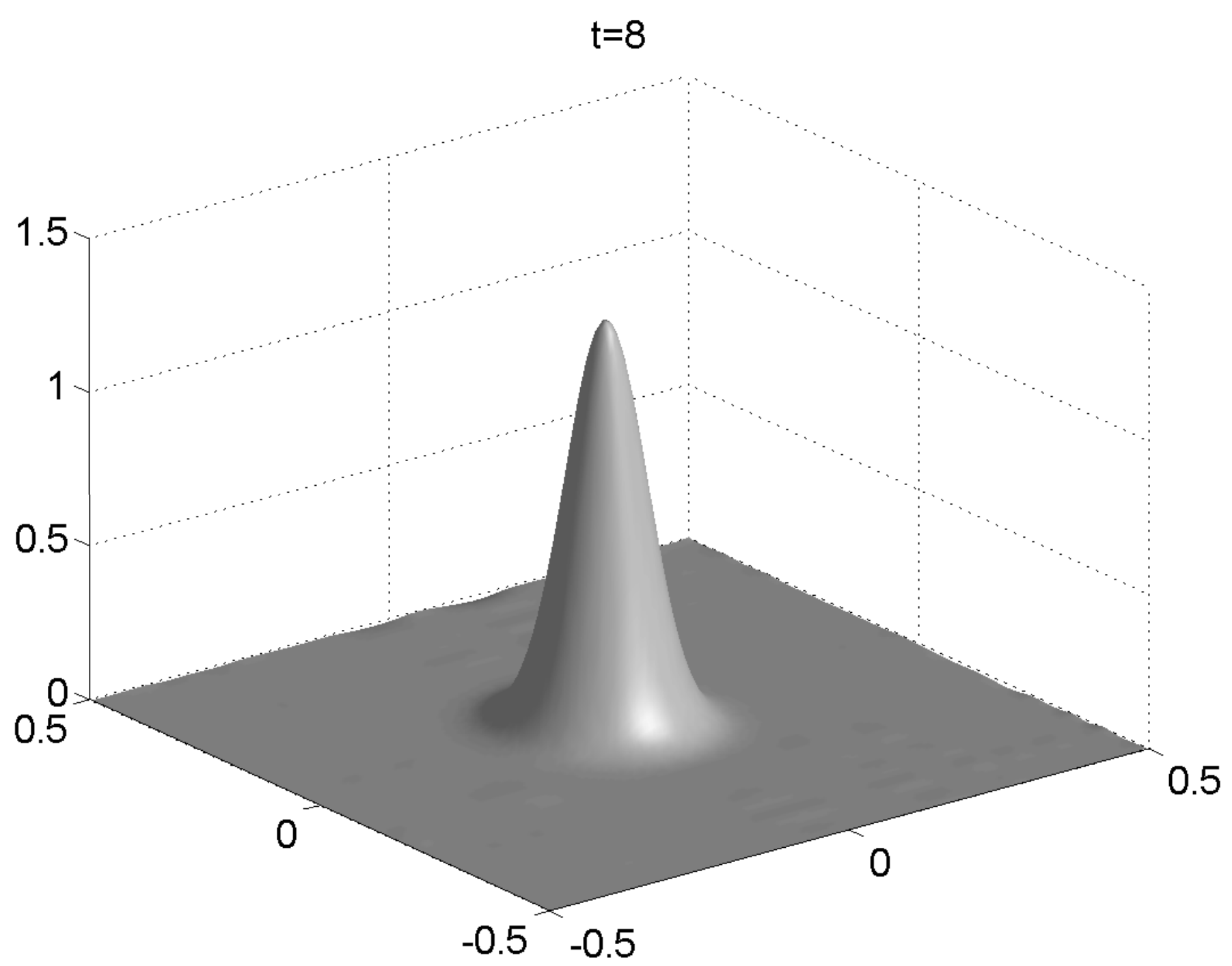}}
\end{center}
\caption{Numerical results of the density at different time for example \ref{exha}.
Here $\e=0.01$, $\tg x=1/32, \tg t=1/32$.} \label{fig41}
\end{figure}
\end{example}

\begin{remark} In analogy to the spectral-splitting method for the Schr\"odinger equation analyzed in
\cite{BJM1}, we find that
$\tg x = \O(\e)$ and $\tg t = \O(1)$, as $\e \rightarrow 0$, is sufficient to guarantee
well-approximated \emph{observable}
of the MD system. A more refined grid in temporal direction is necessary to obtain a
good approximation for the reps. components of
the spinor field itself, typically $\tg t = \mathcal O(\e^2)$ is needed.
\end{remark}

\section{The non-relativistic regime}\label{sec:nr}
Finally we shall also consider the non-relativistic regime for \eqref{dmnr} as $\de \rightarrow 0$. Again we shall
first describe the formal asymptotics and then discuss numerical examples.
\subsection{Formal description of the asymptotic problem}\label{snr}
To describe the non-relativistic  limit of the MD system we first define two pseudo-differential operators
$\Pi_{e/p}^\delta(D)$ via their symbols
\begin{align}
\label{frpr1}
\Pi^\delta_{e/p}(\xi):=\frac{1}{2}\left(\Id{4}\pm \frac{1}{\lambda_{0}(\delta\xi )}\mathcal D_{0}( \delta \xi )\right),
\end{align}
where $\lambda_0(\xi)$, $\mathcal D_0(\xi)$ are given by \eqref{lam}, \eqref{fdo}.
We then define the \emph{(non-relativistic) electronic} and the \emph{(non-relativistic) positronic component}
$\psi^\delta_{e}$, $\psi^\delta_{p}$ by
\begin{equation}
\psi^\delta_{e}(t,\xb):=e^{it/\delta^2}\Pi^\delta_{e}(D)\psi^\delta(t,\xb), \quad
\psi^\delta_{p}(t,\xb):=e^{-it/\delta^2}\Pi^\delta_{p}(D)\psi^\delta(t,\xb),
\end{equation}
where $\psi^\delta$ solves the non-relativistically scaled MD system \eqref{dmnr}. Note the difference in
sign of the phase-factors. This corresponds to subtracting the \emph{rest energy}, which is
positive for electrons but negative for positrons, \cf \cite{BMP, BMS2, Na}.
The above given definition of electronic/positronic wave functions
\emph{should not be confused with the one obtained in the semi-classical regime}, since both definitions are
adapted to the particular scaling of the resp. system under consideration. We remark that up to now
there is \emph{no} satisfactory interpretation in terms of electrons and positrons for the solution of the
\emph{full} MD system \eqref{dm0}, \eqref{dm1}. Indeed there is no such interpretation even for the linear Dirac equation
with external fields, see \eg \,\cite{Sc}.
\begin{remark}
It is easy to see that the formal limit $\delta \rightarrow 0$ of the operators $\Pi_{e/p}^\delta(D)$ yields,
\begin{equation}
\Pi_{e}^0=
\begin{pmatrix}
\Id{2} & 0\\
0 & 0
\end{pmatrix}
,\quad
\Pi_{p}^0=
\begin{pmatrix}
0 & 0\\
0 & \Id{2}
\end{pmatrix}
.
\end{equation}
This explains the interpretation of electrons (resp. positrons) as the \emph{upper} (resp. \emph{lower}) components
of the $4$-vector $\psi^\delta$ for small values of $\delta$, \cf \cite{Sc}.
\end{remark}
It is then shown in \cite{BMS2} (see also \cite{BMP} for easier accessible proofs in the linear case) that
\begin{equation}\label{eq:nrconverge}
\psi^\delta_{e/p}(t,\xb)\stackrel{\delta\rightarrow 0}{\longrightarrow} \varphi_{e/p}(t,\xb),
\quad \mbox{in $H^1(\R^3)\otimes \C^4$},
\end{equation}
where $\varphi_{e}, \varphi_{p}$
solve the mixed electronic/positronic \emph{Schr\"odinger-Poisson system}:
\begin{equation}\label{eq:schp}
\left\{
\begin{aligned}
i \partial_t \varphi_e = & \, - \frac{\Delta}{2} \varphi_e + (V+V^{ex}) \varphi_e,\\
i \partial_t \varphi_p = & \, + \frac{\Delta}{2} \varphi_p + (V+V^{ex}) \varphi_p,\\
-\Delta V = & \, |\varphi_p|^2+|\varphi_e|^2,
\end{aligned}
\right.
\end{equation}
In contrast to the asymptotic problem obtained in the semi-classical limit, this system is \emph{globally}
well posed. The appearance of the Poisson equation can be motivated by performing
a naive Hilbert expansion in the self-consistent fields, \cf \cite{MaMa}, \ie
\begin{equation}
V^\delta = V + \delta \widetilde V + \O(\delta^2), \quad
\Ab^\delta = \Ab + \delta \widetilde \Ab + \O(\delta^2).
\end{equation}
Plugging this into \eqref{dmnr}, comparing equal powers in $\delta$, and having in mind
that $\Jb^\delta \sim \O(1)$ \cite{BMP} gives \eqref{eq:schp}. In \cite{BMS2} the electric potential is proved
to converge in $H^1(\R^3)$ as $\delta \rightarrow 0$, whereas the convergence of the magnetic fields is
not studied in detail. Indeed, it is shown in \cite{BMS2} that if one only aims for a derivation of the
Schr\"odinger-Poisson system, one can even allow for initial data $\Ab^\delta(0,\xb)$, $\partial_t \Ab^\delta(0,\xb)$
which do \emph{not} converge as $\de \rightarrow 0$.
\begin{remark}
If we would, in addition, consider terms of order $\O(\delta)$ too, we (formally) would obtain a \emph{Pauli equation}
for $\varphi_{e/p}$, including the matrix-valued magnetic field term $\sum \sigma_k B_k$, \ie the, so called,
\emph{Pauli-Poiswell system}, \cf \cite{BMP, MaMa}. Moreover we remark that the authors in \cite{BMS2} considered the
MD system in Coulomb gauge, \ie $\diverg \Ab =0$, instead of the Lorentz gauge condition imposed in this work
\eqref{lor}. The reason is rather technical and it is not clear yet if a generalization of their work
to the Lorentz gauged system is possible.
\end{remark}
As before we shall use a time-splitting spectral method \cite{BJM2} to solve the
coupled system of Schr\"odinger-Poisson equations \eqref{eq:schp}:
\newpar
\textbf{Step 1.} First, we solve the following problem:
\begin{equation}\label{eq:schp1}
\left\{
\begin{aligned}
i \partial_t \varphi_e = & \, - \frac{\Delta}{2} \varphi_e,\\
i \partial_t \varphi_p = & \, + \frac{\Delta}{2} \varphi_p,\\
-\Delta V = & \, |\varphi_p|^2+|\varphi_e|^2,
\end{aligned}
\right.
\end{equation}
\textbf{Step 2.} Then we solve the coupled equations
\begin{equation}\label{eq:schp2}
\left\{
\begin{aligned}
i \partial_t \varphi_e = & \, (V+V^{ex}) \varphi_e,\\
i \partial_t \varphi_p = & \, (V+V^{ex}) \varphi_p,
\end{aligned}
\right.
\end{equation}
In step 1, we again use the pseudo-spectral method. In step 2, we can get the exact
solution for this linear ODE system in
time, since $|\varphi_p|^2$ and $|\varphi_e|^2$, resp., are kept invariant by step 2.
\begin{remark}
Let us fix $\varepsilon =1$ and consider $\dt\to0$ in the algorithm given in section
\ref{sec:num}. Based on the expansion of \eqref{eq:numst14}--\eqref{eq:numst16}, we obtain
\be\label{eq:nonr1}
\hat\Phi^{n+1}= \exp\left(\Ld(t-t_n)\right)\hat\Psi^n+\mathcal{O}(\dt),\ee
where in the limit $\dt\to0$ the matrix $\Lambda\in \C^{4\times4}$ simplifies to
$$\Ld=\mbox{diag}[\ld,\ld,-\ld,-\ld], \quad \ld=-i(\dt^{-2}+|\xi|^2/2) .$$
We also have
\bea\label{eq:nonr2}
|\xi|^2\left(\hat V ^{n}+\hat V ^{n+1}\right)
&=&\widehat{\ |\Psi^{n}|^2}+\widehat{\ |\Phi^{n+1}|^2}+\mathcal{O}(\dt)
\eea
and
\bea
|\xi|^2\left(\hat \Ab ^{n}+\hat {\mathbf A}^{n+1} \right)
&=&\mathcal{O}(\dt), \label{eq:nonr3}
\eea
because $ \left<\Phi^{n+1},\ap^k\Phi^{n+1}\right>=\mathcal{O}(\dt)$. If we denote
the \emph{upper} (resp. \emph{lower}) components of the $4$-vector $\Psi$
by $\Psi_e$ (resp. $\Psi_p$), we obtain
\bea\label{eq:nonr4}
\partial_t\left(e^{\pm it/\dt^2}\hat\Phi_{e/p}\right)&=&\mp i\f{|\xi|^2}{2}
\left(e^{\pm it/\dt^2}\hat\Phi_{e/p}\right)+\mathcal{O}(\dt),
\eea
and from \eqref{eq:numst22}, we find
\be\label{eq:nonr6}
\Psi^{n+1}=\exp\left(-i V \tg t\right)\Phi^{n+1}_e+\mathcal{O}(\dt).\ee
Combining  the equations \eqref{eq:nonr2}--\eqref{eq:nonr6}, we conclude that the numerical solutions of our algorithm,
given in section \ref{sec:num}, uniformly converge to the numerical solutions of the above algorithm.
This analysis, previously done for a time-splitting
spectral method for the Zakharov system \cite{JMZ}, shows that one can choose
$h, \tg t$ independent of $\delta$.
\end{remark}
\subsection {Numerical examples for the non-relativistic regime}
\setcounter{lemma}{0}
\begin{example}[\textbf{Purely self-consistent motion II}]
\label{exnr}
Here we consider the MD system \eqref{dmnr}
in a unit cubic with periodic boundary conditions, zero external fields, and initial data
\begin{equation}\label{eq:e51}
\left \{
\begin{aligned}
& \psi^{\delta}(\xb)\big|_{t=0}\equiv \psi^{(0)}(\xb)=
\dpm \chi \exp\left(-\f{|\xb|^2}{4d^2}\right) ,
\quad \chi= (1,1,1,1),\ d=\f{1}{16},\\
& -\Dt V^{(0)}=|\psi^{(0)}|^2,
\quad V^{(1)}(\xb) =0, \\
& -\Dt A_k^{(0)}=\langle \psi^{(0)} , \alpha^k \psi^{(0)} \rangle_{\C^4},
\quad \Ab^{(1)} (\xb)=0.
\end{aligned}
\right.
\end{equation}
Note that the above choice of initial data for $V$ and $A_k$ is
done to avoid initial layers. The impact of this choice on the numerical resolution, \ie the mesh strategy etc.,
is analogous to the \emph{Zakharov system} discussed in \cite{JMZ}.
We also consider the Sch\"odinger-Poisson problem \eqref{eq:schp}
with the initial data
\begin{equation}\label{eq:schp_i}
\varphi_{e}(t,\xb)|_{t=0}=\Pi^\delta_{e}(D)\psi^{(0)}(\xb), \quad
\varphi_{p}(t,\xb)|_{t=0}=\Pi^\delta_{p}(D)\psi^{(0)}(\xb),
\end{equation}
We compare the solution of the MD system
with the (coupled) Sch\"odinger-Poisson problem, \cf Figure \ref{fig51} and Figure \ref{fig52}.
Table \ref{tbnr}, Figure \ref{fig51} and Figure \ref{fig52} illustrates the validity of \eqref{eq:nrconverge}.
The Figures \ref{fig51}--\ref{fig521} also show that $|\Ab^\dt|=\O(\dt)|V^\dt|$, as $\dt\to0$.
\begin{table}[h]
\begin{center}
\caption{Convergence test for example \ref{exnr}: (here $\tg
t= 1/128$, $\tg x=1/64$)}\label{tbnr}
\begin{tabular}{cccc}\hline
$\dt$ &    0.01 & 0.1 &    1.0\\ \hline
$\ba{c}\vspace{-4mm}\\ \dpm\sup_{0\le t\le
1/4}\left|\psi_{e}^\dt -\varphi_{e}\right|^2+
|\psi_{p}^\dt -\varphi_{p}|^2\ea$
&0.101 &   0.345 & 2.407  \vspace*{1mm}\\ \hline
\end{tabular}
\end{center}
\end{table}
\begin{figure} 
\begin{center}
\resizebox{1.4in}{!} {\includegraphics{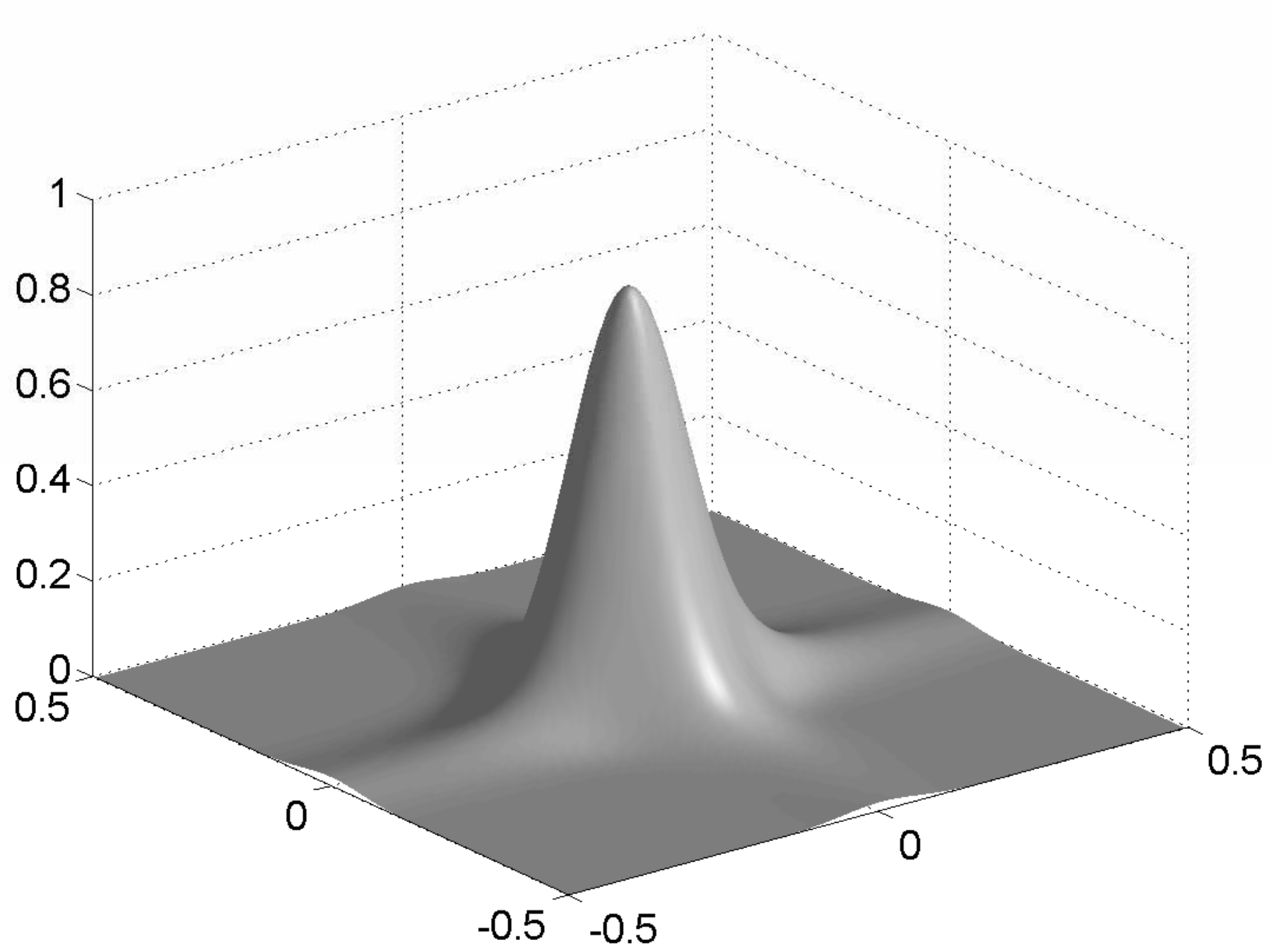}}
\resizebox{1.4in}{!} {\includegraphics{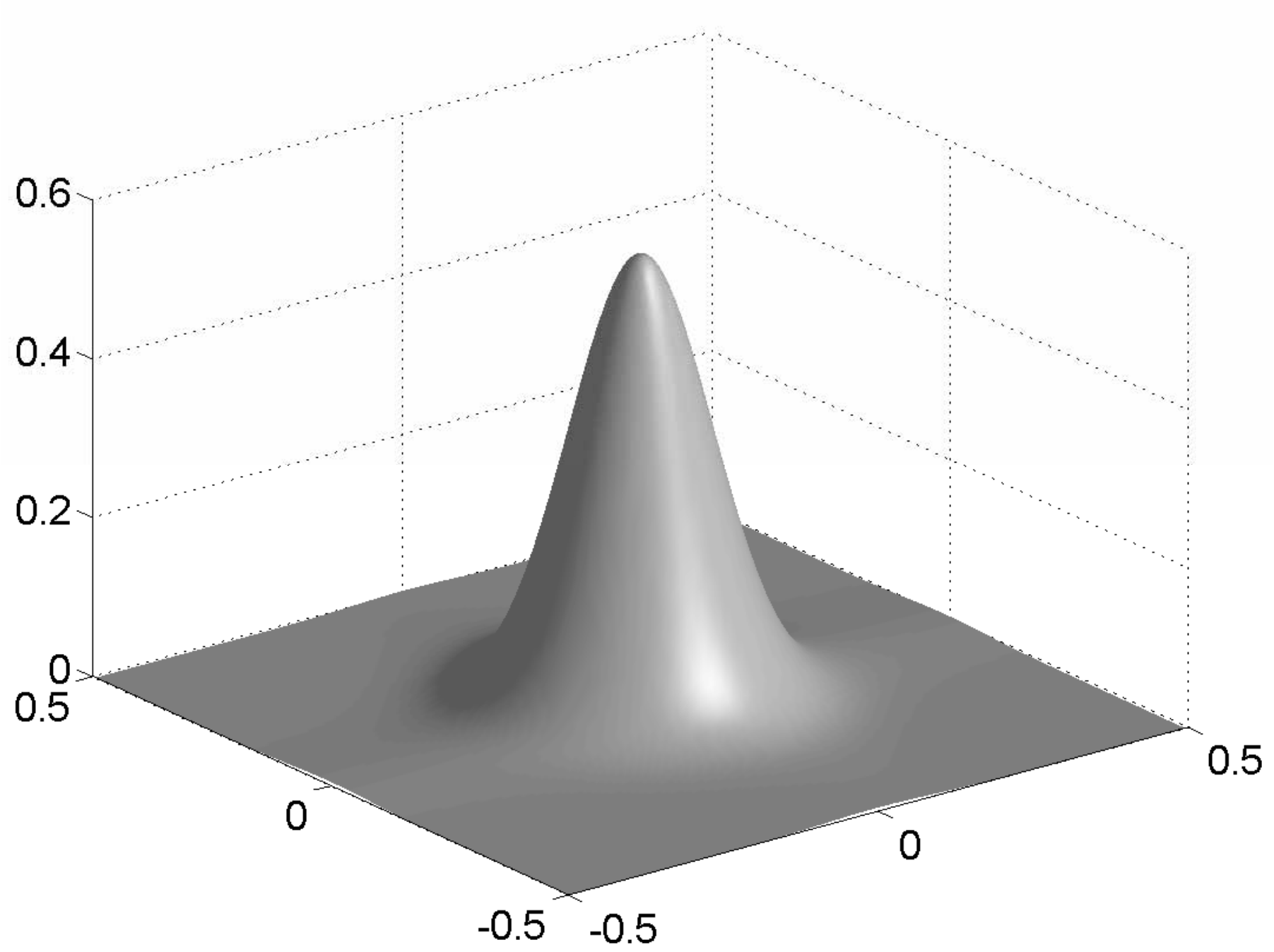}}
\resizebox{1.4in}{!} {\includegraphics{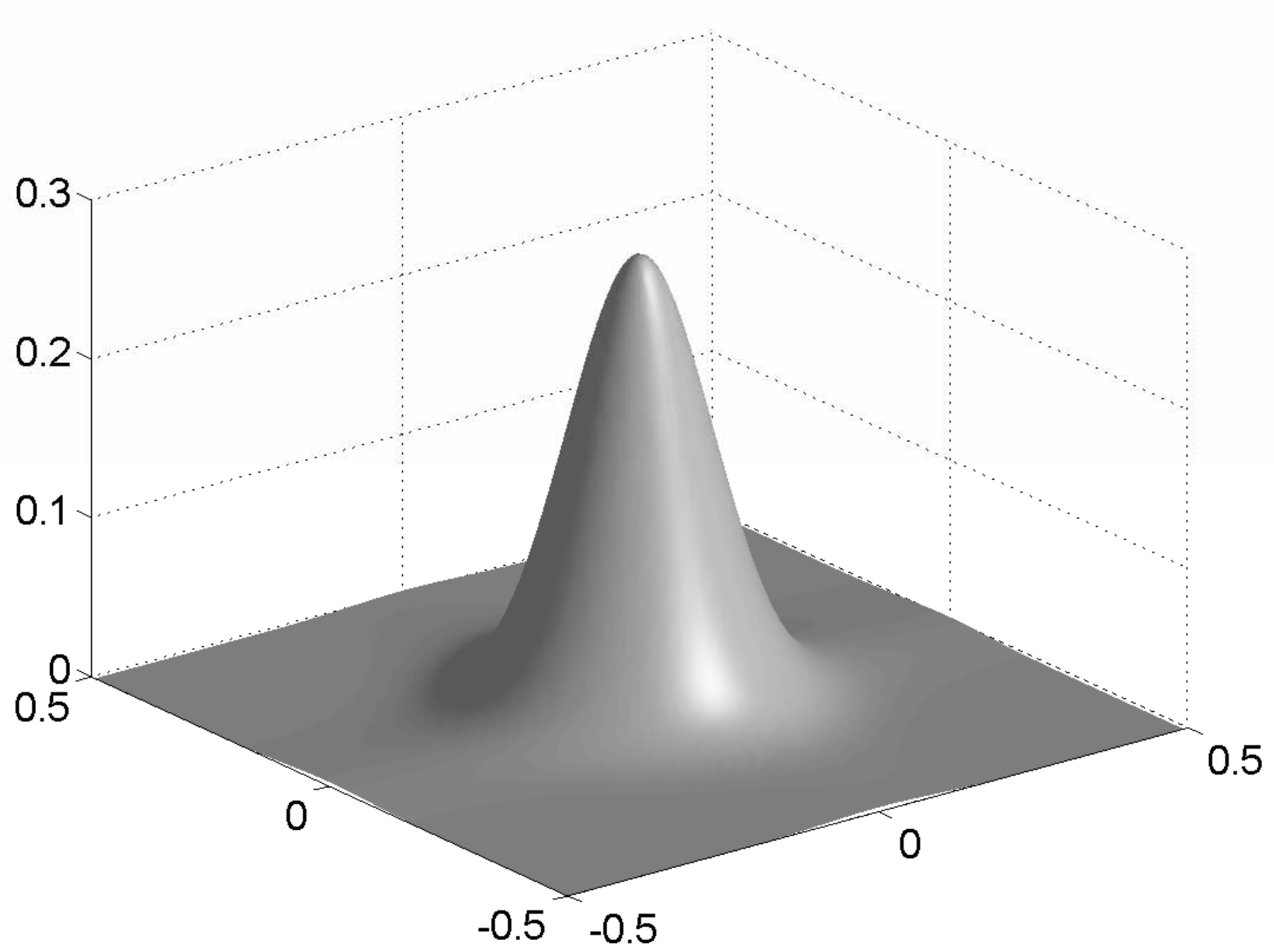}}
\resizebox{1.4in}{!} {\includegraphics{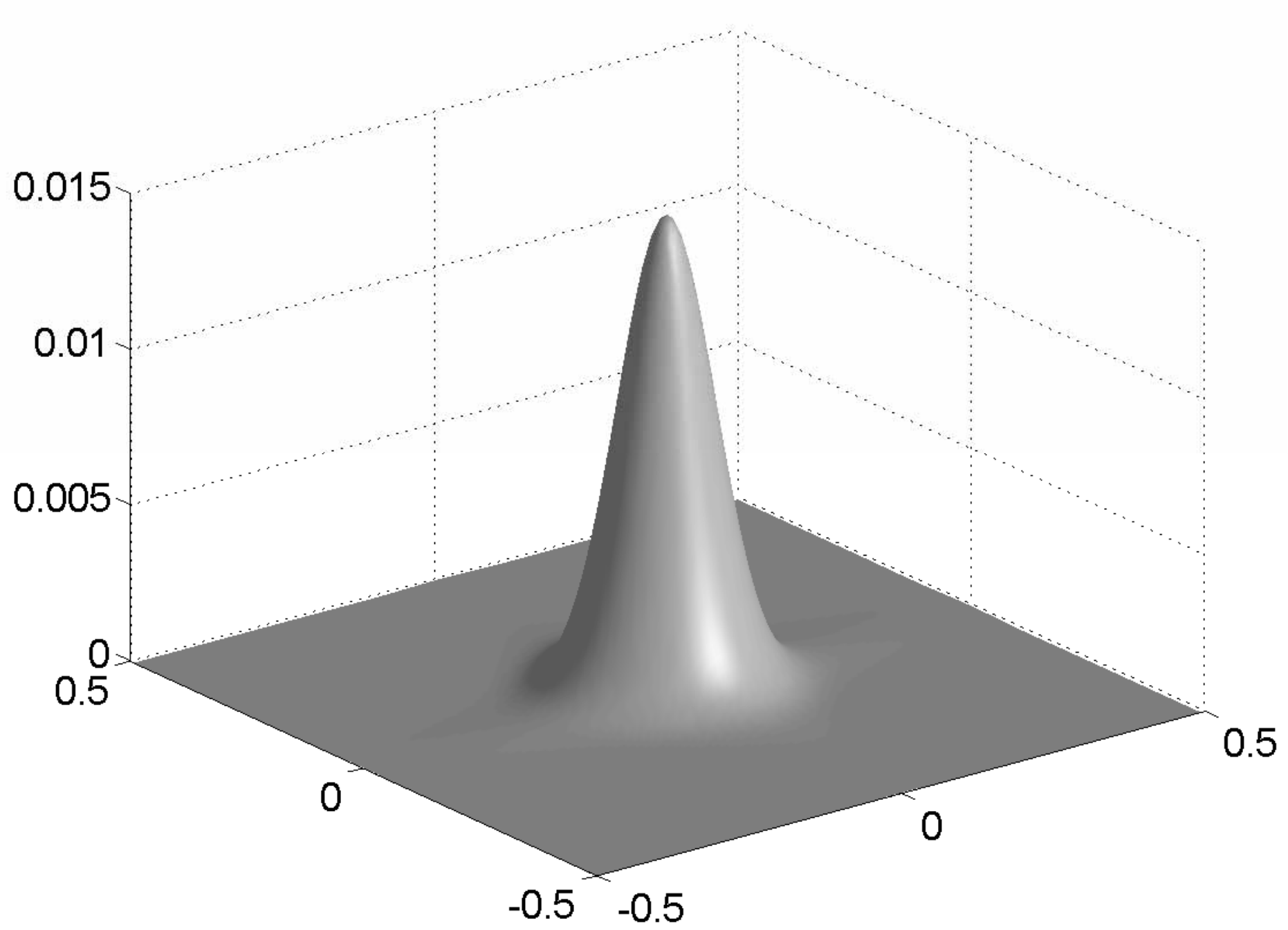}}

{$|\psi^\dt(t,\xb)|^2\big|_{x_3=0}$, $\re\big(\psi_{e,1}^\dt(t,\xb)\big)\big|_{x_3=0}$,
$\im\big(\psi_{e,1}^\dt(t,\xb)\big)\big|_{x_3=0}$ and $V^\dt(t,\xb)\big|_{x_3=0}$ for $\dt=1.0$.}\vspace{1mm}

\resizebox{1.4in}{!} {\includegraphics{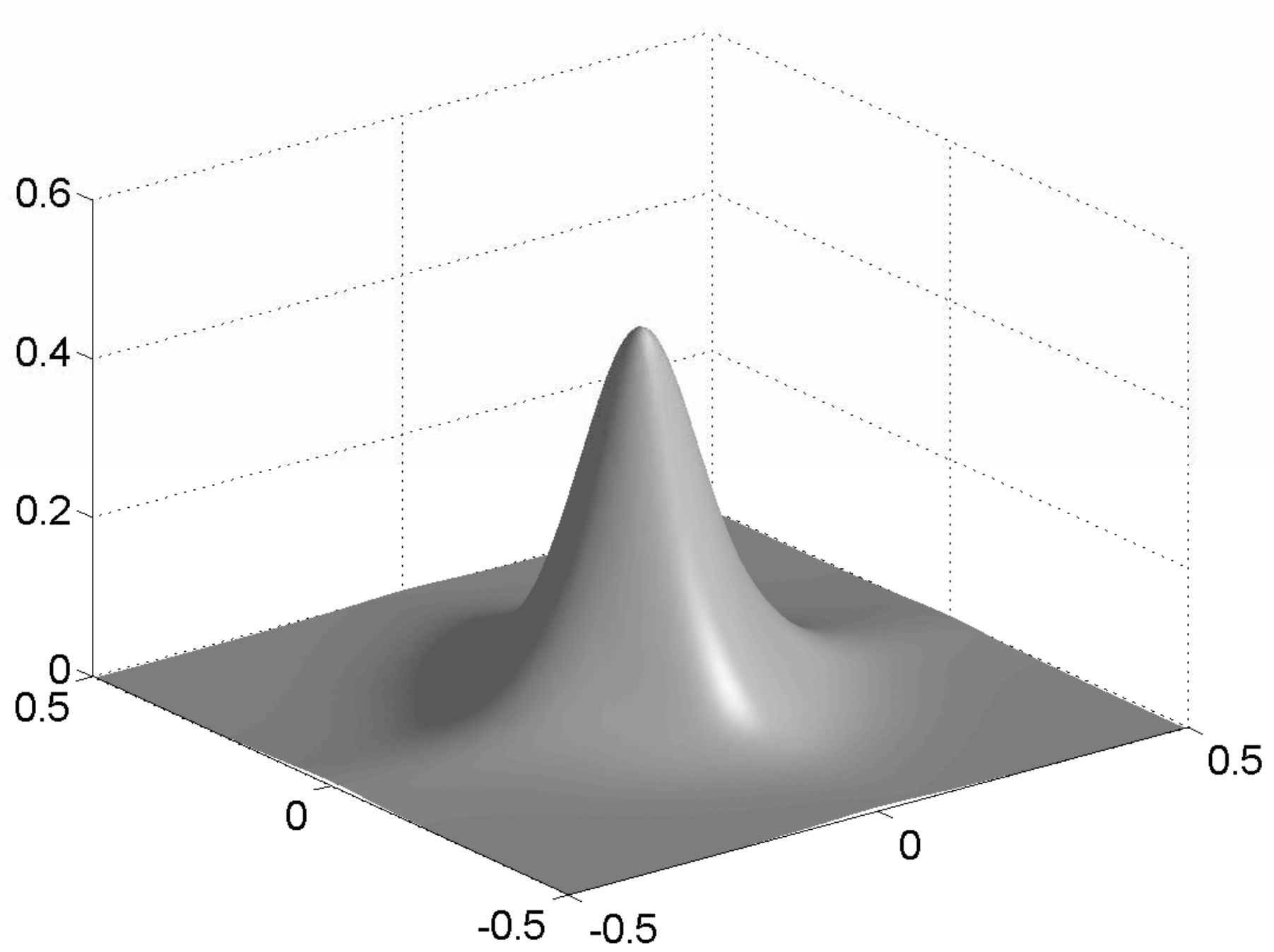}}
\resizebox{1.4in}{!} {\includegraphics{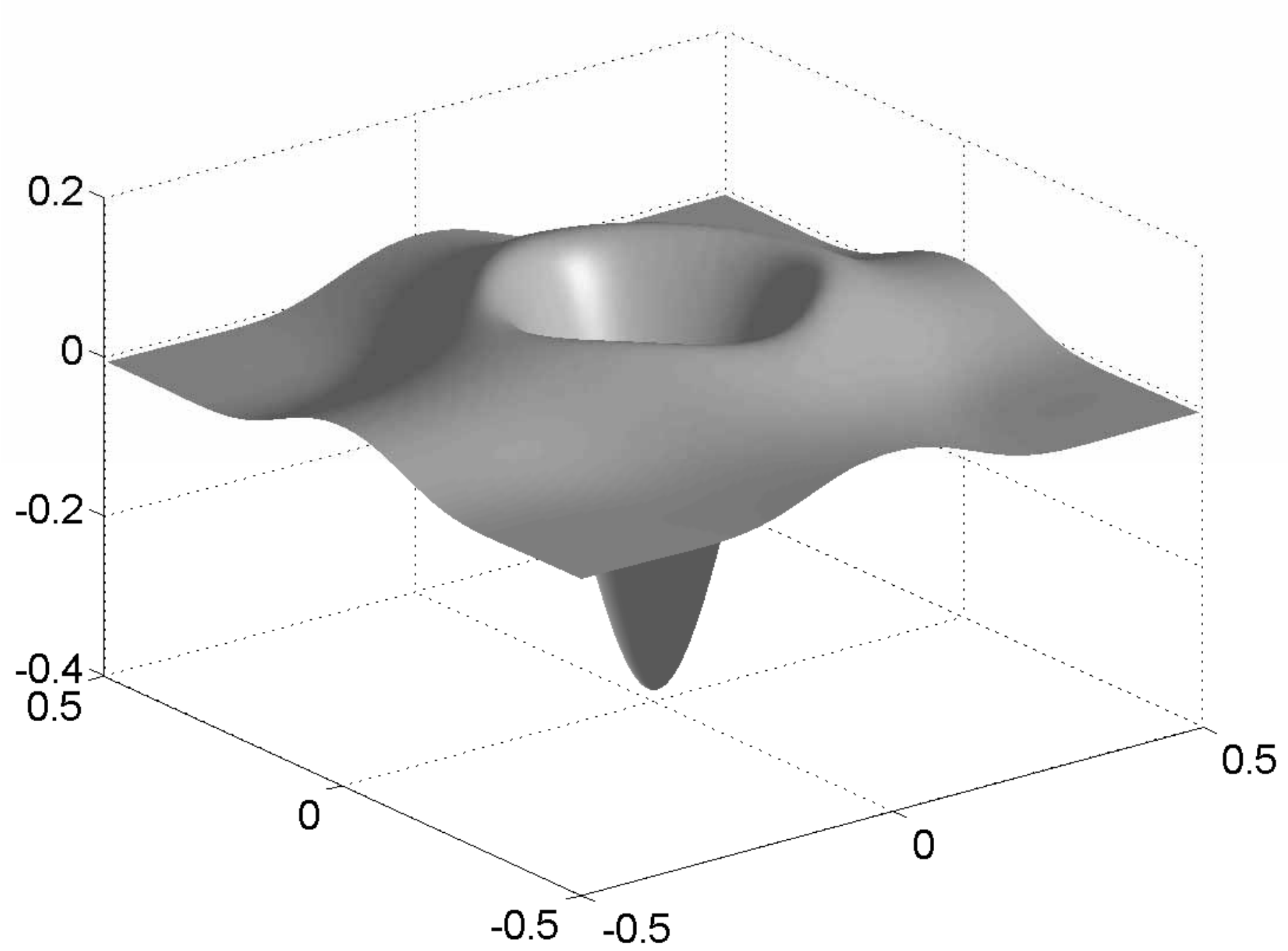}}
\resizebox{1.4in}{!} {\includegraphics{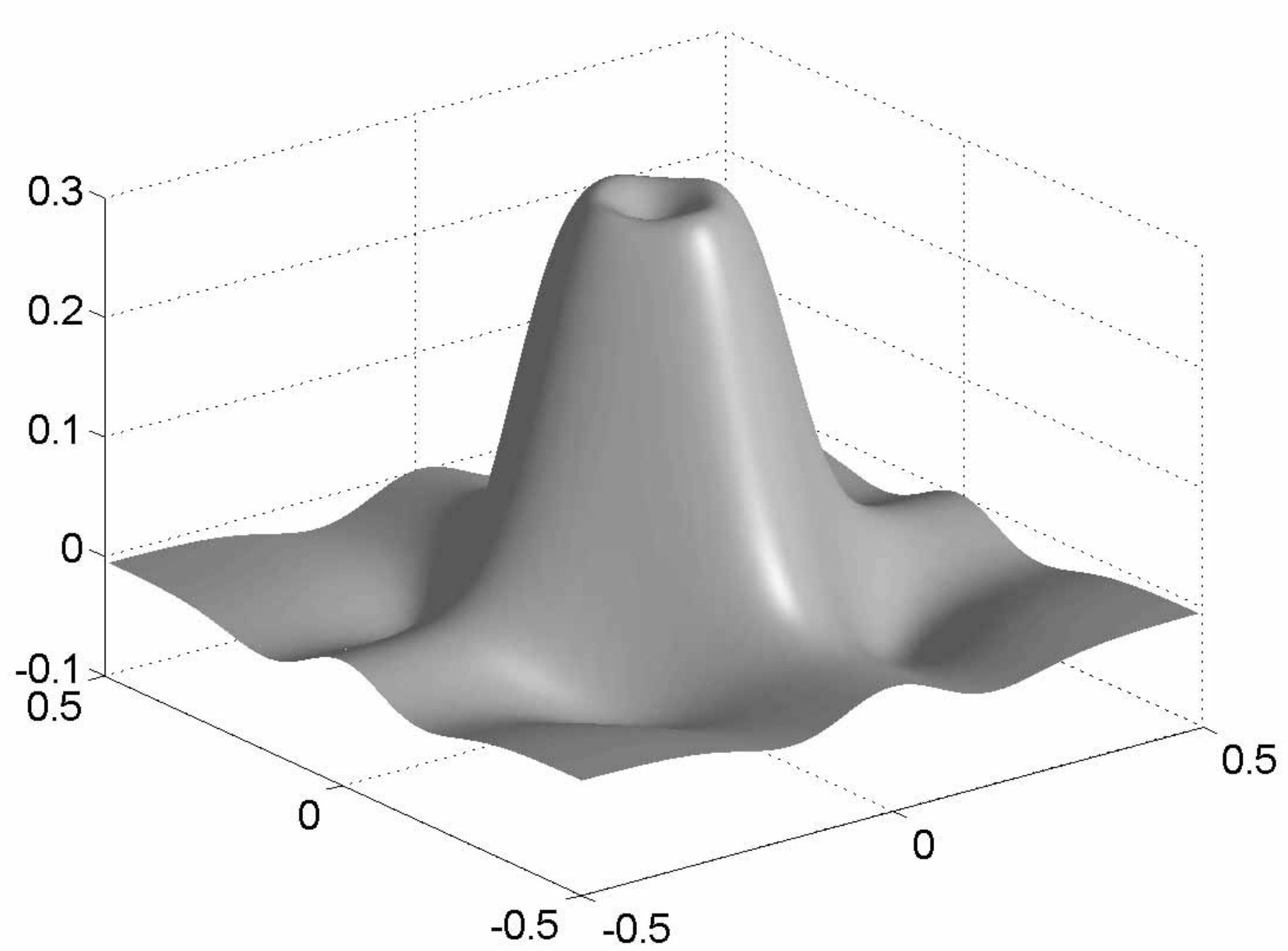}}
\resizebox{1.4in}{!} {\includegraphics{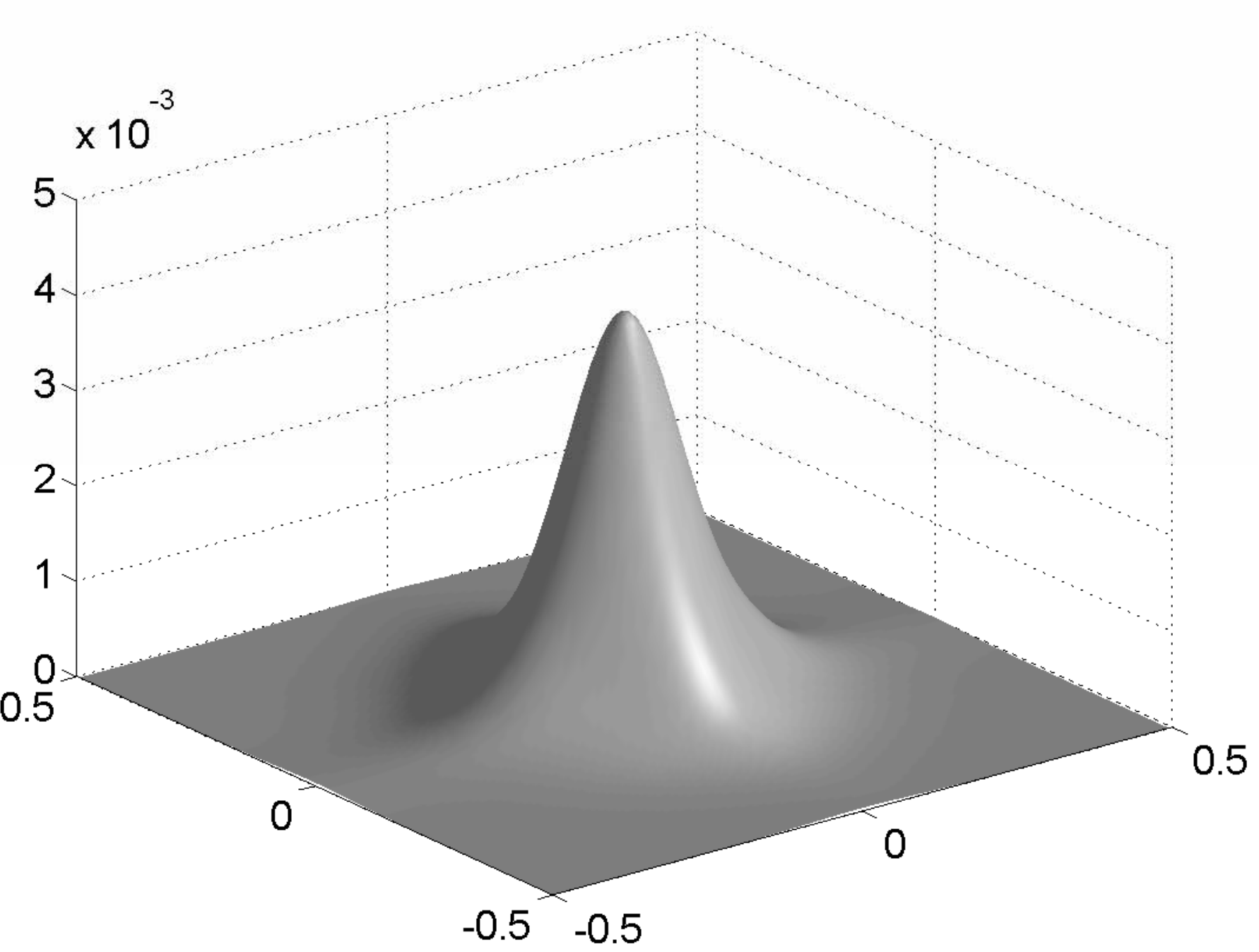}}

{$|\psi^\dt(t,\xb)|^2\big|_{x_3=0}$, $\re\big(\psi_{e,1}^\dt(t,\xb)\big)\big|_{x_3=0}$,
$\im\big(\psi_{e,1}^\dt(t,\xb)\big)\big|_{x_3=0}$ and $V^\dt(t,\xb)\big|_{x_3=0}$ for $\dt=0.01$.}\vspace{1mm}

\resizebox{1.4in}{!} {\includegraphics{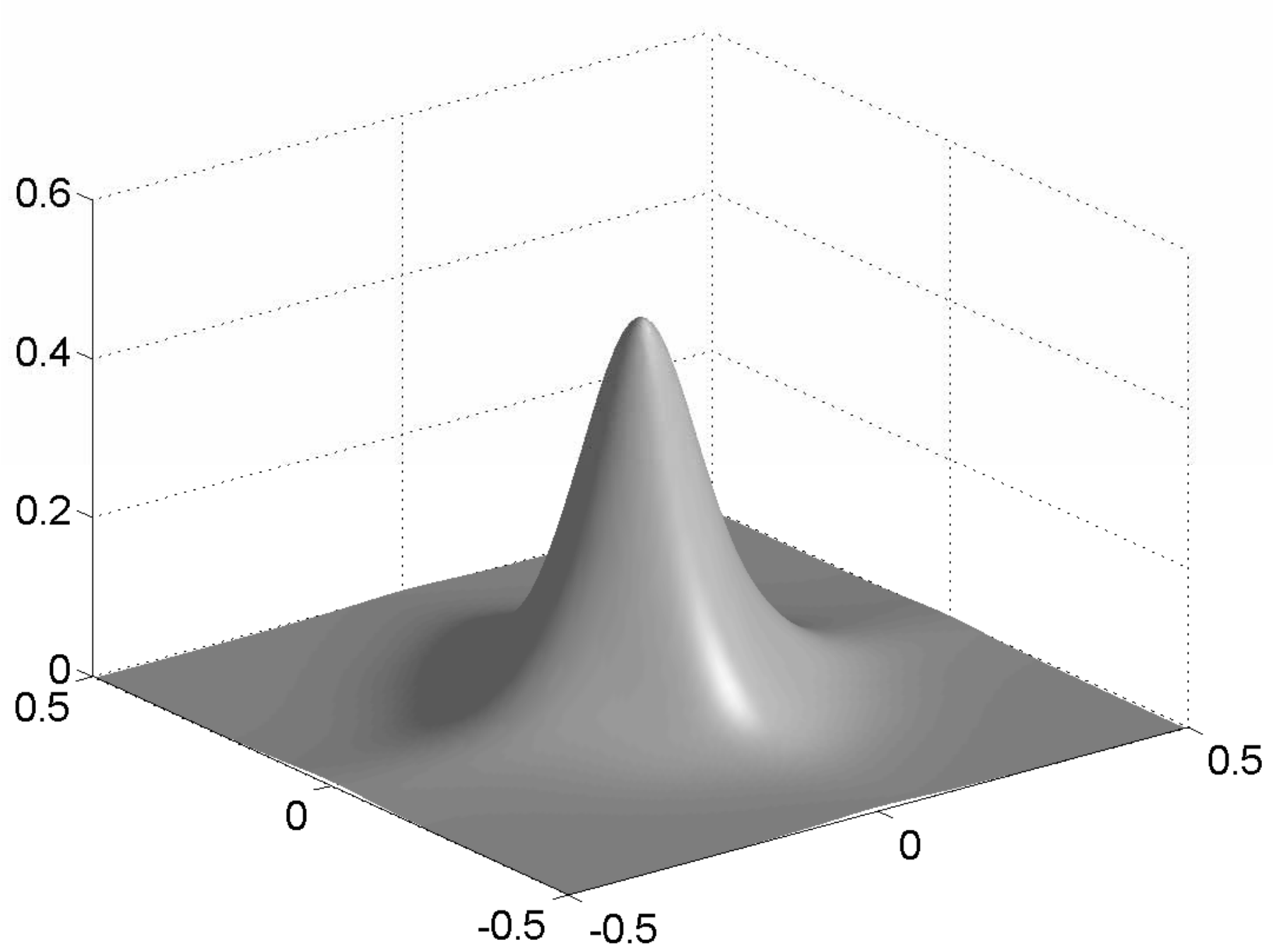}}
\resizebox{1.4in}{!} {\includegraphics{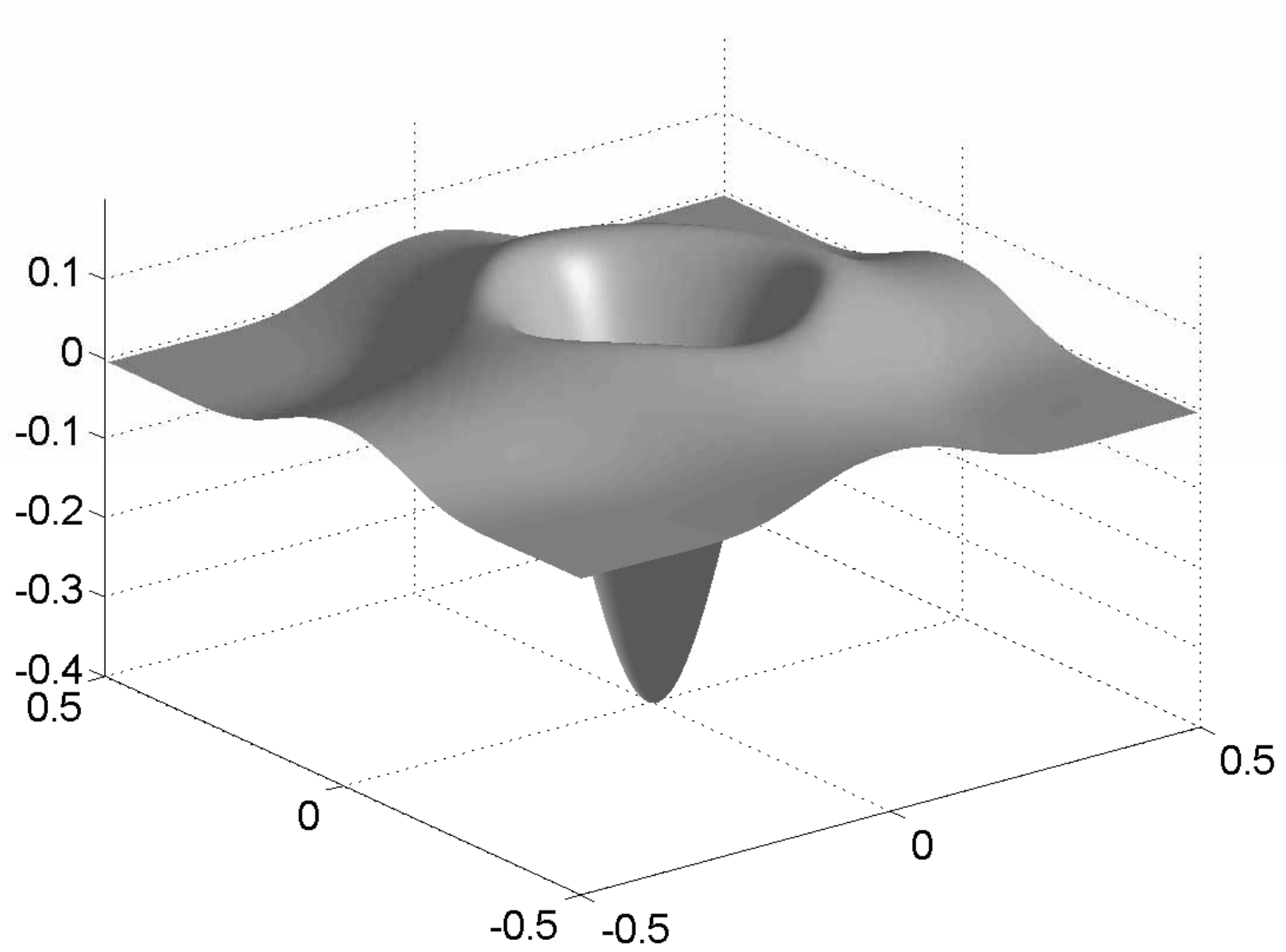}}
\resizebox{1.4in}{!} {\includegraphics{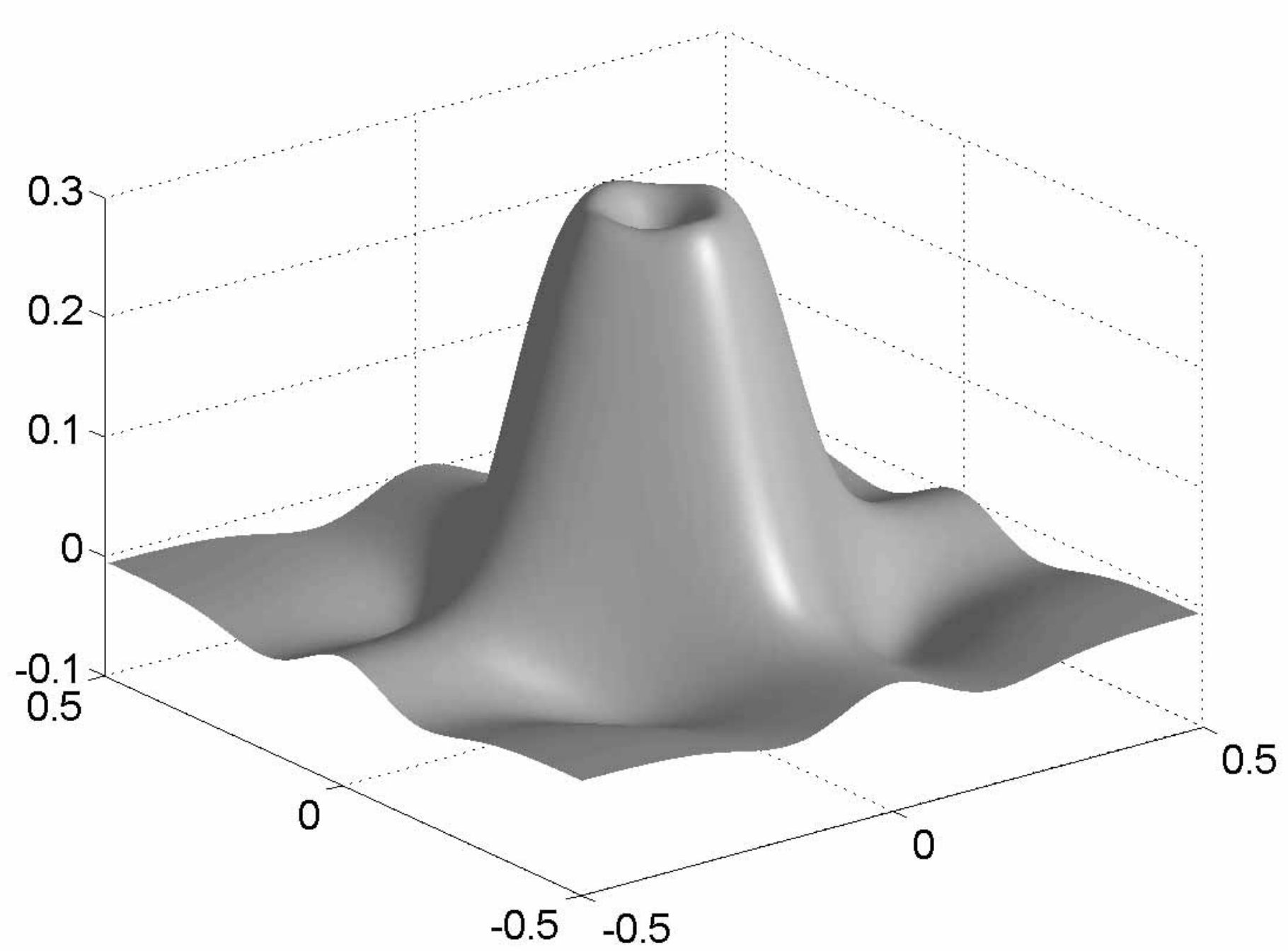}}
\resizebox{1.4in}{!} {\includegraphics{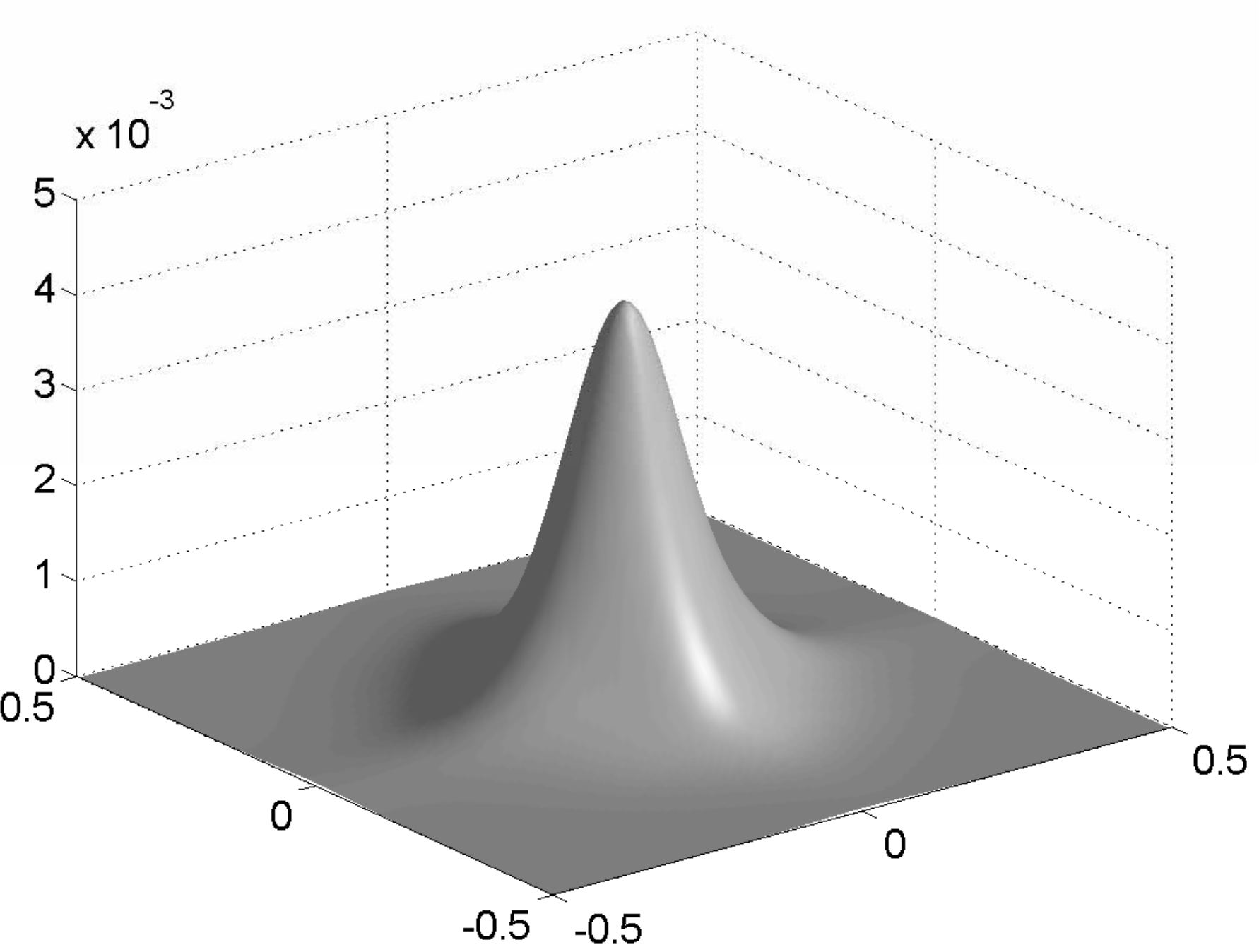}}

{$\left(|\vp_{e}(t,\xb)|^2+|\vp_{p}(t,\xb)|^2\right)\big|_{x_3=0}$, $\re\big(\vp_{e,1}(t,\xb)\big)\big|_{x_3=0}$,
$\im\big(\vp_{e,1}(t,\xb)\big)\big|_{x_3=0}$ and $V(t,\xb)\big|_{x_3=0}$.}
\end{center}
\caption{Numerical results for example \ref{exnr} at t=0.5.
The first row is the solution of the MD system with $\dt=1.0$,
whereas the second row is the solution of the MD system with $\dt=0.01$,
the third line is the solution of the Schr\"odinger-Poisson system.
}\label{fig51}
\begin{center}
\resizebox{1.4in}{!} {\includegraphics{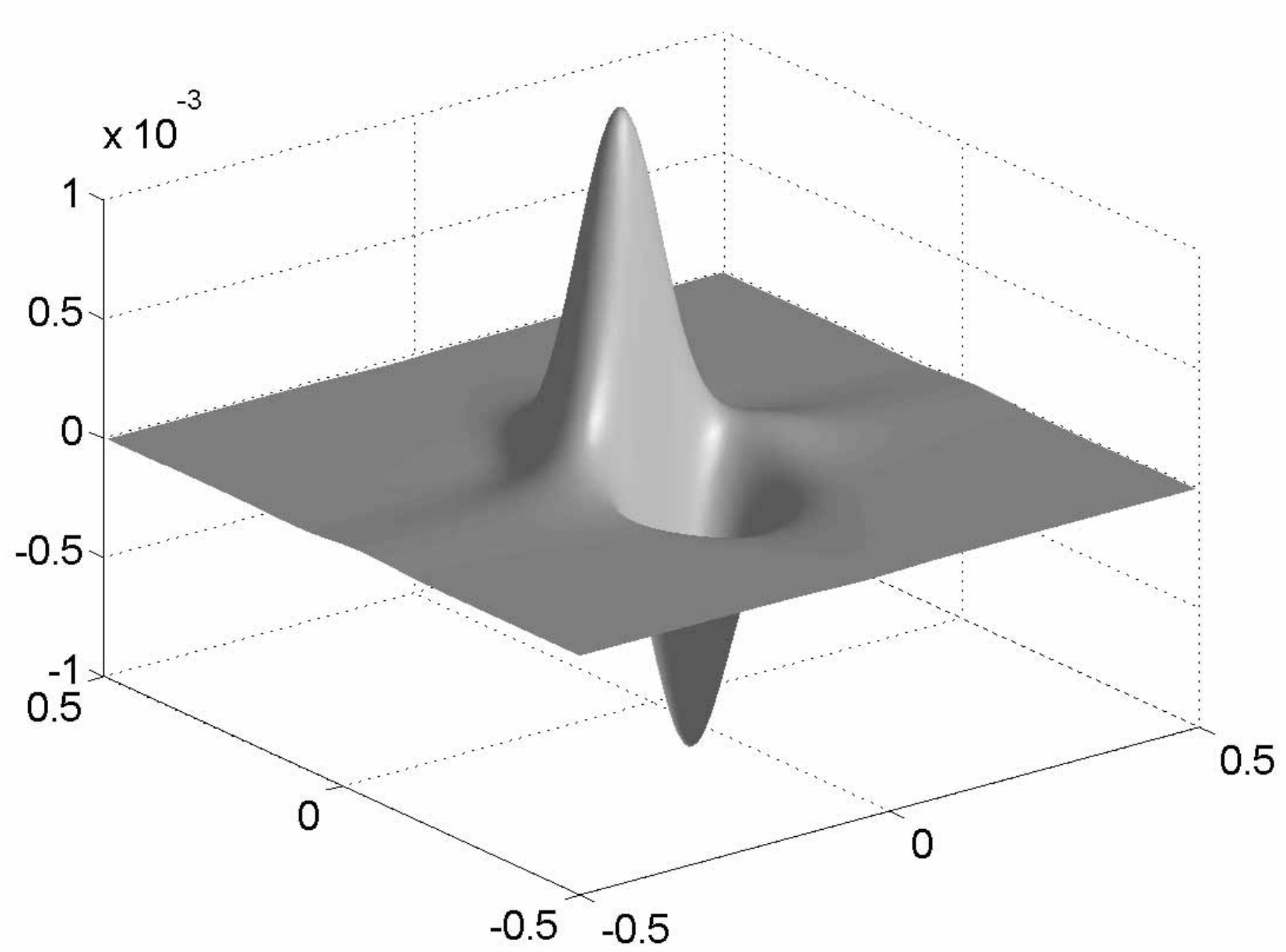}}
\resizebox{1.4in}{!} {\includegraphics{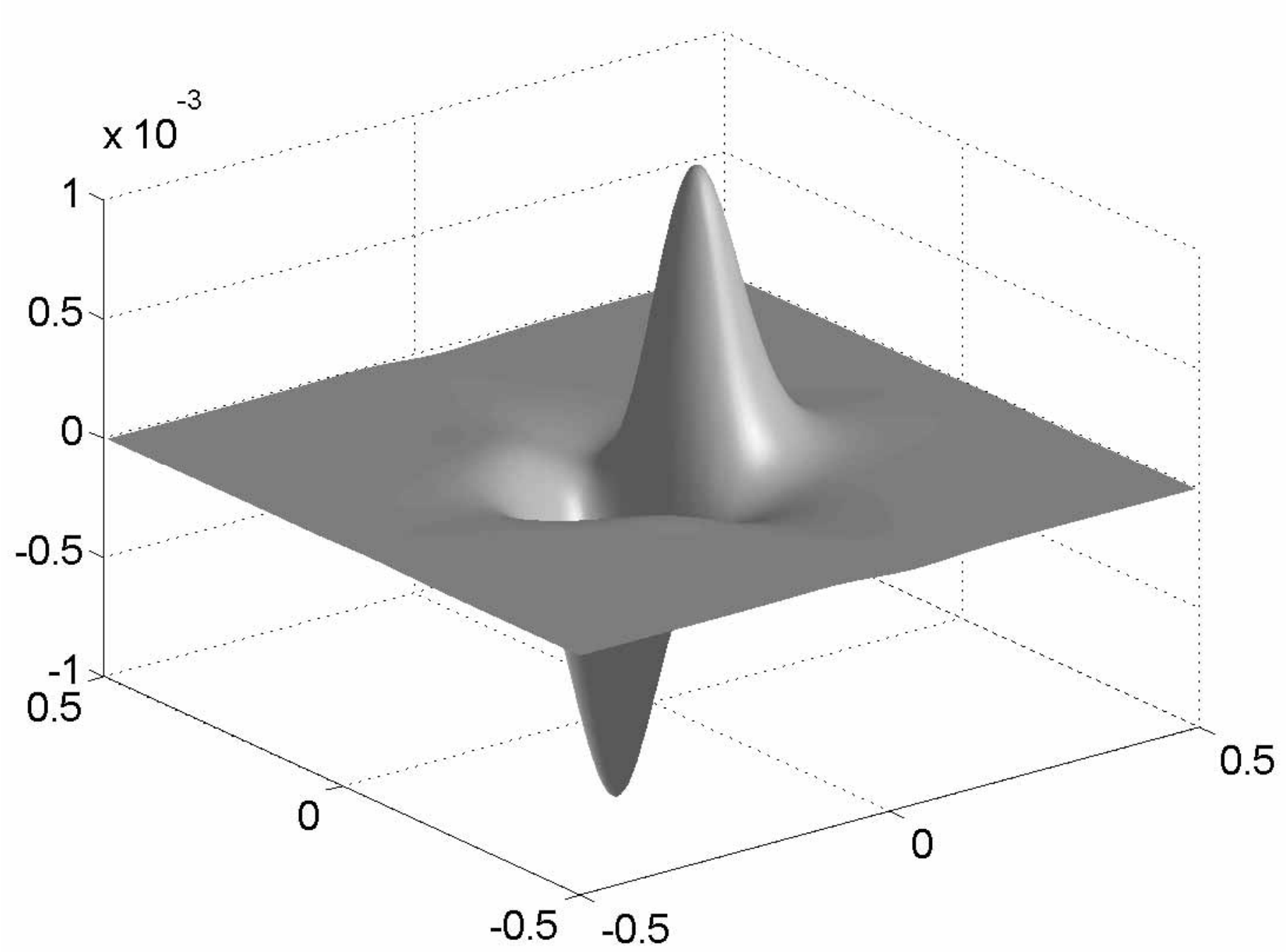}}
\resizebox{1.4in}{!} {\includegraphics{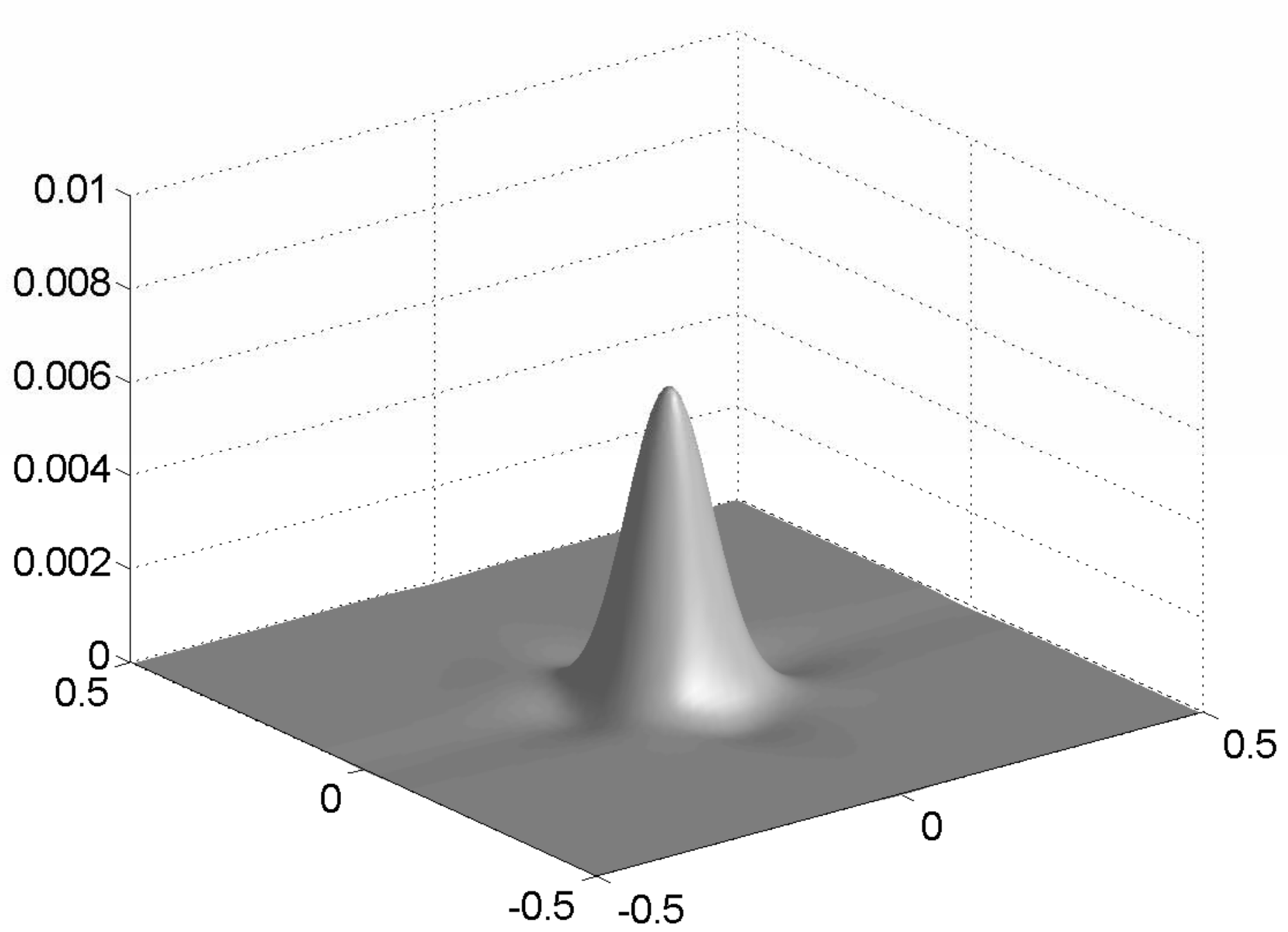}}

{$A_1^\dt(t,\xb)\big|_{x_3=0}$, $A_2^\dt(t,\xb)\big|_{x_3=0}$, $A_3^\dt(t,\xb)\big|_{x_3=0}$
for $\dt=1.0$.}\vspace{1mm}

\resizebox{1.4in}{!} {\includegraphics{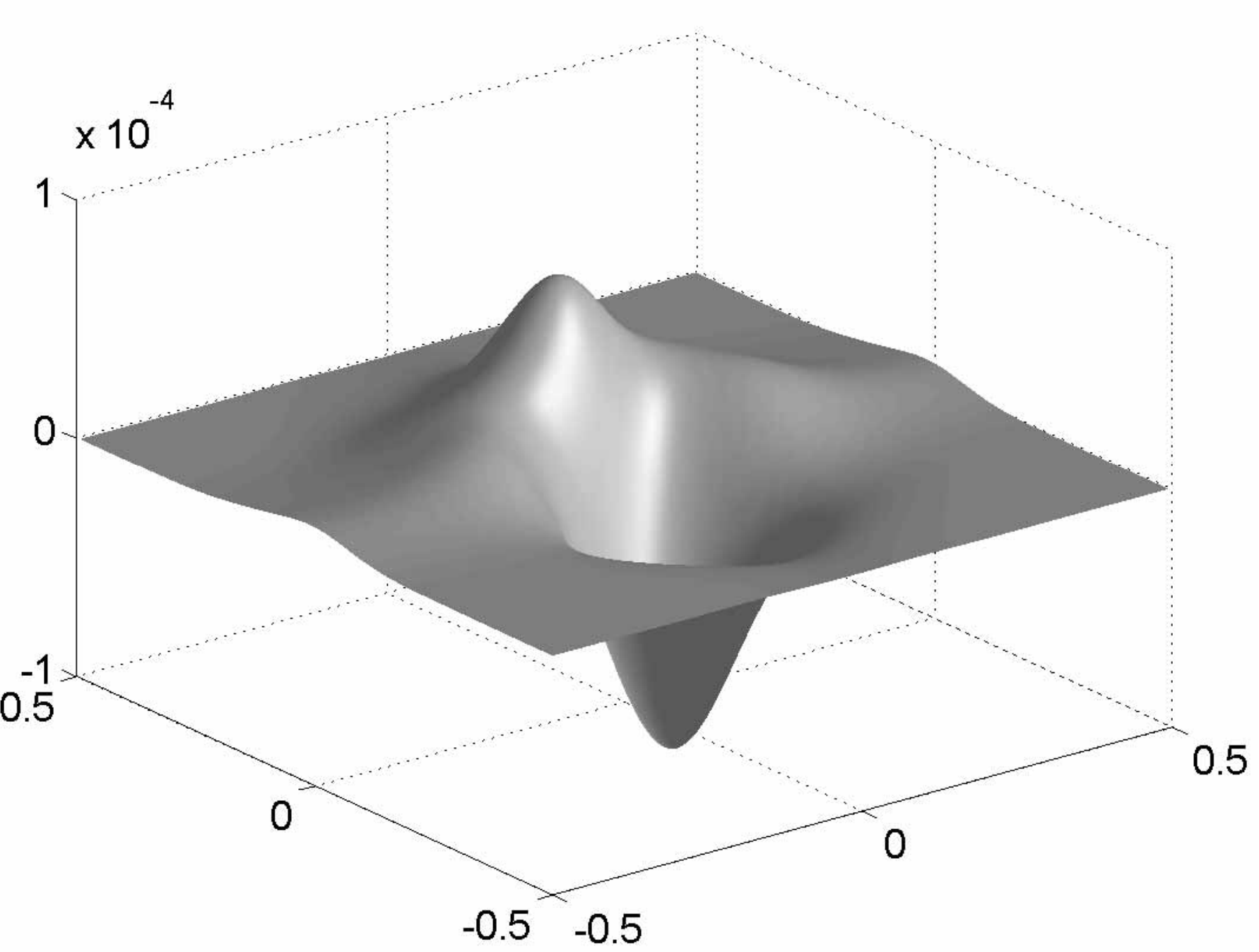}}
\resizebox{1.4in}{!} {\includegraphics{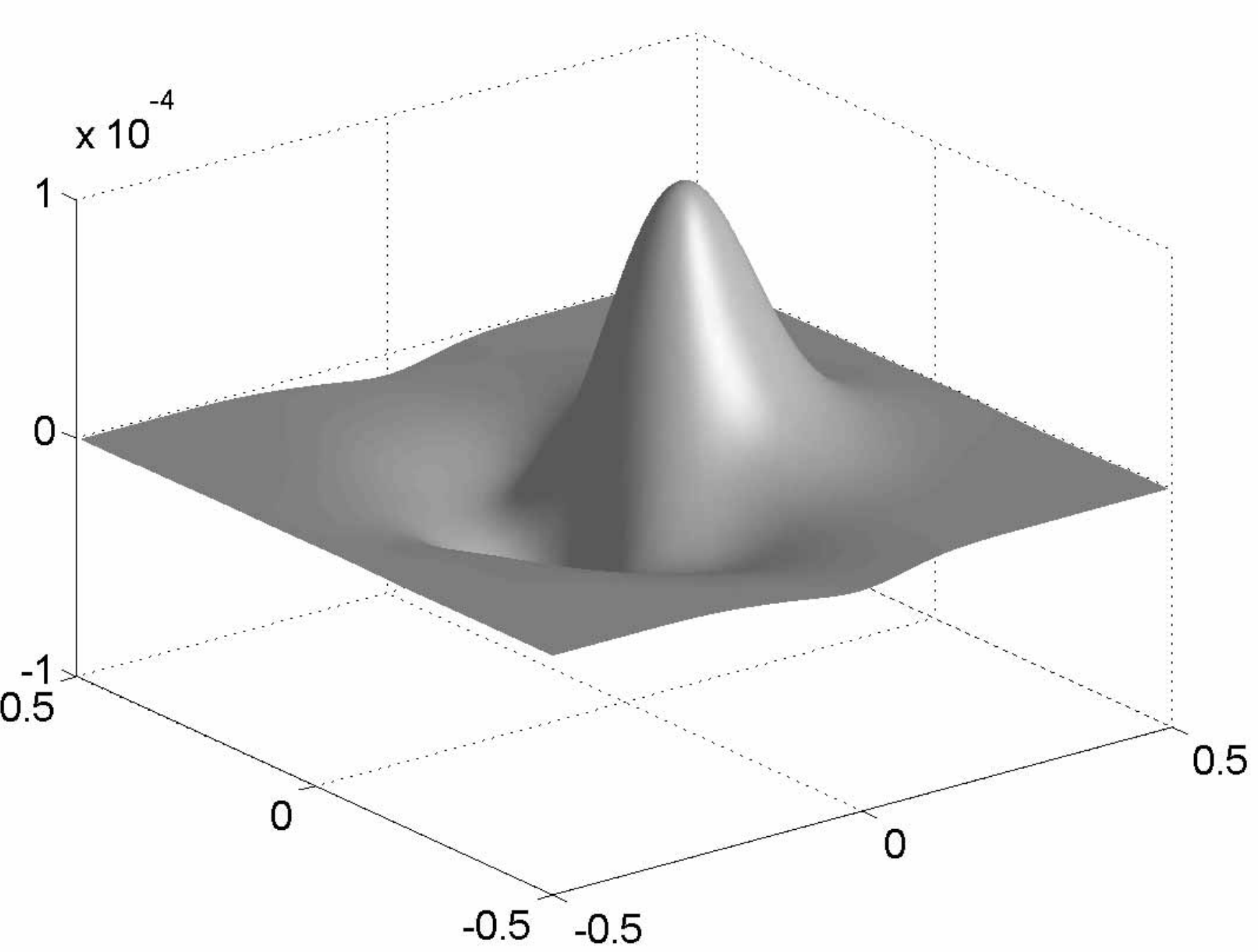}}
\resizebox{1.4in}{!} {\includegraphics{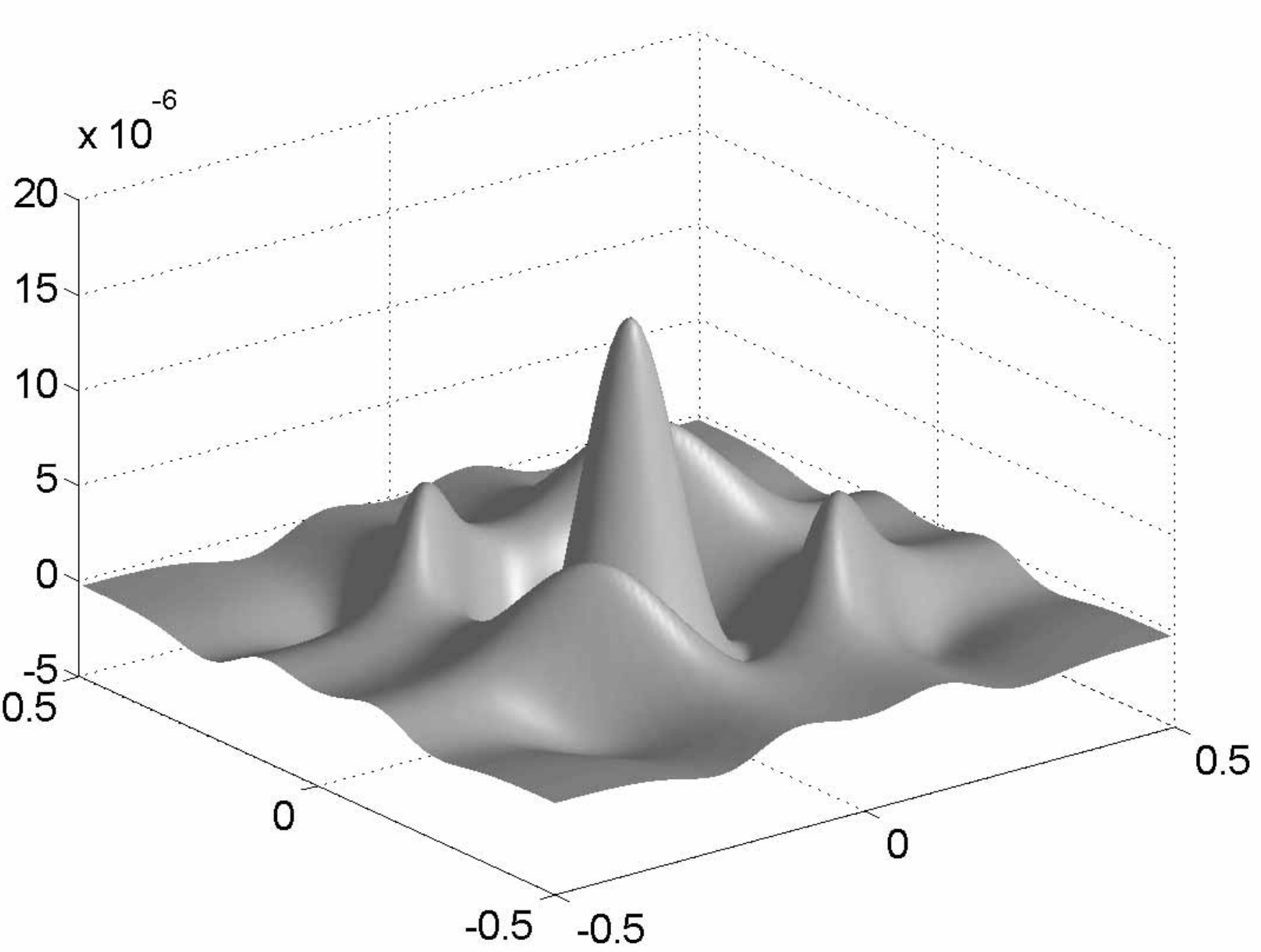}}

{$A_1^\dt(t,\xb)\big|_{x_3=0}$, $A_2^\dt(t,\xb)\big|_{x_3=0}$, $A_3^\dt(t,\xb)\big|_{x_3=0}$
for $\dt=0.01$.}
\end{center}
\caption{Numerical results of the magnetic fields for example \ref{exnr} at t=0.5.
The first row is the solution of the MD system with $\dt=1.0$, whereas
the second row is the solution of the MD system with $\dt=0.01$.
}\label{fig511}
\end{figure}
\begin{figure}
\begin{center}
\resizebox{1.4in}{!} {\includegraphics{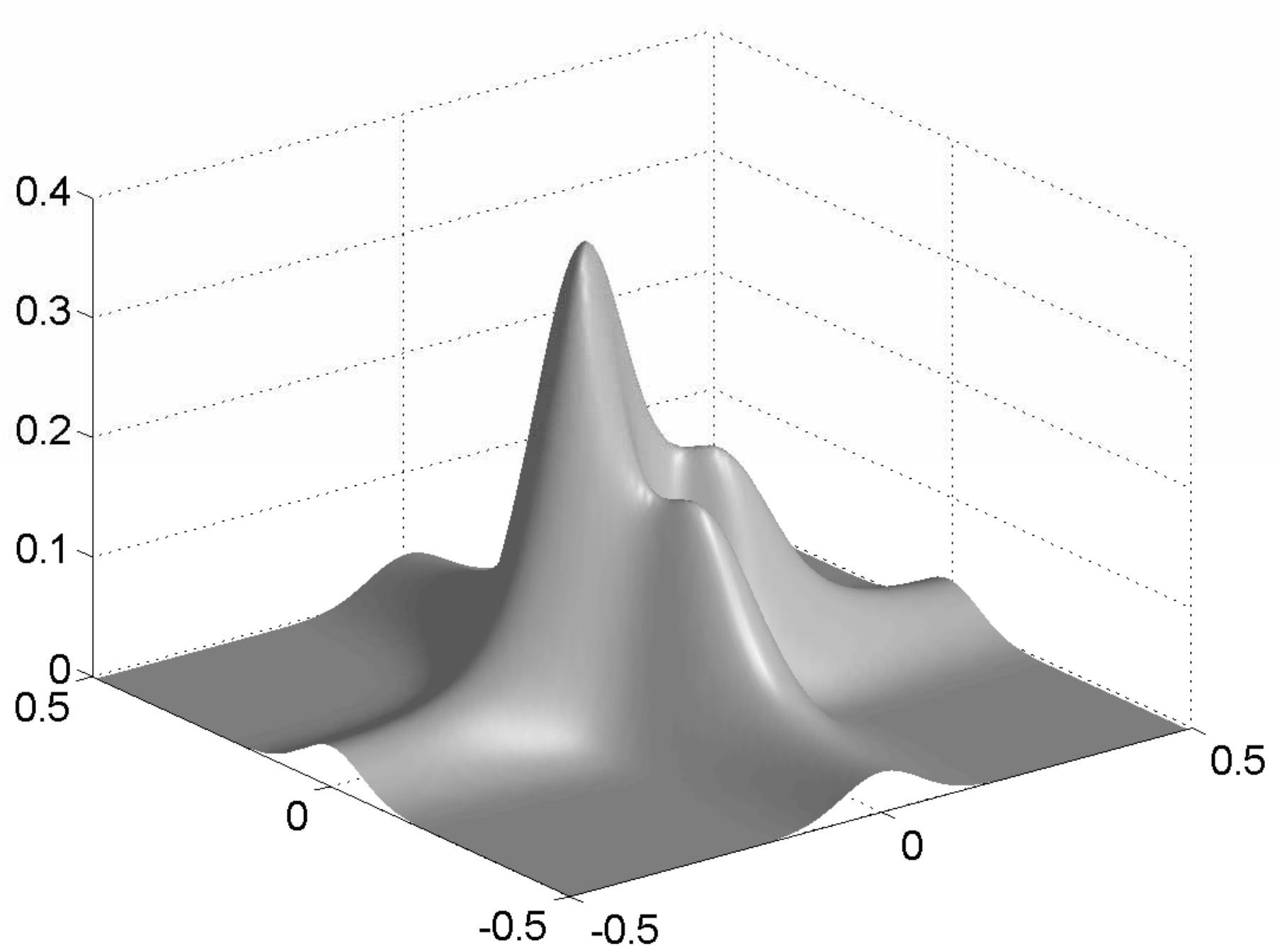}}
\resizebox{1.4in}{!} {\includegraphics{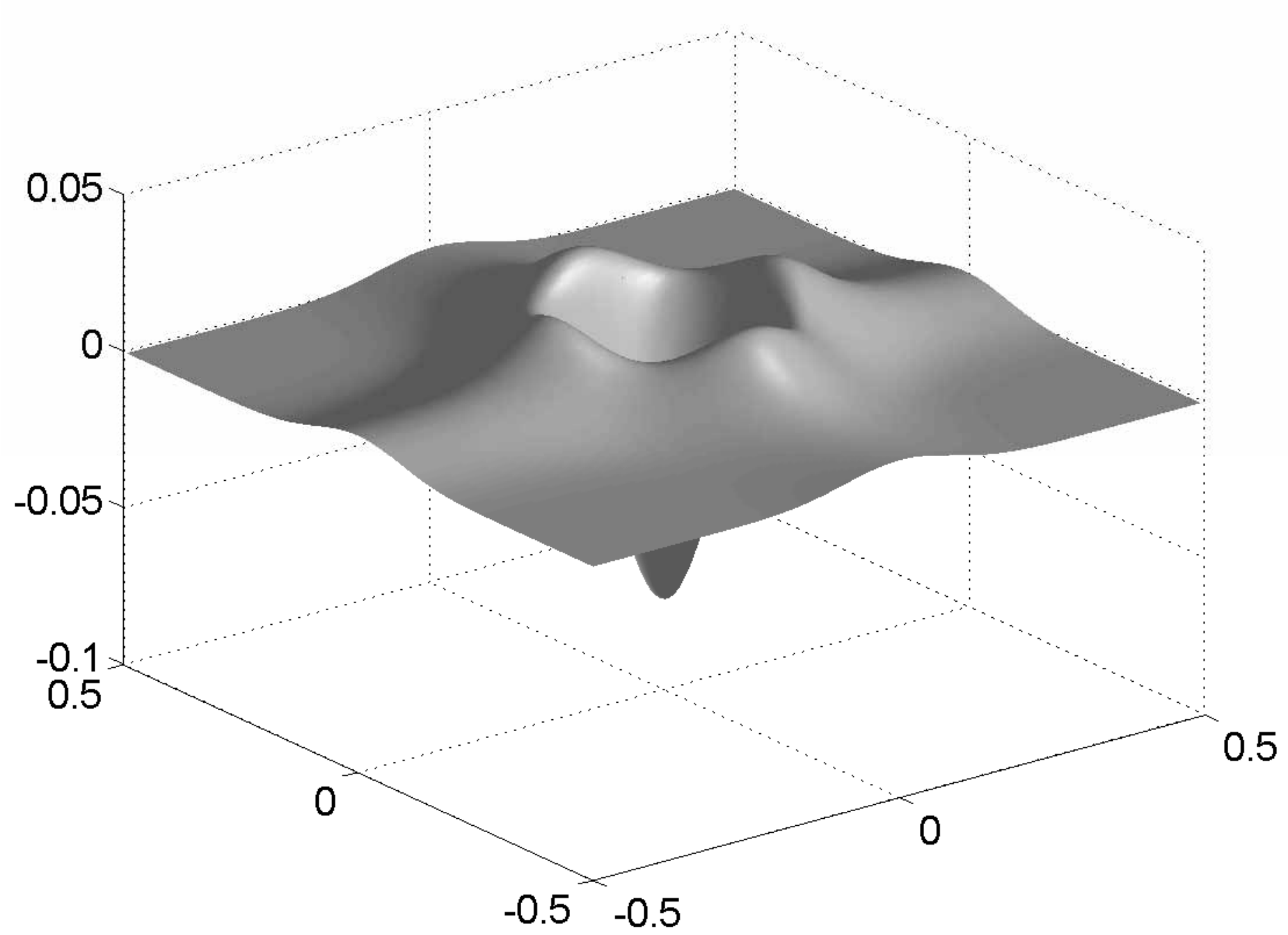}}
\resizebox{1.4in}{!} {\includegraphics{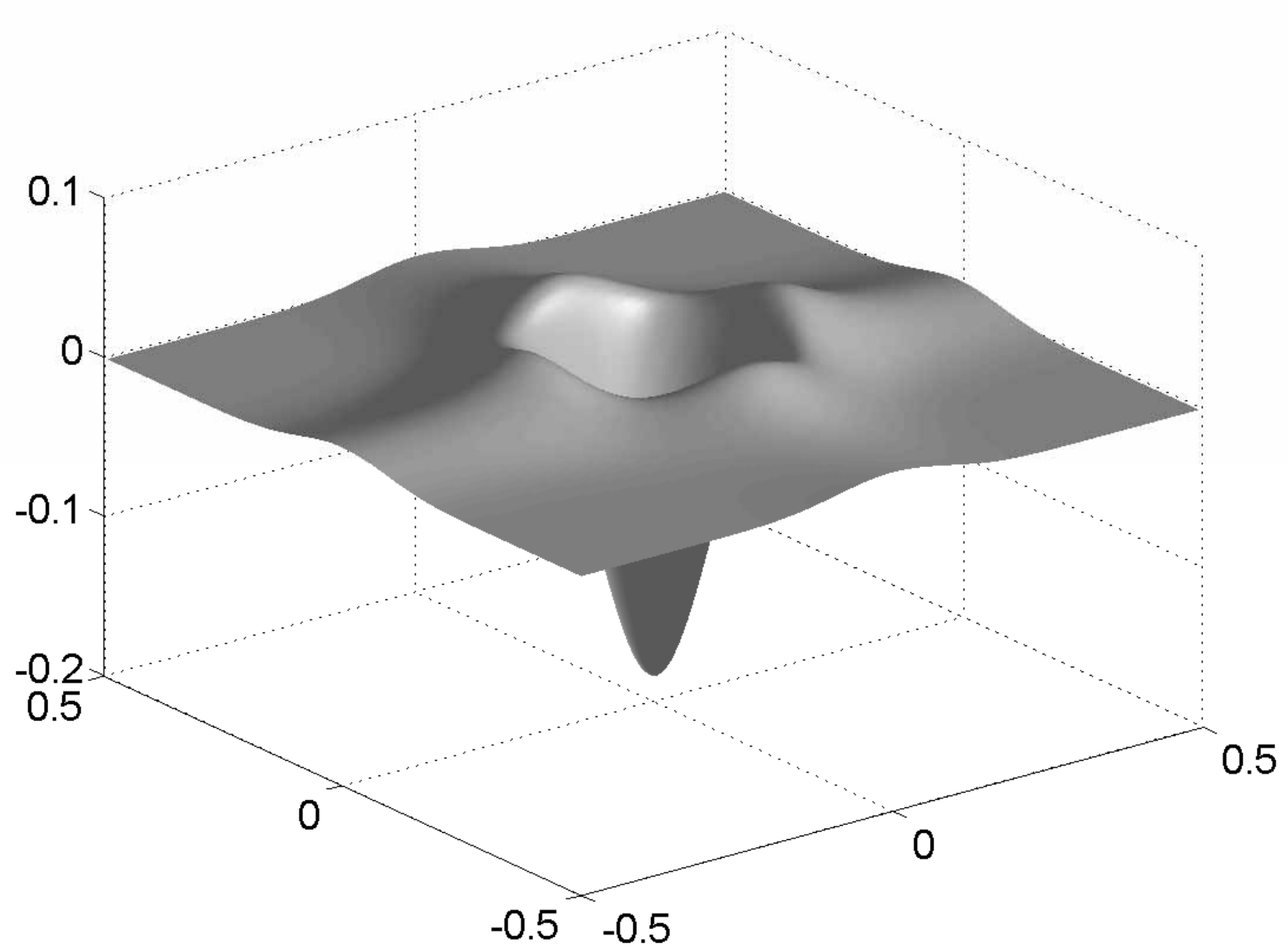}}
\resizebox{1.4in}{!} {\includegraphics{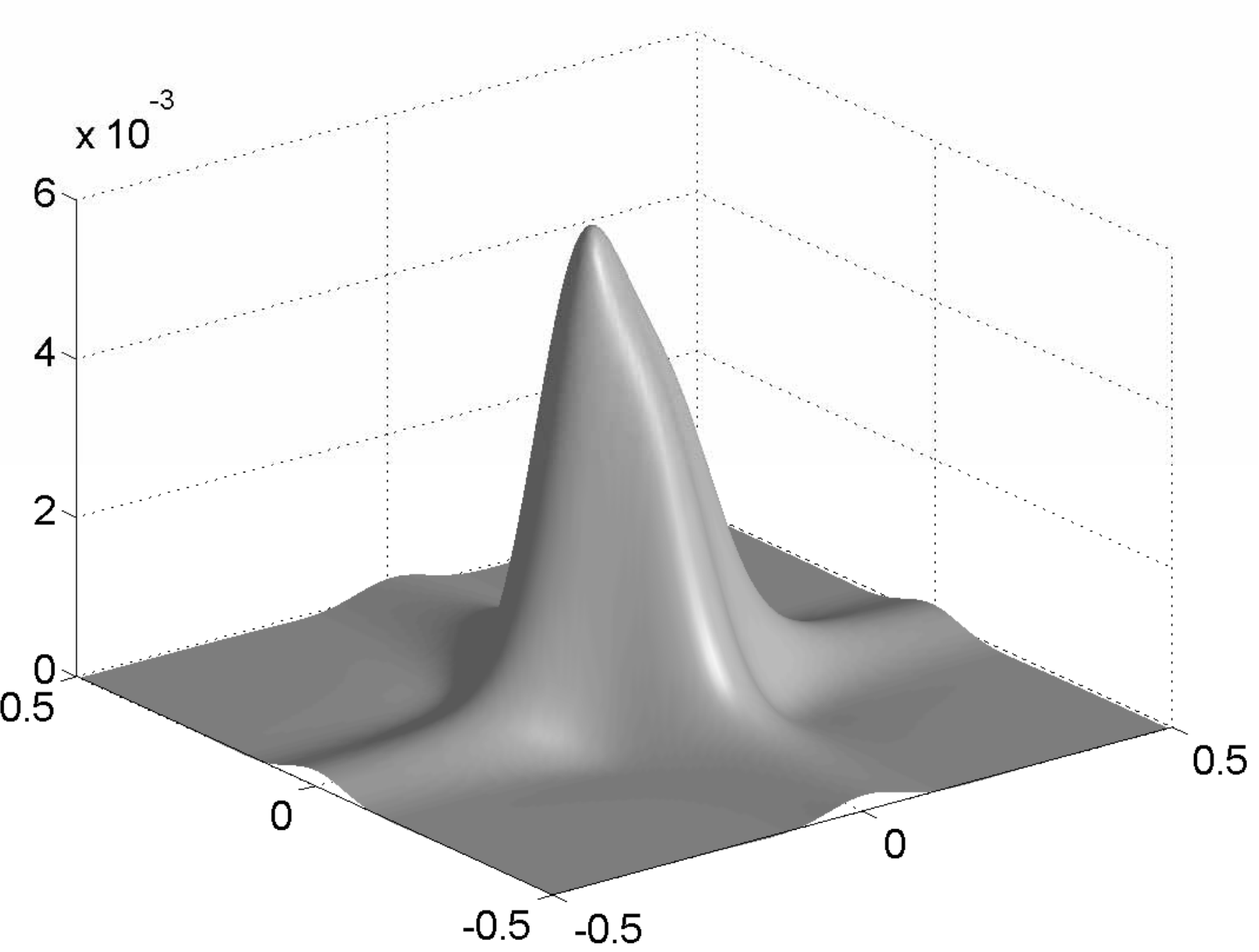}}

{$|\psi^\dt(t,\xb)|^2\big|_{x_3=0}$, $\re\big(\psi_{e,1}^\dt(t,\xb)\big)\big|_{x_3=0}$,
$\im\big(\psi_{e,1}^\dt(t,\xb)\big)\big|_{x_3=0}$ and $V^\dt(t,\xb)\big|_{x_3=0}$ for $\dt=1.0$.}\vspace{1mm}

\resizebox{1.4in}{!} {\includegraphics{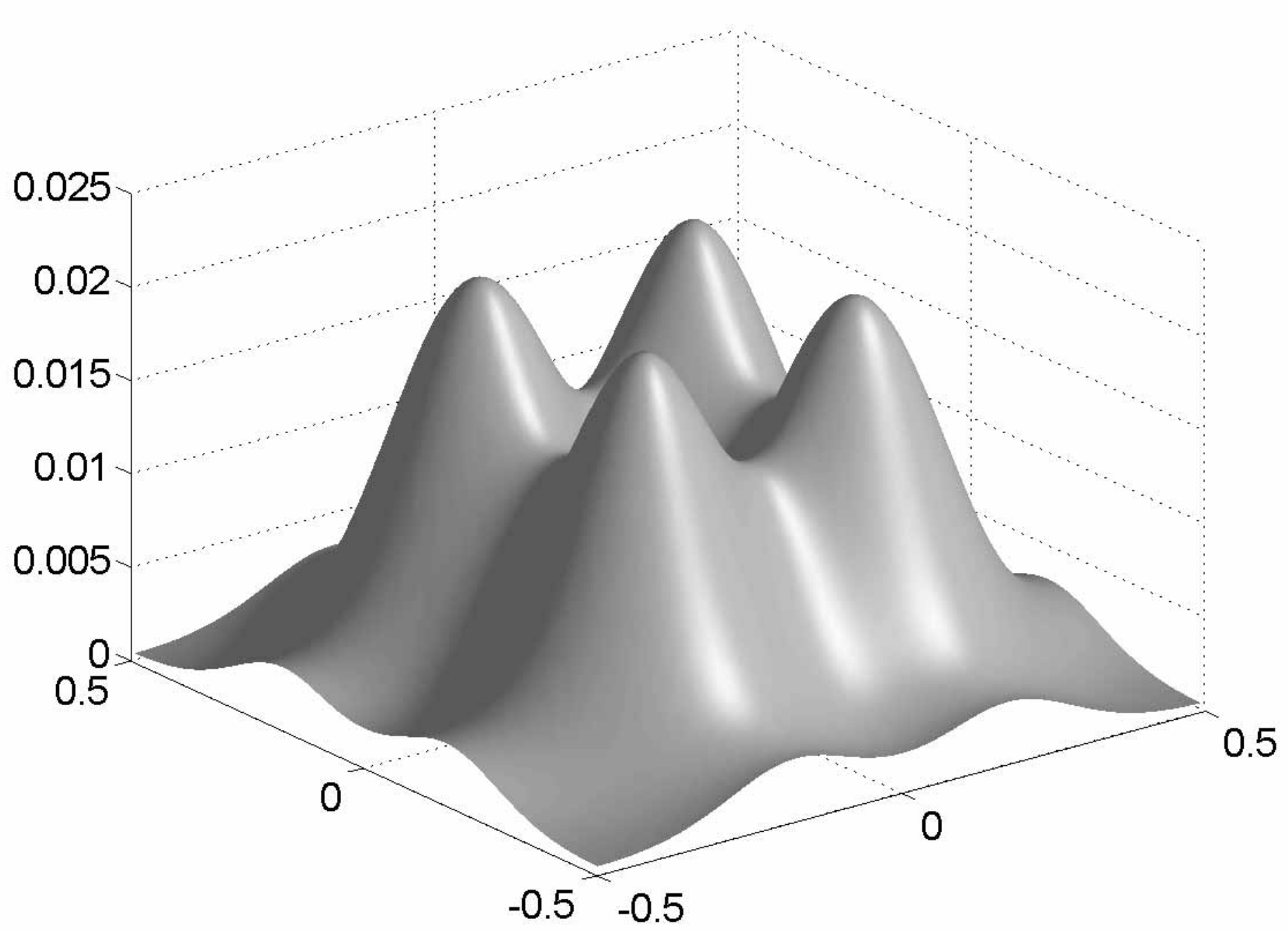}}
\resizebox{1.4in}{!} {\includegraphics{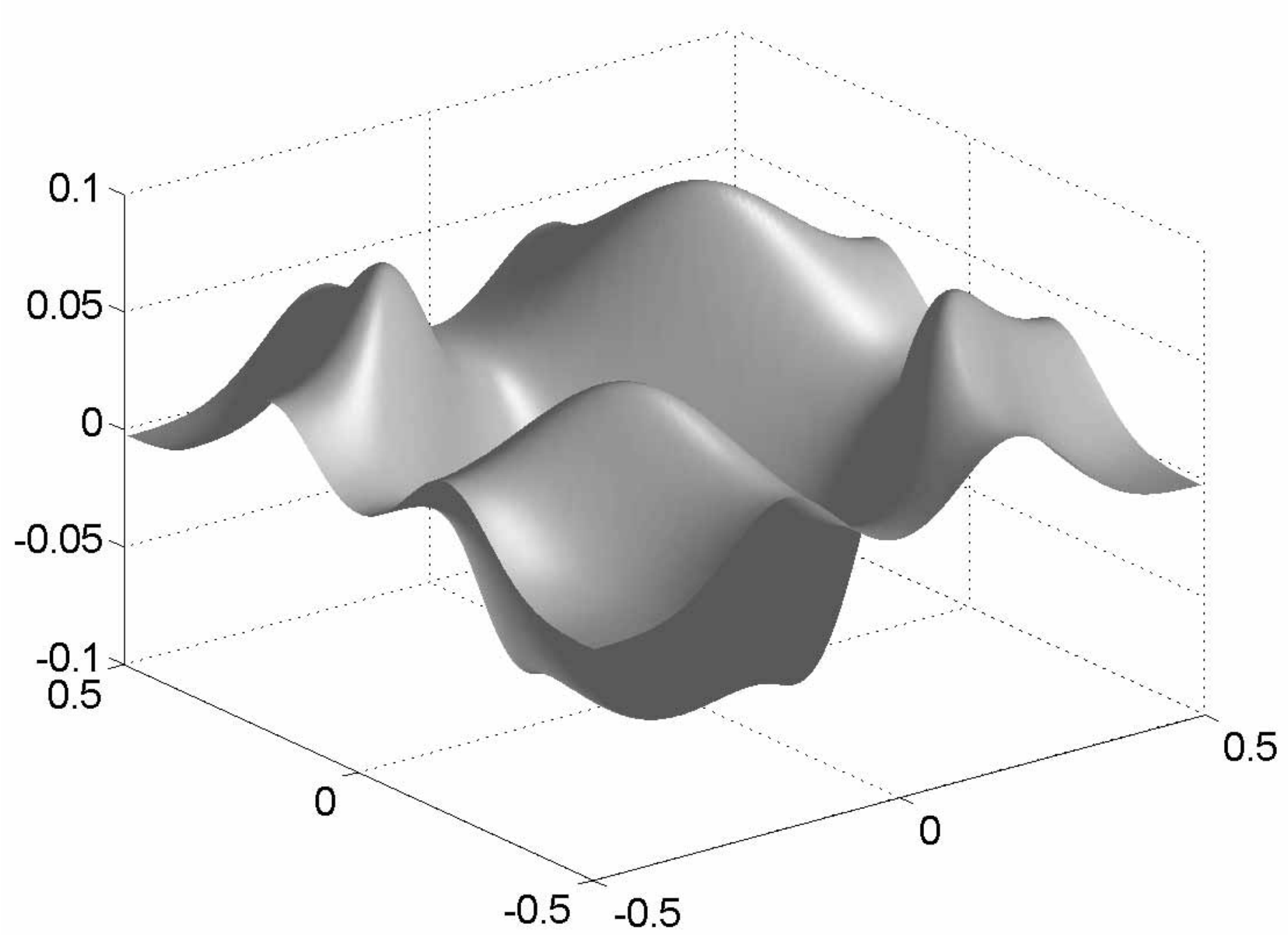}}
\resizebox{1.4in}{!} {\includegraphics{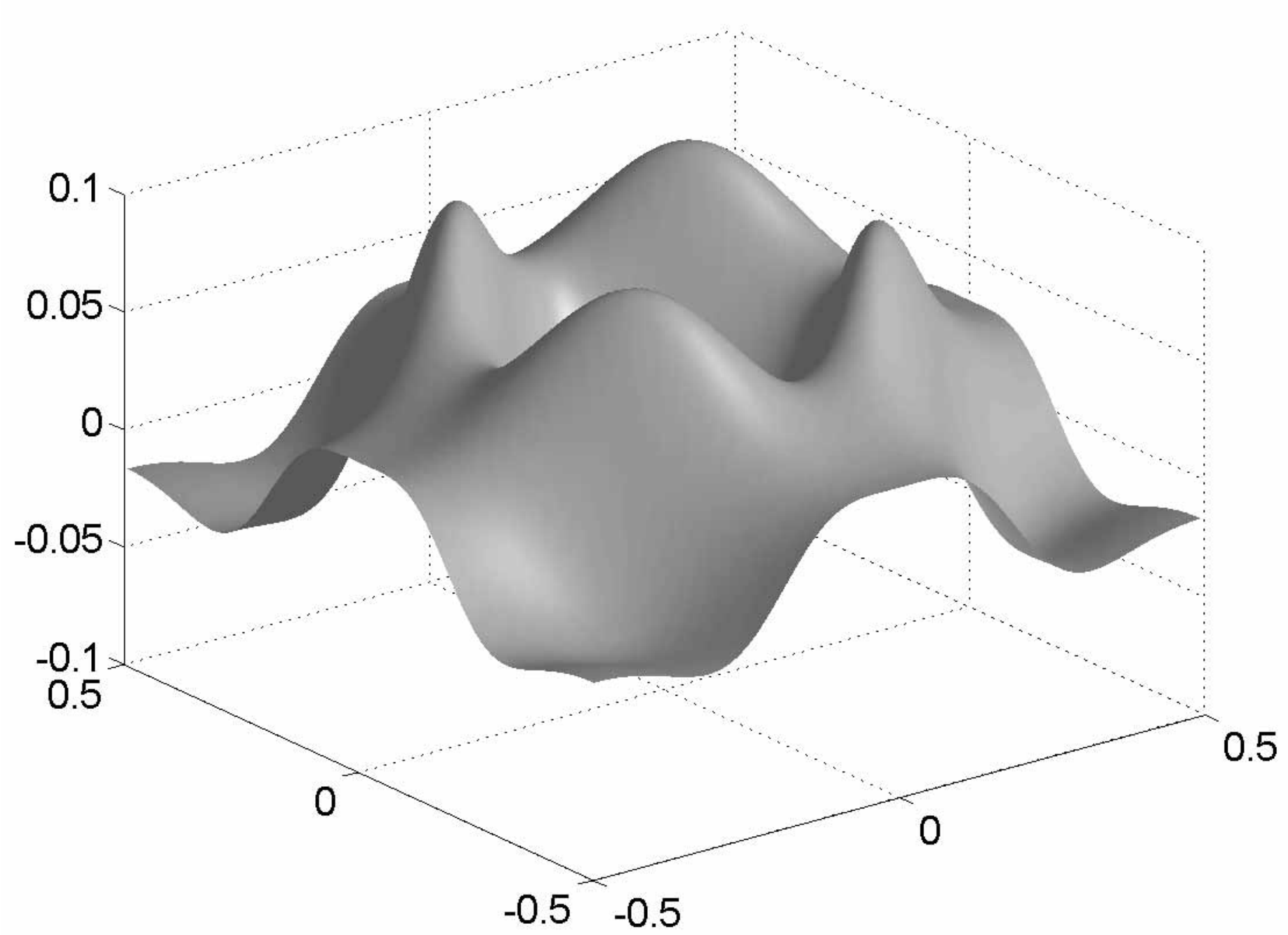}}
\resizebox{1.4in}{!} {\includegraphics{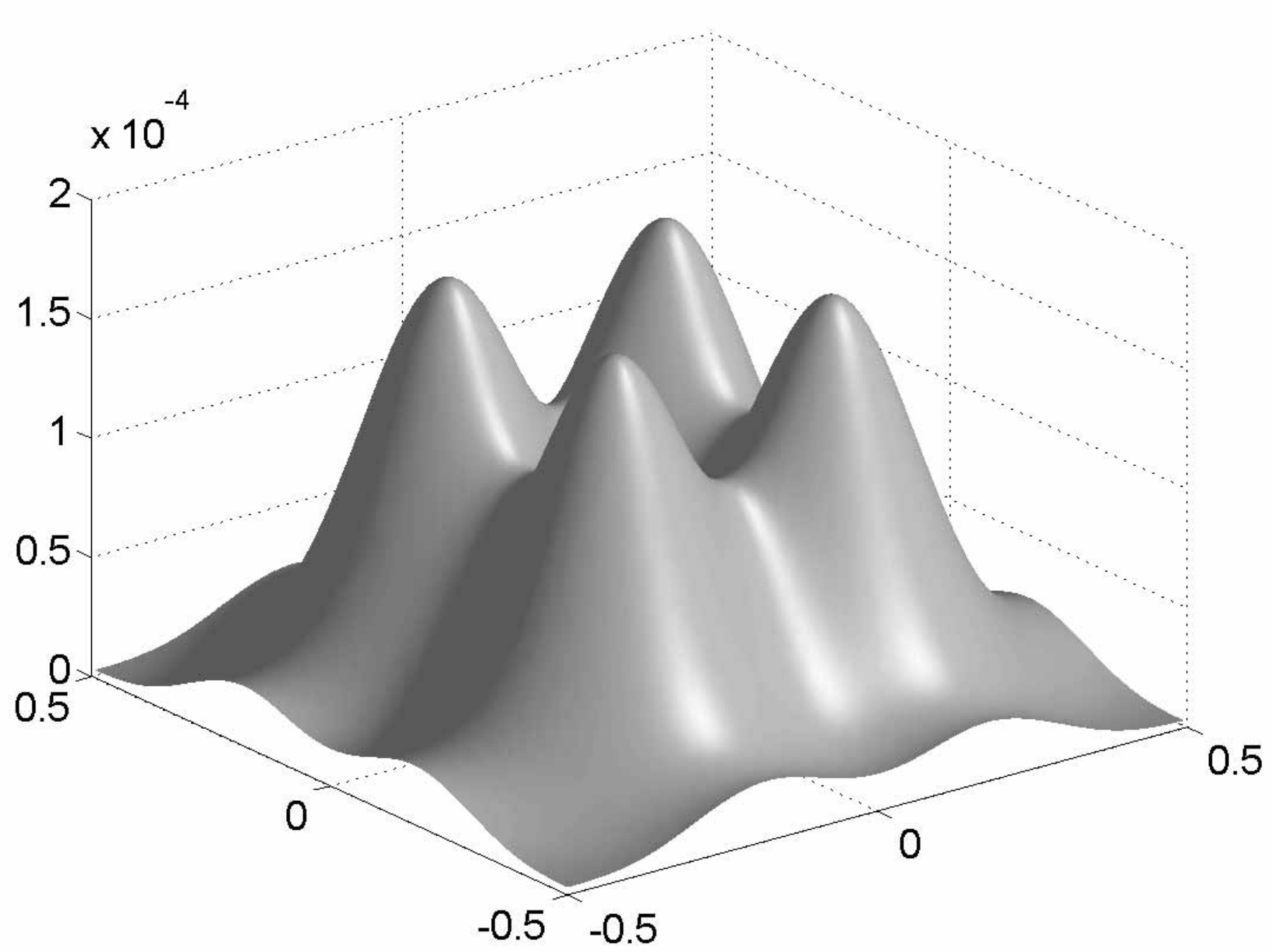}}

{$|\psi^\dt(t,\xb)|^2\big|_{x_3=0}$, $\re\big(\psi_{e,1}^\dt(t,\xb)\big)\big|_{x_3=0}$,
$\im\big(\psi_{e,1}^\dt(t,\xb)\big)\big|_{x_3=0}$ and $V^\dt(t,\xb)\big|_{x_3=0}$ for $\dt=0.01$.}\vspace{1mm}

\resizebox{1.4in}{!} {\includegraphics{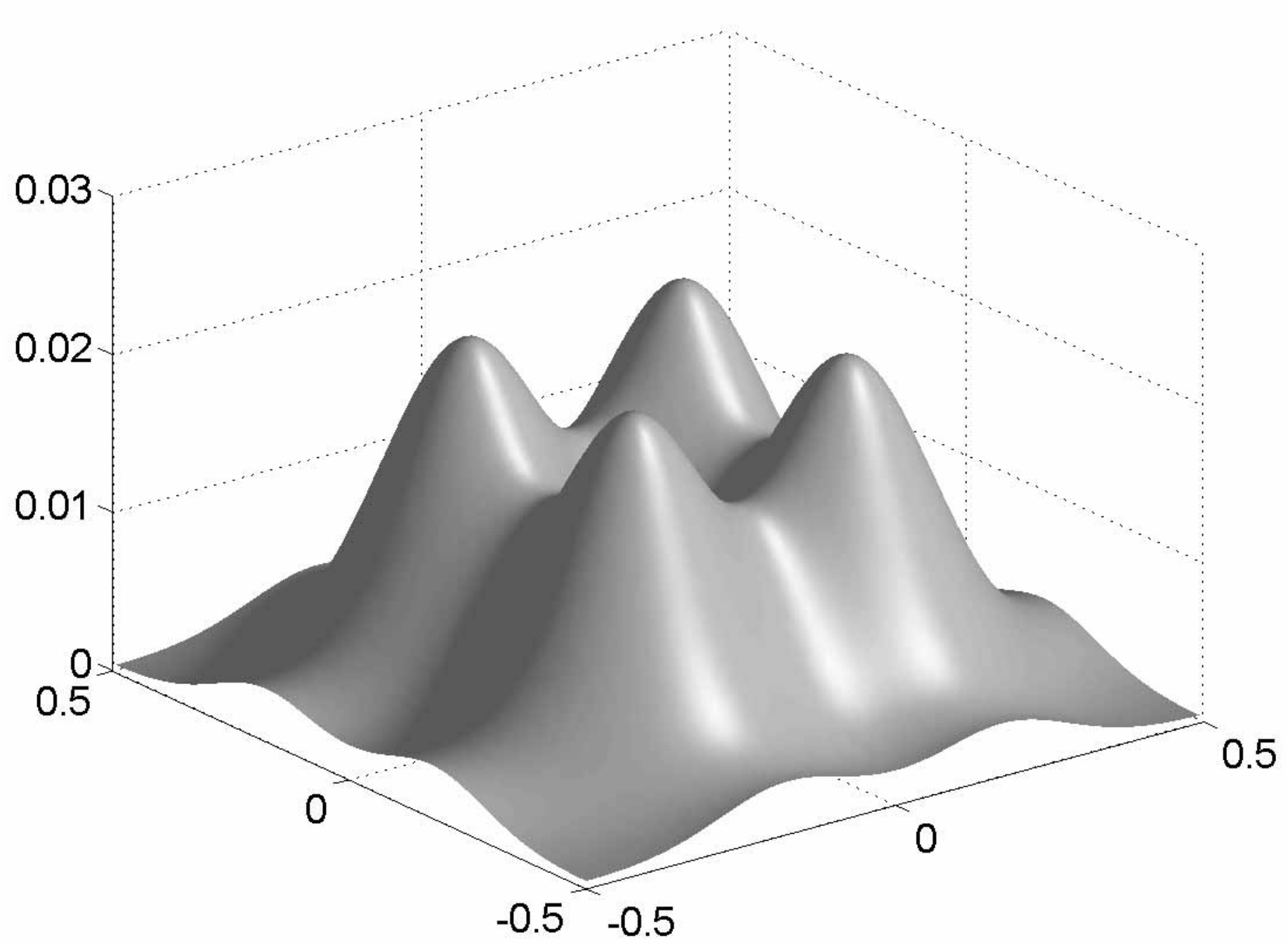}}
\resizebox{1.4in}{!} {\includegraphics{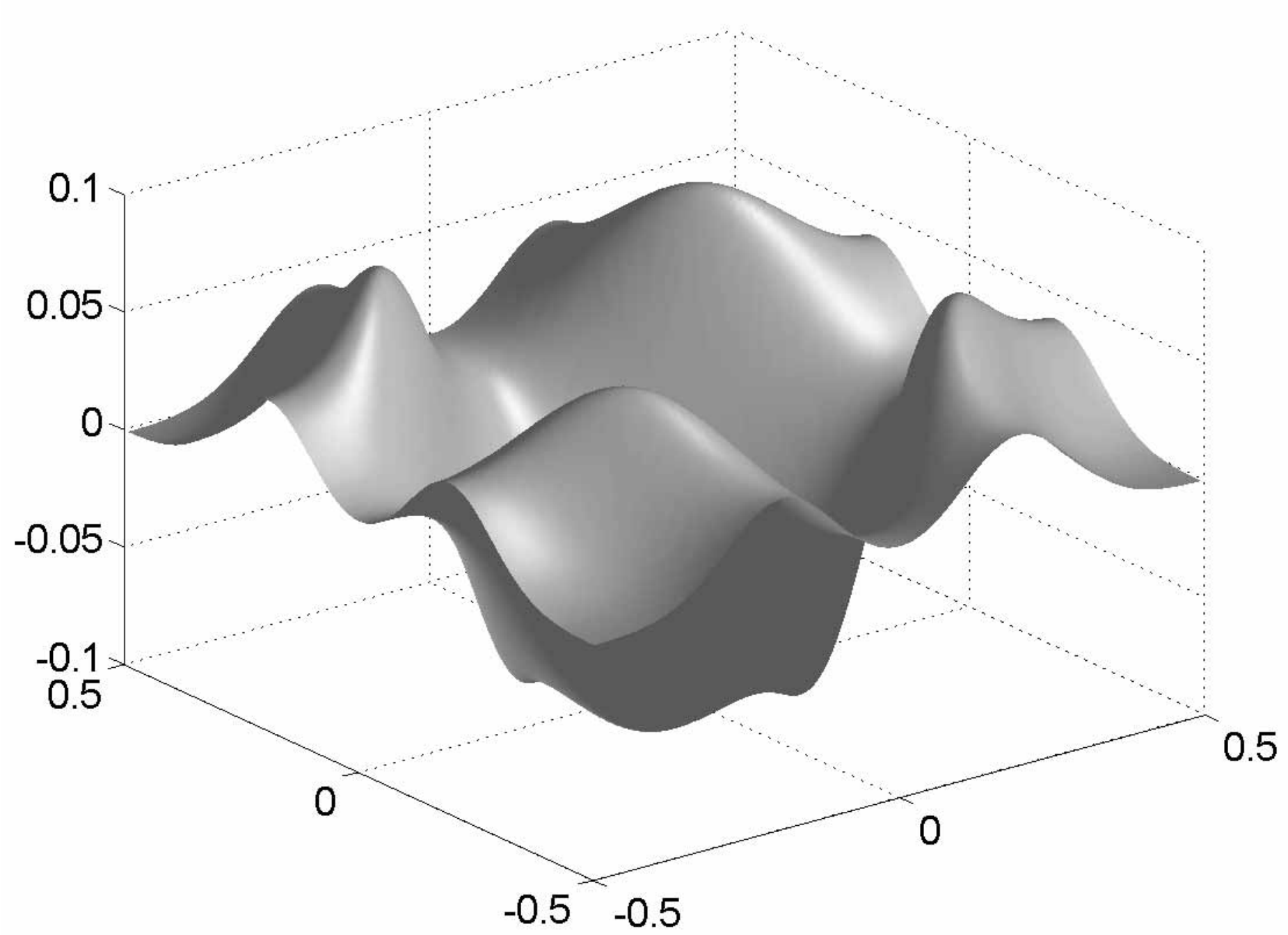}}
\resizebox{1.4in}{!} {\includegraphics{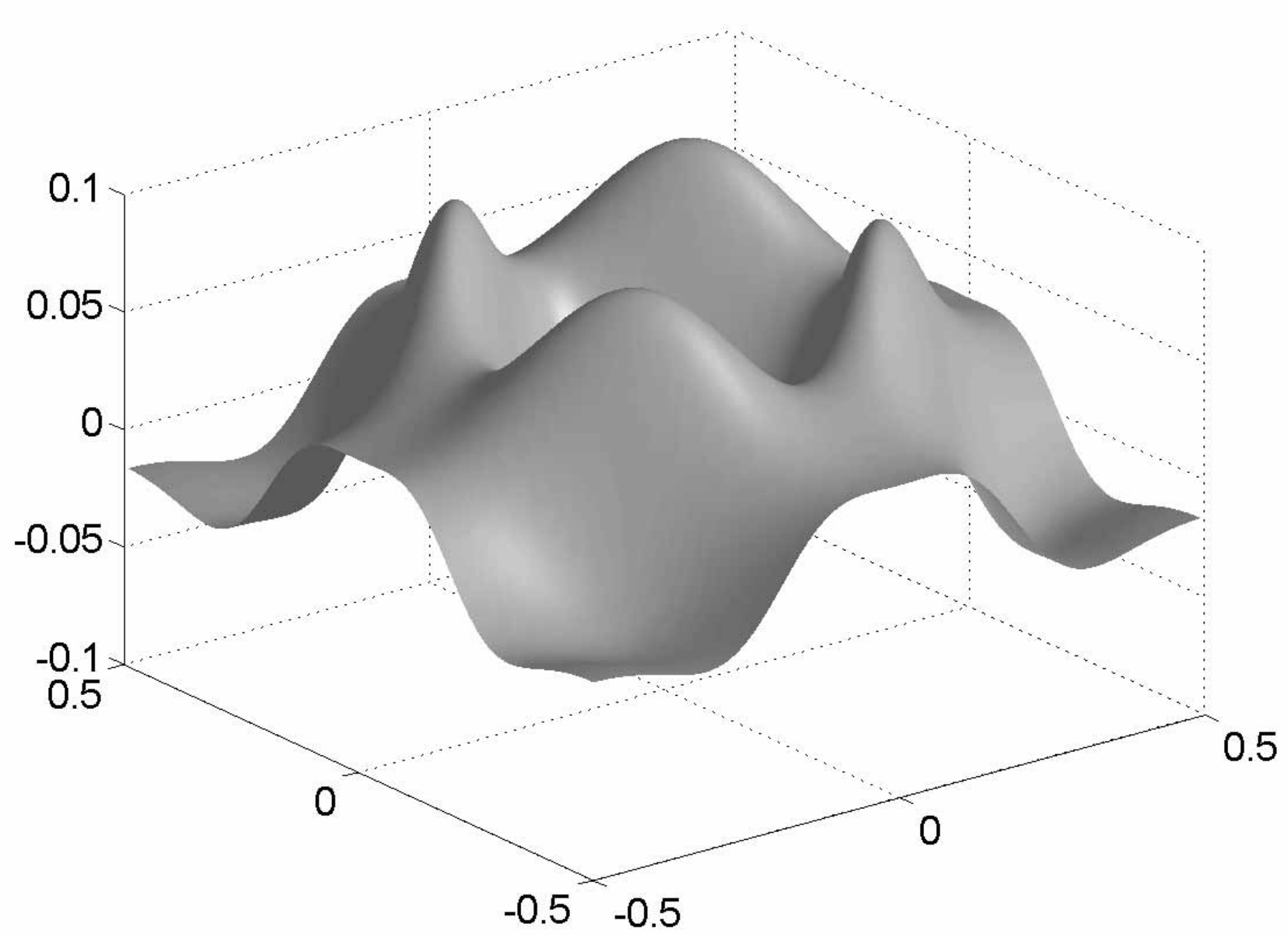}}
\resizebox{1.4in}{!} {\includegraphics{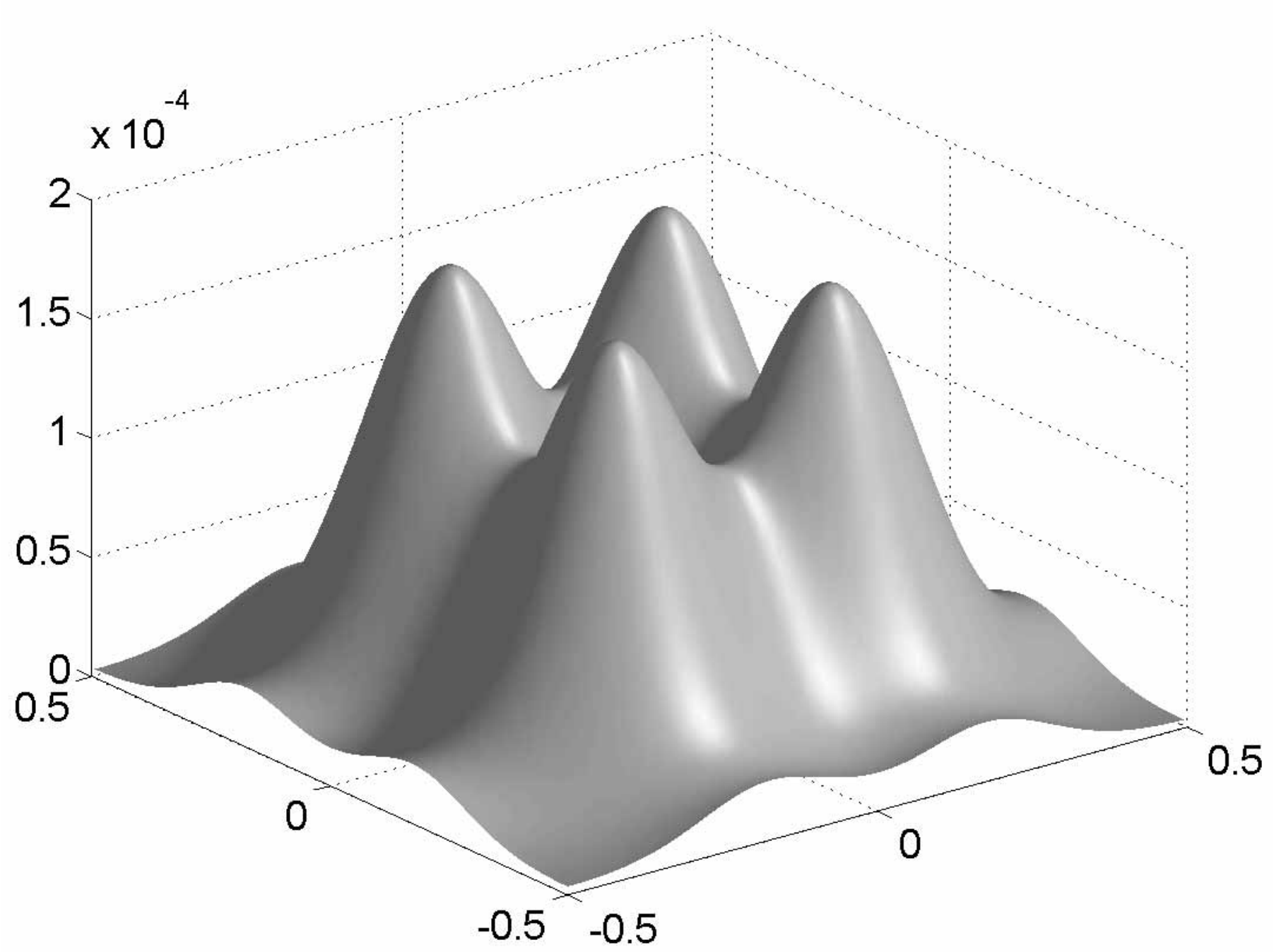}}

{$\left(|\vp_{e}(t,\xb)|^2+|\vp_{p}(t,\xb)|^2\right)\big|_{x_3=0}$, $\re\big(\vp_{e,1}(t,\xb)\big)\big|_{x_3=0}$,
$\im\big(\vp_{e,1}(t,\xb)\big)\big|_{x_3=0}$ and $V(t,\xb)\big|_{x_3=0}$.}
\end{center}
\caption{Numerical results for example \ref{exnr} at t=1.0.
The first row is the solution of the MD system with $\dt=1.0$,
the second row is the solution of the MD system with $\dt=0.01$,
and the third row is the solution of the Schr\"odinger-Poisson problem.
} \label{fig52}
\begin{center}
\resizebox{1.4in}{!} {\includegraphics{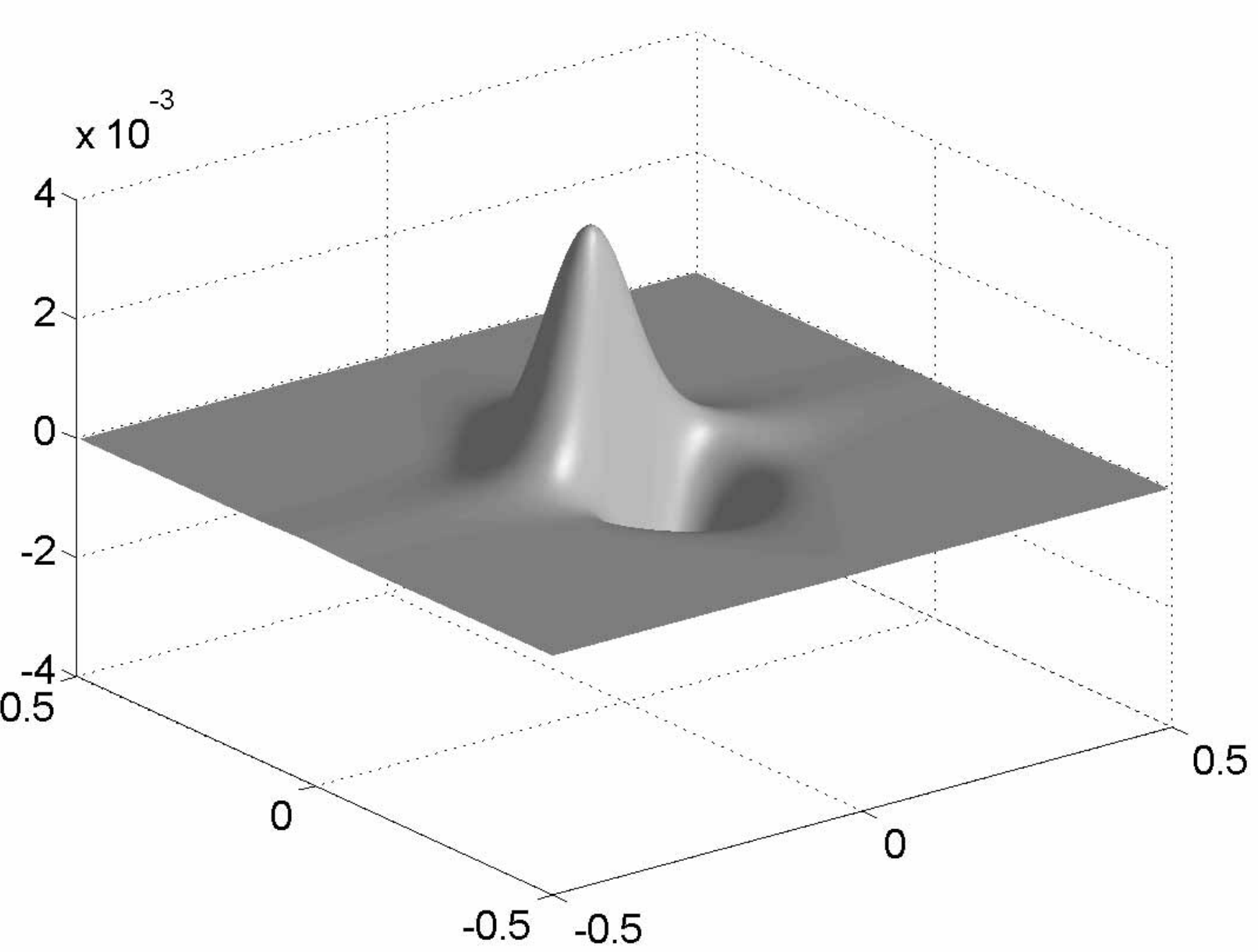}}
\resizebox{1.4in}{!} {\includegraphics{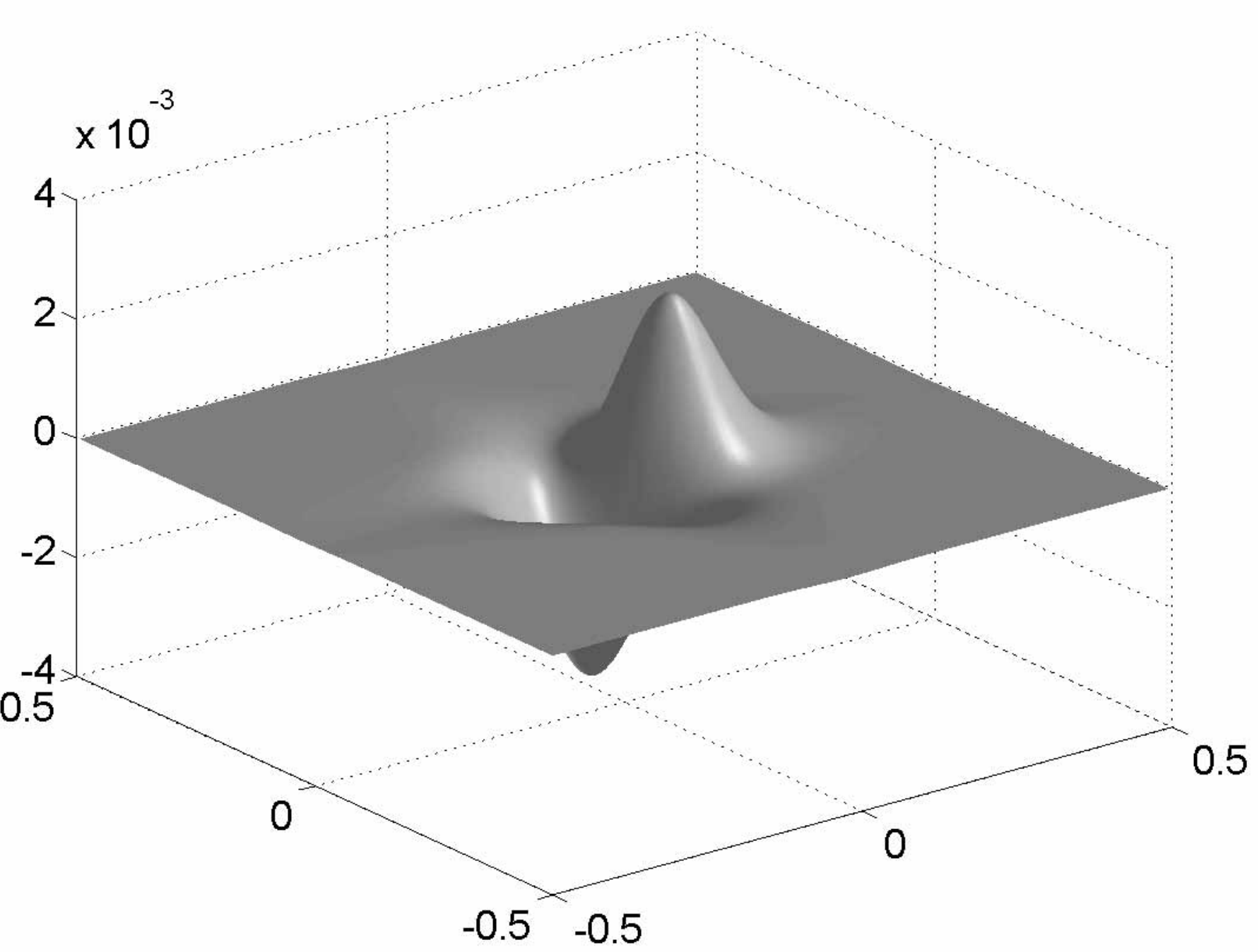}}
\resizebox{1.4in}{!} {\includegraphics{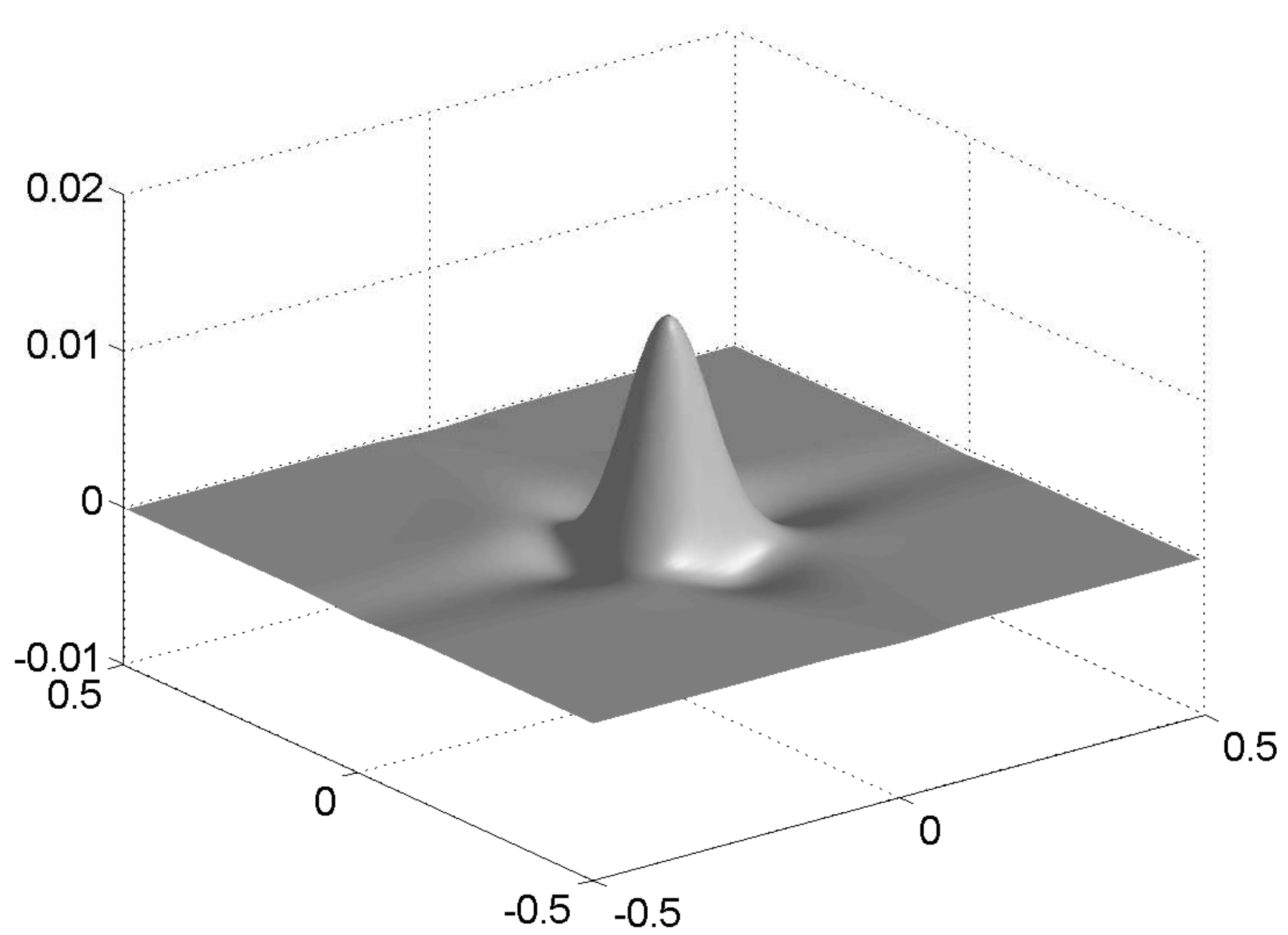}}

{$A_1^\dt(t,\xb)\big|_{x_3=0}$, $A_2^\dt(t,\xb)\big|_{x_3=0}$, $A_3^\dt(t,\xb)\big|_{x_3=0}$
for $\dt=1.0$.}\vspace{1mm}

\resizebox{1.4in}{!} {\includegraphics{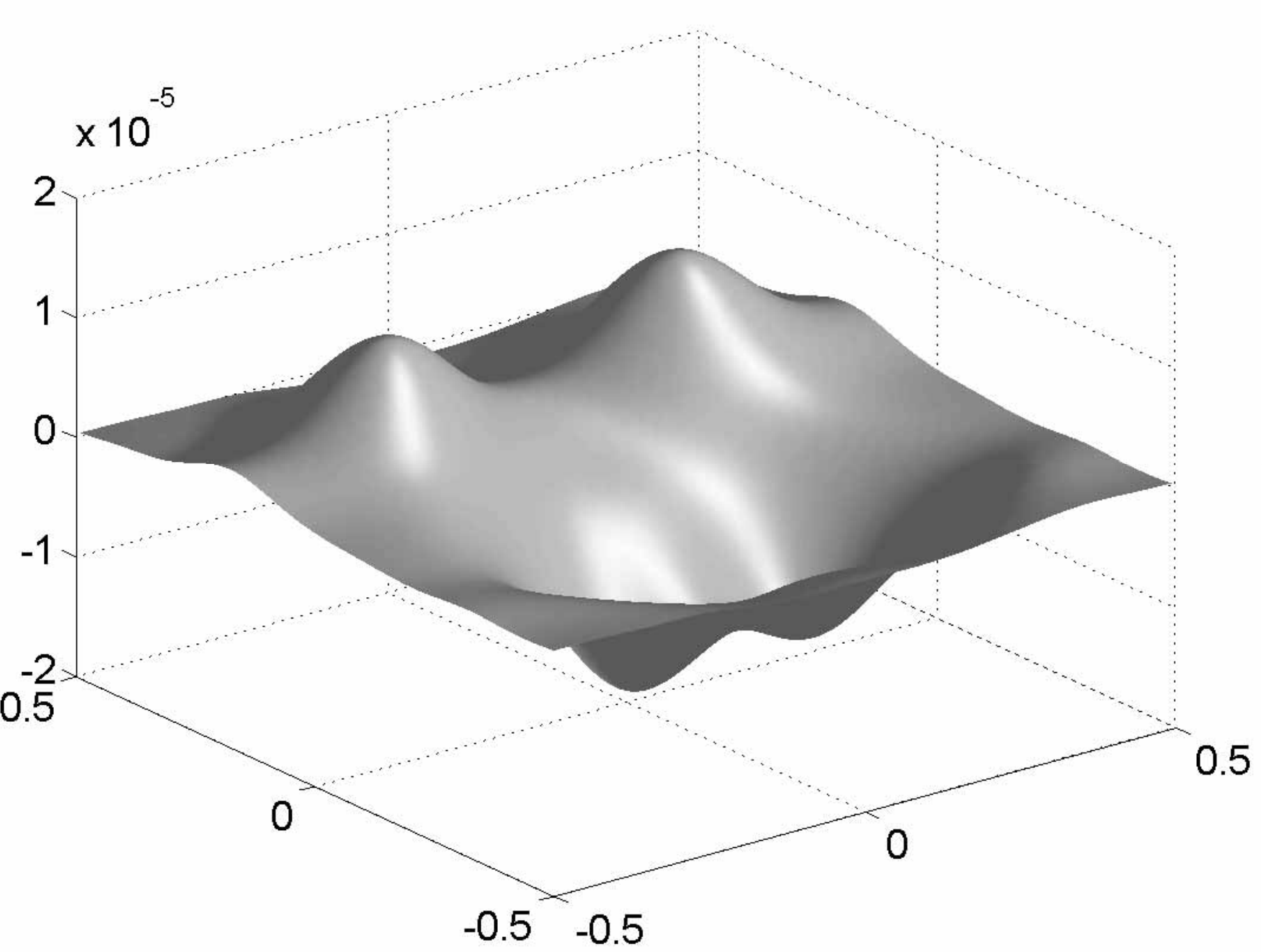}}
\resizebox{1.4in}{!} {\includegraphics{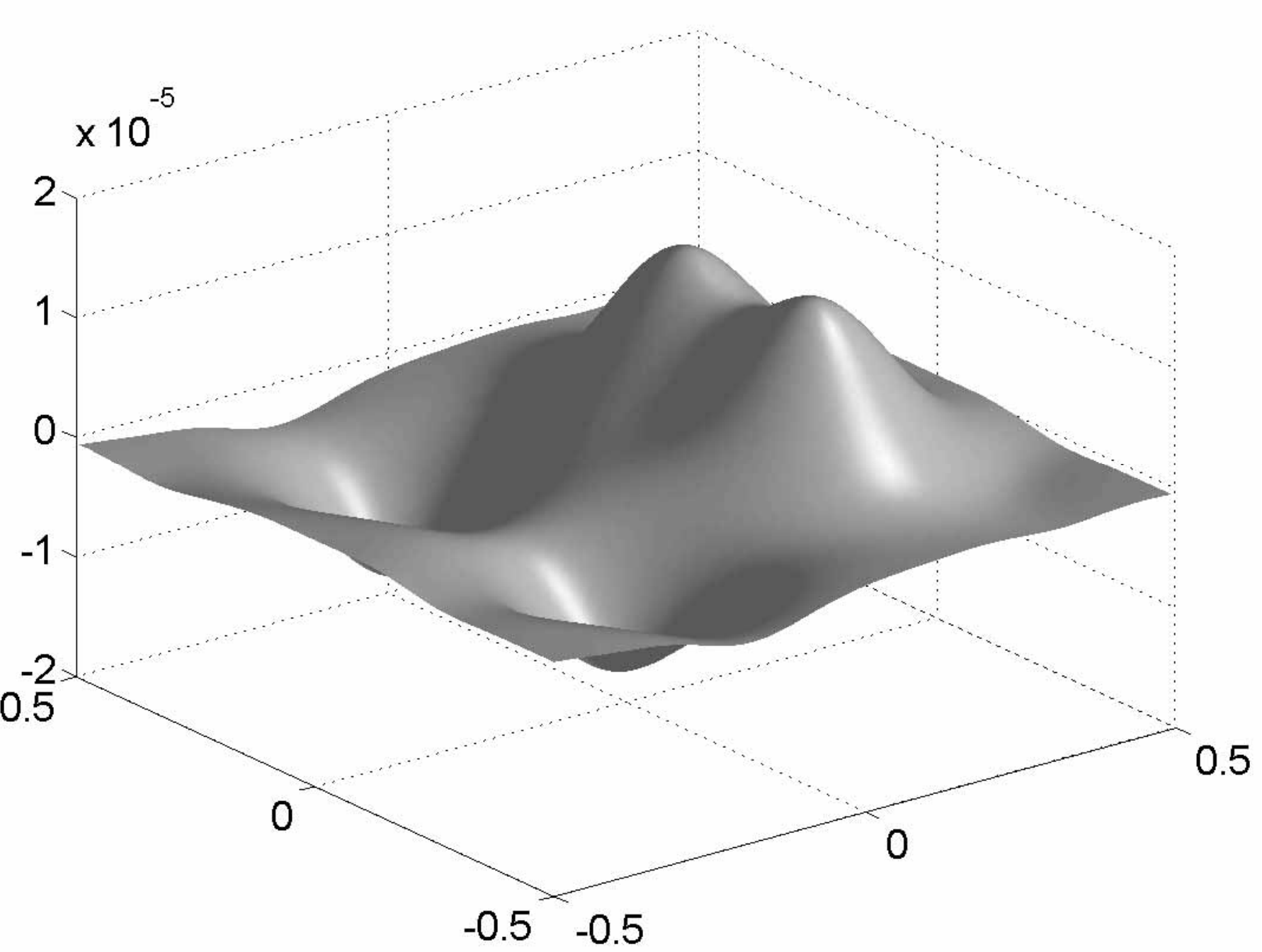}}
\resizebox{1.4in}{!} {\includegraphics{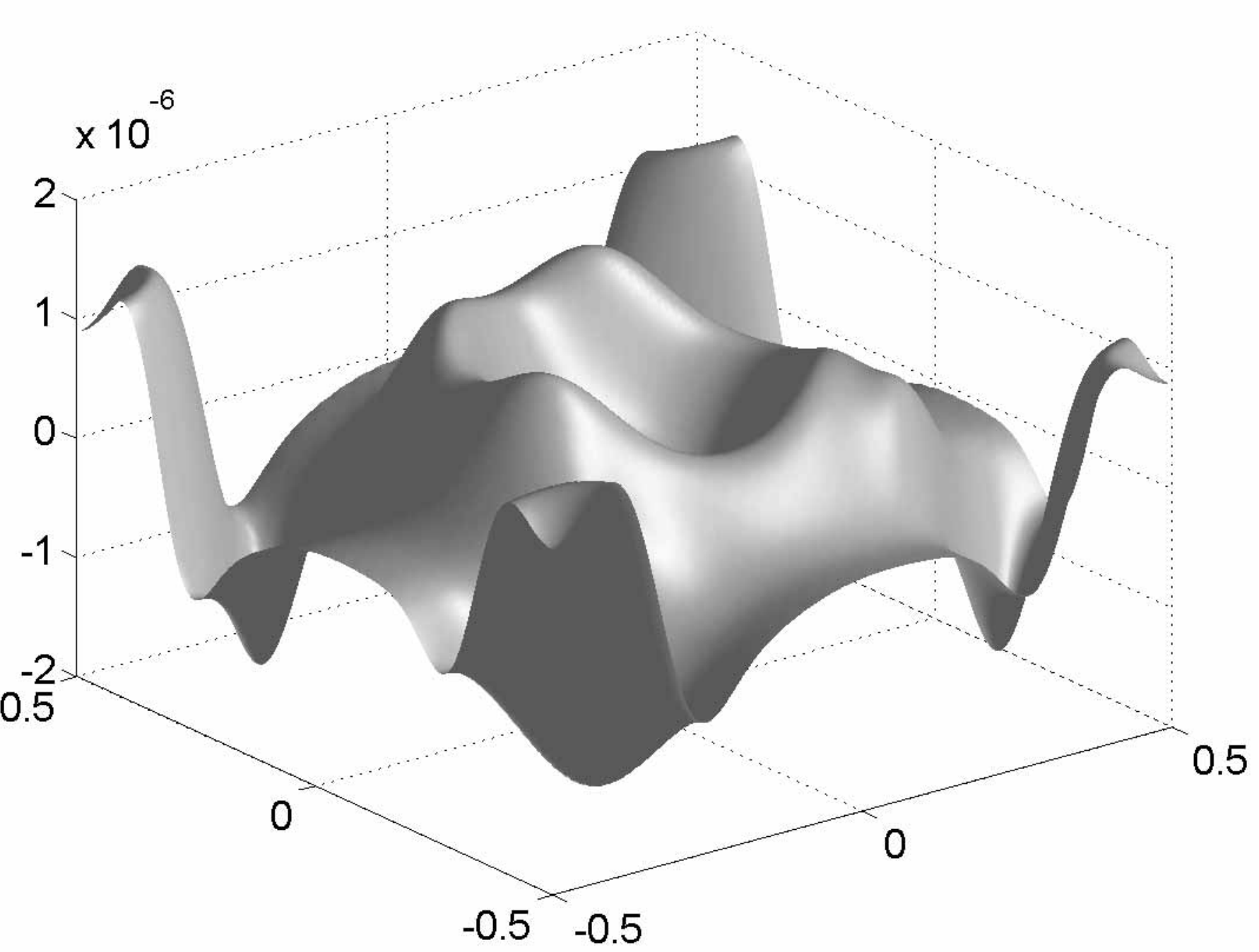}}

{$A_1^\dt(t,\xb)\big|_{x_3=0}$, $A_2^\dt(t,\xb)\big|_{x_3=0}$, $A_3^\dt(t,\xb)\big|_{x_3=0}$
for $\dt=0.01$.}
\end{center}
\caption{Numerical results of the magnetic fields for example \ref{exnr} at t=1.0.
The first row is the solution of the MD system with $\dt=1.0$
and the second row is the solution of the MD system with $\dt=0.01$.
}\label{fig521}
\end{figure}
\end{example}
\begin{example}[\textbf{Harmonic oscillator II}]
\label{exnr2}
Finally, we choose $\Ab^{ex}(\xb)=0$ but include a confining electric
potential of harmonic oscillator type, \ie $V^{ex}(\xb)=C|\xb|^2$. To compete with the effect of
the diffusion term $\btu \psi^\dt$, we choose the large constant $C=100$.
Let us consider the system (\ref{dmnr}) with initial condition
\begin{equation}
\psi^{\dt}\big|_{t=0}= \chi \, \exp{\left(-\frac{(x_1-0.1)^2+(x_2+0.1)^2+x_3^2}{4d^2}\right)},\quad
\chi=(1,0,1,0),\ d=1/16,
\end{equation}
In this case we choose $\dt=10^{-2}$, $\tg t=1/128$, $\tg x=1/64$.
The numerical results are shown in Figure \ref{fig53}.
We see that the wave packet moves in circles due to its interaction with the harmonic potential
and the diffusion term $\btu \psi^\dt$. Note that agreement with the non-relativistic results
is very good also for this test.

\begin{figure} 
\begin{center}
\resizebox{1.4in}{!} {\includegraphics{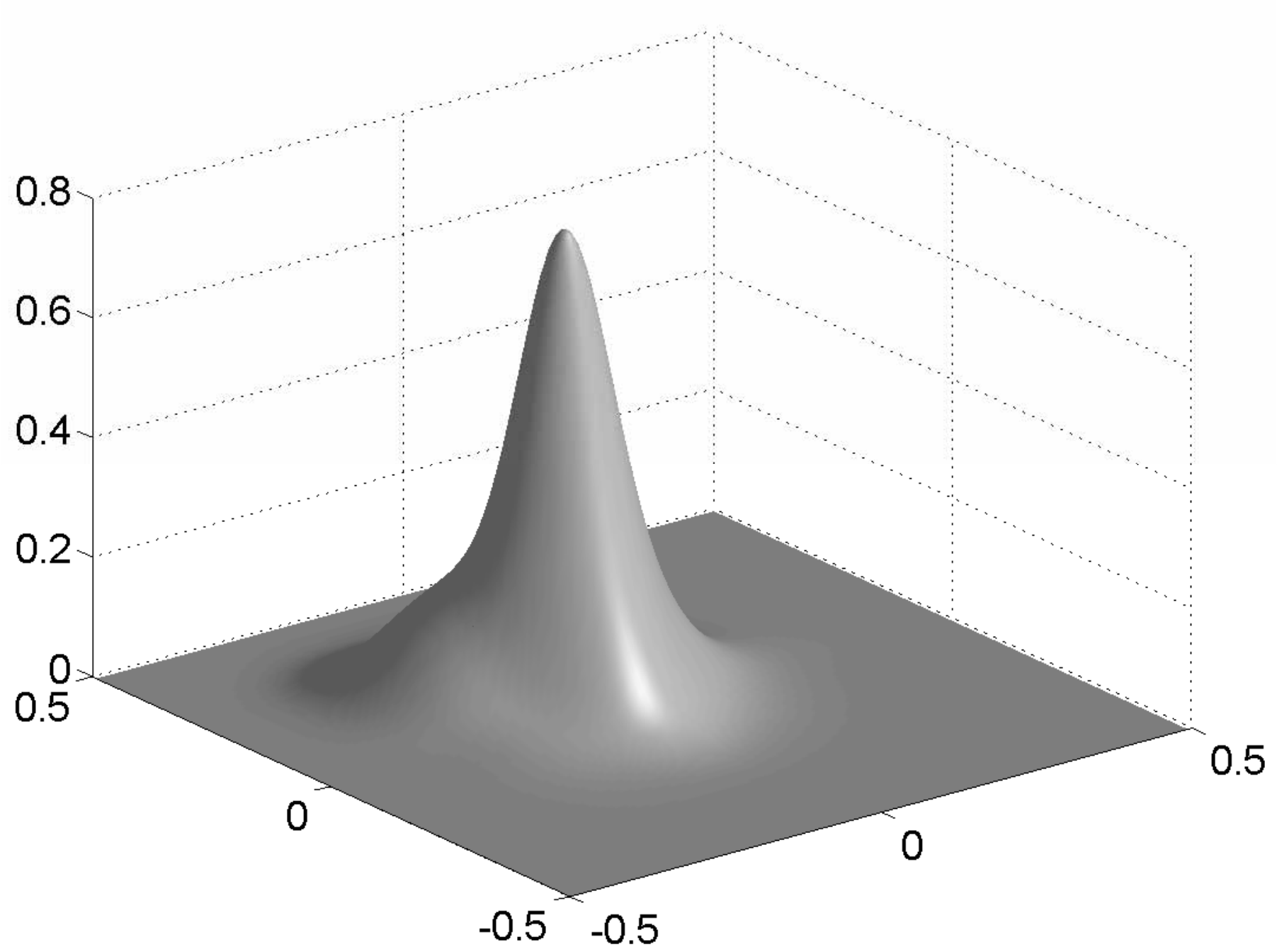}}
\resizebox{1.4in}{!} {\includegraphics{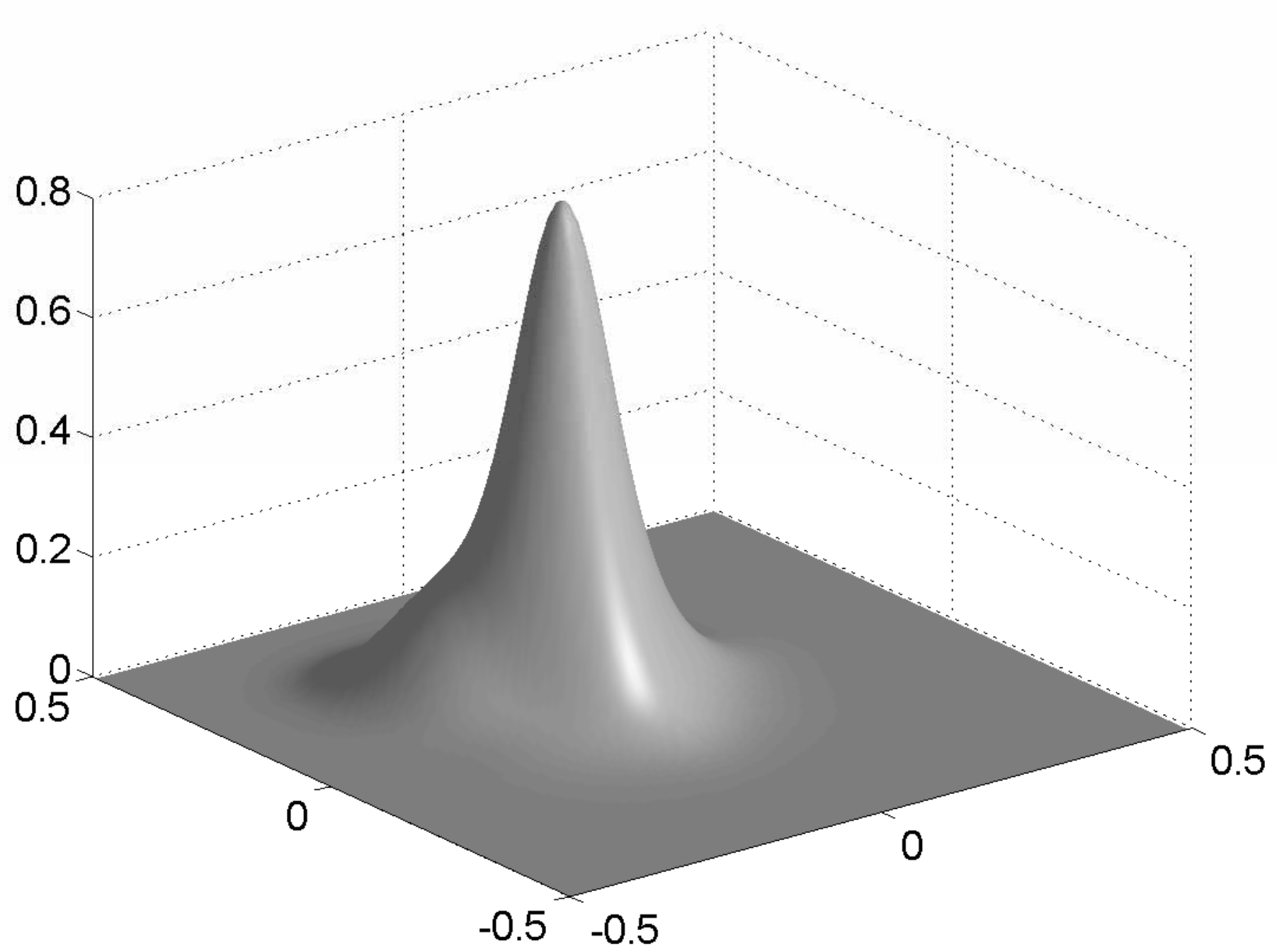}}
\resizebox{1.4in}{!} {\includegraphics{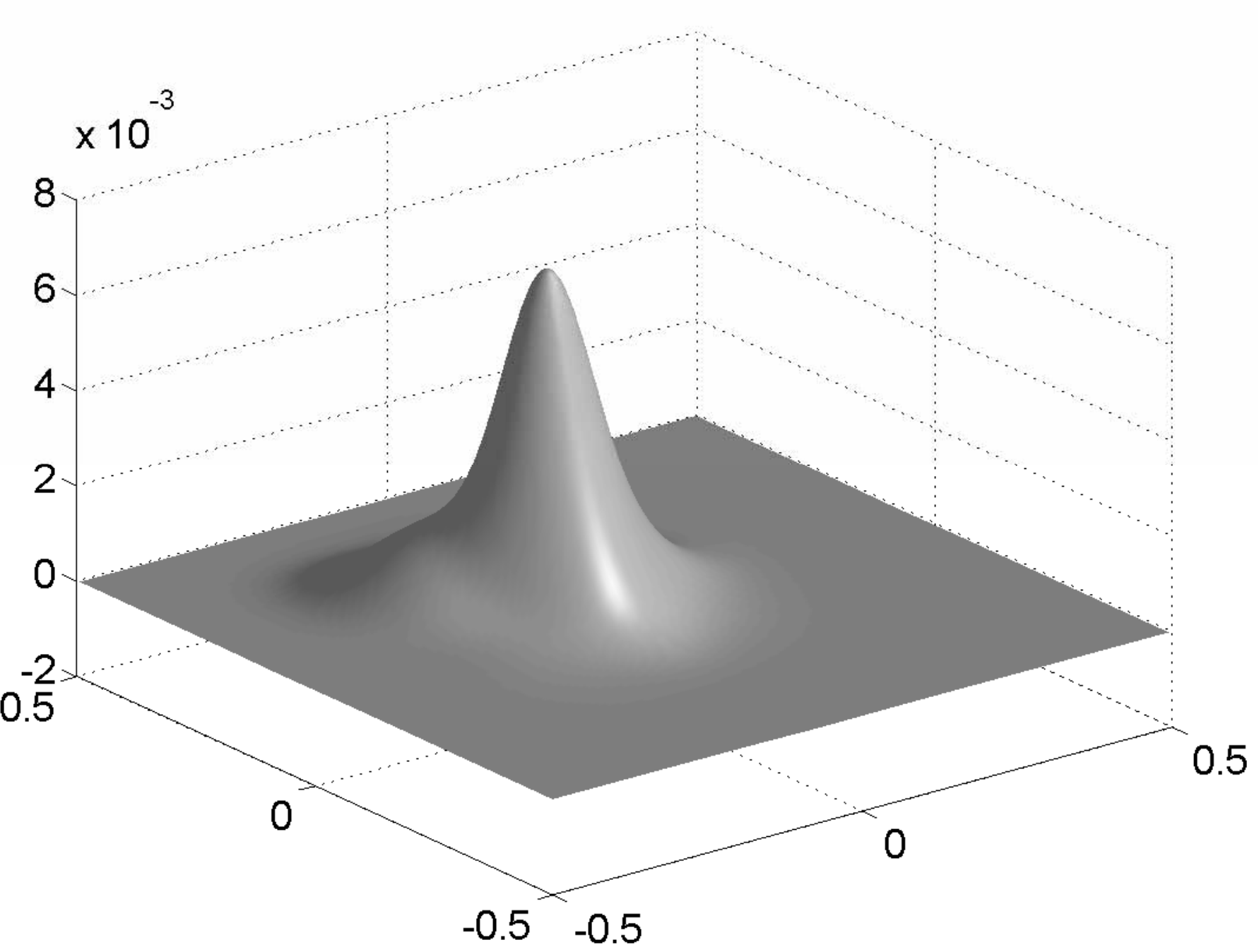}}
\resizebox{1.4in}{!} {\includegraphics{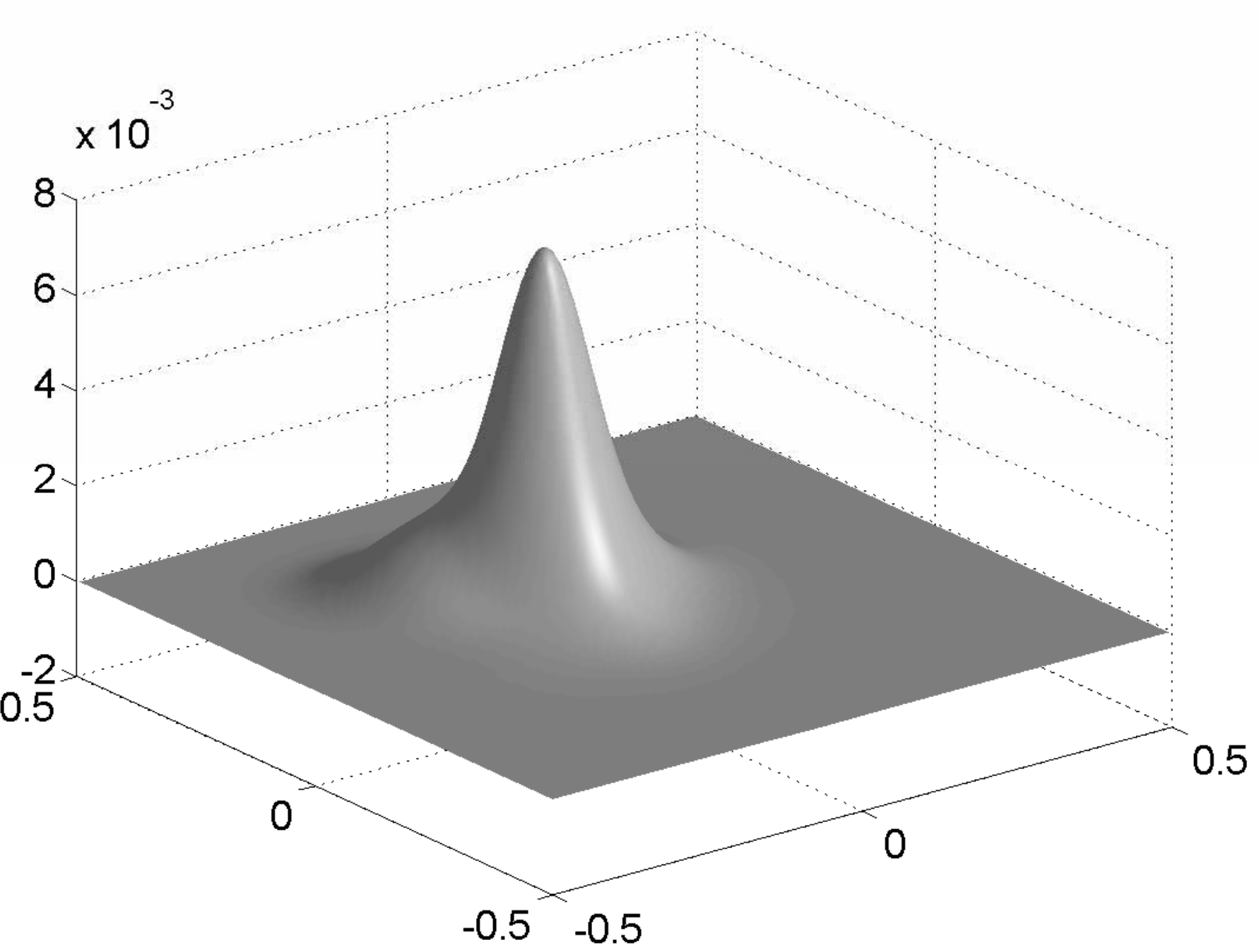}}

{$|\psi^\dt(t,\xb)|^2\big|_{x_3=0}$,
$\left(|\vp_{e}(t,\xb)|^2+|\vp_{p}(t,\xb)|^2\right)\big|_{x_3=0}$,
$V^\dt(t,\xb)\big|_{x_3=0}$ and $V(t,\xb)\big|_{x_3=0}$ at t=0.5.}\vspace{2mm}

\resizebox{1.4in}{!} {\includegraphics{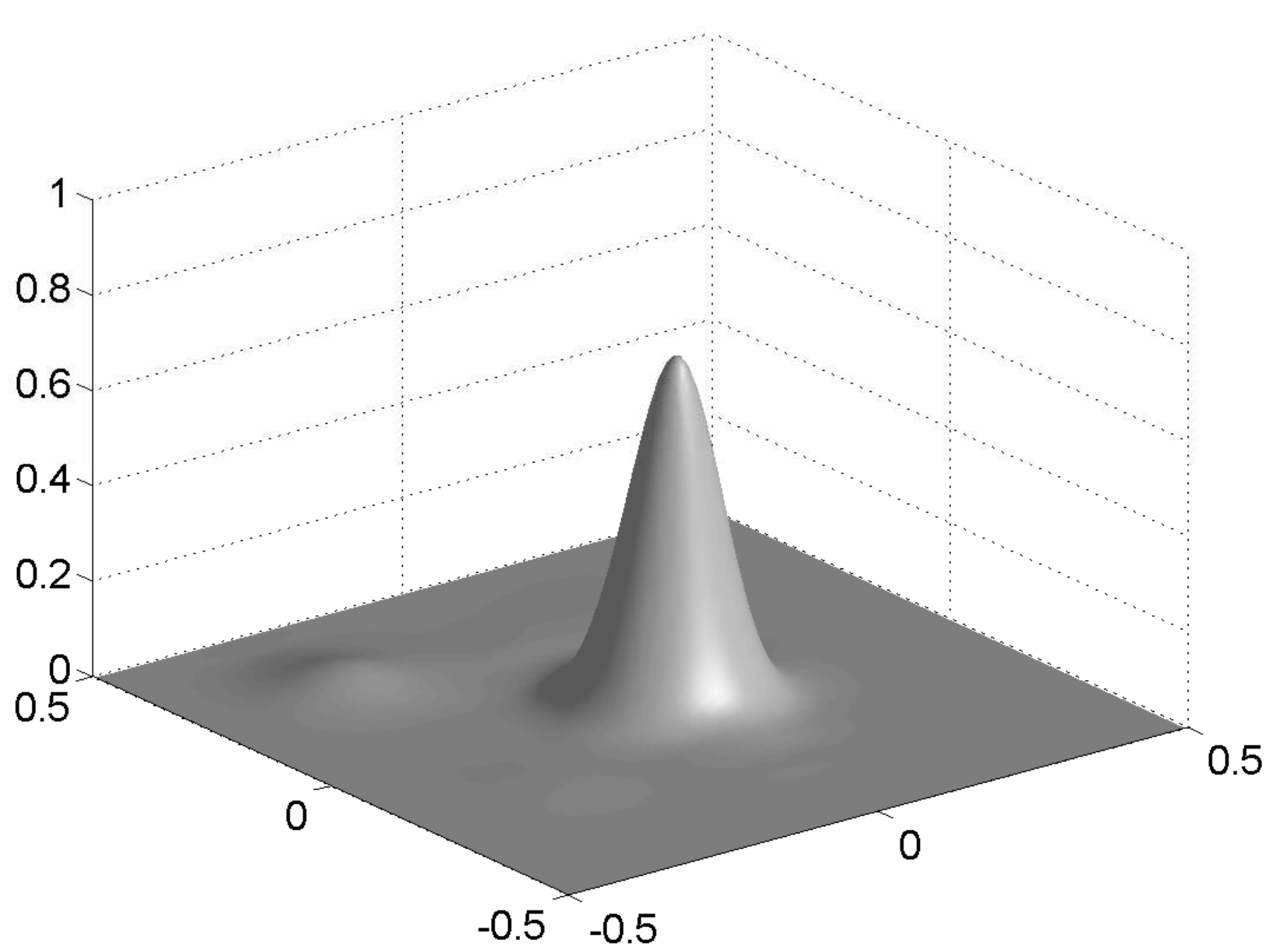}}
\resizebox{1.4in}{!} {\includegraphics{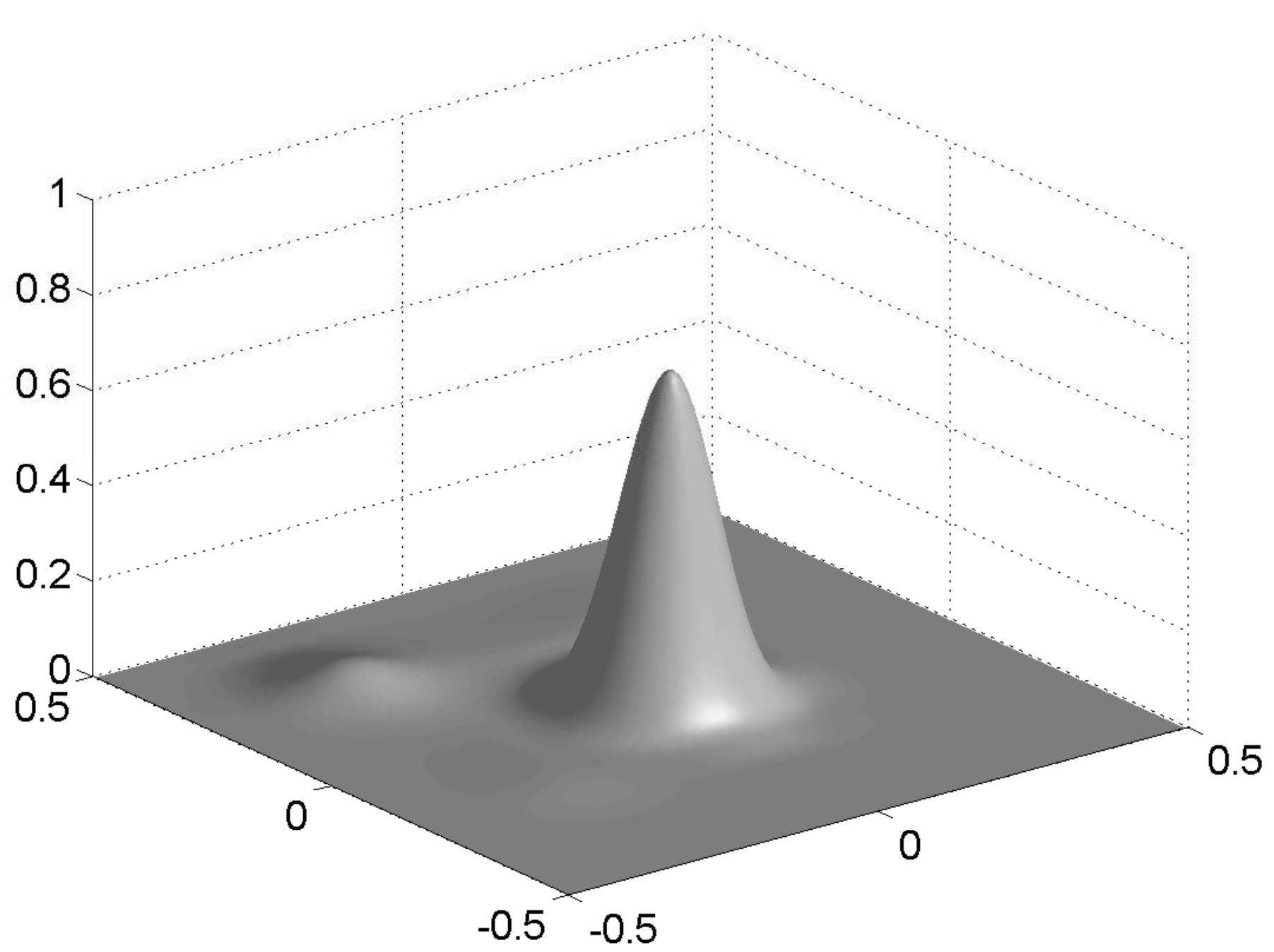}}
\resizebox{1.4in}{!} {\includegraphics{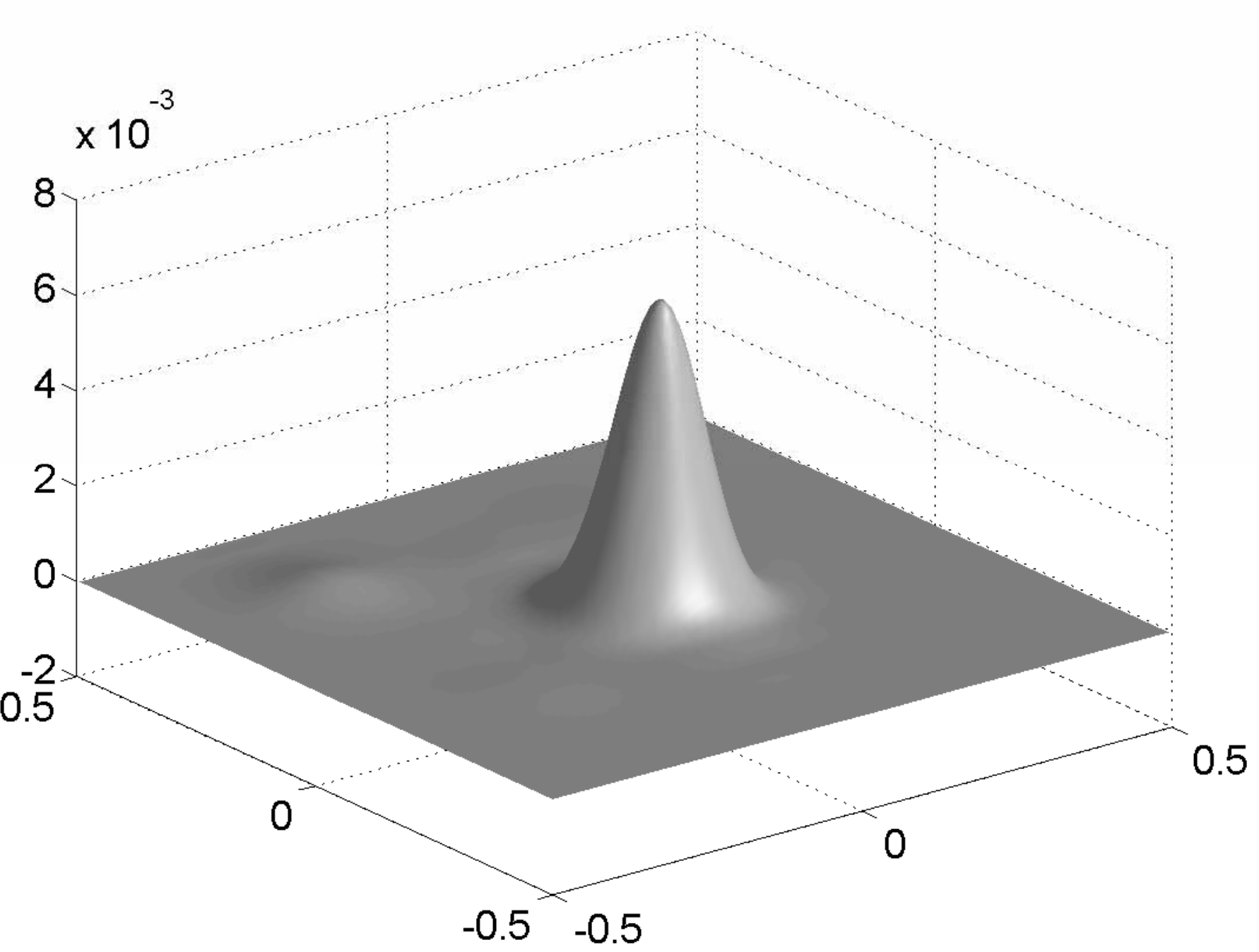}}
\resizebox{1.4in}{!} {\includegraphics{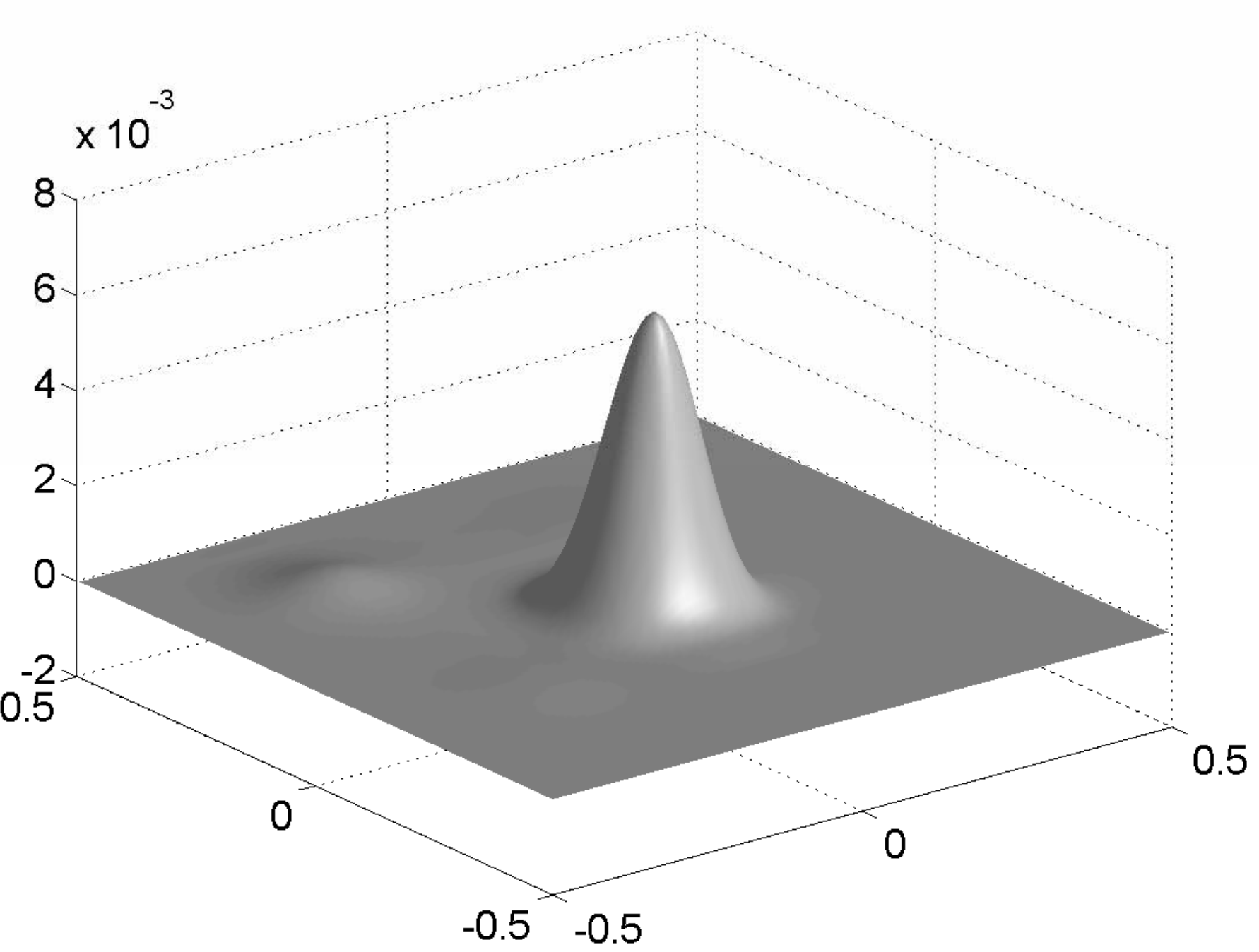}}

{$|\psi^\dt(t,\xb)|^2\big|_{x_3=0}$,
$\left(|\vp_{e}(t,\xb)|^2+|\vp_{p}(t,\xb)|^2\right)\big|_{x_3=0}$,
$V^\dt(t,\xb)\big|_{x_3=0}$ and $V(t,\xb)\big|_{x_3=0}$ at t=1.0.}\vspace{2mm}

\resizebox{1.4in}{!} {\includegraphics{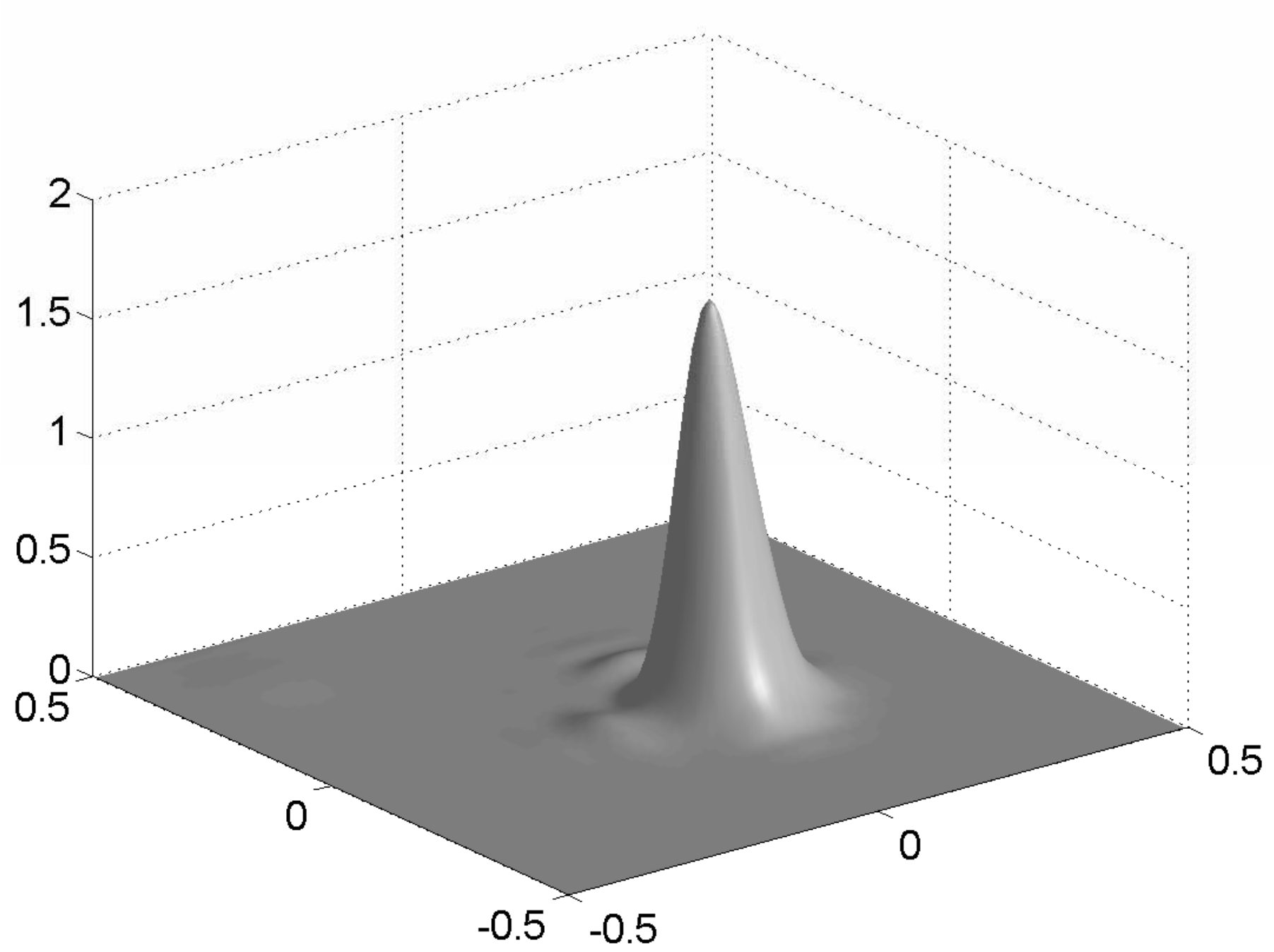}}
\resizebox{1.4in}{!} {\includegraphics{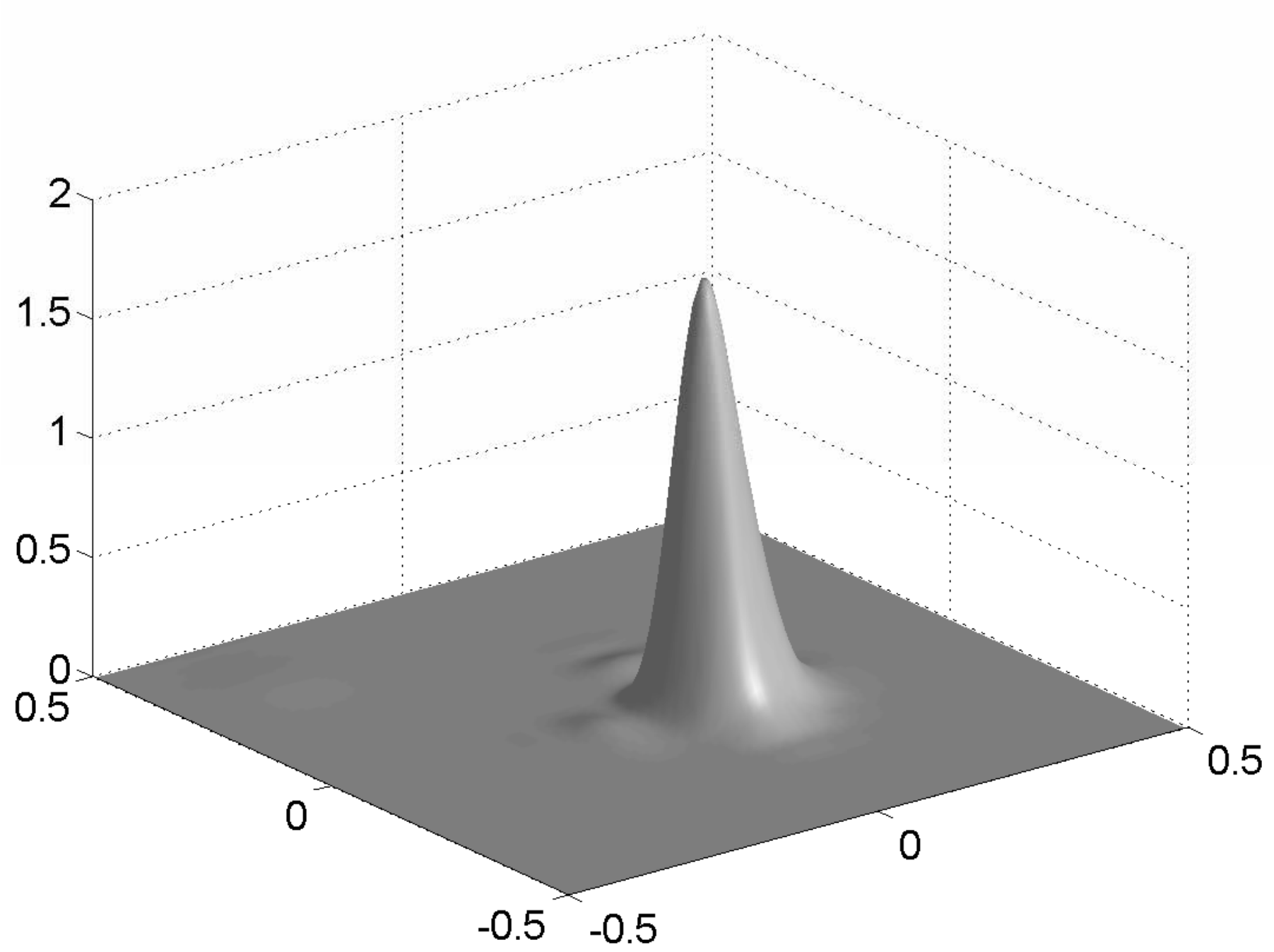}}
\resizebox{1.4in}{!} {\includegraphics{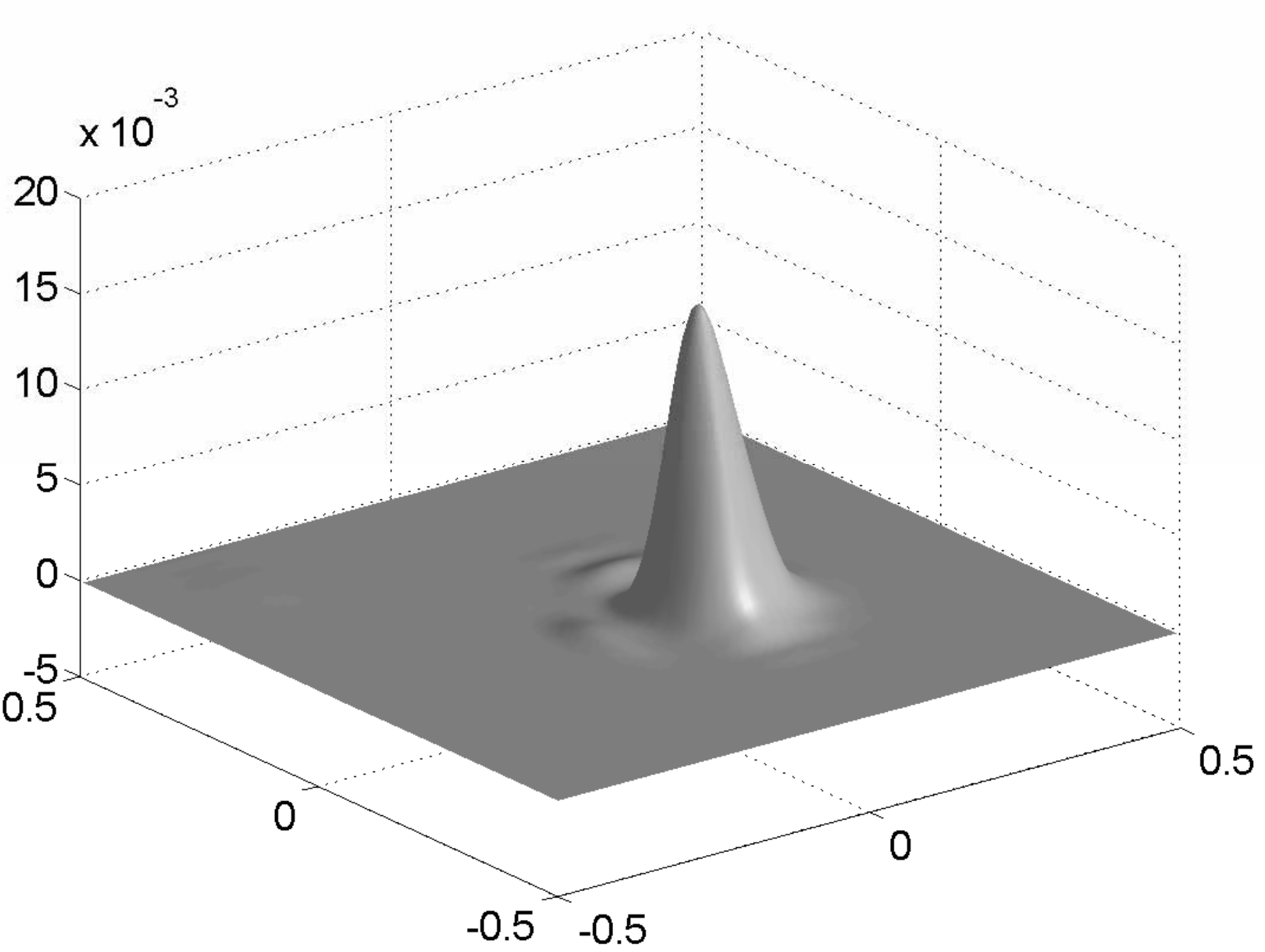}}
\resizebox{1.4in}{!} {\includegraphics{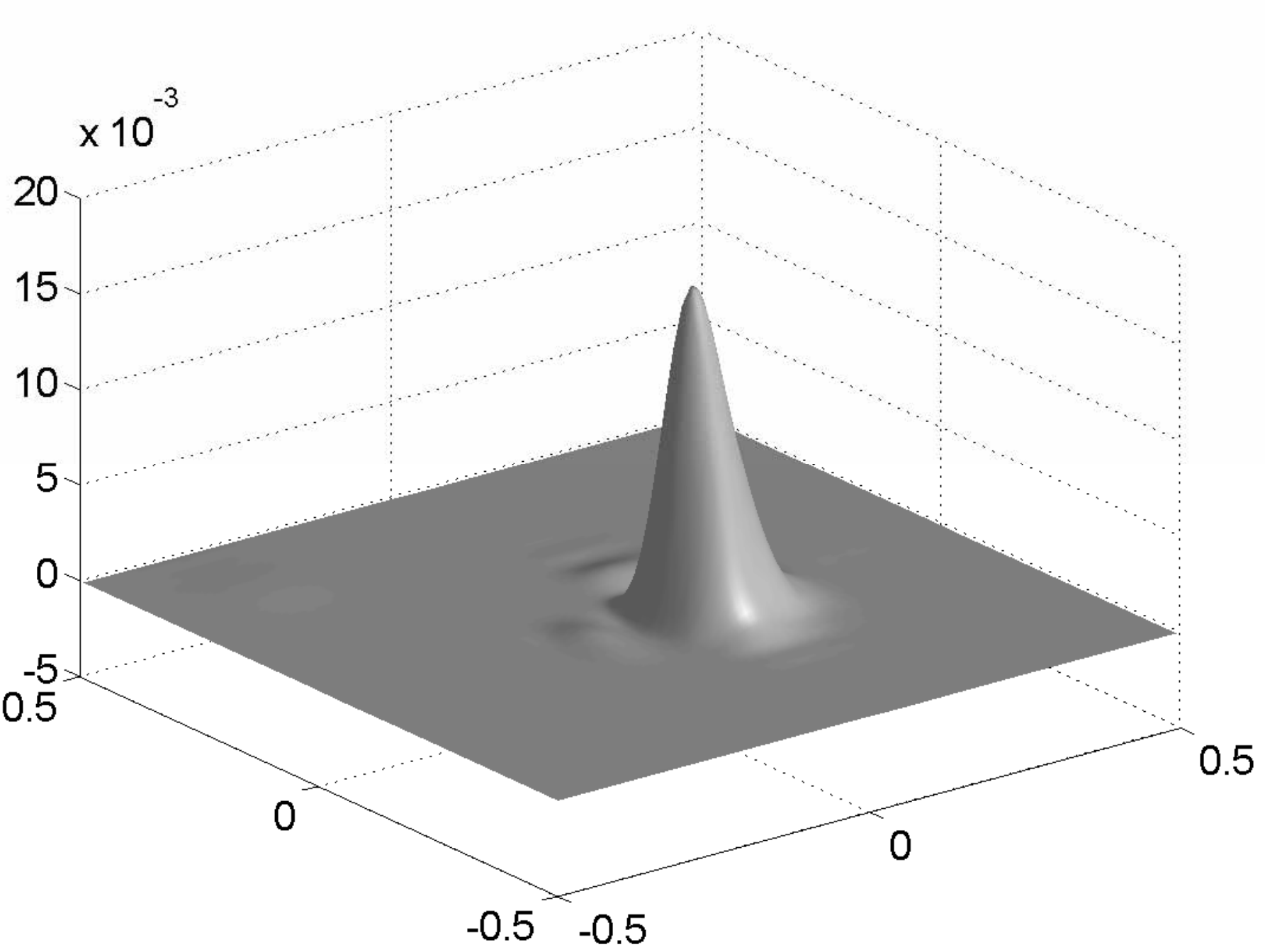}}

{$|\psi^\dt(t,\xb)|^2\big|_{x_3=0}$,
$\left(|\vp_{e}(t,\xb)|^2+|\vp_{p}(t,\xb)|^2\right)\big|_{x_3=0}$,
$V^\dt(t,\xb)\big|_{x_3=0}$ and $V(t,\xb)\big|_{x_3=0}$ at t=1.5.}\vspace{2mm}

\resizebox{1.4in}{!} {\includegraphics{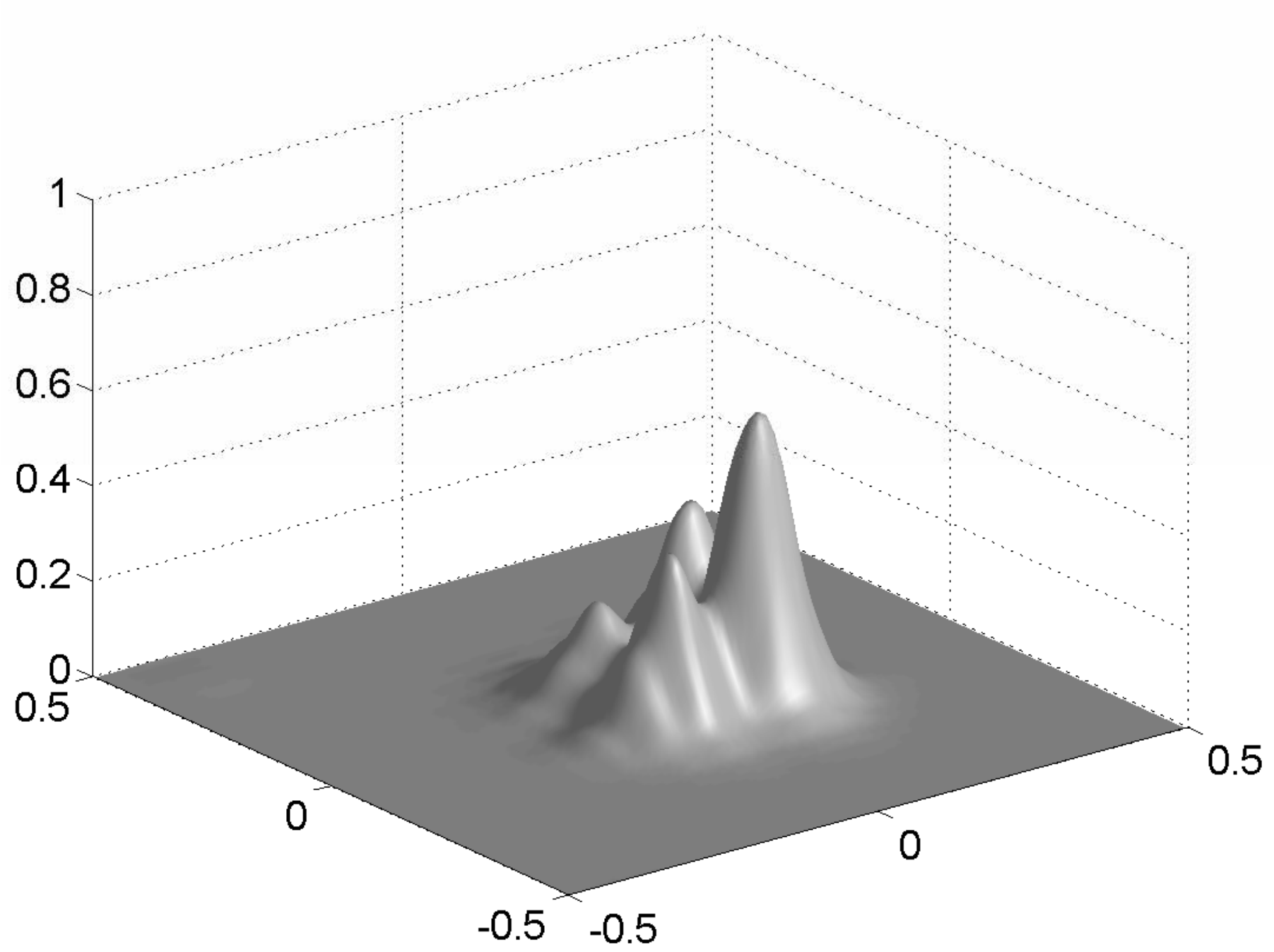}}
\resizebox{1.4in}{!} {\includegraphics{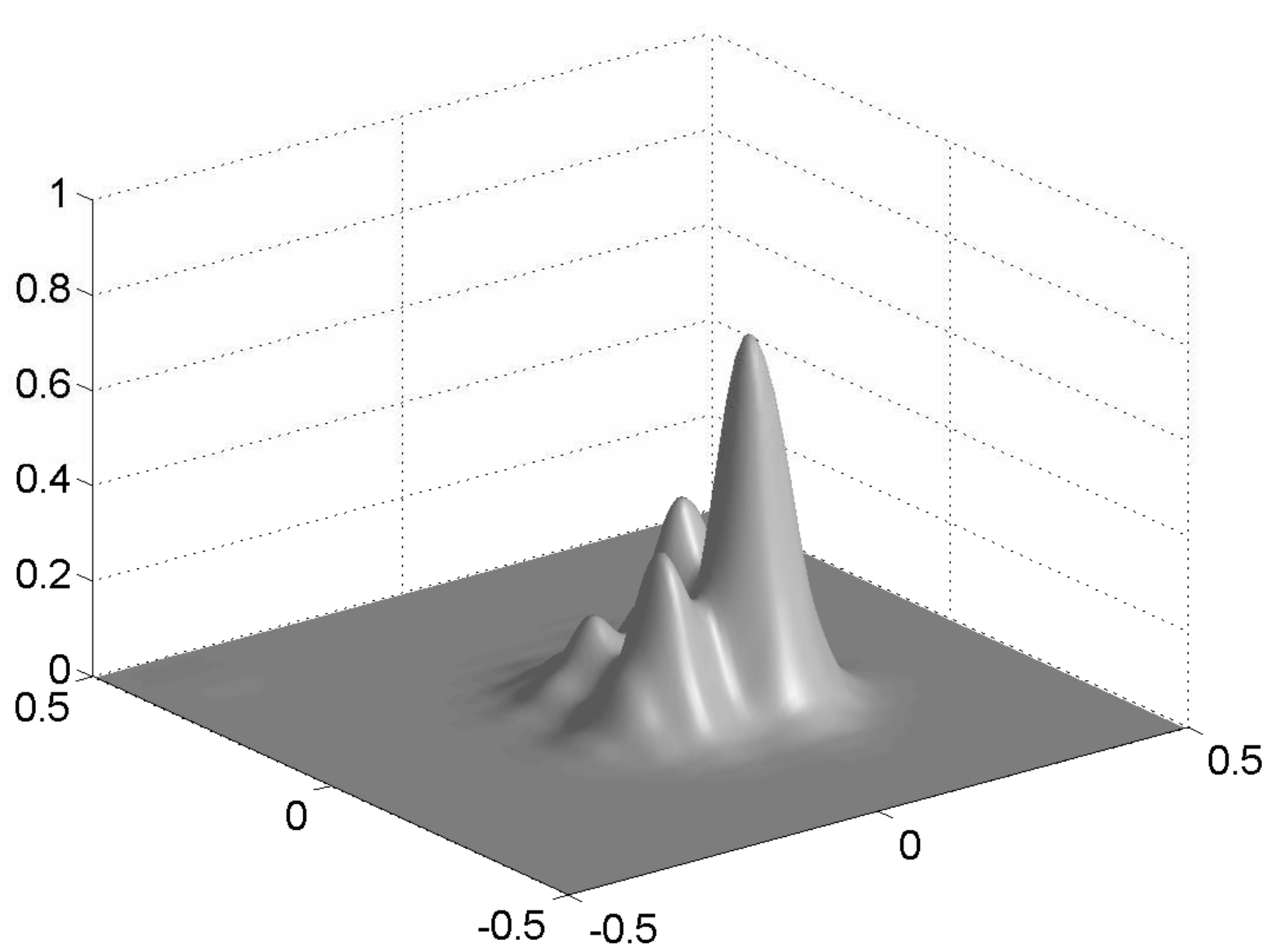}}
\resizebox{1.4in}{!} {\includegraphics{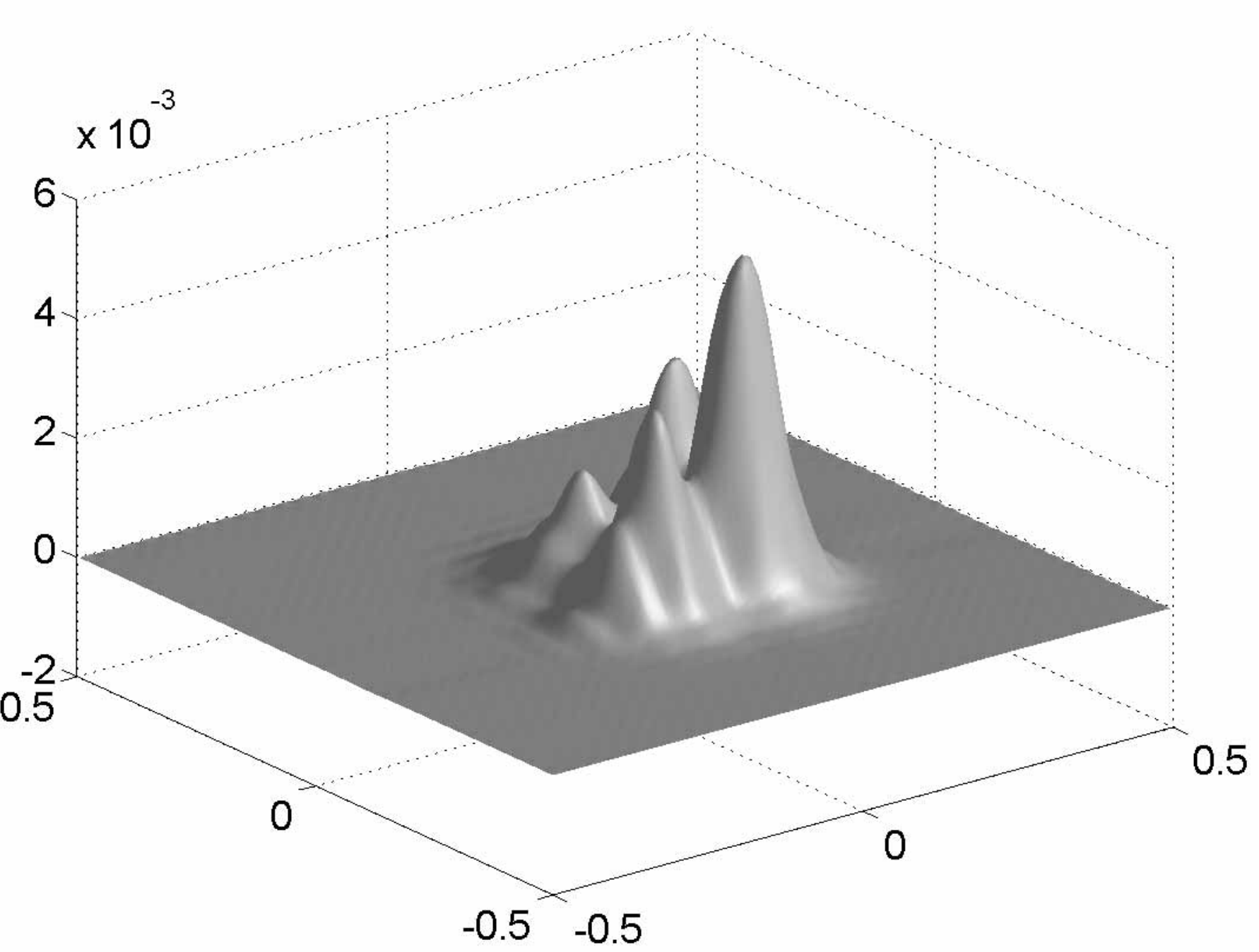}}
\resizebox{1.4in}{!} {\includegraphics{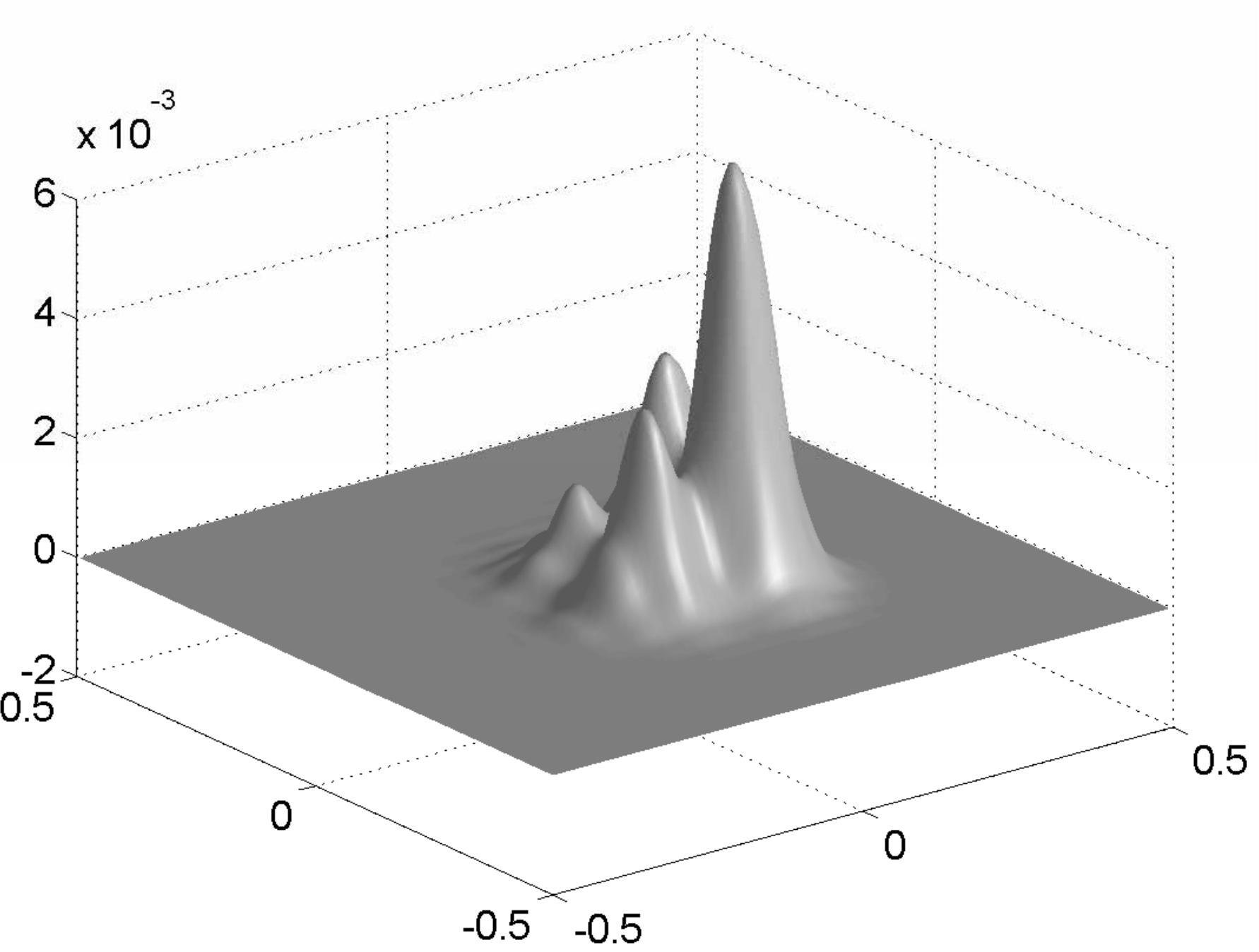}}

{$|\psi^\dt(t,\xb)|^2\big|_{x_3=0}$,
$\left(|\vp_{e}(t,\xb)|^2+|\vp_{p}(t,\xb)|^2\right)\big|_{x_3=0}$,
$V^\dt(t,\xb)\big|_{x_3=0}$ and $V(t,\xb)\big|_{x_3=0}$ at t=2.0.}
\end{center}
\caption{Numerical results of the density for example \ref{exnr2}.
The first and third column are $|\psi^\dt(t,\xb)|^2\big|_{x_3=0}$ and $V^\dt(t,\xb)\big|_{x_3=0}$
for MD system, respectively.
The second and fourth column are $\left(|\vp_{e}(t,\xb)|^2+|\vp_{p}(t,\xb)|^2\right)\big|_{x_3=0}$
and $V(t,\xb)\big|_{x_3=0}$ for Schr\"odinger-Poisson equation, respectively.
Here $\dt=0.01$, $\tg x=1/64, \tg t=1/128$.} \label{fig53}
\end{figure}
\end{example}

\section{Conclusion}
In this work, we presented a time-splitting spectral scheme for the MD system and
similar time-splitting methods for the corresponding asymptotic problems
in the (weakly nonlinear) semi-classical and in the non-relativistic regime.
The proposed scheme conserves the Lorentz gauge condition, is unconditionally stable
and highly efficient as our numerical examples show.
In particular, we presented numerical studies for the creation of positronic modes in
the semi-classical regime as well as numerical evidence for the smallness of the magnetic fields in the
considered non-relativistic scaling. A distinct feature of our time-splitting spectral method, not shared by
previous methods (using the time-splitting spectral approach), is that in the non-relativistic limit,
the scheme exhibits a uniform convergence in the small parameter $\dt$.
\newpar
We finally remark that there are several open questions that deserve further exploration. For example,
it would be an interesting project to derive a better numerical method
for the system of eiconal and transport equations,
describing the semi-classical limit,
which consequently would allow for a more accurate comparison between the solution of the MD system
and limiting WKB-description. A second step then should be the numerical study of the semi-classical
MD equations with stronger nonlinearities,
in particular $\O(1)$-nonlinearities, a so far completely open problem, even from an analytical point of view.


\bibliographystyle{amsplain}


\end{document}